\documentclass{amsart}
\usepackage{etex}
\usepackage{fixltx2e}

\RequirePackage{amsmath}
\RequirePackage{amssymb}
\RequirePackage{amsxtra}
\RequirePackage{amsfonts}
\RequirePackage{latexsym}
\RequirePackage{euscript}
\RequirePackage{amscd}
\RequirePackage{amsthm}
\RequirePackage{xypic}
\usepackage{stmaryrd}
\usepackage{mathtools}
\usepackage[bookmarksopen,bookmarksdepth=2]{hyperref}
\usepackage{mathrsfs}

\usepackage{tikz}
\usepackage{tikz-cd}
\usetikzlibrary{arrows, matrix}
\usepackage[shortlabels]{enumitem}\usetikzlibrary{decorations.pathmorphing}

\xyoption{all}

\usepackage{scalerel}
\newcommand\reallywidehat[1]{\arraycolsep=0pt\relax\begin{array}{c}
\stretchto{
  \scaleto{
    \scalerel*[\widthof{\ensuremath{#1}}]{\kern-.5pt\bigwedge\kern-.5pt}
    {\rule[-\textheight/2]{1ex}{\textheight}} }{\textheight} }{0.5ex}\\           #1\\                 \rule{-1ex}{0ex}
\end{array}
}

\newcommand{\OEA}{\cO_{\mathcal{E},A}}
\newcommand{\OEB}{\cO_{\mathcal{E},B}}
\def\free{\mathrm{free}}
\def\Proj{\mathrm{Proj}}

\newcommand{\Der}{\operatorname{Der}}

\DeclareMathOperator{\Idem}{Idem}

\newcommand{\colim}{\varinjlim}
\newcommand{\twocolim}{\text{2\,-} \colim}

\newcommand{\Ob}{\operatorname{Ob}}

\newcommand{\stimage}{\cZ}

\newcommand{\Kbar}{\overline{K}}

\newcommand{\Sets}{\underline{{Sets}}}

\usepackage{dsfont}

\let\emptyset\varnothing

\def\F{\mathbb F}

\def\Q{\mathbb{Q}}

\def\Z{\mathbb{Z}}

\def\m{\mathfrak m}

\def\Gr{\mathrm{Gr}}

\def\Bbar{\overline{B}}

\def\id{\mathrm{id}}

\def\red{\mathrm{red}}

\def\cont{\mathrm{cont}}

\def\GL{\operatorname{GL}}

\def\Gal{\mathrm{Gal}}

\def\End{\mathrm{End}}

\def\Hom{\mathop{\mathrm{Hom}}\nolimits}

\def\Spec{\mathop{\mathrm{Spec}}\nolimits}
\def\Spf{\mathop{\mathrm{Spf}}\nolimits}

\def\cotimes{\widehat{\otimes}}

\def\m{\mathfrak{m}}

\def\iso{\buildrel \sim \over \longrightarrow}

\swapnumbers

\newcommand{\onto}{\twoheadrightarrow}

\newcommand{\into}{\hookrightarrow}

\newcommand{\To}{\longrightarrow}
\newcommand{\isoto}{\stackrel{\sim}{\To}}

\newlength{\ownl}

\newcommand{\pro}{{\operatorname{pro\, -}}}

\newcommand{\Mor}{{\operatorname{Mor}}}

\newcommand{\AffS}{{\operatorname{Aff}_{/S}}}
\newcommand{\AfffpS}{{\operatorname{Aff}_{\mathrm{pf}/S}}}

\newcommand{\Id}{{\operatorname{Id}}}
\renewcommand{\Im}{{\operatorname{Im}\,}}

\newcommand{\Isom}{{\operatorname{Isom}\,}}

\newcommand{\M}{{\mathcal{M}}}

\newcommand{\V}{{\mathbb{V}}}

\newcommand{\cC}{\mathcal{C}}

\newcommand{\cE}{\mathcal{E}}
\newcommand{\cF}{\mathcal{F}}
\newcommand{\cG}{\mathcal{G}}
\renewcommand{\cH}{\mathcal{H}}

\newcommand{\cO}{\mathcal{O}}

\renewcommand{\cR}{\mathcal{R}}

\newcommand{\cT}{\mathcal{T}}
\newcommand{\cU}{\mathcal{U}}
\newcommand{\cV}{\mathcal{V}}

\newcommand{\cX}{\mathcal{X}}
\newcommand{\cY}{\mathcal{Y}}
\newcommand{\cZ}{\mathcal{Z}}

\newcommand{\gM}{{\mathfrak{M}}}
\newcommand{\gN}{{\mathfrak{N}}}

\newcommand{\gS}{{\mathfrak{S}}}

\newcommand{\gm}{{\mathfrak{m}}}

\newcommand{\gq}{{\mathfrak{q}}}

\newcommand{\tY}{\widetilde{{Y}}}
\newcommand{\tZ}{\widetilde{{Z}}}

\newcommand{\tv}{{\widetilde{{v}}}}

 \newcommand{\xibar    }{\overline{\xi}}

\def\RCS$#1: #2 ${\expandafter\def\csname RCS#1\endcsname{#2}}
\RCS$Revision: 85 $
\RCS$Date: 2011-12-16 18:54:17 -0600 (Fri, 16 Dec 2011) $

 \newcommand{\p}{\mathfrak{p}}

\newcommand{\mf}{\mathfrak}

\newcommand{\rbar}{{\bar{r}}}

 \newcommand{\Qp}{{\Q_p}}

\newcommand{\Zp}{{\Z_p}}

\newcommand{\Fpbar}{{\overline{\F}_p}}

\newcommand{\Fp}{{\F_p}}

\newtheorem{theorem}[subsubsection]{Theorem}
\newtheorem{thm}[subsubsection]{Theorem}
\newtheorem{lemma}[subsubsection]{Lemma}
\newtheorem{lem}[subsubsection]{Lemma}

\newtheorem{cor}[subsubsection]{Corollary}

\newtheorem{prop}[subsubsection]{Proposition}

\theoremstyle{definition}
\newtheorem{df}[subsubsection]{Definition}
\newtheorem{defn}[subsubsection]{Definition}

\theoremstyle{remark}
\newtheorem{remark}[subsubsection]{Remark}
\newtheorem{rem}[subsubsection]{Remark}

\newtheorem{example}[subsubsection]{Example}

\def\numequation{\addtocounter{subsubsection}{1}\begin{equation}}
\def\nummultline{\addtocounter{subsubsection}{1}\begin{multline}}
\def\anumequation{\addtocounter{subsection}{1}\begin{equation}}

\title[Scheme-theoretic images of morphisms of stacks]{``Scheme-theoretic images'' of morphisms of stacks}

\author[M. Emerton]{Matthew Emerton}\email{emerton@math.uchicago.edu}
\address{Department of Mathematics, University of Chicago,
5734 S.\ University Ave., Chicago, IL 60637}

\author[T. Gee]{Toby Gee} \email{toby.gee@imperial.ac.uk} \address{Department of
  Mathematics, Imperial College London,
  London SW7 2AZ, UK}\thanks{The first author was supported in part by the
  NSF grants DMS-1303450, DMS-1601871, and DMS-1902307. The second author was 
  supported in part by a Leverhulme Prize, EPSRC grant EP/L025485/1, Marie Curie Career
  Integration Grant 303605, and by
  ERC Starting Grant 306326.}

\begin{document}
\maketitle
\begin{abstract}
We give criteria for certain morphisms from an algebraic stack to a (not necessarily algebraic)
stack to admit an (appropriately defined) scheme-theoretic image. We apply our
criteria to show
that certain natural moduli stacks of local Galois
representations are algebraic (or Ind-algebraic) stacks.

 \end{abstract}
\setcounter{tocdepth}{1}
\tableofcontents

\section{Introduction}

The goal of this paper is to prove an existence theorem for ``scheme-theoretic 
images'' of certain morphisms of stacks.  We have put scheme-theoretic images
in quotes here because, generally, the objects whose existence we prove will
be certain algebraic spaces or algebraic stacks, rather than schemes.   Like
the usual scheme-theoretic images of morphisms of schemes, though, they will
be closed substacks of the target,
minimal with respect to the property that the given morphism factors through
them. This explains our terminology.

In the case of morphisms of algebraic stacks (satisfying appropriate mild
finiteness conditions), the existence of a scheme-theoretic
image in the preceding sense follows directly from the basic results
about scheme-theoretic images for morphisms of schemes.  
Our interest will be in more general contexts, namely, those in which
the source of the morphism is assumed to be an algebraic stack,
but the target is not; in particular, we will apply our results in one
such situation to construct moduli stacks of Galois representations.

\subsection{Scheme-theoretic images}
We put ourselves in the setting of
stacks in groupoids defined on the big \'etale site of a locally Noetherian
scheme $S$ all of whose local rings $\mathcal O_{S,s}$ at finite type
points
$s \in S$ are $G$-rings. (See Section~\ref{subsec:notation} for any unfamiliar terminology.)
Recall that in this context,
Artin's representability theorem gives a characterisation
of algebraic stacks which are locally of finite presentation over $S$ among all such stacks: namely,
algebraic stacks which are locally of finite presentation over $S$ are precisely
those stacks in groupoids $\cF$ on the big
\'etale site of $S$ that satisfy:

\begin{enumerate}[label={[\arabic*]}]
\item  $\cF$ is limit preserving; \item  (a) $\cF$ satisfies the Rim--Schlessinger condition (RS),
and\\ (b) $\cF$ admits effective Noetherian versal rings at all finite
type points;
\item the diagonal $\Delta: \cF \to \cF\times_S \cF$ is representable by algebraic spaces;\item openness of versality.
\end{enumerate}
See Section~\ref{sec: artin's axioms and
  stacks in groupoids} below for an explanation of these axioms,
and Theorem~\ref{thm: Artin representability} for Artin's theorem.
(In Subsection~\ref{subsec:axiom 4} we also
introduce two further axioms, labelled~[4a] and~[4b],
that are closely related to~[4].  These are used in discussing the
examples of Section~\ref{sec: examples}, but not in the proof
of our main theorem.)

We will be interested in {\em quasi-compact} morphisms $\xi: \cX \to \cF$ of stacks on the big \'etale site of $S$,
where $\cX$ is algebraic and locally of finite presentation
over $S$,  
and $\cF$ is assumed to have a diagonal that is representable
by algebraic spaces and is locally of finite presentation (i.e.\ $\cF$
satisfies~[3], and a significantly weakened form of~[1]).
In this context, we are able
to define a substack
$\cZ$ of $\cF$ which we call the
{\em scheme-theoretic image}
of $\xi$.  (The reason for assuming that $\xi$ is quasi-compact
is that this seems to be a minimal requirement for the formation
of scheme-theoretic images to be well-behaved even for morphisms of
schemes, e.g.\ to be compatible with \emph{fpqc}, or even
Zariski, localisation.) 
If $\cF$ is in fact an algebraic stack, locally of finite
presentation over $S$, then $\cZ$ will coincide with
the usual scheme-theoretic image of $\xi$.  In general,
the substack $\cZ$ will itself satisfy Axioms~[1] 
and~[3].

Our main result is the following theorem (see Theorem~\ref{thm:main-intro assuming [1]}).

\begin{theorem}
\label{thm:main-intro}
Suppose that $\xi: \cX \to \cF$ is a {\em proper} morphism, where $\cX$
is an algebraic stack, locally of finite presentation over $S$, and $\cF$ is a
stack  over $S$
satisfying~{\em [3]},
and whose diagonal is furthermore locally of finite
presentation.
Let $\cZ$ denote the scheme-theoretic image of $\xi$ as discussed
above, and
suppose that $\cZ$ satisfies~{\em [2]}. Suppose also that $\cF$ admits
{\em (}not necessarily Noetherian{\em )}
versal rings at all finite type points
{\em (}in the sense of Definition~{\em \ref{def:versal ring}}
below{\em )}. 

Then $\cZ$ is an algebraic stack, locally of finite presentation over $S$; 
the inclusion $\cZ \hookrightarrow \cF$
is a closed immersion; and the morphism $\xi$ factors 
through a proper, scheme-theoretically surjective morphism $\overline{\xi}: \cX
\to \cZ$. Furthermore, if $\cF'$ is a substack of $\cF$ for which
the monomorphism $\cF' \hookrightarrow \cF$ is representable by algebraic 
spaces and of finite type
{\em (}e.g.\ a closed substack{\em )} with the property that
$\xi$ factors through $\cF'$, then $\cF'$ contains~$\cZ$.
\end{theorem}

We show that if~$\cF$ satisfies~[2], then the assumption that~ $\cZ$
satisfies~[2] follows from the other assumptions in Theorem~\ref{thm:main-intro}
(see Lemma~\ref{lem:we get a sheaf} below).   
By applying the theorem in the case when $\cZ = \cF$ (in which case we say that
$\xi$ is \emph{scheme-theoretically dominant}),
and taking into account this remark,
we obtain the following corollary.  

\begin{cor}
\label{cor:representability-intro}
If $\cF$ is an \'etale stack in groupoids over $S$,
satisfying~{\em [1]}, ~{\em [2]} and~{\em [3]},
for which there exists a
scheme-theoretically dominant proper morphism 
$\xi: \cX \to \cF$ whose domain is an algebraic stack
locally of finite presentation over $S$,
then $\cF$ is an algebraic stack.
\end{cor}

\subsubsection{Some remarks on Theorem~\ref{thm:main-intro} and its proof}
The proof of Theorem~\ref{thm:main-intro} occupies most of the first
three sections of the paper.
One way for the reader to get an idea of the argument is
to read Lemma~\ref{lem:monomorphism}, and then to turn directly
to the proof of Theorem~\ref{thm:main-intro assuming [1]}, taking the
various results referenced in the argument (including even the
definition of the scheme-theoretic image~$\cZ$) on faith.

It is often the case (and it is the case in the proof
of Theorem~\ref{thm:main-intro assuming [1]}) that
the main problem to be overcome when using Artin's axiomatic
approach to proving that a certain stack is algebraic is the verification of axiom
[4] (openness of versality).  Lemma~\ref{lem:monomorphism} (which is a stacky version
of~\cite[Lem.~3.10]{MR0260746})
shows that one can automatically get a slightly weaker version of [4],
where ``smooth'' is replaced by ``unramified'', if the  axioms [1], [2],
and [3] are satisfied.
Morally, this suggests that anything satisfying axioms [1], [2], and
[3] is already close to an Ind-algebraic stack, because it admits
unramified maps (which are at least morally quite close to immersions) from algebraic stacks that are even formally smooth at any
particular point.  So to prove Theorem~\ref{thm:main-intro}, one has
to build on this idea, and then use the extra hypothesis (namely, that there is a proper surjection from an algebraic stack onto the
stack~$\cZ$) to prove axiom~[4].  The argument is ultimately
topological, using properness to eliminate the possibility of having
more and more components building up Zariski locally around a point.

\subsection{Moduli of finite height Galois representations}\label{subsec: intro on
the picture}
\label{subsec: finite flat}
The results in this paper were developed with a view
to applications to the theory of Galois representations, and in particular to
constructing moduli stacks of mod~$p$ and $p$-adic representations of the
absolute Galois group $\Gal(\Kbar/K)$ of a finite extension $K/\Qp$. These
applications will be developed more fully in the
papers~\cite{emertongeepicture,CEGSKisinwithdd}, and we refer the interested
reader to those papers for a fuller discussion of our results and
motivations. 

Let $\rbar:\Gal(\Kbar/K)\to\GL_n(\Fp)$ be a continuous representation. The
theory of deformations of $\rbar$ --- that is, liftings of $\rbar$ to continuous
representations $r:\Gal(\Kbar/K)\to\GL_d(A)$, where $A$ is a complete local
ring with residue field $\Fp$ --- is extremely important in the Langlands
program, and in particular is crucial for proving automorphy lifting theorems
via the Taylor--Wiles method. Proving such theorems often comes
down to studying the moduli spaces of those
deformations which satisfy various $p$-adic Hodge-theoretic conditions; see for
example~\cite{KisinModularity,KisinFM}. 

From the point of view of
algebraic geometry, it seems unnatural to only study ``formal'' deformations of
this kind, and Kisin observed about ten years ago that results on the reduction
modulo $p$ of two-dimensional crystalline representations suggested that there
should be moduli spaces of $p$-adic representations (satisfying certain $p$-adic
Hodge theoretic conditions, for example finite flatness) in which the residual
representations $\rbar$ should be allowed to vary; in particular, the special
fibres of these moduli spaces would be moduli spaces of (for example) finite
flat representations of $\Gal(\Kbar/K)$. Unfortunately, there does not seem to be any simple way of directly constructing
such moduli spaces, and until now their existence has remained a mystery. (We
refer the reader to the introduction to~\cite{emertongeepicture} for a further
discussion of the difficulties of directly constructing moduli spaces of mod~$p$
representations of $\Gal(\Kbar/K)$.) 

Mod~$p$ and $p$-adic Galois representations are studied via integral $p$-adic
Hodge theory; for example, the
theories of $(\varphi,\Gamma)$-modules and Breuil--Kisin modules.
Typically, one begins by analysing $p$-adic representations of the absolute
Galois group~$\Gal(\Kbar/K_\infty)$ of some highly ramified infinite
extension~$K_\infty/K$. (In the theory of $(\varphi,\Gamma)$-modules
this extension is the cyclotomic extension, but in the theory of Breuil--Kisin
modules, it is a non-Galois extension obtained by
extracting $p$-power roots of a uniformiser.)
Having classified these
representations in terms of semilinear algebra (modules over some
ring, equipped with a Frobenius), one then separately
considers the problem of descending the classification to
representations of~$\Gal(\Kbar/K)$. 

More precisely, by the theory of~\cite{MR1106901}, continuous
mod~$p^a$ representations of $\Gal(\Kbar/K_\infty)$ are classified by
\'etale $\varphi$-modules, which are modules over a Laurent series
ring, equipped with a Frobenius.  Following the paper~\cite{MR2562795}
of Pappas and Rapoport,
we consider a moduli stack $\cR$ of \'etale $\varphi$-modules, which,
for appropriate choices of the Frobenius on the Laurent series ring  can be thought of informally as a moduli stack
classifying $\Gal(\overline{K}/K_{\infty})$-representations. (To keep this paper at a
reasonable length, we do not discuss the problem of descending our
results to representations
of~$\Gal(\Kbar/K)$; this is addressed in the papers~\cite{CEGSKisinwithdd,emertongeepicture}.) Pappas and Rapoport prove various properties of the stack~$\cR$
(including that it \emph{is} a stack, which they deduce from deep
results of Drinfeld~\cite{MR2181808} on the \emph{fpqc} locality of
the notion of a projective module over a Laurent series ring),
including that its diagonal is representable by algebraic spaces.
In the present paper
we prove the following theorem about the geometry of $\cR$ (see Theorem~\ref{thm: R
  is an ind algebraic stack and satisfies 1} below).

\begin{theorem} 
\label{thm: existence of Z^a, introduction version}
The stack $\cR$ is an Ind-algebraic stack.  More precisely,
we may write $\cR \cong \varinjlim_{F} \cR_F,$
where each $\cR_F$ is a finite type algebraic stack over $\Z/p^a\Z$,
and where each transition morphism in the inductive limit is a closed
immersion.
\end{theorem}

The theorem is proved by applying Theorem~\ref{thm:main-intro} to
certain morphisms $\cC_F \to \cR$ (whose sources~$\cC_F$ are algebraic
stacks) so as to prove that their scheme-theoretic images $\cR_F$ are
algebraic stacks; we then show that $\cR$ is naturally identified with
the inductive limit of the $\cR_F$, and thus that it is an
Ind-algebraic stack. (The index~$F$ is a certain element of a power
series ring; replacing the indexing set with a cofinal subset given by
powers of a fixed~$F$, one can write~$\cR$ as an Ind-algebraic stack
with the inductive limit being taken over the natural numbers.)
 
The stacks~$\cC_F$ were defined by Pappas and Rapoport
in~\cite{MR2562795}, where it is proved that they are algebraic,
and that the morphisms $\cC_F \to \cR$ are representable by algebraic spaces and
proper. The definition of the stacks~$\cC_F$
is motivated by the
papers~\cite{KisinModularity,KisinCrys,MR2373358}, in which Kisin (following
work of Breuil and Fontaine) studied certain
integral structures on \'etale $\varphi$-modules, in particular (what are
now called)
Breuil--Kisin modules of height at most $F$, where~$F$ is a power of
an Eisenstein polynomial for a finite extension~$K/\Qp$.
A Breuil--Kisin module is a module with
Frobenius over a power series ring, satisfying a condition
(depending on~$F$) on the cokernel of the Frobenius. Inverting the
formal variable~$u$ in the power series ring gives a functor from the
category of Breuil--Kisin modules of height at most~$h$ to the category of
\'etale $\varphi$-modules. We say that the Galois representations
corresponding to the \'etale $\varphi$-modules in the essential image of
this functor have \emph{height at most $F$}.

With an eye to future applications, we work in a general context
in this paper, and in particular we allow considerable flexibility in
the choice of the polynomial~$F$ and the Frobenius on the
coefficient rings. In particular our \'etale $\varphi$-modules do
not obviously correspond to representations of
some~$\Gal(\overline{K}/K_\infty)$ (but the case that they do is the
main motivation for our constructions and theorems). A key point in the argument (since it is one of the hypotheses of
Theorem~\ref{thm:main-intro}) is that $\cR$ admits versal rings at all
finite type points, and that the scheme-theoretic images $\cR_F$
satisfy~[2]. In the setting of representations
of~$\Gal(\overline{K}/K_\infty)$,  versal rings for~$\cR$ are given
be the
framed deformation rings associated to continuous mod $p$ representations
of $\Gal(\overline{K}/K_{\infty})$, as we show
in~\cite{emertongeepicture,CEGSKisinwithdd}; then the basic input to
the verification of~[2] for~$\cR_F$ is the theory
of finite height framed Galois deformation rings (which are proved to be Noetherian by Kim in~\cite{MR2782840}, reflecting the fact that the scheme-theoretic
images  $\cR_F$ turn out to be finite-type algebraic stacks). In our
more general setting we work instead with lifting rings for $\varphi$-modules.

\subsection{Further remarks on the contents and organisation of the paper}
In the remainder of Section~1,
we describe our notation and conventions (Subsection~\ref{subsec:notation}),
and also record some simple lemmas in local algebra that will be needed
later in the paper (Subsection~\ref{subsec:local}).

In
Section~\ref{sec: artin's axioms and stacks in groupoids},
we explain Artin's axioms (as listed above) in some detail,
and present many related definitions and results.
We do not expect the reader experienced in the theory of stacks
to find much that is novel in this section, and indeed, many of
the results that we have included
are simple variants of results that are already
in the literature.  
In light of this, it might be worthwhile to
offer some justification for the length of this section.
Primarily, we have been guided by the demands of the arguments
presented in Section~\ref{sec:images};  these
demands have largely dictated the choice of material
presented in 
Section~\ref{sec: artin's axioms and stacks in groupoids},
and its organisation. 
Additionally,
we anticipate that the typical reader of this paper interested in
the application of our results to the moduli of Galois representations will not already be completely familiar with
the foundational results discussed in this section,
and so we have made an effort to include a more careful
discussion of these results, as well as more references to the literature, than
might strictly be required for the typical reader interested only in
Theorem~\ref{thm:main-intro} and Corollary~\ref{cor:representability-intro}.

Moreover, the basic idea of our argument came from a careful reading
of~\cite[\S 3]{MR0260746}, especially the proof of Theorem~3.9 therein,
which provided a model argument for deducing Axiom~[4] (openness of versality)
from a purely geometric assumption on the object to be represented.
We also found the several (counter)examples that Artin presents
in~\cite[\S 5]{MR0262237} to be illuminating.
For these reasons, among others, we have chosen to discuss Artin's 
representability theorem in terms that are as close as possible to 
Artin's treatment in~\cite[\S 3]{MR0260746}
and~\cite[\S 5]{MR0262237}, making the minimal changes necessary to adapt the
statement of the axioms, and of the theorem, to the stacky situation.
Of course, such adaptations have been presented by many authors,
including Artin himself, but these works have tended to focus
on developments of the theory (such as the use of obstruction theories
to verify openness of versality) which are unnecessary for our purposes.
Ultimately, we found the treatment of Artin representability in the
Stacks Project~\cite[\href{http://stacks.math.columbia.edu/tag/07SZ}{Tag 07SZ}]{stacks-project}  
to be closest in spirit to the approach we wanted to take,
and it forms the basis for our treatment of the theorem here.
However, for the reason discussed above, of wanting to follow
as closely as possible Artin's original treatment, we have 
phrased the axioms in different terms to the way
they are phrased in the Stack Project, terms which are closer
to Artin's original phrasing. 

One technical reason for preferring Artin's phrasing is the emphasis
that it places  on the role of pro-representability (or equivalently, versality).
As already noted, the main intended application of Theorem~\ref{thm:main-intro}
is to the construction of moduli of Galois representations, 
and phrasing the axioms in a way which emphasises pro-representability
makes it easy to incorporate the formal deformation theory of Galois
representations into our arguments
(one of the main outputs of that formal deformation theory
being various pro-representability statements of the kind that
Theorem~\ref{thm:main-intro} requires as one of its inputs. In fact,
in the interests of generality we work with $\varphi$-modules that do
not evidently correspond to Galois representations, so we do not
directly invoke results from Galois deformation theory, but rather
adapt some arguments from that theory to our more general setting.)

On a few occasions it has seemed sensible to us to state and prove
a result in its natural level of generality, even if this
level of generality is not strictly required for 
the particular application we have in mind.     We have also developed
the analogue of Artin's axioms [4a] and [4b] 
of~\cite[\S 3]{MR0260746} (referred to as [4] and [5] 
in~\cite[\S 5]{MR0262237}) in the stacky context; while not
necessary for the proof of Theorem~\ref{thm:main-intro}, thinking in terms
of these axioms helps to clarify some of the foundational results
of Section~\ref{sec: artin's axioms and stacks in groupoids}
(e.g.\ the extent to which the unramifiedness condition
of Lemma~\ref{lem:monomorphism} can be promoted to the condition of being
a monomorphism, as in Corollary~\ref{cor:monomorphism}).

Just to inventory the contents of 
Section~\ref{sec: artin's axioms and stacks in groupoids}
a little more precisely: in Subsections~\ref{subsec:axiom 1}
through~\ref{subsec:axiom 4} we discuss each of Artin's axioms in turn.
In Subsection~\ref{subsec:[1]-ification} we develop a partial analogue
of~\cite[Prop.~3.11]{MR0260747}, which allows us to construct 
stacks satisfying~[1] by defining their values on 
algebras of finite presentation over the base and then taking
appropriate limits.
In Subsections~\ref{subsec: 1 and 3} and~\ref{subsec: 2 and 3}
we develop various further technical consequences of Artin's axioms.
Of particular importance is Lemma~\ref{lem:monomorphism},
which is a generalisation to the stacky context of
one of the steps appearing in the proof of~\cite[Lem.~3.10]{MR0260746}:
it provides unramified algebraic approximations to stacks satisfying Axioms~[1], [2],
and~[3], and so is the key to establishing openness of versality
(i.e.~Axiom~[4]) in certain contexts in which it is not assumed to hold
{\em a priori}. In Subsection~\ref{subsec: Artin representability}
we explain how our particular formulation of Artin's axioms does indeed
imply his representability theorem for algebraic stacks.

In Section~\ref{sec:images}, 
after a preliminary discussion in Subsection~\ref{subsec:images part one}
of the theory of scheme-theoretic images
in the context of morphisms of algebraic stacks,
in Subsection~\ref{subsec:scheme images two},
we present our main definitions, and give the proof of
Theorem~\ref{thm:main-intro}. In Subsection~\ref{subsec:relatively representable}, we investigate the behaviour of our
constructions with respect to base change (both of the target
stack~$\cF$, and of the base scheme~$S$); as well as being of
intrinsic interest, this will be important in our applications to
Galois representations in~\cite{CEGSKisinwithdd,emertongeepicture}.

In Section~\ref{sec: examples}, we give various examples of stacks
and Ind-stacks, which illustrate the results of
Section~\ref{sec:images}, and the roles of the various hypotheses of Section~\ref{sec: artin's axioms and
  stacks in groupoids} in the proofs of our main results. We also
prove some basic results about Ind-stacks which are used in the proof
of Theorem~\ref{thm: existence of Z^a, introduction version}.

The paper concludes with Section~\ref{sec: finite flat},
in which we define the moduli stacks of \'etale $\varphi$-modules 
and prove (via an application of Theorem~\ref{thm:main-intro})
that they are Ind-algebraic stacks.

\subsection{Acknowledgements}We would like to thank Alexander
Beilinson, Laurent Berger, Brian
Conrad, Johan de Jong, Vladimir Drinfeld, Ashwin Iyengar,  Wansu Kim, Mark Kisin, 
Jacob Lurie, Martin Olsson, Mike Roth,
Nick Rozenblyum and Evan Warner for
helpful conversations and correspondence.
We are particularly grateful to Mark Kisin, whose has generously
shared his ideas and suggestions about the moduli
spaces of Galois representations over the years. 
We are also indebted to Brian Conrad, Michael Harris, and Florian Herzig for their
close readings of, and many
helpful comments on and corrections to,
various versions of this manuscript.
The first author would like to thank the members of his working group
at the Stacks Project Workshop --- Eric Ahlqvist, Daniel Bragg, Atticus Christensen, Ariyan Javanpeykar, Julian Rosen --- for their helpful comments.
Finally, we thank the anonymous referee for their careful reading 
of the manuscript, and their many corrections and suggestions for
improvement.

\subsection{Notation and conventions}\label{subsec:notation}We follow the conventions
of~\cite{stacks-project} except where explicitly noted, and we refer to this
reference for background material. We note that the references
to~\cite{stacks-project} in the electronic version of this paper are clickable,
and will take the reader directly to the web page of the corresponding entry. We use Roman letters $X,Y,\dots$ for schemes
and algebraic spaces, and calligraphic letters $\cX,\cY,\dots$ for (possibly
non-algebraic) stacks (and more generally, for categories fibred in
groupoids). We elide the difference between commutative diagrams and
2-commutative diagrams, referring to both as simply commutative diagrams.

Since we follow~\cite{stacks-project}, we do not assume (unless otherwise
stated) that our algebraic spaces and algebraic stacks have quasi-compact or quasi-separated
diagonals. This is in contrast to references such
as~\cite{MR0302647,MR0399094,MR1771927}, and occasionally requires us to make
some simple additional arguments; the reader interested only in our applications
to moduli stacks of Galois representations should feel free to impose the
additional hypotheses on the diagonal that are common in the stacks literature,
and will lose nothing by doing so.

\subsubsection{Choice of site}
One minor difference between our approach and that taken
in~\cite{stacks-project} is that we prefer to only assume that the stacks that
we work with are stacks in groupoids for the \'etale topology, rather than the
\emph{fppf} topology. This ultimately makes no difference,
as the definition of an algebraic stack can be made using either the \'etale
or \emph{fppf} topologies
\cite[\href{http://stacks.math.columbia.edu/tag/076U}{Tag
  076U}]{stacks-project}. In practice, this means that we will sometimes cite
results from~\cite{stacks-project} that apply to stacks in groupoids for the
\emph{fppf} topology, but apply them to stacks in groupoids for the \'etale
topology. In each such case, the proof goes over unchanged to this setting.

To ease terminology, from now on we will
refer to a stack in groupoids for the \'etale topology (on some
given base scheme $S$) simply as a stack (or a stack over~$S$). 
(On a few occasions in the manuscript,
we will work with stacks in sites other than
the \'etale site, in which case we will be careful to signal this
explicitly.)

\subsubsection{Finite type points}
If $S$ is a scheme, and $s\in S$ is a point of $S$, we let $\kappa(s)$ denote
the residue field of $s$. A \emph{finite type point} $s\in S$ is a point such
that the morphism $\Spec \kappa(s)\to S$ is of finite type. By
\cite[\href{http://stacks.math.columbia.edu/tag/01TA}{Tag
  01TA}]{stacks-project}, a morphism $f:\Spec k\to S$ is of finite type if and
only if there is an affine open $U\subseteq S$ such that the image of $f$ is a
closed point $u\in U$, and $k/\kappa(u)$ is a finite extension. In a Jacobson
scheme, the finite type points are precisely the closed points; more
generally, the finite type points of any scheme $S$ are dense in every closed
subset of $S$ by~\cite[\href{http://stacks.math.columbia.edu/tag/02J4}{Tag
  02J4}]{stacks-project}.
If $X \to S$ is a finite type morphism, then a morphism $\Spec k \to X$
is of finite type if and only if the composite $\Spec k \to S$ is of finite type,
and so in particular a point $x \in X$ is of finite type
if and only if the composite $\Spec \kappa(x) \to X \to S$ is of finite type.

\subsubsection{Points of categories fibred in groupoids}
If $\cX$ is a category fibred in groupoids, then a \emph{point} of $\cX$ is an equivalence class
of morphisms from spectra of fields $\Spec K\to\cX$, where we say that $\Spec
K\to\cX$ and $\Spec L\to\cX$ are equivalent if there is a field $M$ and a
commutative diagram \[\xymatrix{\Spec M \ar[d]\ar[r]&\Spec L\ar[d]\\ \Spec
  K\ar[r]&\cX.}\](This is an equivalence relation
by~\cite[\href{http://stacks.math.columbia.edu/tag/04XF}{Tag
  04XF}]{stacks-project}; strictly speaking, this proves the claim in the case that $\cX$ is an
algebraic stack, but the proof goes over identically to the general case that
$\cX$ is a category fibred in groupoids.) If $\cX$ is furthermore an algebraic
stack, then the set of points of $\cX$ is denoted $|\cX|$;
by~\cite[\href{http://stacks.math.columbia.edu/tag/04XL}{Tag
  04XL}]{stacks-project} there is a natural topology on $|\cX|$, which has in
particular the property that if $\cX\to\cY$ is a morphism of algebraic stacks, then the
induced map $|\cX|\to|\cY|$ is continuous.

If $\cX$ is a category fibred in groupoids, then a \emph{finite type point} of $\cX$ is a point that can be represented by a
morphism $\Spec k\to\cX$ which is locally of finite type. If $\cX$ is an
algebraic stack, then by~\cite[\href{http://stacks.math.columbia.edu/tag/06FX}{Tag
  06FX}]{stacks-project}, a point $x\in |\cX|$ is of finite type if and only if
there is a scheme $U$, a smooth
morphism $\varphi:U\to\cX$
and a finite type point $u\in U$ such that $\varphi(u)=x.$ The
set of finite type points of an algebraic stack $\cX$ is dense in any closed subset of $|\cX|$ by~\cite[\href{http://stacks.math.columbia.edu/tag/06G2}{Tag 06G2}]{stacks-project}.

If $X$ is an {\em algebraic space} which is locally of finite type
over a locally Noetherian base scheme $S$,
then any finite type point of $X$ may be
represented by a {\em monomorphism} $\Spec k \to X$ which is locally
of finite type;
this representative is unique up to unique isomorphism,
and any other morphism $\Spec K \to X$ representing $x$ factors
through this one.  (See Lemma~\ref{lem:monomorphisms at points} below.)

\subsubsection{pro-categories}\label{subsubsec:procats}We will make
several uses of the formal pro-category $\pro\cC$ associated to a
category~$\cC$, in the sense of~\cite{MR1603480}. Recall that
an object of $\pro\cC$ is a projective system   $(\xi_i)_{i\in I}$ of
objects of~$\cC$, and  the morphisms between two pro-objects $\xi=(\xi_i)_{i\in I}$ and
$\nu=(\eta_j)_{j\in J}$ are by
definition \[\Mor(\xi,\eta)=\varprojlim_{j\in J}\varinjlim_{i\in I}\Mor(\xi_i,\eta_j).\]

We will apply this definition in particular to categories of Artinian local rings
with fixed residue fields in Section~\ref{subsec:axiom 2}, and to the
category of affine schemes locally of finite presentation over a fixed
base scheme in Section~\ref{subsec:[1]-ification}, as well as to
categories (co)fibred in groupoids over these categories.
\subsubsection{$G$-rings}\label{subsubsec:Grings}
Recall (\cite[\href{http://stacks.math.columbia.edu/tag/07GH}{Tag 07GH}]{stacks-project}) that a Noetherian ring $R$ is a $G$-ring if for every prime $\p$ of $R$,
the (flat) map $R_\p\to(R_\p)^\wedge$ is regular.
By~\cite[\href{http://stacks.math.columbia.edu/tag/07PN}{Tag
  07PN}]{stacks-project}, this is equivalent to demanding that for every pair of primes $\mathfrak q \subseteq \mathfrak p \subset R$
the algebra
$(R/\mathfrak q)_\mathfrak p^\wedge \otimes_{R/\mathfrak q} \kappa(\mathfrak q)
$
is geometrically regular over $\kappa(\mathfrak q)$ (where $\kappa(\mf q)$ denotes the residue field of
$\gq$; recall~\cite[\href{http://stacks.math.columbia.edu/tag/0382}{Tag 0382}]{stacks-project} that if $k$ is a field, a Noetherian $k$-algebra~$A$ is
geometrically regular if and only if $A\otimes_kk'$ is regular for
every finitely generated field extension $k'/k$). Excellent rings are $G$-rings by definition.

In our main results we will assume that our base
scheme $S$ is locally Noetherian, and that its local rings $\cO_{S,s}$ at all finite type points
$s\in S$ are $G$-rings. This is a replacement
of Artin's assumption that $S$ be of finite type over a
field or an excellent DVR; this more general setting is permitted
by improvements in Artin approximation, due essentially to
Popescu (\cite{MR818160,MR868439,MR1060701}; see also~\cite{MR1935511}
and~\cite[\href{http://stacks.math.columbia.edu/tag/07GB}{Tag
  07GB}]{stacks-project}). However, since this assumption will not always be in force, we will
indicate when it is assumed to hold.

\subsubsection{Groupoids}
We will make use of groupoids in algebraic spaces, and we will use the notation
for them which is introduced
in~\cite[\href{http://stacks.math.columbia.edu/tag/043V}{Tag
  043V}]{stacks-project}, which we now recall. A groupoid in algebraic spaces
over a base algebraic space $B$ is
a tuple $(U,R, s,t,c)$ where $U$ and $R$ are algebraic spaces over $B$, and
$s,t:R\to U$ and $c:R\times_{s,U,t}R\to R$ are morphisms of algebraic spaces
over $B$ whose $T$-points form a groupoid for any scheme $T\to B$. (The maps
$s,t,c$ give the source, target and composition laws for the arrows of the
groupoid.) Given such a groupoid in algebraic spaces, there are unique morphisms
 $e:U\to R$ and $i:R\to R$ of algebraic spaces over $B$ which give the identity
 and inverse maps of the groupoid, and we sometimes denote the groupoid in
 algebraic spaces by the tuple $(U,R, s,t,c,e,i)$.

\subsubsection{Properties of morphisms}
In most cases, we follow the terminology and conventions for 
properties 
of morphisms of stacks introduced in~\cite{stacks-project}.  We recall
some of the general framework of those conventions here.

An important concept, defined for morphisms of categories fibred in
groupoids, and so in particular for morphisms of stacks,
is that of being {\em representable by algebraic spaces}.
Following~\cite[\href{http://stacks.math.columbia.edu/tag/04SX}{Tag
04SX}]{stacks-project}, we say that a morphism $\cX \to \cY$
of categories fibred in groupoids is representable by algebraic
spaces if for any morphism $T \to \cY$ with $T$ a scheme,
the fibre product $\cX\times_{\cY} T$ is (representable by)
an algebraic space.  (This condition then continues to hold
whenever $T$ is an algebraic
space~\cite[\href{http://stacks.math.columbia.edu/tag/0300}{Tag
0300}]{stacks-project}.)
A morphism of algebraic stacks is representable by algebraic
spaces if and only if the associated diagonal morphism is a
monomorphism~\cite[\href{http://stacks.math.columbia.edu/tag/0AHJ}{Tag 0AHJ}]
{stacks-project}. 

If $P$ is a property of morphisms of algebraic spaces which
is preserved under arbitrary base-change, and which is {\em fppf}
local on the target,
then~\cite[\href{http://stacks.math.columbia.edu/tag/03YK}{Tag
03YK}]{stacks-project}
provides a general mechanism for defining the property $P$
for morphisms of categories fibred in groupoids that are representable
by algebraic spaces: namely,
such a morphism $f:\cX \to \cY$ is defined to have property $P$ if and only
if for any morphism $T \to \cY$ with $T$ a scheme,
the base-changed morphism $\cX\times_{\cY} T \to T$ (which is a morphism
of algebraic spaces, by assumption) has property $P$
(and it is equivalent to impose the same condition with
$T$ being merely an algebraic space, because an algebraic space by definition
has an \'etale (and therefore {\em fppf}) cover by a scheme, and $P$ is {\em
  fppf} local on the target by assumption).

If $P$ is a property of morphisms of algebraic spaces
which is smooth local on the source-and-target,
then~\cite[\href{http://stacks.math.columbia.edu/tag/06FN}{Tag 06FN}]
{stacks-project} extends the definition of $P$ to arbitrary morphisms
of algebraic stacks (in particular, to morphisms that are not necessarily
representable by algebraic spaces): a morphism $f: \cX \to \cY$
is defined to have property $P$ if it can be lifted to a morphism
$h: U \to V$ having property $P$, where $U$ is a smooth cover of $\cX$
and $V$ is a smooth cover of $\cY$.  If $P$ is furthermore preserved under arbitrary base-change and {\em fppf} local on the target,
so that the definition
of~\cite[\href{http://stacks.math.columbia.edu/tag/03YK}{Tag
03YK}]{stacks-project} applies,
then the two definitions coincide in the case
of morphisms that are representable by algebraic
spaces~\cite[\href{http://stacks.math.columbia.edu/tag/06FM}{Tag
06FM}]{stacks-project}.

Many additional properties of morphisms of algebraic
stacks are defined 
in~\cite[\href{http://stacks.math.columbia.edu/tag/04XM}{Tag 04XM}]
{stacks-project}.  In Subsection~\ref{subsec:axiom 3} below,
we further extend many of these definitions to the case
of morphisms of stacks whose source is assumed to be algebraic,
but whose target is assumed only to satisfy condition~[3] of Artin's
axioms.

\subsection{Some local algebra}
\label{subsec:local}
In this subsection, we state and prove some results from local 
algebra which we will need in what follows.  

\begin{lemma}\label{lem:fibre products of local rings}
If $B \to A$ and $C\to A$
are local morphisms from a pair of complete Noetherian local rings
to an Artinian local ring, $C \to A$ is surjective, and the residue field of~$A$ is finite over
the residue field of~$B$, then the fibre product $B\times_A C$ is a complete Noetherian local ring,
and the natural morphism $B\times_A C\to A$ is local.
\end{lemma}
\begin{proof}
Write $R := B\times_A C$, and 
$\mathfrak m_R := \mathfrak m_B \times_{\mathfrak m_A} \mathfrak m_C$.
Since $B\to A$ and $C \to A$ are local morphisms of local rings,
we see that if $(b,c) \in R$, 
with both $b$ and $c$ lying over the element $a \in A$,
then if $a \in \mathfrak m_A$, we have $(b,c) \in \mathfrak m_R$,
while if $a \not\in \mathfrak m_A$, then $(b,c) \not\in \mathfrak m_R$.
In the latter case, we find that $(b,c)$ is furthermore
a unit in $R$.  Thus
$R$ is a local ring with maximal ideal $\mathfrak m_R$,
and the natural morphism $R \to A$ is local.

If we choose $r\geq 0$ so that $\mathfrak m_A^r = 0$ (which is
possible, since $A$ is Artinian), then each of $\mathfrak m_B^r$
and $\mathfrak m_C^r$ has vanishing image in $A$,
and so we see that 
\numequation
\label{eqn:maximal ideal comparison}
\mathfrak m_R^r \subseteq \mathfrak m_B^r \times \mathfrak m_C^r
\subseteq \mathfrak m_R.
\end{equation}
From this inclusion, and the fact that $B\times C$ is $\mathfrak m_B\times
\mathfrak m_C$-adically complete, it follows that $R$ is
$\mathfrak m_R$-adically complete.  (Indeed, we see that
$R$ is open and closed as a topological subgroup of $B\times C$,
and that the induced topology on $R$ coincides with the $\mathfrak m_R$-adic
topology.)

Finally, to see that~$R$ is Noetherian, we use the hypothesis that
$C\to A$ is surjective, which implies that the
residue fields~$k_C$ and~$k_A$ of~$C$ and~$A$ are isomorphic, as are the residue
fields~$k_R$ and~$k_B$ of~$R$ and~$B$, which are subfields of
$k_C=k_A$.  

Then the inclusion $\m_R=\m_B\times_{\m_A}\m_C\into\m_B\times\m_C$
induces an
inclusion \[\m_R/\m_R^2\into\m_B/\m_B^2\times\m_C/\m_C^2,\]and
since~$B$ and~$C$ are Noetherian, and the extension degree $[k_R:k_B]=[k_C:k_A]$ is
finite, the target of the inclusion is a finite-dimensional
$k_R$-vector space. 
It follows
that~$\m_R/\m_R^2$ is also finite-dimensional, and therefore that~$R$
is Noetherian, as required.\end{proof}

\begin{lemma}\label{lem: Cohen rings lemma}
Let $B \to A$ be a local morphism from a complete Noetherian local ring
to an Artinian local ring, which induces a finite extension of residue
fields.
Then this morphism admits a factorisation $B \to B' \to A,$
where $B \to B'$ is a faithfully flat local morphism
of complete local Noetherian rings, and $B' \to A$ is surjective
{\em (}and so in particular induces an isomorphism on residue fields{\em )}.
\end{lemma}
\begin{proof}
Let $k_B \subseteq k_A$ be the embedding of residue fields induced by
the given morphism $B \to A$.   Let $\Lambda_{k_A}$ denote
a Cohen ring with residue field $k_A$, and choose (as we may)
a surjection $\Lambda_{k_A}[[x_1,\ldots,x_d]] \onto A$
(for some appropriate value of~$d$).
Let $\Bbar$ denote the image of $B$ in $A$,
and let $A'$ denote the fibre product $A' := \Bbar \times_A \Lambda_{k_A}[[x_1,
\ldots, x_d]];$
then Lemma~\ref{lem:fibre products of local rings} 
shows that $A'$ is a complete Noetherian local subring of $\Lambda_{k_A}[[x_1,
\ldots, x_d]]$ with residue field $k_B$.  If $\Lambda_{k_B}$ 
denotes a Cohen ring with residue field $k_B$, then we may find
a local morphism $\Lambda_{k_B} \to A'$ inducing the identity on
residue fields.  The composite
\numequation
\label{eqn:cohen ring map}
\Lambda_{k_B} \to A' \subseteq \Lambda_{k_A} [[x_1,\ldots,x_d]]
\end{equation}
 is flat.

By~\cite[Thm.\ 0.19.8.6(i)]{MR0173675}, the composite morphism $\Lambda_{k_B} \to A' \to \Bbar$
(the second arrow being the projection)
may be lifted to a morphism $\Lambda_{k_B} \to B$. Now define $B' :=
B \cotimes_{\Lambda_{k_B}} \Lambda_{k_A}[[x_1,\ldots,x_d]]$
(the completed tensor product). By~\cite[Lem.\ 0.19.7.1.2]{MR0173675}, $B'$ is a
complete local Noetherian ring, flat over $B$.

The given morphisms $B \to A$ and $\Lambda_{k_A}[[x_1,\ldots,x_d]]
\to A$ induce a surjection $B' \to A$,
and $B \to B' \to A$ is the required factorisation of our given
morphism $B \to A$.
(Note that flat local morphisms of 
local rings are automatically faithfully flat.)
\end{proof}

\section{Stacks in groupoids and Artin's axioms}
\label{sec: artin's axioms and stacks in groupoids}
Since Artin first introduced his axioms characterising algebraic spaces~\cite{MR0260746}, many versions of these axioms have appeared in the works of various authors.
In this paper we have tried to follow Artin's original treatment closely, and the labelling of our four axioms is chosen to match
the labelling in~\cite{MR0260746}. 

In this section we will discuss each of the four axioms,
explain why they imply representability (essentially,
by relating them to the axioms given
in~\cite[\href{http://stacks.math.columbia.edu/tag/07SZ}{Tag 07SZ}]{stacks-project})
and 
also discuss some related foundational material.

As noted in the introduction, our basic setting will
be that of stacks in groupoids on the big \'etale site
of a scheme $S$. A general reference for the basic definitions and properties of
such stacks is~\cite{stacks-project}.
As remarked in the introduction, from now on we will
refer to such a stack in groupoids simply as a stack.
At times we will furthermore assume that $S$ is locally Noetherian,
and in Subsection~\ref{subsec: Artin representability},
where we present Artin's representability theorem,
we will additionally assume
that the local rings $\mathcal O_{S,s}$ are $G$-rings,
for each finite type point $s \in S$.

\subsection{Remarks on Axiom [1]}
\label{subsec:axiom 1}
We begin by recalling the definition of limit preserving.

\begin{df}
\label{def:limit preserving}
A category fibred in groupoids (e.g.\ an algebraic stack)
$\cF$ over $S$ is {\em limit preserving} if,
whenever we have a projective limit $T = \varprojlim T_i$
of affine $S$-schemes, the induced functor 
\numequation
\label{eqn:limit functor}
\twocolim\cF(T_i)\to\cF(T)
\end{equation}
is an equivalence of categories. 

More concretely, as in~\cite[\href{http://stacks.math.columbia.edu/tag/07XK}{Tag 07XK}]{stacks-project} this means that each object of
$\cF(T)$ is isomorphic to the restriction to $T$ of an object of $\cF(T_i)$ for
some $i$; that for any two objects $x,y$ of $\cF(T_i)$, any morphism between the
restrictions of $x,y$ to $T$ is the restriction of a morphism between the
restrictions of $x,y$ to $T_{i'}$ for some $i'\ge i$;
and that for any two objects $x,y$ of $\cF(T_i)$, if two morphisms $x \rightrightarrows y$
coincide after restricting to $T$, than they coincide after restricting to 
$T_{i'}$ for some $i' \geq i$.  (Since we are considering categories
fibred in groupoids, it suffices to check this last condition when one
of the morphisms is the identity.)

\end{df}

We have the following related definition~\cite[\href{http://stacks.math.columbia.edu/tag/06CT}
{Tag 06CT}]{stacks-project}.

\begin{df}
\label{def:limit preserving on objects}
A morphism $\cF \to \cG$ of categories fibred in groupoids
(e.g.\ of algebraic stacks) over $S$ is said to be {\em limit 
preserving on objects} if for any affine $S$-scheme $T$,
written as a projective limit of affine $S$-schemes $T_i$,
and any morphism $T \to \cF$ for which the composite morphism
$T \to \cF \to \cG$ factors through $T_i$ for some $i$,
there is a compatible factorisation of the morphism $T \to \cF$ 
through $T_{i'}$, for some $i' \geq i$.

Somewhat more precisely, given a commutative diagram
$$\xymatrix{T \ar[r] \ar[d]  & \cF\ar[d]
\\ T_i \ar[r] & \cG}
$$
we may factor it in the following manner:
$$\xymatrix{T \ar[rd] \ar[d]  & \\
T_{i'} \ar[r]\ar[d] & \cF\ar[d] \\
T_i \ar[r] & \cG}
$$
\end{df}

\medskip

We also make the following variation on the preceding definition.

\begin{df}
We say that 
a morphism $\cF \to \cG$ of categories fibred in groupoids
(e.g.\ of algebraic stacks) over $S$ is {\em \'etale locally limit 
preserving on objects} if for any affine $S$-scheme $T$,
written as a projective limit of affine $S$-schemes $T_i$,
and any morphism $T \to \cF$ for which the composite morphism
$T \to \cF \to \cG$ factors through $T_i$ for some $i$,
then there is an affine \'etale surjection $T'_{i'} \to T_{i'}$,
for some $i' \geq i$, and a morphism $T'_{i'} \to \cF$,
such that, if we write $T' := T_i'\times_{T_i} T,$
then the resulting  diagram 
$$\xymatrix{T' \ar[rd]\ar[d]  &  \\ T_{i'}' \ar[d]\ar[r]  & \cF \ar[d] \\ 
T_i  \ar[r] & \cG }
$$
commutes. \end{df}

The following lemma relates Definitions~\ref{def:limit preserving} and~\ref{def:limit
preserving on objects}.

\begin{lem}
\label{lem:limit preserving vs. limit preserving on objects}
If $\cF$ is a category fibred in groupoids over $S$,
then the following are equivalent:
\begin{enumerate}
\item $\cF$ is limit preserving.
\item Each of the morphisms $\cF \to S$, $\Delta:\cF \to \cF\times_S \cF$,
and $\Delta_{\Delta}: \cF \to \cF\times_{\cF\times_S \cF} \cF$ is limit preserving on objects.
\end{enumerate}
\end{lem}
\begin{proof}
	This is just a matter of working through the definitions.
Indeed, the morphism $\cF \to S$ being limit preserving on objects
is equivalent to the functor~(\ref{eqn:limit functor}) being essentially
surjective (for all choices of $\displaystyle T = \varprojlim_{i}T_i$),
the diagonal morphism being limit preserving on objects
is equivalent to this functor being full, and the double diagonal being
limit preserving on objects is equivalent to this functor being
faithful.

More precisely, by definition the morphism $\cF\to S$ is limit
preserving on objects if and only if for every
$\displaystyle T = \varprojlim_{i}T_i$ as above, any morphism
$T\to\cF$ factors through some~$T_i$; equivalently, if and only if
every object of~$\cF(T)$ is isomorphic to the restriction to~$T$ of an
object of~$\cF(T_i)$ for some~$i$; equivalently,
if and only if
the functor~(\ref{eqn:limit functor}) is
essentially surjective.
Similarly,
the morphism  $\Delta:\cF \to \cF\times_S
\cF$ is limit preserving on objects if and only if
for any pair of objects of $\cF(T_i)$ (for some value of $i$),
a morphism between their images in $\cF(T)$ arises as the restriction
of a morphism between their images in $\cF(T_{i'})$, for some $i' \geq i$;
or equivalently,
if and only if 
the functor~(\ref{eqn:limit functor}) is full.
The claim about the double diagonal is
similar, and is left to the interested reader.
\end{proof}

We can strengthen the preceding lemma
when the category fibred in groupoids involved
is actually a stack.

\begin{lemma}
\label{lem:limit preserving vs. locally limit preserving on objects}
If $\cF$ is a stack over $S$,
then the following are equivalent:
\begin{enumerate}
\item $\cF$ is limit preserving.
\item The morphism $\cF \to S$ is \'etale locally limit preserving
on objects, while each of $\Delta:\cF \to \cF\times_S \cF$ 
and $\Delta_{\Delta}: \cF \to \cF\times_{\cF\times_S \cF} \cF$ is limit preserving on objects.
\end{enumerate}
\end{lemma}
\begin{proof}
Taking into account Lemma~\ref{lem:limit preserving vs. limit preserving
on objects}, we see that the lemma will follow if we show that the 
assumptions of~(2) imply that $\cF \to S$ is limit preserving on objects.
Thus we put ourselves in the situation described in
Definition~\ref{def:limit preserving on objects} (taking $\cG = S$),
namely
we give ourselves an affine $S$-scheme $T$, written as a projective limit
$T = \varprojlim_i T_i$ of affine $S$-schemes, and we suppose given
a morphism $T \to \cF$; we must show that this morphism factors through
$T_i$ for some $i$.  Applying the assumption that $\cF \to S$
is \'etale locally limit preserving on objects, we find that,
for some $i$, we may find an \'etale cover $T_i'$ of $T_i$ 
and a morphism $T'_i \to \cF$, for some value of $i$,
through which the composite $T' := T_i'\times_{T_i} T \to T \to \cF$
factors;
our goal is then to show that, for some $i' \geq i,$ we 
may find a morphism $T_{i'} \to \cF$ through which the morphism
$T \to \cF$ itself factors.

For any $i' \geq i,$ we let $T_{i'}' := T_i'\times_{T_i} T_{i'}.$
In order to find the desired morphism $T_{i'} \to \cF$,
it suffices to equip the composite $T'_{i'} \to T'_i \to \cF$
with descent data to $T_{i'}$, in a manner compatible with 
the canonical descent data of the composite $T' \to T \to \cF$
to $T$.  That this is possible follows easily from the assumptions
on the diagonal and double diagonal of $\cF$ (\emph{cf.} the proof of
Lemma~\ref{lem:restriction to finitely presented algebras}~(2) below).
\end{proof}

We will have use for the following finiteness results.

\begin{lem}
\label{lem:morphism between limit preserving}
If $\cF \to \cG \to \cH$ are morphisms between categories fibred
in groupoids over $S$, 
and if both the composite morphism $\cF \to \cH$ and
the diagonal morphism $\Delta :\cG \to \cG \times_{\cH} \cG$
are limit preserving on objects,
then the morphism $\cF \to \cG$ is also limit preserving on objects.
\end{lem}
\begin{proof}
Let $T=\varprojlim T_i$ be a projective limit of affine $S$-schemes, and
suppose that we are given a morphism $T\to\cF$ such that the composite
$T\to\cF\to\cG$ factors through $T_i$ for some $i$. We must show that there is
a compatible factorisation of
$T\to\cF$  through $T_{i'}$ for some $i'\ge i$. 

Since the composite $\cF\to \cH$ is limit preserving on objects,
we may factor $T\to\cF$ through some
$T_j$, in such a way that the composites $T \to T_j \to \cF \to \cH$ and
$T \to T_i \to \cG \to \cH$ coincide. Replacing $i,j$ by some common $i''\ge i,j$, we have two morphisms
$T_{i''}\to\cG$ (one coming from the given morphism $T_i\to\cG$, and one from
the composite $T_j\to\cF\to\cG$) which induce the same morphism to $\cH$, and which 
agree when pulled-back to $T$. Since $\Delta$ is limit
preserving on objects, they agree over some $T_{i'}$, for some $i'\ge i''$, as required.
\end{proof}

\begin{cor}
  \label{cor: morphism between limit preserving}If $\cF$ and $\cG$ are categories
  fibred in groupoids over $S$, both of which are limit preserving, then any morphism $\cF\to\cG$ 
  is limit preserving on objects.
\end{cor}
\begin{proof}This follows directly from Lemmas~\ref{lem:limit preserving vs. limit preserving on objects}
and~\ref{lem:morphism between limit preserving} (taking $\cH=S$ in the latter).
\end{proof}

The next
lemma (which is essentially~\cite[Prop.\ 4.15(i)]{MR1771927}) explains why
the condition of being limit preserving is sometimes referred to as being locally
of finite presentation. 

\begin{lem}\label{lem: limit preserving iff locally of finite presentation}If $\cF$ is
  an algebraic stack over $S$, then the following are equivalent:
  \begin{enumerate}
  \item $\cF$ is limit preserving.
  \item $\cF \to S$ is limit preserving on objects.\item $\cF$ is locally of finite presentation over $S$.
  \end{enumerate}

\end{lem}
\begin{proof}
  Note that (1)$\implies$(2) by definition (since (2) is just the condition that
  the functor~(\ref{eqn:limit functor}) be
essentially surjective),
while the equivalence of~(2) and~(3) is a special case of Lemma~\ref{lem:lfp equals limit preserving
on objects} below. (The reader may easily check that the present lemma 
is not used in the proof of that result, and so there is no circularity in
appealing to it.)

It remains to show that (3)$\implies$(1). 
Assuming that~(3) holds,
we claim that the diagonal
$\Delta:\cF\to\cF\times_S\cF$ is locally of finite presentation. To see this,
choose a smooth surjection $U \to \cF$ whose source is a scheme (which exists because $\cF$ is assumed algebraic). Since
$U\to\cF$ and $\cF\to S$ are locally of finite presentation by assumption, we
see that $U\to S$ is locally of finite presentation
by~\cite[\href{http://stacks.math.columbia.edu/tag/06Q3}{Tag
  06Q3}]{stacks-project}. Then the diagonal $U\to U\times_S U$ is locally of
finite presentation by~\cite[\href{http://stacks.math.columbia.edu/tag/0464}{Tag
  0464}]{stacks-project}
and~\cite[\href{http://stacks.math.columbia.edu/tag/084P}{Tag
  084P}]{stacks-project}, and therefore the composite \[U\to U\times_S U\to
\cF\times_S U\to \cF\times_S \cF\] is locally of finite presentation (note that
the last two morphisms are base changes of the smooth morphism
$U\to\cF$). Factoring this morphism
as \[U\to\cF\stackrel{\Delta}{\to}\cF\times_S\cF,\] (where $\Delta$ is
representable by algebraic spaces, because~$\cF$ is algebraic), we see
from~\cite[\href{http://stacks.math.columbia.edu/tag/06Q9}{Tag
  06Q9}]{stacks-project} that $\Delta$
is locally of finite
presentation, as claimed. A similar argument shows that the double diagonal $\Delta_{\Delta}:\cF\to \cF \times_{\cF\times_S\cF} \cF$
is locally of finite presentation.   
Applying Lemma~\ref{lem:lfp equals limit preserving on objects} below
(or~\cite[\href{http://stacks.math.columbia.edu/tag/06CX}{Tag 06CX}]{stacks-project}),
we find that each of $\Delta$, $\Delta_{\Delta}$ and the structure map $\cF\to S$
is limit preserving on objects. 
It now follows from Lemma~\ref{lem:limit preserving vs. limit preserving on
  objects}
that $\cF$ is limit preserving, as required.
\end{proof}

\subsection{Remarks on Axiom [2]}
\label{subsec:axiom 2}
Throughout this subsection, we assume that $S$ is locally Noetherian, and that
$\cF$ is a category fibred in groupoids over $S$, which is a stack for the
Zariski topology. We denote by $\widehat{\cF}$ the restriction of $\cF$ to the
category of finite type Artinian local $S$-schemes.

We begin by discussing Axiom~[2](a), which is the Rim--Schlessinger condition (RS).
Consider pushout diagrams
$$\xymatrix{Y \ar[d]\ar[r] & Y' \ar[d] \\ Z \ar[r] & Z' }$$
of $S$-schemes, 
with the horizontal arrows being closed immersions. We have an induced
functor \numequation\label{eqn: RS functor}\cF(Z') \to \cF(Y') \times_{\cF(Y)}
\cF(Z).\end{equation} 

\begin{defn}
  We say that $\cF$ satisfies condition (RS) if the functor~(\ref{eqn: RS
    functor}) is an equivalence of categories whenever $Y,Y',Z,Z'$ are
  finite type local Artinian $S$-schemes.
\end{defn}
\begin{lem}
  \label{lem: algebraic stacks satisfy RS}If $\cF$ is an algebraic stack, then
  $\cF$ satisfies \emph{(RS)}.
\end{lem}
\begin{proof}
  This is immediate from~\cite[\href{http://stacks.math.columbia.edu/tag/07WN}{Tag 07WN}]{stacks-project}.
\end{proof}
The same condition appears under a different name, and with slightly
different phrasing, in~\cite[Lem.\ 1.2]{hall2013artin}. We recall
this, and some
closely related notions that we will occasionally use.

\begin{defn}
\label{df:homogeneity}
Following~\cite{hall2013artin},
we say that $\cF$ is~$\textbf{Art}^{\textbf{fin}}$-homogeneous (resp.\
$\textbf{Art}^{\textbf{triv}}$-homogeneous,
resp.\ $\textbf{Art}^{\textbf{sep}}$-homogeneous,
resp.\ $\textbf{Art}^{\textbf{insep}}$-homogeneous)
if the functor~(\ref{eqn: RS
    functor}) is an equivalence of categories whenever $Y,Z$ are local Artinian
  $S$-schemes of finite type over $S$ (resp.\ with the induced extension of
  residue fields being trivial, resp.\ separable,
resp.\ purely inseparable), and $Y\to Y'$ is a nilpotent
  closed immersion.
\end{defn}

\begin{lem}
  \label{lem: RS is Art-fin homogeneous}$\cF$ satisfies {\em (RS)} if and only if it
  is~\emph{$\textbf{Art}^{\textbf{fin}}$}-homogeneous.
\end{lem}
\begin{proof}
This is just a matter of comparing the two definitions.  Precisely:
  a closed immersion of local Artinian schemes is automatically
  nilpotent. Conversely, a finite type nilpotent thickening of a local Artinian scheme is
  local Artinian, and the pushout of local Artinian schemes of finite type over $S$ is
  also local Artinian of finite type over $S$. (Recall from
  \cite[\href{http://stacks.math.columbia.edu/tag/07RT}{Tag
    07RT}]{stacks-project} that if above we write $Y=\Spec A$, $Y'=\Spec A'$,
  $Z=\Spec B$, then the pushout $Z'=\Spec B'$ is just given by $B'=B\times_A A'$.)\end{proof}

The following lemma relates the various conditions of
Definition~\ref{df:homogeneity}.
\begin{lem}
\label{lem:homogeneity conditions}
The condition of~$\mathbf{Art}^{\mathbf{fin}}$-homogeneity
of Definition~{\em \ref{df:homogeneity}}
implies
each of~$\mathbf{Art}^{\mathbf{sep}}$-homogeneity
and~$\mathbf{Art}^{\mathbf{insep}}$-homogeneity, and these conditions in turn imply
$\mathbf{Art}^{\mathbf{triv}}$-homogeneity. 
If $\cF$ is a stack
for the \'etale site, then conversely,
$\mathbf{Art}^{\mathbf{triv}}$-homogeneity
implies
$\mathbf{Art}^{\mathbf{sep}}$-homogeneity, 
while
$\mathbf{Art}^{\mathbf{insep}}$-homogeneity
implies
$\mathbf{Art}^{\mathbf{fin}}$-homogeneity. 
If $\cF$ is furthermore a stack for the fppf site, then
$\mathbf{Art}^{\mathbf{triv}}$-homogeneity 
implies
$\mathbf{Art}^{\mathbf{fin}}$-homogeneity. 
\end{lem} 
\begin{proof}
This follows immediately from~\cite[Lem.\ 1.6]{hall2013artin}
and its proof (see also the proof of~\cite[Lem.\ 2.6]{hall2013artin}).
\end{proof}

We now discuss Axiom~[2](b). To begin, we recall some definitions
from~\cite[\href{http://stacks.math.columbia.edu/tag/06G7}{Tag
  06G7}]{stacks-project}.

Fix  a Noetherian ring $\Lambda$, and a finite ring map $\Lambda\to k$,
whose target is a field. Let the kernel of this map be
$\m_\Lambda$ (a maximal ideal of $\Lambda$). We let
$\cC_\Lambda$ be the category whose objects are pairs $(A,\phi)$ consisting of an
Artinian local $\Lambda$-algebra $A$ and a $\Lambda$-algebra
isomorphism $\phi:A/\m_A\to k$, and whose morphisms are given by local $\Lambda$-algebra homomorphisms
compatible with $\phi$. Note that any such $A$ is finite over $\Lambda$, and
that the morphism $\Lambda\to A$ factors through $\Lambda_{\m_\Lambda}$, so that
we have $\cC_\Lambda=\cC_{(\Lambda_{\m_\Lambda})}$ in an evident sense.

There are some additional categories, closely related to $\cC_{\Lambda}$,
that we will also consider.  We let $\widehat{\cC}_{\Lambda}$ denote
the category of complete Noetherian local $\Lambda$-algebras $A$
equipped with a $\Lambda$-algebra isomorphism $A/\mathfrak m_A \iso k$,
while we let $\pro\cC_{\Lambda}$ denote the category of formal pro-objects
from $\cC_{\Lambda}$ in the sense of Section~\ref{subsubsec:procats}.  If $(A_i)_{i\in I}$ is an object of $\pro\cC_{\Lambda}$,
then we form the actual projective limit $A:= \varprojlim_{i \in I}
A_i,$ thought of as a topological ring (endowed with the projective
limit topology, each $A_i$ being endowed with its discrete topology).
In this manner we obtain an equivalence of categories between
$\pro\cC_{\Lambda}$ and the category of topological pro-(discrete
Artinian) local $\Lambda$-algebras equipped with a $\Lambda$-algebra isomorphism between
their residue fields and $k$ \cite[\S A.5]{MR1603480}.  We will frequently identify an
object of $\pro\cC_{\Lambda}$ with the associated topological local
$\Lambda$-algebra $A$.
There is a fully faithful embedding of $\widehat{\cC}_{\Lambda}$
into $\pro\cC_{\Lambda}$ given by associating to any object $A$
of the former category the pro-object $(A/\mathfrak m_A^i)_{i \geq 1}.$ 
In terms of topological $\Lambda$-algebras, this amounts to regarding
$A$ as a topological $\Lambda$-algebra by equipping it with its
$\mathfrak m_A$-adic topology.

\begin{remark}
\label{rem:pseudo-compact}
We note that objects of $\pro\cC_{\Lambda},$ when regarded as topological rings,
are examples of {\em pseudo-compact} rings, in the sense of \cite{MR0232821}.  
In particular, any morphism of such rings has closed image, and induces a topological
quotient map from its source onto its image; consequently,
a homomorphism $A \to B$ of such rings is surjective if and only if it is induced
by a compatible collection of surjective morphisms $A_i \to B_i$ for projective
systems $(A_i)_{i \in I}$ and $(B_i)_{i \in I}$ of objects in $\cC_{\Lambda}$.
(See the discussion beginning on \cite[p.~390]{MR0232821}, especially the statement
and proof of Lem.~1 and of~Thm.~3.)
\end{remark}

We have the usual notion of a category cofibred in groupoids over
$\cC_\Lambda$, for which
see~\cite[\href{http://stacks.math.columbia.edu/tag/06GJ}{Tag
  06GJ}]{stacks-project}. The particular choices of $\cC_\Lambda$ and categories cofibred in groupoids
over $\cC_\Lambda$ that we are interested in arise as follows
(see~\cite[\href{http://stacks.math.columbia.edu/tag/07T2}{Tag
  07T2}]{stacks-project} for more details). Let $\cF$ be a
category fibred in groupoids over $S$, let $k$ be a field and let $\Spec k\to
S$ be a morphism of finite type.

Let $x$ be an object of $\cF$ lying over
$\Spec k$, and let $\Spec\Lambda\subseteq S$ be an affine open so that $\Spec k\to
S$ factors as $\Spec k\to\Spec\Lambda\to S$, where $\Lambda\to k$ is
finite. (For the existence of such a $\Lambda$, and the independence of
$\cC_\Lambda$ of the choice of $\Lambda$ up to canonical equivalence,
see~\cite[\href{http://stacks.math.columbia.edu/tag/07T2}{Tag
  07T2}]{stacks-project}.) Write $p:\cF\to S$ for the tautological morphism.
We then let $\widehat{\cF}_{x}$ be the category whose:
\begin{enumerate}
\item objects are morphisms $x\to x'$ in $\cF$ such that $p(x')=\Spec A$, with
  $A$ an Artinian local ring, and the morphism $\Spec k\to\Spec A$
  given by $p(x)\to p(x')$ corresponds to a ring homomorphism $A\to k$
  identifying $k$ with the residue field of $A$, and
\item morphisms $(x\to x')\to(x\to x'')$ are commutative diagrams in $\cF$ of
  the form \[\xymatrix{x'&& x''\ar[ll]\\& x\ar[ul]\ar[ur]&}\]
\end{enumerate}
Note that the ring $A$ in~(1) is an object of $\cC_\Lambda$. Under the assumption that $\cF$ satisfies (RS), $\widehat{\cF}_{x}$ is a
deformation category
by~\cite[\href{http://stacks.math.columbia.edu/tag/07WU}{Tag
  07WU}]{stacks-project}. By definition, this means that
$\widehat{\cF}_{x}(\Spec k)$ is equivalent to a category with a single object and morphism,
and that $\widehat{\cF}_{x}$ is cofibred in groupoids and satisfies a natural analogue
of (RS) (more precisely, an analogue of $\mathbf{Art}^{\mathbf{triv}}$-homogeneity).

The category $\widehat{\cF}_{x}$ naturally extends to its
completion, which by definition is the pro-category $\pro\widehat{\cF}_{x}$, which is a category cofibred in groupoids over
$\pro{\cC}_{\Lambda}$. (Note that this is a more general definition
than that of~\cite[\href{http://stacks.math.columbia.edu/tag/06H3}{Tag
  06H3}]{stacks-project}, which restricts the definition to $\widehat{\cC}_\Lambda$.) 
There is a fully faithful embedding of $\widehat{\cF}_{x}$ into its completion,
which attaches to any object of $\widehat{\cF}_{x}$ lying over an Artinian $\Lambda$-algebra
the corresponding pro-object, and this embedding induces an equivalence between
$\widehat{\cF}_{x}$ and the restriction of its completion to
$\cC_{\Lambda}$. We therefore also denote the completion of $\widehat{\cF}_{x}$ by $\widehat{\cF}_{x}$. 
If $A$ is an object of $\pro\cC_{\Lambda}$, then we will usually
denote an object of the completion of $\widehat{\cF}_{x}$ lying over
$A$ by a morphism $\Spf A \to \widehat{\cF}_{x}$.

We also introduce the notation $\cF_x$ to denote the following category
cofibred in groupoids over $\pro{\cC}_{\Lambda}$: if $A$ is an object
of $\pro{\cC}_{\Lambda}$, then $\cF_x(A)$ denotes the groupoid
consisting of morphisms $\Spec A \to \cF$, together with an isomorphism
between the restriction of this morphism to the closed point of $\Spec A$
and the given morphism $x:\Spec k \to \cF$.   
If $A$ is Artinian,
then $\cF_x(A) = \widehat{\cF}_x(A)$.  In general, there is a natural
functor $\cF_x(A) \to \widehat{\cF}_x(A)$ (the functor
of {\em formal completion}); a morphism $\Spf A \to \cF$
lying in the essential image of this functor is said to be {\em
  effective}.

\begin{remark}\label{rem:change of view-point for formal objects}
If $A$ is an object of $\pro\cC_{\Lambda}$, then we may consider the
formal scheme $\Spf A$ as defining a sheaf of sets on the \'etale site of $S$,
via the definition $\Spf A := \varinjlim_{i} \Spec A_i$ (writing $A$ as the 
projective limit of its discrete Artinian quotients $A_i$, and taking
the inductive limit in the category of \'etale sheaves; this is a special
case of the Ind-constructions considered in
Subsection~\ref{subsec:non-representable} below).   Of course,
we may also regard the resulting sheaf $\Spf A$ as a stack (in setoids).

Giving a morphism $\Spf A \to \widehat{\cF}_x$ in the sense described above
is then equivalent to giving a morphism of stacks $\Spf A \to
\cF$ which induces the given morphism $x: \Spec k \to \cF$ when composed
with the natural morphism $\Spec k \to \Spf A$.  This view-point is useful
on occasion; for example, we say that $\Spf A \to \widehat{\cF}_x$ is
a {\em formal monomorphism} if the corresponding morphism of stacks
$\Spf A \to \cF$ is a monomorphism.  (Concretely, this amounts to the
requirement that the induced morphism $\Spec A_i \to \cF$ is a monomorphism
for each discrete Artinian quotient $A_i$ of $A$.)\end{remark}

We now introduce the notion of a versal ring at the
morphism $x$, which will be used in the definition of Axiom~[2](b). (See Remark~\ref{rem:independence of morphism for 2b}
for a discussion of why we speak of a versal ring at a morphism,
rather than at a point.) As above, fix the finite type morphism $\Spec k \to S$, 
an affine open subset $\Spec \Lambda \to S$ through which 
this morphism factors,
and the lift of this morphism to a morphism
$x: \Spec k \to \cF$.

\begin{df}
\label{def:versal ring}

Let $A_x$ be a topological local $\Lambda$-algebra corresponding 
to an element of $\pro\cC_{\Lambda}$.  We say that a morphism
$\Spf A_x \to \widehat{\cF}_x$ is {\em versal} if it is smooth, 
in the sense 
  of~\cite[\href{http://stacks.math.columbia.edu/tag/06HR}{Tag
    06HR}]{stacks-project},
    i.e.\ satisfies the infinitesimal lifting property with
  respect to morphisms in $\cC_{\Lambda}$. More precisely, given a
  commutative diagram
$$\xymatrix{\Spec A \ar[r] \ar[d] & \Spec B \ar[d] \ar@{-->}[ld]\\
\Spf A_{x} \ar[r] & \widehat{\cF}_{x}}$$
in which the upper arrow is the closed immersion corresponding
to a surjection $B \to A$ in $\cC_{\Lambda}$, 
and the left hand
vertical arrow corresponds to a morphism in
$\pro\cC_\Lambda$ (equivalently, it is continuous when~$A_x$ is given
its projective limit topology, and~$A$ is given the discrete topology)  we can fill in the
dotted arrow (with a morphism coming from a morphism in $\pro\cC_\Lambda$) so that the diagram remains commutative.

We refer to $A_x$ as a {\em versal ring} to $\cF$ at the morphism
$x$. We say that $A_x$ is an \emph{effective versal ring} to $\cF$ at the morphism
$x$ if the morphism $\Spf A_x \to \widehat{\cF}_x$ is effective.

We say that~$\cF$ \emph{admits versal rings at all finite type
  points} if there is a versal ring for every morphism $x:\Spec
k\to\cF$ whose source is a finite type $\cO_S$-field. We say
that~$\cF$ \emph{admits effective versal rings at all finite type
  points} if  there is an effective versal ring for every morphism $x:\Spec
k\to\cF$ whose source is a finite type $\cO_S$-field.
\end{df}
Then Axiom~[2](b) is by definition the assertion that $\cF$ admits
Noetherian effective versal rings at all finite type
  points.

\begin{remark}
\label{rem:independence of morphism for 2b}
One complication in verifying Axiom~[2](b)
is that, at least {\em a priori}, it does not depend simply 
on the finite type point of $\cF$ represented by a given morphism
$x: \Spec k \to \cF$ (for a field $k$ of finite type over $\cO_S$),
but on the particular morphism.  More concretely, if we are given
a finite type morphism $\Spec l \to \Spec k$, and if
we let $x'$ denote the composite $\Spec l \to \Spec k \to \cF$,
then it is not 
obvious that validity of~[2](b) for either of $x$ or $x'$ implies the validity
of~[2](b) for the other.   One of the roles of Axiom~[2](a) in the theory
is to bridge the gap between different choices of field defining the
same finite type point of~$\cF$, and indeed one can 
show that in the presence of~[2](a), the property of $\widehat{\cF}_x$
admitting a Noetherian versal ring depends only on the finite type point of $\cF$ underlying $x$;
see part~(1) of
Lemma~\ref{lem:independence of morphism for 2b} below.

The problem of showing that effectivity is independent of the choice
of representative of the underlying finite type point of $\cF$ seems slightly
more subtle,
and to obtain a definitive result we have to make additional assumptions
on $\cF$, and possibly on $S$.
This is the subject of the parts (2), (3), and (4) of 
Lemma~\ref{lem:independence of morphism for 2b}.

\end{remark}
In the case when the morphism $\Spec l \to \Spec k$ is separable,
it is possible to make a softer argument to pass the existence of a versal
ring from from $x$ to $x'$, without making any Noetherianness assumptions.
We do this in Lemma~\ref{lem:change of versal ring under separable
  extension}; we firstly recall 
from~\cite[\href{http://stacks.math.columbia.edu/tag/07W7}{Tag
  07W7},\href{http://stacks.math.columbia.edu/tag/07WW}{Tag
  07WW}]{stacks-project} a useful formalism for considering
such a change of residue field.

\begin{rem}\label{rem: material on change of residue field}
  Let $\cF$ be a category fibred in groupoids over~$S$, and fix a
  morphism $x: \Spec k \to \cF$, where $k$ is a finite
  type~$\cO_S$-field.  Suppose that we are given a finite type
  morphism $\Spec l \to \Spec k$, so that $l/k$ is a finite extension
  of fields, and let $x'$ denote the composite
  $\Spec l \to \Spec k \to \cF$. Write $\cC_{\Lambda,k}$ for the
  category of local Artinian $\Lambda$-algebras with residue field $k$,
  and similarly write $\cC_{\Lambda,l}$ for the category of local
  Artinian $\Lambda$-algebras with residue field $l$.

  We let $(\widehat{\cF}_x)_{l/k}$ denote the category cofibred in
  groupoids over $\cC_{\Lambda,l}$ defined by setting
  $(\widehat{\cF}_{x})_{l/k}(B) := \widehat{\cF}_x(B\times_l k)$, for
  any object $B$ of $\cC_{\Lambda,l}$.  If~$\cF$ satisfies~[2](a), it
  follows
  from~\cite[\href{http://stacks.math.columbia.edu/tag/07WX}{Tag
    07WX}]{stacks-project} that there is a natural equivalence of
  categories cofibred in groupoids
  $(\widehat{\cF}_{x})_{l/k} \iso \widehat{\cF}_{x'}$; an examination
  of the proof shows that if~$l/k$ is separable, the same conclusion
  holds if~$\cF$ is only assumed to be
  $\mathbf{Art}^{\mathbf{sep}}$-homogeneous.
\end{rem}

\begin{rem}\label{rem: minimal versal rings}
  Recall \cite[\href{http://stacks.math.columbia.edu/tag/06T4}{Tag
    06T4}]{stacks-project} that we say that a Noetherian versal
  morphism $\Spf A_x\to\widehat{\cF}_x$ is \emph{minimal} if whenever
  we can factor this morphism through a morphism
  $\Spf A\to\widehat{\cF}_x$, the underlying map $A\to A_x$ is
  surjective. This notion is closely related to conditions on the
  tangent space of $\widehat{\cF}_x$, in the following way.

By definition, the \emph{tangent space} $T\widehat{\cF}_x$
of~$\widehat{\cF}_x$ is the $k$-vector space
$\widehat{\cF}_x(k[\epsilon])$. As explained
in~\cite[\href{http://stacks.math.columbia.edu/tag/06I1}{Tag
  06I1}]{stacks-project}, there is a natural action
of~$\Der_\Lambda(k,k)$ on $T\widehat{\cF}_x$, and the versal morphism
$\Spf A_x\to\widehat{\cF}_x$ gives rise to a
$\Der_\Lambda(k,k)$-equivariant morphism $d:\Der_{\Lambda}(A_x,k)\to T\widehat{\cF}_x$.

By~\cite[\href{http://stacks.math.columbia.edu/tag/06IR}{Tag
  06IR}]{stacks-project}, a Noetherian versal
  morphism $\Spf A_x\to\widehat{\cF}_x$ is minimal provided that the
  morphism~$d$ is bijective on~$\Der_\Lambda(k,k)$-orbits. Conversely,
  if~$\cF$ is $\mathbf{Art}^{\mathbf{triv}}$-homogeneous, then it
  follows from the proof
  of~\cite[\href{http://stacks.math.columbia.edu/tag/06J7}{Tag
    06J7}]{stacks-project} that~$\widehat{\cF}_x$ satisfies the
  condition (S2)
  of~\cite[\href{http://stacks.math.columbia.edu/tag/06HW}{Tag
    06HW}]{stacks-project}, and it then follows
  from~\cite[\href{http://stacks.math.columbia.edu/tag/06T8}{Tag
    06T8}]{stacks-project} that if
$\Spf A_x\to\widehat{\cF}_x$ is minimal,
then $d$ is bijective
  on~$\Der_\Lambda(k,k)$-orbits. 
\end{rem}

\begin{lem}
\label{lem:change of versal ring under separable extension}
Let $\cF$ be a category fibred in groupoids over $S$ which
is $\mathbf{Art}^{\mathbf{sep}}$-homogeneous.
Suppose given $x:\Spec k \to \cF$, with $k$ 
a finite type $\cO_S$-field, let $l/k$ be a finite separable
extension,  and let $x'$ denote the composite $\Spec l \to \Spec k
\stackrel{x}\longrightarrow \cF$.
Suppose also that $\Spf A_x \to \widehat{\cF}_x$ is a versal ring to~$\cF$ at the
morphism~$x$,
so that in particular $A_x$ is a pro-Artinian
local $\cO_S$-algebra with residue field $k$.
Let $A_{x'}/A_x$ denote the finite \'etale extension of $A_x$
corresponding {\em (}via the topological invariance of the \'etale site{\em )}
to the finite extension $l$ of $k$, so that in particular $A_{x'}$ is a 
pro-Artinian local $\cO_S$-algebra  with residue field $l$.

Then the induced morphism $\Spf A_{x'} \to \widehat{\cF}_{x'}$ realises
$A_{x'}$ as a versal ring to~$\cF$ at the morphism~$x'$. If $A_x$ is
Noetherian, and the morphism $\Spf
A_x\to\widehat{\cF}_x$ is minimal, then so is the
morphism $\Spf A_{x'} \to \widehat{\cF}_{x'}$.\end{lem}
\begin{proof}
By Remark~\ref{rem: material on change of residue field}, since $\cF$ is
$\textbf{Art}^{\textbf{sep}}$-homogeneous and~$l/k$ is separable, we
have a natural equivalence of groupoids $(\widehat{\cF}_{x})_{l/k}
\iso \widehat{\cF}_{x'}$. Suppose given a commutative diagram
$$\xymatrix{\Spec A \ar[r] \ar[d] & \Spec B \ar[d] \ar@{-->}[ld]\\
\Spf A_{x'} \ar[r] & \widehat{\cF}_{x'}}$$
in which the upper arrow is the closed immersion corresponding
to a surjection $B \to A$ in $\cC_{\Lambda,l}$;
we wish to show that we can fill in the dotted arrow in such a way
that resulting diagram remains commutative.
Passing to the pushout with $k$ over~ $l$, and noting that the morphism $A_x \to A_{x'}$ factors through $A_{x'} \times_l k,$ we obtain a commutative
diagram
 $$\xymatrix{\Spec A\times_l k \ar[r] \ar[d] & \Spec B\times_l k \ar@{-->}[d] \ar[rd]
\\
\Spf A_{x'}\times_l k \ar[r] & \Spf A_{x} \ar[r] & \widehat{\cF}_x}$$
where the dotted arrow exists by the versality of the morphism  $\Spf A_x \to \widehat{\cF}_x.$
Now consider the diagram
$$\xymatrix{\Spec A \ar[r]\ar[dd] & \Spec B \ar[d]\ar@{-->}[ldd] \\
& \Spec B\times_l k \ar[d] \\
\Spf A_{x'} \ar[r] & \Spf A_x.}
$$
Since $\Spf A_{x'} \to \Spf A_x$ is formally \'etale, we may fill
in the dotted arrow so as to make the resulting diagram commutative.  
This dotted arrow also makes the original diagram commutative, 
so $A_{x'}$ is a versal ring to~$\cF$ at the morphism~$x'$.

Finally, suppose that the versal morphism $\Spf A_x\to\widehat{\cF}_x$
is minimal. By Remark~\ref{rem: minimal versal rings}, the natural
morphism $\Der_{\Lambda}(A_x,k)\to T\widehat{\cF}_x$ is a bijection
on $\Der_{\Lambda}(k,k)$-orbits, and we need to show that the natural
morphism $\Der_{\Lambda}(A_{x'},l)\to T\widehat{\cF}_{x'}$ is a
bijection on $\Der_\Lambda(l,l)$-orbits.

By~\cite[\href{http://stacks.math.columbia.edu/tag/06I0}{Tag
  06I0},\href{http://stacks.math.columbia.edu/tag/07WB}{Tag
07WB}]{stacks-project} there is a natural isomorphism of $l$-vector spaces
$T\widehat{\cF}_x\otimes_k l\isoto T\widehat{\cF}_{x'}$. Since
$A_{x'}/A_x$ is \'etale (and $l/k$ is separable), restriction induces
isomorphisms
$\Der_\Lambda(A_{x'},k)\cong\Der_\Lambda(A_x,k)$ and
$\Der_\Lambda(l,k)\cong\Der_\Lambda(k,k)$, and thus we also have 
natural isomorphisms of $l$-vector spaces $\Der_\Lambda(A_x,k)\otimes_k
l\isoto\Der_\Lambda(A_{x'},l)$ and
$\Der_\Lambda(l,l)\isoto \Der_\Lambda(k,k)\otimes_kl$.   One easily 
checks that the base-change by $l$ over $k$ of 
the morphism $d: \Der_{\Lambda}(A_x,k) \to T\widehat{\cF}_x$
coincides, with respect to these identifications,
with the morphism $d: \Der_{\Lambda}(A_{x'}, l) \to T\widehat{\cF}_{x'}$,
in a manner which identifies, under the isomorphism $\Der_{\Lambda}
(k,k)\otimes_k l \cong \Der_{\Lambda}(l,l)$, the action of
$\Der_{\Lambda}(k,k)\otimes_kl$ 
with the action of $\Der_{\Lambda}(l,l)$.  The result
follows.\end{proof}

We now engage in a slight digression; namely, we use the theory of
versal rings in order to define
the notion of the complete local
ring at a finite type point of an algebraic space locally of finite type
over $S$. 

In order to motivate this concept,
we recall first
that quasi-separated algebraic spaces are decent, in the sense of~\cite[\href{http://stacks.math.columbia.edu/tag/03I8}{Tag 03I8}]{stacks-project},
which is to say that 
any point of such
an algebraic space $X$ is represented by a quasi-compact
monomorphism $x: \Spec k
\to X$.  Given a point $x$ of a quasi-separated algebraic space,
thought of as such a quasi-compact monomorphism,
Artin and Knutson (\cite[Defn.\ 2.5]{MR0262237}, \cite[Thm.\ 6.4]{MR0302647}) define a Henselian local ring of $X$ 
at $x$; completing this local ring gives our desired complete local ring
in this case.   Rather than imposing a quasi-separatedness hypothesis at
this point and appealing to these results,
we adopt a slightly different approach, which  will allow us to
define a complete local ring
at any finite type 
point $x$ of a locally finite type algebraic space $X$ over~$S$.
(Note  that throughout the discussion,  we will
maintain our assumption that $S$ is locally Noetherian.)

To begin with, we note that finite type points always admit
representatives that are monomorphisms (regardless of any separatedness
hypothesis); indeed, we have the following lemma.

\begin{lemma}
\label{lem:monomorphisms at points}
Any finite type point of an algebraic space~$X$,
locally of finite type over the locally Noetherian scheme $S$,
admits a representative
$\Spec k \to X$ which a monomorphism.  This representative is unique up
to unique isomorphism, the field $k$ is a finite type $\cO_S$-field,
and any other representative $\Spec K \to X$ of the given point
factors through this monomorphic representative in a unique fashion.
\end{lemma}
\begin{proof}
We apply the criterion of part~(1) of
\cite[\href{http://stacks.math.columbia.edu/tag/03JU}{Tag 03JU}]{stacks-project}.  More precisely, we choose an \'etale morphism $U \to X$ with $U$ a scheme over $S$, which will again be locally of finite type over $S$ (since smooth morphisms are locally of finite type, and $X$ is assumed to be locally of finite type over $S$).  We choose a pair of finite type points $u,u' \in U$ lying
over the given finite type point of $X$; as noted in the proof of the 
lemma just cited, we must verify that the underlying topological 
space of the scheme $u\times_X u'$ is finite.  But since the diagonal
map $X \to X \times_S X$ is a monomorphism (as $X$ is an algebraic space),
the fibre product $u \times_X u'$ maps via a monomorphism into
$u \times_ S u'$.   This latter scheme has a finite underlying topological
space (since $u$ and $u'$ are each the Spec of some finite type $\cO_S$-field).
Since  monomorphisms induce embeddings on underlying topological spaces,
we see that $u\times_X u'$ also has a finite underlying topological space.
By the above cited lemma,
this implies the existence of the desired monomorphism $\Spec k \to X$
representing the given point.  Denote this monomorphism by $x$.

If $x': \Spec K \to X$ is any other representative of the same point,
then  we may consider the base-changed morphism
$\Spec k \times_X \Spec K \to \Spec K,$ which is again a monomorphism.
Since its source is non-empty (as $x$ and $x'$ represent 
the same point of $X$), and as its target is the Spec of a field,
it must be an isomorphism; equivalently, the morphism $x'$ 
must factor through $x$.  Since $x$ is a monomorphism, this factorisation
is unique.  Since $K$ can be chosen to be a finite type $\cO_S$-field,
we see that $k$ must in particular be a finite type $\cO_S$-field.
If $x'$ is also a monomorphism, then we may reverse the roles of $x$
and $x'$, and so conclude that $x$ is determined uniquely
up to unique isomorphism.
This completes the proof of the lemma.
\end{proof}

The following proposition then constructs complete local rings
at finite type points of locally finite type algebraic spaces over $S$.

\begin{prop}
\label{prop:complete local rings}
If $X$ is an algebraic space,
locally of finite type over the locally Noetherian scheme~$S$,
and if $x: \Spec k \to X$ is a monomorphism, for some field $k$
of finite type over $\mathcal O_S$, then there is an effective
Noetherian versal
ring $A_x$ to $X$ at the morphism $x$ with the property that 
the corresponding morphism
$\Spf A_x \to X$ 
is a formal monomorphism.  Furthermore, the ring $A_x$, equipped with
its morphism $\Spec A_x \to X$ inducing the given morphism $x$, is unique up to unique isomorphism.
Finally,
if $A$ is any object of $\cC_{\Lambda}$, then any morphism $\Spec A \to X$
factors uniquely through the morphism $\Spec A_x \to X$.
\end{prop}
\begin{proof}
Let $A_x$ be a minimal versal ring to $X$ at the morphism $x,$
in the sense 
of~\cite[\href{http://stacks.math.columbia.edu/tag/06T4}{Tag 06T4}]{stacks-project}
(which exists, and is Noetherian, by virtue
of~\cite[\href{http://stacks.math.columbia.edu/tag/06T5}{Tag 06T5}]{stacks-project}
and the fact that $X$, being an algebraic space, 
admits a Noetherian versal ring at the morphism~$x$; see e.g.\ 
  Theorem~\ref{thm: Artin representability} below).
We will show that the morphism $\Spf A_x \to X$ is a formal monomorphism.
To this end, we choose an \'etale surjective morphism
$U \to X$ whose source is a scheme
(such a morphism exists, since $X$ is an algebraic space).   
It suffices
(by~\cite[\href{http://stacks.math.columbia.edu/tag/042Q}{Tag
  042Q}]{stacks-project}, and the definition of a formal
monomorphism in Remark~\ref{rem:change of view-point for formal objects}) to show that the base-changed morphism
$U\times_X \Spf A_x \to U$ is a formal monomorphism. We begin by describing this morphism more explicitly.

Abusing notation slightly, we write $x$ to denote the point $\Spec k$,
as well its monomorphism into $X$.
The pull-back of $U$ over  the monomorphism $x \to X$ is then an \'etale
morphism $U_x \to x$ with non-empty source.  We may write $U_x$ as
a disjoint union of points $u_i$, each of which is of the form
$u_i = \Spec l_i,$
for some finite separable extension $l_i$ of $k$. 
Since $U_x \to U$ is a monomorphism (being the base-change of a monomorphism),
each of the composites $u_i \to U_x \to U$ is also a monomorphism;
in other words, for each $i$,
the field $l_i$ is also the residue field of the image of $u_i$ in $U$;
in light of this, we identify $u_i$ with its image in $U$.

For each $i$, let $A_i$ denote the finite \'etale extension of $A_x$
corresponding, via the topological invariance of the \'etale site,
to the finite extension $l_i/k$.
The projection $U\times_{X} \Spf A_x
\to \Spf A_x$ is formally \'etale
(being the pull-back of an \'etale morphism),
and thus admits a natural identification with the morphism
$\coprod_{i \in I} \Spf A_{i} \to \Spf A_x$.
Since the points $u_i$ of $U$ are distinct for distinct values of $i$,
in order to show that $\coprod_{i\in I} \Spf A_i \cong U\times_X \Spf A_x
\to U$ is a formal monomorphism, it suffices to show that each
of the individual morphisms $\Spf A_i \to U$ is a formal monomorphism;
and we now turn to doing this.

For each $i$, we let $\cC_{\Lambda, l_i}$ denote the category of local
Artinian $\cO_S$-algebras with residue field $l_i$.
We write $x_i$ to denote the composite $\Spec l_i = u_i \to U \to X$. The infinitesimal lifting property for the \'etale morphism
$U \to X$ shows that the induced morphism $\widehat{U}_{u_i} \to 
\widehat{X}_{x_i}$ is an equivalence of categories cofibred in
groupoids over $\cC_{\Lambda,l_i}$.
Lemma~\ref{lem:change of versal ring under separable extension}
shows that $A_{i}$ is a versal ring to $X$ at $x_i$,
and, by the noted equivalence, it is thus also a versal ring to $U$
at $u_i$.   In fact, since $A_{x}$ is a minimal versal ring at $x$,
each of the rings $A_{i}$ is a minimal versal ring at $u_i$.

On the other hand, since $u_i$ is a point of the scheme $U$, 
the complete local ring
$\widehat{\cO}_{U,u_i}$ is a minimal versal ring to $U$ at $u_i$;
thus we may identify $A_{i}$ with $\widehat{\cO}_{U,u_i}$,
so that the morphism $\Spf A_i \to U$ is identified with the canonical
morphism $\Spf \widehat{\cO}_{U,u_i} \to U$.  This latter morphism
is a formal monomorphism, and thus so is the former.

The fact that the morphism $\Spf A_x \to X$ is effective
is a consequence of $X$ being an algebraic space;
see~\cite[\href{http://stacks.math.columbia.edu/tag/07X8}{Tag 07X8}]{stacks-project}.
The uniqueness of this morphism, up to unique isomorphism, follows
from Lemma~\ref{lem: formal versal monomorphism is unique up to unique
    isomorphism.} below; the fact that its effectivisation is unique up to unique isomorphism
again follows
from~\cite[\href{http://stacks.math.columbia.edu/tag/07X8}{Tag
  07X8}]{stacks-project}. 

Since $A_x$ is versal to $X$ at $x$,
any morphism $\Spec A \to X$, for $A$ an object
of $\cC_{\Lambda}$, factors through the morphism $\Spf A_x \to X$.  Of course,
it then factors through the induced morphism $\Spec B \to X$, for
some Artinian quotient $B$ of $A_x$, and hence also through the morphism
$\Spec A_x \to X$.  The uniqueness of this factorisation again follows
from Lemma~\ref{lem: formal versal monomorphism is unique up to unique
    isomorphism.} below.\end{proof}

\begin{lem}
  \label{lem: formal versal monomorphism is unique up to unique
    isomorphism.}Let $\cF$ be a category fibred in groupoids, and
  suppose that $\Spf
  A_x\to\widehat{\cF}_x$ is a versal morphism at the morphism  $x:\Spec
  k\to\cF$, where $k$ is a finite type $\cO_S$-field. Suppose also that $\Spf
  A_x\to\widehat{\cF}_x$ is a formal monomorphism.

Then if~$A$ is any object of $\cC_\Lambda$, any morphism $\Spec
A\to\widehat{\cF}_x$ factors uniquely through the morphism $\Spf
  A_x\to\widehat{\cF}_x$. Furthermore, the ring $A_x$, together with the
  morphism $\Spf
  A_x\to\widehat{\cF}_x$, is uniquely determined up to unique
  isomorphism by the property of being a versal formal monomorphism.
\end{lem}
\begin{proof}
Since the morphism $\Spf A_x \to \widehat{\cF}_x$ is versal by assumption,
  any morphism  $\Spec
A\to\widehat{\cF}_x$ factors through this morphism;
that this factorisation is unique is
  immediate from the definition of a formal monomorphism. If $\Spf A'_x\to\widehat{\cF}_x$ is
  another versal formal monomorphism, then applying this property to
  the discrete Artinian quotients of $A'_x$, and then reversing the
  roles of $A_x$ and $A'_x$, we find the required unique isomorphisms.
\end{proof}
\begin{df}\label{df:complete local rings}
We refer to the ring $A_x$ of Proposition~\ref{prop:complete local rings} as the \emph{complete local ring of $X$ at the point $x$}.
\end{df}

In Subsection~\ref{subsec:non-representable} below,
we generalise the notion of the complete local ring at a point
to certain Ind-algebraic spaces;
see Definition~\ref{df:Ind complete local rings}.
We now state and prove
a result which will be used in Section~\ref{subsec:relatively
  representable}, and which uses this generalisation.

\begin{lem}\label{lem:versal rings and base change}
  If $\cF$ admits versal rings at all finite type points
{\em (}in the sense of Definition~{\em \ref{def:versal ring})}, and if $\cF'\to\cF$ is
  representable by algebraic spaces and locally of finite presentation,
then $\cF'$ admits versal rings at all
  finite type points. In fact, if $x':\Spec k \to \cF'$
is a morphism from a finite type $\cO_S$-field to $\cF',$
inducing the morphism $x:\Spec k \to \cF$,
and if $\Spf A \to \cF$ is a versal ring
to $\cF$ at $x$, 
then $X:= \cF'\times_{\cF} \Spf A$ is defined as
an Ind-locally finite type algebraic space over $S$ {\em (}in the
sense of Definition~{\em \ref{df:Ind-lft alg space}} below{\em )},
the morphism~$x'$ induces a lift of $x$ to $X$,
and the complete local ring of $X$ at $x$ is a versal ring to~ $\cF'$ at $x'$.
\end{lem}
\begin{proof}If we write $A \cong \varprojlim_{i \in I} A_i$ as a projective limit
of finite type local Artinian $\cO_S$-algebras,
then we define $X := \varinjlim_{i \in I} \cF' \times_{\cF} \Spec A_i$;
thus $X$ is an Ind-locally finite type algebraic space over $S$,
which is clearly well-defined
(as a sheaf of setoids on the \'etale site of $S$)
independently of the choice
of description of $A$ as a projective limit. The only claim, then, that is not immediate from the definitions
is that the composite
  morphism $\Spf \cO_{X,x} \to X\to\cF'$ is versal. We leave this as an easy
  exercise for the reader; it is essentially
  immediate from the versality of~$\Spf A$.\end{proof}

We now introduce the notion of a presentation of a deformation
category by an effectively 
Noetherianly pro-representable smooth groupoid in functors, which is
closely related to the notion of admitting an effective Noetherian
versal ring. Our reason
for introducing this notion is to prove Lemma~\ref{lem:versal rings give Art-triv if the relation is
  Art-triv} and Corollary~\ref{cor: 3 and 1 and versal rings implies
  arttriv.}; under an appropriate hypothesis on the diagonal of~$\cF$,
these will enable us to deduce Axiom [2](a) from [2](b).  

We say that a set-valued functor on $\cC_{\Lambda}$ is {\em pro-representable}
if it representable by an object of $\pro\cC_{\Lambda}$.  If $A$ is
the associated topological ring to the representing pro-object,
then we will frequently denote this functor by $\Spf A$.  We say
that a functor is {\em Noetherianly pro-representable} if it 
is pro-representable 
by an object $A$ of $\widehat{\cC}_{\Lambda}$.\footnote{In~\cite{stacks-project}, what we call {\em Noetherian pro-representability}
is called simply {\em pro-representability};
see~\cite[\href{http://stacks.math.columbia.edu/tag/06GX}{Tag 06GX}]{stacks-project}.
However, we will need to consider more general pro-representable functors,
and so we need to draw a distinction between 
the general case and the case of pro-representability
by an object of $\widehat{\cC}_{\Lambda}$.}

We refer to~\cite[\href{http://stacks.math.columbia.edu/tag/06K3}{Tag 06K3}]{stacks-project} for the definition of
a groupoid in functors over~ $\cC_\Lambda$, and then make the following
related definitions.

\begin{df}
	\label{def:groupoid defs}
	(1)
We say that a groupoid in functors over $\cC_{\Lambda}$, say
$(U,R,s,t,c)$, is \emph{smooth} if $s,t:R\to U$ are
smooth\footnote{The term {\em smooth} is used here in the 
sense of~\cite[\href{http://stacks.math.columbia.edu/tag/06HG}{Tag 06HG}]{stacks-project}; i.e.\ we require the infinitesimal
lifting property with respect to morphisms in $\cC_{\Lambda}$.  Other
authors might use the term {\em versal} here, because the residue 
field is being held fixed.}
; equivalently, if the quotient morphism $U\to [U/R]$ is smooth. 

(2)
We say that
$(U,R,s,t,c)$ is ({\em Noetherianly}) \emph{pro-representable} 
if $U$ and $R$ are each ({\em Noetherianly}) \emph{pro-representable}.
\end{df}

A \emph{presentation} of $\widehat{\cF}_{x}$ is an equivalence
$[U/R]\isoto\widehat{\cF}_{x}$ of categories cofibred in groupoids over
$\cC_\Lambda$,
where $(U,R,s,t,c)$ is a groupoid in functors over $\cC_\Lambda$.
Suppose given such a presentation
by a groupoid in functors that is Noetherianly pro-representable,
in the sense of Definition~\ref{def:groupoid defs}, and let
$A_{x}\in\Ob(\widehat{\cC}_\Lambda)$ be an object that pro-represents $U$. We then obtain
an induced  morphism
\numequation
\label{eqn:morphism from U to F}
\Spf A_x = U\to [U/R]\to\widehat{\cF}_{x}. 
\end{equation}

\begin{df}
	\label{def:effectively pro-representable}
We say that the given presentation is \emph{effectively}
Noetherianly pro-representable if the morphism~(\ref{eqn:morphism
       from U to F}) is effective,
i.e.\ arises as the formal completion of a morphism $\Spec
A_{x}\to\cF$.
\end{df}

The existence of an effectively Noetherianly
pro-representable presentation by a smooth groupoid in functors is
closely related to the property of having effective versal rings, as
we will now see.

\begin{lem}\label{lem:relationship of versal rings to pro-representability}
If $[U/R] \iso \widehat{\cF}_x$ is a presentation
of $\widehat{\cF}_x$ by a smooth groupoid in functors,
for which $U$ is pro-representable by a topological local $\Lambda$-algebra
$A$, then the morphism $\Spf A = U \to \widehat{\cF}_x$ is versal.
Conversely,
if $A$ is the topological local $\Lambda$-algebra
corresponding to some element of $\pro\cC_{\Lambda}$,
and if $\Spf A \to \widehat{\cF}_x$ is versal,
then if we write $U := \Spf A$ and $R = U\times_{\widehat{\cF}_x} U,$
and $s$, $t$ for the two projections $R\to U$, 
then $(U,R,s,t)$ is a smooth groupoid in functors, and the natural morphism
$[U/R] \to \widehat{\cF}_x$ is an equivalence,
and thus equips $\widehat{\cF}_x$ with a presentation by a smooth groupoid
in functors.
\end{lem}
\begin{proof}
Essentially by definition, if $[U/R] \iso \widehat{\cF}_x$ is a presentation
of $\widehat{\cF}_x$ by a smooth groupoid in functors,
then the induced morphism $U \to [U/R] \iso \widehat{\cF}_x$ is smooth
(in the sense 
  of~\cite[\href{http://stacks.math.columbia.edu/tag/06HR}{Tag
    06HR}]{stacks-project}),
and so by definition is versal.
The converse statement follows from~\cite[\href{http://stacks.math.columbia.edu/tag/06L1}{Tag
    06L1}]{stacks-project}).
\end{proof}

\begin{rem}\label{rem:relationship of versal rings to pro-representability}
Suppose that $\widehat{\cF}_x$ admits a presentation by a smooth Noetherianly
pro-representable groupoid in functors.
Then, if $U = \Spf A \to \widehat{\cF}_x$ is a versal morphism with $A$
an object of $\widehat{\cC}_{\Lambda}$ (so that $U$ is Noetherianly
pro-representable),
one finds that $R := U\times_{\widehat{\cF}_x} U$ is also
Noetherianly pro-representable 
(\emph{cf.}\ the proof of~\cite[\href{http://stacks.math.columbia.edu/tag/06L8}{Tag
    06L8}]{stacks-project}), 
and so the equivalence $[U/R] \iso \widehat{\cF}_x$ of Lemma~\ref{lem:relationship of versal rings to pro-representability}
gives a particular smooth Noetherianly pro-representable presentation
of $\widehat{\cF}_x$.

In Lemma~\ref{lem:versal rings give Art-triv if we know [3]}
below we will show that if the diagonal 
of $\cF$ satisfies an appropriate hypothesis,
and if we are given a versal morphism
$U = \Spf A \to \widehat{\cF}_x$ from
a (not necessarily Noetherianly) pro-representable functor $U$,
then (without any {\em a priori} hypothesis that
$\widehat{\cF}_x$ admits a presentation by a smooth pro-representable
groupoid in functors)
the fibre product
$R := U \times_{\widehat{\cF}_x} U$ is also pro-representable,
and thus (taking into account the isomorphism $[U/R] \iso \widehat{\cF}_x$
of Lemma~\ref{lem:relationship of versal rings to pro-representability})
the existence of a versal ring to $\cF$ at the morphism $x$ will
imply that $\widehat{\cF}_x$ in fact admits a presentation
by a pro-representable smooth groupoid in functors.
\end{rem}

We close our discussion of Axiom~[2] by stating a lemma
that relates the existence of presentations by smooth pro-representable
groupoids in functors
to the conditions of Definition~\ref{df:homogeneity}.

\begin{lem}\label{lem:versal rings give Art-triv if the relation is
  Art-triv}Suppose, for every morphism $x: \Spec k \to \cF$,
with $k$ a finite type $\cO_S$-field,
that $\widehat{\cF}_x$ admits a presentation
by a pro-representable smooth groupoid in functors.
Then $\cF$ is $\mathbf{Art}^{\mathbf{triv}}$-homogeneous
{\em (}in the sense
of Definition~{\em \ref{df:homogeneity}}{\em )}.\end{lem}
\begin{proof}
By~\cite[\href{http://stacks.math.columbia.edu/tag/06KT}{Tag
    06KT}]{stacks-project},
we need only check that a pro-representable functor is
  $\mathbf{Art}^{\mathbf{triv}}$-homogeneous (or in the language
  of~\cite{stacks-project}, satisfies (RS)).
In the case of Noetherianly pro-representable functors,
this is~\cite[\href{http://stacks.math.columbia.edu/tag/06JB}{Tag
    06JB}]{stacks-project}, and the proof in the general case is identical.
\end{proof}

\subsection{Remarks on Axiom [3]} 
\label{subsec:axiom 3}

Recall that a category fibred in groupoids $\cF$ satisfies Axiom~[3] if
and only if the diagonal
$\Delta:\cF \to \cF \times_S \cF$ is representable by algebraic spaces;
equivalently (by~\cite[\href{http://stacks.math.columbia.edu/tag/045G}{Tag
  045G}]{stacks-project}), if and only if $X\times_{\cF}Y$ is an algebraic space
whenever $X\to\cF,Y\to\cF$ are morphisms from algebraic spaces $X,Y$. We begin with the following lemma.

\begin{lemma}
\label{lem:fibre products}
Let $\cF$ be a category fibred in groupoids satisfying Axiom~{\em [3]}.
 If $\cX$ and $\cY$ are categories fibred in groupoids satisfying~{\em [3]},
then for any morphisms of categories fibred in groupoids $\cX, \cY \to \cF$,
the fibre product $\cX \times_{\cF} \cY$ is again a category fibred in groupoids satisfying~{\em [3]}.
If $\cX$ and $\cY$ are furthermore {\em (}algebraic{\em )} stacks,
then the fibre product is also an {\em (}algebraic{\em )} stack.
\end{lemma}
\begin{proof}
The claim for algebraic stacks is proved
in~\cite[\href{http://stacks.math.columbia.edu/tag/04TF}{Tag 04TF}]{stacks-project}.
(Note that, as stated, that result actually deals with stacks in the {\em fppf} topology;
here, as explained in Section~\ref{subsec:notation}, we are applying the
analogous result for the \'etale topology.) An examination of the proof of that result also gives the claim for categories
fibred in groupoids satisfying~[3]. The claim for stacks then follows from~\cite[\href{http://stacks.math.columbia.edu/tag/02ZL}{Tag 02ZL}]{stacks-project}.
\end{proof}

Our next goal in this section is to extend the definition of
certain properties of morphisms of algebraic stacks to morphisms
of stacks $\cX \to \cF$ whose source is assumed algebraic, 
but whose target is merely assumed to satisfy~[3].
To this end, we first note the following lemma.

\begin{lemma}
\label{lem:base-changing morphisms}
Let $\cX \to \cY \to \cF$ be morphisms of stacks, with $\cX$ and $\cY$
algebraic stacks, and with $\cF$ assumed to satisfy~{\em [3]}.  If $P$
is a property of morphisms of algebraic stacks that is preserved under
arbitrary base-change \emph{(}by morphisms of algebraic stacks\emph{)},
 then the morphism
$\cX \to \cY$ has the property $P$ if and only if, 
for every morphism of stacks $\cZ \to \cF$ with $\cZ$ being algebraic,
the base-changed morphism
$\cX\times_{\cF} \cZ \to \cY \times_{\cF} \cZ$ has property $P$.
\end{lemma}
\begin{proof}
The indicated base-change can be rewritten as the
base-change of the morphism $\cX \to \cY$ via the morphism
$(\cY\times_{\cF} \cZ) \to \cY.$
Since $P$ is assumed to be preserved under arbitrary base-changes,
we see that if $\cX \to \cY$ has property $P$,
so does the base-change
$\cX\times_{\cF} \cZ \to \cY \times_{\cF} \cZ$.

Conversely, suppose that all such base-changes have
property $P$; then, in particular, the morphism 
$\cX \times_{\cF} \cY \to \cY\times_{\cF} \cY$ has property $P$.
Thus so does the pull-back of this morphism via the diagonal
$\Delta: \cY \to \cY\times_{\cF} \cY.$
This pull-back may be described by the usual ``graph'' Cartesian diagram
(letting $f$ denote the given morphism $\cX \to \cY$)
$$\xymatrix{\cX \ar[d]_-{\Gamma_f := \id_{\cX} \times f} \ar^-f[r] &
\cY \ar^-{\Delta}[d] \\
\cX \times_{\cF} \cY \ar[r]^-{f\times \id_{\cY} } &
\cY\times_{\cF} \cY }
$$
from which we deduce that the original morphism $f$
has property $P$.
\end{proof}

\begin{example}
\label{ex:diagonals}
If $\cX \to \cF$ is a morphism of stacks, with $\cX$ being algebraic
and $\cF$ satisfying~[3], then we may apply the preceding
lemma to the diagonal morphism $\Delta: \cX \to \cX \times_{\cF} \cX$.
(Note that the source of this morphism is an algebraic stack
by assumption, and the target is an algebraic stack 
by Lemma~\ref{lem:fibre products}.) If $\cZ \to \cF$ is 
any morphism of stacks with $\cZ$ being algebraic,
then (since the formation of diagonals is compatible with base-change),
the base-change of the diagonal may be naturally identified
with diagonal of the base-change $$\Delta: (\cX\times_{\cF} \cZ)
\to (\cX\times_{\cF} \cZ) \times_{\cZ} (\cX\times_{\cF} \cZ).$$
In particular, since the the properties of being representable by algebraic
spaces, and of being locally of finite type, are preserved under any base-change,
and hold for the diagonal of any morphism of algebraic
stacks~\cite[\href{http://stacks.math.columbia.edu/tag/04XS}{Tag 04XS}]{stacks-project}),
we see that $\Delta: \cX \to \cX\times_{\cF} \cX$ 
is representable by algebraic spaces, and is locally of finite type.
(Proposition~\ref{prop:morphisms between stacks with repble diagonals have
repble diagonals}
below generalises the first of these
statements to the case when $\cX$ is also assumed only to satisfy~[3].)
\end{example}

We now make the following definition.

\begin{df}
\label{def:properties defined by base-change}
Assume that $\cF$ is a stack satisfying Axiom~[3]. Given an algebraic stack $\cX$ and a morphism $\cX \to \cF$,
and a property $P$ of morphisms of algebraic stacks that is preserved
under arbitrary base-change, then
we say that 
$\cX \to \cF$ has property $P$,
if and only if,
for any algebraic stack $\cY$ equipped with a morphism $\cY \to \cF$,
the base-changed morphism $\cX \times_{\cF} \cY \to \cY$ (which by Lemma~\ref{lem:fibre products} is a morphism of algebraic stacks) has property $P$. 
\end{df}

\begin{remark}
	In the spirit
	of~\cite[\href{http://stacks.math.columbia.edu/tag/03YK}{Tag 03YK}]{stacks-project},
	it might be better to restrict this definition to properties
	that are furthermore {\em fppf} local on the target. 
	Since, in any case,
	we only apply the definition to such properties,
	we don't worry about this subtlety here.
\end{remark}

\begin{remark}\label{rem: various properties of morphisms of algebraic stacks}
We will apply the preceding definition to the following properties of morphisms
of algebraic stacks:\begin{itemize}
\item Locally of finite presentation (this property is smooth local on the
  source-and-target,
so is defined by~\cite[\href{http://stacks.math.columbia.edu/tag/06FN}{Tag 06FN}]{stacks-project}).
\item
 Locally of finite type (again, this is smooth local on the
  source-and-target, so is defined by~\cite[\href{http://stacks.math.columbia.edu/tag/06FN}{Tag 06FN}]{stacks-project}).
\item Quasi-compact, \cite[\href{http://stacks.math.columbia.edu/tag/050U}{Tag 050U}]{stacks-project}.
\item Finite type, which as usual is defined to be locally of finite type and
  quasi-compact, \cite[\href{http://stacks.math.columbia.edu/tag/06FS}{Tag 06FS}]{stacks-project}.
\item
 Universally closed, \cite[\href{http://stacks.math.columbia.edu/tag/0513}{Tag 0513}]{stacks-project}.

\item
 Surjective, \cite[\href{http://stacks.math.columbia.edu/tag/04ZS}{Tag 04ZS}]{stacks-project}.
\item Separated, i.e.\ having proper diagonal
\cite[\href{http://stacks.math.columbia.edu/tag/04YW}{Tag 04YW}]{stacks-project}.
\item Proper, which as usual we define to be separated, finite type, and universally
  closed.
\item Representable by algebraic spaces, which
is equivalent to the condition that the diagonal
morphism be a
monomorphism~\cite[\href{http://stacks.math.columbia.edu/tag/0AHJ}{Tag 0AHJ}]
{stacks-project}. 
\item Monomorphism
\cite[\href{http://stacks.math.columbia.edu/tag/04ZW}{Tag 04ZW}]{stacks-project},
which is equivalent to the condition that the diagonal be an isomorphism,
or that the morphism, thought of as a functor between categories fibred
in groupoids, is fully faithful
\cite[\href{http://stacks.math.columbia.edu/tag/04ZZ}{Tag 04ZZ}]{stacks-project}. 
Note that monomorphisms are necessarily representable by algebraic spaces.
\item
Closed immersion,
\cite[\href{http://stacks.math.columbia.edu/tag/04YL}{Tag 04YL}]{stacks-project}.
(Closed immersions also admit an alternative characterisation as being 
the proper monomorphisms. To see this, note that
by~\cite[\href{http://stacks.math.columbia.edu/tag/045F}{Tag
  045F}]{stacks-project} it is enough to prove the same statement for morphisms
of algebraic spaces. Since closed immersions of algebraic spaces are
representable by definition, and proper monomorphisms are representable
by~\cite[\href{http://stacks.math.columbia.edu/tag/0418}{Tag
  0418}]{stacks-project}, we reduce to the case of morphisms of schemes, which is~\cite[18.12.6]{MR0238860}.)
\item Unramified, which is defined to be
  locally of finite type, with \'etale diagonal \cite[Appendix B]{MR2818725}.
(Note that the diagonal morphism is representable by algebraic 
spaces~\cite[\href{http://stacks.math.columbia.edu/tag/04XS}{Tag 04XS}]
{stacks-project}, and so the notion of it being \'etale is defined
by~\cite[\href{http://stacks.math.columbia.edu/tag/04XB}{Tag
04XB}]{stacks-project}.  Note also that an unramified morphism
is representable by algebraic spaces if and only if its diagonal
is an open immersion, since open immersions of algebraic
spaces are precisely the \'etale monomorphisms; see~\cite[17.9.1]{MR0238860}
for this statement
in the context of morphisms of schemes, which implies the statement for
algebraic spaces, because open immersions of algebraic spaces are representable
by definition, and \'etale monomorphisms of algebraic spaces are representable
by~\cite[\href{http://stacks.math.columbia.edu/tag/0418}{Tag 0418}]{stacks-project}.)\end{itemize}

\end{remark}

\begin{remark}
Some of the preceding properties, when interpreted via
the mechanism of Definition~\ref{def:properties defined by base-change},
also admit a more direct interpretation.   In particular,
if $\cX \to \cF$ is a morphism of stacks with $\cX$ being algebraic
and $\cF$ satisfying~[3], then the diagonal morphism
$\Delta: \cX \to \cX \times_{\cF} \cX$ {\em is} a morphism
of algebraic stacks, which
is furthermore representable by algebraic spaces
(as noted in Example~\ref{ex:diagonals}),
and so we know what it means for it
to be proper, or \'etale (for
example).  The following lemma incorporates this, and some similar
observations, to give more direct interpretations of some
of the preceding properties.
\end{remark}

\begin{lemma}
\label{lem:characterizing morphisms}
Let $f:\cX \to \cF$ be a morphism of stacks, with $\cX$ being algebraic
and $\cF$ satisfying~{\em [3]}.
\begin{enumerate}
\item The morphism $f$ is quasi-compact, in the
sense of Definition~{\em \ref{def:properties defined by base-change}},
if and only if, for every morphism
$\cY \to \cF$ with $\cY$ a quasi-compact algebraic stack,
the algebraic stack $\cX \times_{\cF} \cY$ is quasi-compact.
\item The morphism $f$ is universally closed,
in the sense of Definition~{\em \ref{def:properties defined by base-change}},
if and only if
for every morphism $\cY \to \cF$ with $\cY$ an algebraic stack,
the induced morphism
$|\cX\times_{\cF} \cY| \to |\cY|$ is closed.
\item If $P$ is a property which is preserved under arbitrary base-change,
which is smooth local on the source-and-target, and which is fppf
local on the target,
then $f$ satisfies $P$,
in the sense of Definition~{\em \ref{def:properties defined by base-change}},
if and only if for some {\em (}or, equivalently,
any{\em )} smooth cover $U \to \cX$ of $\cX$ by a scheme,
the composite morphism $U \to \cF$ {\em (}which is representable by
algebraic spaces, since $\cF$ satisfies {\em [3])} satisfies
condition~$P$ in the sense
of~{\em \cite[\href{http://stacks.math.columbia.edu/tag/03YK}{Tag
03YK}]{stacks-project}}.\footnote{The assumption that $P$ be \emph{fppf}
local on the target is included purely in order for the definition
of~\cite[\href{http://stacks.math.columbia.edu/tag/03YK}{Tag
03YK}]{stacks-project} to apply to~$P$.}
\item Let $P'$ be a property of morphisms of algebraic stacks
that are representable by algebraic spaces. Assume further that $P'$ is preserved
under arbitrary base-change,
and let $P$ be the property of morphisms of algebraic stacks 
defined by
the requirement that the corresponding diagonal morphism should satisfy~$P'$.
Then $f$ satisfies $P$,
in the sense of Definition~{\em \ref{def:properties defined by base-change}},
if and only if the diagonal morphism
$\Delta_f: \cX \to \cX\times_{\cF} \cX$ satisfies~$P'$.
\item The morphism $f$ is representable by algebraic spaces,
in the sense of Definition~{\em \ref{def:properties defined by base-change}} and
Remark~{\em \ref{rem: various properties of morphisms of algebraic stacks}},
if and only if it is representable by algebraic spaces in the usual
sense, i.e.\ if and only if for any morphism $T \to \cF$ with
$T$ a scheme, the base-change $\cX\times_{\cF} T$ is an algebraic space;
and these conditions are equivalent in turn to the condition
that the diagonal $\Delta_f$ be a monomorphism.
\item If $f: \cX \to \cF$ is locally of finite type,
in the sense of Definition~{\em \ref{def:properties defined by base-change}},
then $\Delta_f:\cX \to \cX\times_{\cF} \cX$ is locally of
finite presentation.
\end{enumerate}
\end{lemma}
\begin{proof}
Suppose that $f$ is quasi-compact,
in the sense of Definition~\ref{def:properties defined by base-change}.
If $\cZ \to \cF$ is a morphism
of stacks, with $\cZ$ being algebraic and quasi-compact,
then by definition $\cX\times_{\cF} \cZ \to \cZ$
is a quasi-compact morphism. Since $\cZ$ is quasi-compact,
it again follows by definition that $\cX\times_{\cF} \cZ$
is a quasi-compact algebraic stack.  Conversely,
suppose that $\cX\times_{\cF} \cZ$ is quasi-compact,
for every morphism $\cZ \to \cF$ with $\cZ$ being a quasi-compact
algebraic stack.  Let $\cY \to \cF$ be any morphism of stacks
with $\cY$ algebraic, and let $\cZ \to \cY$ be a morphism
of algebraic stacks with $\cZ$ quasi-compact.  Then
the base-change
$$(\cX \times_{\cF} \cY)\times_{\cY} \cZ \to \cY\times_{\cY} \cZ \iso \cZ$$
may be naturally identified with the base-change $\cX\times_{\cF} \cZ
\to \cZ;$
in particular, $(\cX\times_{\cF} \cY) \times_{\cY} \cZ$ is quasi-compact.
The base-changed morphism $\cX\times_{\cF} \cY \to \cY$ is thus
quasi-compact by definition.  Since $\cY \to \cF$ was arbitrary,
we conclude that $f$ is quasi-compact,
in the sense of Definition~\ref{def:properties defined by base-change};
this proves~(1).

The proofs of~(2) and of the first claim of (5) proceed
along identical lines to the proof of~(1).
To prove~(3), suppose first that $f$ has property $P$,
in the sense of Definition~\ref{def:properties defined by base-change}.
If $T \to \cF$ is any morphism from a scheme to $\cF$,
and $U \to \cX$ is any smooth surjection,
then the base-changed morphism $U\times_{\cF} T \to \cX \times_{\cF} T$
(which is naturally identified with the base-change of the morphism
$\cU \to \cX$ by the morphism of algebraic stacks $\cX\times_{\cF} T \to \cX$)
is smooth, while the morphism $\cX \times_{\cF} T \to T$ satisfies
$P$ by assumption.  Thus the composite $U\times_{\cF} T \to T$
satisfies $P$ (as $P$ is assumed to be smooth local on the source-and-target),
and so the morphism $U \to \cF$ satisfies $P$ in the sense
of~\cite[\href{http://stacks.math.columbia.edu/tag/03YK}{Tag
03YK}]{stacks-project}.  
Conversely, if this latter condition holds, and if $T \to \cF$
is a morphism whose source is a scheme, then 
the morphism $U\times_{\cF} T \to T$ satisfies $P$.  Since
$U\times_{\cF} T \to \cX\times_{\cF} T$ is smooth, and $P$ is assumed
to be smooth local on the source-and-target, we find that  $\cX\times_{\cF} T
\to T$ satisfies $P$. 
Now suppose that $\cY \to \cF$ is any morphism of stacks whose source
is algebraic, and let $T \to \cY$ be a smooth surjection.
We have just seen that $\cX\times_{\cF} T \to T$ satisfies~$P$.
Since $P$ is smooth local on the source-and-target,
and since each of the morphisms $\cX\times_{\cF} T \to \cX\times_{\cF} \cY$
and $T \to \cY$ are smooth and surjective
(the latter by assumption and the former because it is naturally
identified with the base-change of the latter by the projection
$\cX\times_{\cF} \cY \to \cY$), we find that
$\cX\times_{\cF} \cY \to \cY$ satisfies $P$, as required.
Thus, by definition, the morphism $f$ satisfies~$P$,
in the sense of Definition~\ref{def:properties defined by base-change}.

To prove~(4), note that the formation of diagonals is compatible
with base-change, after which~(4) follows from
Lemma~\ref{lem:base-changing morphisms}.
The second claim of~(5) follows from~(4), applied in the 
case when $P'$ is the property of being a monomorphism
(since~\cite[\href{http://stacks.math.columbia.edu/tag/0AHJ}{Tag 0AHJ}]
{stacks-project} shows that the associated property $P$ 
is then precisely the property of being representable by algebraic
spaces). 

Claim~(6) also follows from
Lemma~\ref{lem:base-changing morphisms},
and the fact that the analogous claim holds 
for morphisms of algebraic stacks. (For lack of a reference, we now give a proof
of this property; that is, we show that if   $f:\cX\to\cY$ is a morphism of
algebraic stacks  which is
locally of finite type, then $\Delta_f$ is locally of finite presentation. Choose a smooth surjection $V\to\cY$ from a scheme
$V$. By~\cite[\href{http://stacks.math.columbia.edu/tag/06Q8}{Tag
  06Q8}]{stacks-project}, it is enough to prove that the base change of
$\Delta_f$ by $V\to\cY$ is locally of finite presentation. This base change is
the diagonal of the base change $\cX\times_\cY V\to V$, so we may
replace $\cX$ by $\cX\times_\cY V$ and $\cY$ by $V$, and therefore reduce to the
case of a morphism
$\cX\to V$.

Now choose a smooth surjection from a scheme
$U\to\cX$. By~\cite[\href{http://stacks.math.columbia.edu/tag/06Q9}{Tag
  06Q9}]{stacks-project} it suffices to show that the composite
$U\to\cX\to\cX\times_V\cX$ is locally of finite presentation. Factoring this
as the composite \[U\to U\times_VU\to U\times_V\cX\to \cX\times_V\cX\]we see that
it suffices to show that $U\to U\times_V U$ is locally of finite presentation;
that is, we have reduced to the case of schemes, which is~\cite[\href{http://stacks.math.columbia.edu/tag/0818}{Tag 0818}]{stacks-project}.)
\end{proof}

\begin{example}
\label{ex:proper morphisms}
Lemma~\ref{lem:characterizing morphisms} shows that a morphism of stacks
$\cX \to \cF$, with $\cX$ being algebraic and $\cF$ 
satisfying~[3], 
is proper,
in the sense of Definition~\ref{def:properties defined by base-change},
if and only if the following 
properties hold:
(i) for some (or, equivalently, any)
smooth surjection $U \to \cX$, with $U$ a scheme, 
the composite morphism $U \to \cF$ (which is a morphism
of stacks representable by algebraic spaces) is locally
of finite type in the sense
of~\cite[\href{http://stacks.math.columbia.edu/tag/03YK}{Tag
03YK}]{stacks-project};  
(ii) for any morphism of stacks $\cY \to \cF$ with $\cY$ algebraic,
the base-changed
morphism $\cX\times_{\cF} \cY \to \cY$ induces
a closed morphism
$|\cX\times_{\cF} \cY| \to |\cY|$;
 (iii) in the context of~(ii),
if $\cY$ is furthermore quasi-compact, then the base-changed
algebraic stack $\cX\times_{\cF} \cY$ is quasi-compact;
(iv) the diagonal morphism $\Delta: \cX \to \cX\times_{\cF} \cX$
is proper.
\end{example}

\begin{example}
\label{ex:unramified morphisms}
Lemma~\ref{lem:characterizing morphisms} shows that a morphism of stacks
$\cX \to \cF$, with $\cX$ being algebraic and $\cF$ 
satisfying~[3], 
is unramified,
in the sense of Definition~\ref{def:properties defined by base-change},
if and only if the following 
properties hold:
(i) for some (or, equivalently, any)
smooth surjection $U \to \cX$, with $U$ a scheme, 
the composite morphism $U \to \cF$ (which is a morphism
of stacks representable by algebraic spaces) is locally
of finite type in the sense
of~\cite[\href{http://stacks.math.columbia.edu/tag/03YK}{Tag
03YK}]{stacks-project};  
(ii) the diagonal morphism $\Delta: \cX \to \cX\times_{\cF} \cX$
(which is {\em a priori} locally of finite presentation,
by~(i) and Lemma~\ref{lem:characterizing morphisms}~(6)) is \'etale.
The unramified morphism $\cX \to \cF$ is furthermore 
representable by algebraic spaces if the diagonal
is in fact an open immersion (taking into account~(5) of
the preceding lemma, and the fact that open immersions
are precisely the \'etale monomorphisms).
\end{example}

\begin{example}
Lemma~\ref{lem:characterizing morphisms}~(4) shows 
that a morphism of stacks
$f: \cX \to \cF$, with $\cX$ being algebraic and $\cF$ satisfying~[3], 
is a monomorphism
in the sense of Definition~\ref{def:properties defined by base-change},
if and only if the diagonal
$\Delta_f: \cX \to \cX\times_{\cF} \cX$ is an isomorphism.
In turn, this is equivalent to the condition that $f$,
thought of as a functor between categories fibred in groupoids,
is fully faithful. 
Lemma~\ref{lem:characterizing morphisms}~(5) shows 
that $f$ is then in particular necessarily
representable by algebraic spaces in the usual sense. 
\end{example}

\begin{df}\label{defn:closed substack} We say that a substack $\cF'$ of $\cF$ is a {\em closed} substack if the inclusion morphism $\cF' \hookrightarrow \cF$ is representable by
algebraic spaces, and is a closed immersion in the sense
of~\cite[\href{http://stacks.math.columbia.edu/tag/03YK}{Tag
 03YK}]{stacks-project}, i.e.\ has the property that for any
morphism $T \to \cF$ with $T$ an affine scheme,
the base-changed morphism $T\times_{\cF} \cF' \to T $ is a closed immersion
of schemes. 
\end{df}

\begin{lemma}
\label{lem:closed immersion}
Assume that $S$ is locally Noetherian, and
that $\cF$ is a stack over $S$ which satisfies~{\em [3]},
for which the morphism $\cF \to S$ is limit preserving on objects. If $\cZ \hookrightarrow \cF$ is a closed immersion, then it is limit preserving on objects.
\end{lemma}
\begin{proof}
By definition of what it means for $\cZ \to \cF$ to be a closed immersion,
if $T \to \cF$ is a morphism from an algebraic space, then the induced
morphism $\cZ \times_{\cF} T \to T$ is a closed immersion.  Since $\cF$ satisfies~[3],
the source is an algebraic stack; but if $T$ is a scheme, then in fact the source 
will be a scheme (being a closed substack of a scheme). In particular the morphism $\cZ \to \cF$ is representable by algebraic spaces,
and so to check that it is limit preserving on objects, it suffices to check
that it is locally of
finite presentation~\cite[\href{http://stacks.math.columbia.edu/tag/06CX}{Tag 06CX}]{stacks-project}.
For this, it suffices to check, in the preceding context, 
that if $T$ is an affine scheme then $\cZ \times_{\cF} T \to T$ is locally of finite presentation.

If we write $T = \varprojlim_i T_i$ as the
projective limit of finite type affine $S$-schemes~$T_i$,
then, since $\cF \to S$ is limit preserving on objects and $S$ is locally Noetherian,
we may factor the morphism $T \to \cF$ through one of the $T_i$,
and hence reduce to the case when $T$ is finite type over the locally Noetherian
scheme $S$.  But in this case $T$ itself is Noetherian, and so the
closed immersion $\cZ \times_{\cF} T \to T$ is in fact of finite presentation.
This proves the lemma.
\end{proof}

Our next goal is to state and prove a lemma which is a variant
of~\cite[\href{http://stacks.math.columbia.edu/tag/06CX}{Tag
06CX}]{stacks-project}
(which states that for a morphism
between categories fibred in groupoids that is representable
by algebraic spaces, being locally of finite presentation is equivalent
to being limit preserving on objects).
Before doing this, we recall some
standard facts about descending finitely presented morphisms of affine schemes
over $S$ to morphisms of finitely presented affine schemes over $S$.

\begin{lem}
  \label{lem:descending etale to finite presented affine}Suppose that $T,T'$ are
  affine schemes over~$S$, and that we are given a morphism of finite
  presentation $T'\to T$. Suppose that $T$ may be written as a limit
  $\varprojlim T_i$ of affine schemes of finite presentation over $S$. Suppose also that we are given a morphism $T'\to T''$
  with $T''$ a scheme locally of finite presentation over $S$.

Then for some $j$ we may find a factorisation of the given morphism $T'\to T''$
which fits into a commutative
diagram \[\xymatrix{T'\ar[d]\ar[r]&T'_j\ar[d]\ar[r]&T''\\ T\ar[r]&T_j&}\]in
which the square is Cartesian, $T'_j$ is an affine scheme of finite presentation
over~$S$, and the morphism $T'_j\to T_j$ is of finite presentation. Furthermore,
if $T'\to T$
is surjective, then $T'_j\to T_j$ can be taken to be surjective; and
if $T'\to T$ is assumed to be \'etale {\em (}resp.\ fppf, resp.\ an
open immersion{\em )}, then $T'_j\to T_j$ can also
be taken to be \'etale {\em (}resp.\ fppf, resp.\ an open immersion{\em )}.
\end{lem}
\begin{rem}
  \label{rem:if S quasisep, T is automatically a limit}If $S$ is
  quasi-separated, then the hypothesis that $T$ may be written as a limit of
  affine schemes of finite presentation is automatically satisfied; see
  Theorem~\ref{thm: Thomason and Trobaugh, affine is pro-fp-affine} below.
\end{rem}
\begin{proof}[Proof of Lemma~\ref{lem:descending etale to finite presented affine}]
  By~\cite[\href{http://stacks.math.columbia.edu/tag/01ZM}{Tag
    01ZM}]{stacks-project}, there is some $j_0$ and a morphism of finite
  presentation $T'_{j_0}\to T_{j_0}$ such that the pull-back of this morphism to
  $T$ is the given morphism $T'\to T$. For $j\ge j_0$, set
  $T'_j:=T_j\times_{T_{j_0}}T'_{j_0}$. Since products commute with limits, we
  have $T'=\varprojlim T'_j$.

Since $T_{j_0}$ is affine, and $T'_{j_0}\to T_{j_0}$ is of finite presentation
(hence quasi-compact and quasi-separated),
it follows that $T'_{j_0}$ is quasi-compact and quasi-separated. It then follows
from~\cite[\href{http://stacks.math.columbia.edu/tag/01ZN}{Tag
  01ZN}]{stacks-project} that we may assume that $T'_{j}$ is affine for all
$j\ge j_0$. If $T'\to
T$ is surjective (resp.\ \'etale, resp.\ \emph{fppf}, resp.\ an open immersion), then by~\cite[\href{http://stacks.math.columbia.edu/tag/07RR}{Tag 07RR}]{stacks-project} (resp.\ \cite[\href{http://stacks.math.columbia.edu/tag/07RP}{Tag
  07RP}]{stacks-project}, resp.\
\cite[\href{http://stacks.math.columbia.edu/tag/04AI}{Tag
  04AI},\href{http://stacks.math.columbia.edu/tag/07RR}{Tag
  07RR}]{stacks-project}, resp.\ \cite[\href{http://stacks.math.columbia.edu/tag/07RP}{Tag
  07RP},\href{http://stacks.math.columbia.edu/tag/07RQ}{Tag
  07RQ}]{stacks-project} (recall that open immersions are precisely
the \'etale monomorphisms)) we may assume that $T'_{j}\to T_{j}$ is
also surjective  (resp.\ \'etale, resp.\ \emph{fppf}, resp.\ an open immersion) for all $j\ge j_0$.

Since $T''$ is locally of finite presentation over $S$, it follows
from~\cite[\href{http://stacks.math.columbia.edu/tag/01ZC}{Tag
  01ZC}]{stacks-project} that $T'\to T''$ factors through $T'_j$ for some
$j\ge j_0$, as required.\end{proof}

\begin{lem}
\label{lem:lfp equals limit preserving on objects}
A morphism $\cX \to \cF$, where $\cX$ is an algebraic stack and
$\cF$ is a stack satisfying~{\em [3]},
is locally of finite presentation,
in the sense of Definition~{\em \ref{def:properties defined by base-change}},
if and only if it is limit preserving on objects.
\end{lem}
\begin{proof}
Let $U \to \cX$ be a smooth surjection from a scheme (which exists,
since $\cX$ is an algebraic stack).  This morphism is representable
by algebraic spaces (again, since $\cX$ is an algebraic stack),
and is locally of finite presentation (since it is smooth), and hence
it is limit preserving on objects by~\cite[\href{http://stacks.math.columbia.edu/tag/06CX}{Tag
06CX}]{stacks-project}.  Thus if the morphism $\cX \to \cF$
is limit preserving on objects, so is the composite
$U \to \cF$~\cite[\href{http://stacks.math.columbia.edu/tag/06CW}{Tag 06CW}]{stacks-project}.
Since $\cF$ satisfies~[3], this morphism is representable by algebraic
spaces, and so we may apply~\cite[\href{http://stacks.math.columbia.edu/tag/06CX}{Tag 06CX}]{stacks-project}
again to deduce that $U \to \cF$ is locally of finite presentation.
Hence the same is true of the morphism $\cX \to \cF$,
by part~(3) of Lemma~\ref{lem:characterizing morphisms}.

Conversely, suppose that $\cX \to \cF$
is locally of finite presentation, 
in the sense of Definition~\ref{def:properties defined by base-change}.
Suppose given a morphism $T \to \cX$,
for an affine $S$-scheme $T$,
written as a projective limit $T  =\varprojlim_i T_i$ of affine $S$-schemes $T_i$,
and suppose further that the composite $T \to \cX \to \cF$ factors through
one of the $T_i$.
We must show that  our given morphism $T \to \cX$ factors in a compatible manner
through $T_{i'} \to \cX$ for some sufficiently large $i' \geq i$. 
Replacing $\cX$ with $\cX\times_{\cF} T_i$ (which is an algebraic stack,
since $\cF$ satisfies~[3]) and $\cF$ by $T_i$, we may in fact assume
that $\cF = T_i$,
which we do from now on.

Let $U \to \cX$ be a
smooth surjection,
and write $R := U \times_{\cX} U$,
so that there is a natural isomorphism $[U/R] \iso \cX$.
The assumption that $\cX \to T_i$ is locally of finite presentation
is by definition equivalent to supposing that $U \to T_i$ is locally of finite presentation.

Let $p_1,p_2: R \rightrightarrows U$ denote the two projections. 
Consider the pull-backs $U_T$ 
and $R_T$.  There is a natural identification $R_T \iso U_T \times_T U_T,$
via which the base-changes of the natural projections $p_i$ become identified
with the natural projections $U_T \times_T U_T \rightrightarrows U_T.$
There is a natural isomorphism of stacks $[U_T/R_T] \iso T.$

We may find (for example by~\cite[\href{http://stacks.math.columbia.edu/tag/055V}{Tag 055V}]{stacks-project}) an \'etale slice $T'$ of $U_T$, i.e.\ a morphism $T' \to U_T$ for which
the composite $T' \to U_T \to T$ is \'etale and surjective. We then define $T'' := T'\times_{T} T'.$
We have the composite morphisms $T' \to U_T \to U$
and $T'' \to R_T \to R,$ with respect to which the two projections $T'' \rightrightarrows
T'$ are compatible with the two projections $p_i: R \to U.$
Thus there is an induced morphism $[T'/T''] \to [U/R],$ which is naturally identified
with our original morphism $T \to \cX$.

Since \'etale morphisms are open,
and since $T$ is quasi-compact (being affine), we may replace $T'$ by a quasi-compact
open subscheme for which the induced morphism to $T$ remains surjective.  Finally,
replacing $T'$ by the disjoint union of the members of a finite affine open cover,
we may in fact assume that $T'$ is also affine.  The morphism $T' \to T$ is then
an affine morphism that is locally of finite presentation (being \'etale),
and hence is actually of finite presentation.   
By Lemma~\ref{lem:descending etale to finite presented affine} we may write the affine \'etale morphism $T' \to T$
as the projective limit of affine \'etale morphisms $T_{i'}' \to T_{i'}$ (starting from
a sufficiently large value of~$i'$).   We write $T_{i'}'' := T_{i'}'\times_{T_{i'}} T_{i'}'.$

Since $U \to T_i$ is locally of finitely presentation, 
as is $R \to T_i$, by 
applying~\cite[\href{http://stacks.math.columbia.edu/tag/01ZM}{Tag
01ZM}]{stacks-project} and~\cite[\href{http://stacks.math.columbia.edu/tag/07SJ}{Tag 07SJ}]{stacks-project} we find that we may 
factor the morphisms $T' \to U$ and $T'' \to R$ through
$T'_{i'}$ and $T''_{i'}$, for some sufficiently large value of $i'$,
in such a manner
that the projections $T_{i'}'' \rightrightarrows T_{i'}'$
are compatible with the projections $R \rightrightarrows U.$
Thus we obtain a morphism $[T_{i'}'/T_{i'}''] \to [U/R],$
which we may identify with the required morphism $T_{i'} \to \cX$.
\end{proof}

We close this section with some further propositions, 
the first of which generalises one of the observations
of Example~\ref{ex:diagonals} to the case when both source
and target are assumed merely to satisfy~[3].

\begin{prop}
  \label{prop:morphisms between stacks with repble diagonals have repble
    diagonals}
Let $f:\cX\to\cY$ be a morphism between categories fibred in groupoids
  satisfying~\emph{[3]}. Then $\Delta_f:\cX\to\cX\times_{\cY}\cX$ is representable by algebraic spaces.
\end{prop}
\begin{proof}
  The proof is essentially identical to that
  of~\cite[\href{http://stacks.math.columbia.edu/tag/04XS}{Tag
    04XS}]{stacks-project}. Let $T\to\cX\times_{\cY}\cX$ be a morphism from
  a scheme $T$; by definition, this is the data of a triple $(x,x',\alpha)$
  where $x,x'$ are objects of $\cX$ over $T$, and $\alpha:f(x)\to f(x')$ is a
  morphism in the fibre category of $\cY$ over $T$. Since $\cX$, $\cY$
  satisfy~[3], the sheaves $\Isom_{\cX}(x,x')$ and $\Isom_{\cY}(f(x),f(y))$ are
  algebraic spaces over~$T$
  by~\cite[\href{http://stacks.math.columbia.edu/tag/045G}{Tag
    045G}]{stacks-project}. The morphism $\alpha$ corresponds to a section of
  the morphism $\Isom_{\cY}(f(x),f(x'))\to T$.

If $T'\to T$ is a morphism of schemes, then we see that a $T'$-valued point of
$\cX\times_{\cX\times_{\cY}\cX}T$ is by definition an isomorphism $x|_{T'}\isoto
x'|_{T'}$ whose image under $f$ is $\alpha|_{T'}$. Putting this together, we see
that there is a fibre product diagram of sheaves over $T$ \numequation\label{eqn: fibre
  product diagram for diagonals}
\xymatrix{\cX\times_{\cX\times_{\cY}\cX}T\ar[d]\ar[r]&\Isom_{\cX}(x,x')\ar[d]\\
  T\ar^-{\alpha}[r]&\Isom_{\cY}(f(x),f(x'))}\end{equation}Thus
$\cX\times_{\cX\times_{\cY}\cX}T$ is an algebraic space, as required.
\end{proof}

\begin{prop}\label{prop: double diagonals of stacks satisfying [3]}If $\cX$ is a
  category fibred in groupoids satisfying~{\em [3]}, and $\Delta:\cX\to\cX\times_S\cX$ is limit preserving on
  objects, then $\Delta_\Delta:\cX\to\cX\times_{\cX\times_S\cX}\cX$
  is representable by algebraic spaces and limit preserving on objects.
\end{prop}
\begin{proof}By Proposition~\ref{prop:morphisms between stacks with repble diagonals have repble
    diagonals}, $\Delta_\Delta$ is representable by algebraic spaces.
    (It follows from Lemma~\ref{lem:fibre products} that $\cX \times_S \cX$ satisfies~[3],
    since $\cX$ does, so that the proposition does indeed apply.)
    By~\cite[\href{http://stacks.math.columbia.edu/tag/06CX}{Tag
    06CX}]{stacks-project}, we see that $\Delta$ is locally of finite
  presentation, and that to show that $\Delta_\Delta$ is limit preserving on
  objects, it is equivalent to show that it is locally of finite presentation. Consider
  the diagram~(\ref{eqn: fibre
    product diagram for diagonals}) in the case that $f$ is
  $\Delta$; we must show that for all choices of $T$, the left hand vertical arrow
  is locally of finite presentation. 

  It therefore suffices to show that the right hand vertical arrow is locally of
  finite presentation; but this is by definition the diagonal
  map \[\Isom(x,x')\to\Isom(x,x')\times_S \Isom(x,x'),\]which is locally of
  finite presentation
  by~\cite[\href{http://stacks.math.columbia.edu/tag/084P}{Tag
    084P}]{stacks-project} (note that since $\Delta$ is locally of finite presentation, so
  is $\Isom(x,x')$).
\end{proof}

\begin{prop}
  \label{prop: limit preserving and base change}Let $\cX\to\cY$ be a morphism of
  categories fibred in groupoids which is representable by algebraic
  spaces. Assume that~$\cY$ satisfies~\emph{[3]}. Then:
  \begin{enumerate}
  \item $\cX$ satisfies~\emph{[3]}.
  \item If $\cX\to\cY$ and $\Delta_{\cY}:\cY\to\cY\times_S\cY$ are locally of finite
    presentation, then so is $\Delta_{\cX}:\cX\to\cX\times_S\cX$. \item If $\cX\to\cY$ is locally of finite presentation, and $\cY$ satisfies~\emph{[1]}, then so does~$\cX$.
  \end{enumerate}

\end{prop}
\begin{proof} Throughout the proof, we will freely make use of~\cite[\href{http://stacks.math.columbia.edu/tag/06CX}{Tag
    06CX}]{stacks-project},  the
  equivalence of being locally of finite presentation and being limit
  preserving on objects for morphisms representable by algebraic spaces.
 We begin with~(1); we need to show that the morphism $\cX\to \cX\times_S\cX$ is representable by algebraic
  spaces. This may be factored as \[\cX\to\cX\times_\cY\cX\to\cX\times_S\cX;\]
  since the second morphism is a base change of $\cY\to\cY\times_S\cY$, it is
  enough to show that $\cX\to\cX\times_\cY\cX$ is representable by algebraic
  spaces.

To this end, suppose that $T\to\cX\times_{\cY}\cX$ is a morphism whose
source is an algebraic space. This induces a morphism $T\to\cY$, and
we may consider the base change of our whole situation by this
morphism. Writing $\cX_T$ for $\cX\times_{\cY}T$, we may reinterpret
the given morphism $T \to \cX \times_{\cY} \cX$ as a section of the morphism
$\cX_T \times_T \cX_T \to T,$
and it is then
straightforward to see that we in fact
have \numequation\label{eqn:base change by T}\cX\times_{\cX\times_{\cY}\cX} T=
\cX_T\times_{\cX_T\times_T\cX_T} T.\end{equation} Now, since $\cX\to\cY$ is
representable by algebraic spaces, $\cX_T$ is an algebraic space, so
that $\cX_T\times_{\cX_T\times_T\cX_T} T$ is an algebraic space, as required.

 We now consider~(2). We may factor $\Delta_{\cX}$ as the
composite\[\cX\to\cX\times_{\cY}\cX\to\cX\times_S\cX.\]Since the second arrow is
a base change of $\Delta_{\cY}$, which is locally of finite
presentation  by assumption, it follows from
\cite[\href{http://stacks.math.columbia.edu/tag/06CV}{Tag
  06CV},\href{http://stacks.math.columbia.edu/tag/06CW}{Tag
  06CW}]{stacks-project} that it suffices to show that
$\cX\to\cX\times_{\cY}\cX$ is locally of finite presentation. By definition,
we need to show that all the base changes of this morphism by
morphisms $T\to\cX\times_{\cY}\cX$ with source an algebraic space are
locally of finite presentation;
examining~(\ref{eqn:base change by T}), we see that we are reduced to
the case that $\cY=T$ is an algebraic space, which is a special case
of Lemma~\ref{lem:characterizing morphisms}~(6).

To prove (3), by Lemma~\ref{lem:limit preserving vs. limit preserving on
  objects} we need to show that each of $\cX\to S$,
$\Delta_{\cX}:\cX\to\cX\times_S\cX$,
$\Delta_{\Delta}:\cX\to\cX\times_{\cX\times_S\cX}\cX$ is limit preserving on
objects; the same result also shows that both $\cY\to S$ and $\Delta_{\cY}$ are limit preserving on
objects, so that $\Delta_{\cX}$ is limit preserving on objects by (2). By Proposition~\ref{prop: double diagonals of stacks satisfying [3]},
it is then enough to show that $\cX\to S$ is limit preserving on objects; but
this is immediate from~\cite[\href{http://stacks.math.columbia.edu/tag/06CW}{Tag
  06CW}]{stacks-project}, applied to the composite $\cX\to\cY\to
S$.
\end{proof}

\begin{lem}
  \label{lem:3.2 running assumptions descend to monos}If $\cF'\to\cF$
  is a monomorphism of categories fibred in groupoids, and if~
  $\Delta_{\cF}$ is representable by algebraic spaces and locally of
  finite presentation, then~$\Delta_{\cF'}$ is also representable by
  algebraic spaces and locally of finite presentation.
\end{lem}
\begin{proof}
  Since $\cF'\to\cF$ is a monomorphism, $\cF'\to\cF'\times_\cF\cF'$ is
  an equivalence, and therefore~$\Delta_{\cF'}$ is the base change
  of~$\Delta_{\cF}$ via the natural morphism
  $\cF'\times_S\cF'\to\cF\times_S\cF$. The result follows.
\end{proof}
\subsection{Remarks on Axiom~[4]}\label{subsec:axiom 4}
We begin with some preliminary definitions and results related
to the concept of smoothness of morphisms.

\begin{defn}
  \label{defn:etale}As in Section~\ref{subsec:axiom 3}, we will  use the notions of {\em unramified} and {\em \'etale}  for morphisms
between algebraic stacks that are not necessarily representable
by algebraic spaces. We say that a morphism of algebraic stacks is unramified if it is locally of
finite type and has \'etale diagonal, and that it is \'etale if it is unramified, flat, 
and locally of finite presentation \cite[Appendix B]{MR2818725}.
Note that, although the definition of unramified morphism
includes the condition that the diagonal be \'etale,
there is no circularity, because diagonal morphisms are representable
by  algebraic spaces, and in this case the notation of
\'etale is defined
following~\cite[\href{http://stacks.math.columbia.edu/tag/04XB}{Tag
04XB}]{stacks-project}.  By~\cite[Prop.\ B2]{MR2818725} a morphism of algebraic
stacks is \'etale if and only if it is smooth and unramified.
\end{defn}

\begin{defn}
  \label{defn:formally smooth etc for categories in groupoids, or at least
    stacks}Let $f:\cX\to\cY$ be a morphism of stacks. Then we say that $f$ is
  \emph{formally unramified} (resp.\ \emph{formally smooth}, resp.\
  \emph{formally \'etale}) if for every affine $\cY$-scheme $T$, and every closed
  subscheme $T_0\into T$ defined by a nilpotent ideal sheaf, the functor
  $\Hom_\cY(T,\cX)\to\Hom_\cY(T_0,\cX)$ is fully faithful (resp.\ essentially
  surjective, resp.\ an equivalence of categories).
\end{defn}
\begin{prop}
  \label{prop:relationship of smooth to formally smooth etc}
  \begin{enumerate}
  \item A morphism of algebraic stacks is smooth if and only if it is formally smooth and
    locally of finite presentation.
  \item A morphism of algebraic stacks is unramified if and only if it is
    formally unramified and
    locally of finite type.
  \item A morphism of algebraic stacks is \'etale if and only if it is formally \'etale and
    locally of finite presentation.
  \end{enumerate}

\end{prop}
\begin{proof}
  See~\cite[Cor.\ B.9]{MR2818725}.
\end{proof}
\begin{defn}
\label{def:versions of smoothness}
Assume that $S$ is locally Noetherian, let $\cF$ be a category fibred
in groupoids over $S$, let $U$ be a scheme locally of finite type over $S$ equipped with a morphism
$\varphi:U \to \cF$, and let $u\in U$ be a point. \begin{enumerate}
\item
We say that $\varphi$ is {\em versal at} $u$ if for any diagram
$$\xymatrix{Z_0 \ar[d] \ar[r] & U \ar^-{\varphi}[d] \\
Z \ar[r] & \cF,}$$
where $Z_0$ and $Z$ are Artinian local schemes with 
the latter being a nilpotent thickening of the former,
and where the closed point $z\in Z_0$ maps to $u$,
inducing an isomorphism
$\kappa(u) \cong \kappa(z)$, 
we may lift the morphism $Z \to \cF$ to a morphism $Z \to U$.
\item We say that $\varphi$ is {\em formally smooth at} $u$
if for any diagram as in~(1) 
where $Z_0$ and $Z$ are Artinian local schemes with 
the latter being a nilpotent thickening of the former,
and where the closed point $z\in Z_0$ maps to $u$,
inducing a finite extension
$\kappa(u) \hookrightarrow \kappa(z)$, 
we may lift the morphism $Z \to \cF$ to a morphism $Z \to U$.
\item If $\varphi$ is representable by algebraic spaces and locally
of finite presentation,
then we say that $\varphi$ is {\em smooth at} $u$
if for any finite type $S$-scheme $X$ equipped with a morphism
$X \to \cF$ over $S$, there is an open neighbourhood $U'$ of $u$
in $U$ such the base-change morphism $U'\times_{\cF} X \to X$
is smooth.
\item If $\varphi$ is is representable by algebraic spaces
and locally of finite presentation,
then we say that $\varphi$ is {\em smooth in a neighbourhood of} $u$
if there exists a neighbourhood $U'$ of $u$ such that
the restriction $\varphi_{| U'} : U' \to \cF$ (which is again
a morphism representable by algebraic spaces) is smooth.
\end{enumerate}
\end{defn}
\begin{remark}
If $u \in U$ is a finite type point,
then $\varphi$ is {\em versal at} $u$ in the sense of~(1) of the preceding
definition if and only if the induced morphism
$\Spf \widehat{\mathcal O}_{U,u} \to \cF_{\varphi(u)}$
is versal in the sense discussed in Subsection~\ref{subsec:axiom 2} above.
Note that
we follow~\cite[\href{http://stacks.math.columbia.edu/tag/07XF}{Tag 07XF}]{stacks-project}
in using the terminology {\em versal at} $u$,
rather than {\em formally versal at} $u$,
as some other sources (e.g.\ \cite{hall2013artin}) do.

Our definition of formal smoothness at a point follows that of
\cite[Def.~2.1]{hall2013artin}.  
In~\cite[Def.~3.1]{MR0260746},
Artin defines the notion of formally \'etale at a point (for a morphism
from a scheme to a functor), but his definition is not quite the
obvious analogue of the definition of formal smoothness given here,
in that he does not impose any condition on the degree of the
extension of the residue field at the closed point of $Z_0$ 
over the residue field at $u$.
The condition on the residue field that we impose 
(following \cite{hall2013artin}) makes it slightly easier to
verify formal smoothness (see in particular Lemma~\ref{lem:versions of
  smoothness}~(2) below), and is harmless in practice.  (Indeed,
if $\cF$ is a stack satisfying~[1] and~[3],
then part~(3) of Lemma~\ref{lem:versions of smoothness} below shows that at finite type points
the variant definition of formal smoothness, in which
we impose no condition on the extension of residue fields,
is equivalent to the definition given above.)
\end{remark}

\begin{remark}
The notion of smoothness is defined
for morphisms between algebraic stacks,
or, more generally (via Definition~\ref{def:properties defined by base-change}),
for morphisms from an algebraic stack to a stack satisfying~[3].
Thus the notion of {\em smooth at a point} can naturally be extended
to morphisms whose source is an algebraic stack,
and whose target satisfies~[3].

The notion of {\em versality at a point} or {\em formal smoothness
at a point} of an algebraic
stack is slightly more problematic to define, 
since a point of an algebraic stack is defined as an equivalence
class of morphisms from the spectrum of a field, 
and so we can't speak of the residue field at a point of a stack.
In Definition~\ref{def:formal smoothness at a point for stacks},
we {\em will} give a definition of formal smoothness at a point for 
a morphism whose source is an algebraic stack, under slightly
restrictive conditions on the target of the morphism
(which, however, will not be too restrictive for the applications
of this notion that we have in mind).   The key to making
the definition work is Corollary~\ref{cor:formal smoothness is smooth
local}, which shows (under suitable hypotheses) that formal
smoothness at a point (for morphisms from a scheme)
can be detected smooth locally.
\end{remark}

The following lemma,
which is essentially drawn from \cite[\S 2]{hall2013artin},
relates the various notions
of Definition~\ref{def:versions of smoothness}.
We remark that
part~(3) of the lemma provides an analogue, for formal smoothness at a point,
of Artin's~\cite[Lem.~3.3]{MR0260746},
which provides a characterisation of morphisms from a 
scheme to a functor that are formally \'etale at a point.

\begin{lemma}
\label{lem:versions of smoothness}
Suppose that we are in the context of
Definition~{\em \ref{def:versions of smoothness}}, and that $u$ is a finite type
point.
\begin{enumerate}
\item In general, we have that {\em (4)} $\implies$~{\em (3)}
$\implies$~{\em (2)} $\implies$~{\em (1)}.
\item If $\cF$ satisfies ~{\em (RS)},
then~{\em (1)} $\Longleftrightarrow$~{\em (2)}; i.e.\
$\varphi$ is versal at $u$ if and only if it is formally smooth at $u$.
\item  Consider the following conditions:
\begin{enumerate}[(a)]
\item  $\varphi$ is formally smooth at $u$.
\item For any finite type $S$-scheme $X$,
and any morphism $X \to \cF$ over~$S$,
the base-changed morphism $U \times_{\cF} X \to X$ contains
the fibre over $u$ in its smooth locus.
\item For any locally finite type algebraic stack $\cX$ over $S$,
and any morphism $\cX \to \cF$ over~$S$,
the base-changed morphism $U \times_{\cF} \cX \to \cX$ contains
the fibre over $u$ in its smooth locus.
\item For any algebraic stack $\cX$ over $S$,
and any morphism $\cX \to \cF$ over~$S$,
the base-changed morphism $U \times_{\cF} \cX \to \cX$ contains
the fibre over $u$ in its smooth locus.
\end{enumerate}
If $\varphi$ is representable by algebraic spaces and locally
of finite type, then conditions (a), (b), and (c) are equivalent.  If
furthermore $\cF$ is a stack satisfying~{\em [1]} and~{\em [3]},
then all four of these conditions are equivalent.
\item 
If $\cF$ is an algebraic stack which is locally of finite presentation
over $S$, then conditions~{\em (1)}, {\em (2)}, {\em (3)}, and~{\em (4)}
of Definition~{\em \ref{def:versions of smoothness}} are
equivalent.  \end{enumerate}
\end{lemma}
\begin{proof}
That condition~(2) of Definition~\ref{def:versions of smoothness}
implies condition~(1) is immediate, as is the fact that
condition~(4) implies condition~(3).  To see that condition~(3)
implies condition~(2), note that if $A$ is an Artinian local
$\mathcal O_S$-algebra
whose residue field is of finite type over $\mathcal O_S$, 
then $\Spec A$ is of finite type over $S$.  Thus, if condition~(3) holds
and we are in the situation of condition~(2), the base-changed
morphism $U\times_{\cF} Z \to Z$ is smooth in a neighbourhood 
of the image of the induced map $Z_0 \to U\times_{\cF} Z$,
and since smooth morphisms are formally smooth, we may find the
desired lift $Z \to U$.  Thus condition~(3) of the definition
implies condition~(2).
This completes the proof of part~(1) of the present lemma.

Part~(2) is ~\cite[Lem.~2.3]{hall2013artin}; we recall the (short) argument.
As explained in the previous paragraph, $Z_0, Z$ are
automatically of finite type over $S$. Let $W_0$ be the image of $Z_0$ in
$\Spec\cO_{U,u}$, and let $W$ be the pushout of $W_0$ and $Z$ over $Z_0$. Then by~(RS)
we have a commutative diagram $$\xymatrix{Z_0 \ar[d] \ar[r] & W_0 \ar[d] \ar[r]& U \ar^-{\varphi}[d] \\
Z\ar[r]& W \ar[r] & \cF,}$$ and versality at $u$ allows us to lift the morphism
$W\to \cF$ to a morphism $W\to U$. The composite morphism $Z\to W\to U$ gives
the required lifting.

We turn to proving~(3).  
We first note that,
if $\cF$ is a stack satisfying~[1] and~[3],
and $\varphi: U \to \cF$ with $U$ locally of finite type
over $S$, then $\varphi$ is representable by algebraic
spaces, by~[3], and is locally of finite type,
by Lemma~\ref{lem:another finiteness lemma}~(1) below.

Also, if $\cX \to \cF$ is a morphism from an algebraic
stack to $S$, then $\cX$ admits a smooth surjection $X \to \cX$
whose source is a scheme (which is locally of finite type
over $S$ if $\cX$ is, since smooth morphisms are locally of finite
presentation), and, since smoothness can be tested smooth-locally
on the source, we see that (c) or (d) for $\cX$ is equivalent
to the corresponding statement for $X$.  Thus we need only
consider the case when $\cX$ is a scheme $X$ from now on.

If $\cF$ is a stack (on the \'etale site, and hence also on the Zariski site),
and if $X \to \cF$ is a morphism whose source
is a scheme, then condition~(d) may be checked Zariski locally
on $X$, and thus the general case of~(d) follows from the
case when $X$ is affine.  
If $\cF$ furthermore satisfies~[1],
then we may factor $X \to \cF$ through a finite type $S$-scheme,
and thus assume that $X$ is locally of finite type over $S$.

Thus we see that if $\cF$ is a stack satisfying~[1] and~[3], 
then (d) follows from~(c), while clearly~(d) implies~(c).
Again, it is clear that (c) implies (b),
and an argument essentially
identical to the proof of (1) above
shows that (b) implies~(a).  Thus it remains to show
that~(a) implies~(c), for maps $X \to \cF$
where $X$ is a scheme locally of finite type over $S$, under the assumption that
$\varphi$ is representable by algebraic spaces and locally
of finite type.

We must
verify that the smooth locus of $U\times_{\cF} X \to X$
contains the fibre over~$u$. The projection $U \times_{\cF} X \to X$ is  locally of
finite type, and hence $U\times_{\cF} X$ is locally of finite type
over $S$.  As $U$ is also locally of finite type over $S$,
the projection $U\times_{\cF} X \to U$ is locally of finite type
as well. Since $U\times_{\cF} X$ is an algebraic space,
it admits an \'etale cover by a scheme $V$.  
Since smoothness may be checked smooth
locally on the source, 
it suffices to show that the fibre of $V$
over $u$ is contained in the smooth locus of the composite
$V \to U\times_{\cF} X \to X$.
Since $V  \to U$ is locally of finite type,
the fibre of $V$ over $u$ is a scheme locally of finite type over $\kappa(u)$,
and so it suffices to show that every closed point
of this fibre lies in the smooth locus of $V \to X$.
By~\cite[Prop.\ 17.14.2]{MR0238860},  it in fact suffices to show that
this morphism is formally smooth at each of these points.

Let $v$ be such a point, and suppose given a commutative 
diagram 
$$\xymatrix{Z_0 \ar[r] \ar[d] & V \ar[d] \\
Z \ar[r] & X}$$
as in the definition of formal smoothness at $v$.  We may fit this into
the larger diagram
$$\xymatrix{Z_0 \ar[r] \ar[d] & V \ar[r]\ar[d] & U \ar[d] \\
Z \ar[r] & X \ar[r] & \cF}$$
Our assumption that $\varphi$ is formally smooth at $u$ implies
that we may lift the composite of the lower horizontal arrows
to a morphism $Z \to U$ (it is here that we use the assumption
that $v$ is a closed point of the fibre over $u$, so that the residue field
at $v$ is a finite extension of the residue field at $u$),
and hence to a morphism $Z \to U\times_{\cF} X$.  Since $V$ is \'etale,
and so in particular formally smooth,
over $U\times_{\cF} X$, we may then further lift this morphism to
a morphism $Z \to V$, as required.
This completes the proof that (a) implies~(c).

It remains to prove part~(4) of the lemma.
Lemma~\ref{lem: algebraic stacks satisfy RS} and
part~(2) show that conditions~(1) and~(2) of Definition~\ref{def:versions
of smoothness} are equivalent when $\cF$ is an algebraic stack. Taking into account the statement of part~(1), it suffices
to show that condition~(2) of Definition~\ref{def:versions of smoothness}
implies condition~(4). 
Thus we suppose that $\varphi:U \to \cF$ is formally smooth at the
point $u \in U$. The equivalence between conditions (a) and (d) of 
part~(3) of the present lemma (which holds, since $\cF$ is 
an algebraic stack, locally of finite presentation over $S$,
and thus satisfies~[1], by Lemma~\ref{lem: limit preserving iff locally
of finite presentation}, and~[3], by definition)
shows (taking $\cX = \cF$) that
$\varphi$ is smooth in a neighbourhood of $u$, as required.
\end{proof}

As a corollary of Lemma~\ref{lem:versions of smoothness}~(3),
we next show (under mild assumptions on the
morphism $\varphi$) that formal smoothness at a point can
be checked smooth locally on~$U$.

\begin{cor}
\label{cor:formal smoothness is smooth local}
Let $\varphi:U \to \cF$ be a morphism {\em (}over $S${\em )}  whose
source~$U$ is a locally of finite type $S$-scheme, and whose
target~$\cF$ is
 a category fibred in groupoid, and suppose that
$\varphi$ is representable by algebraic spaces,
and is locally of finite type.\footnote{As seen in the 
proof of Lemma~\ref{lem:versions of smoothness}~(3),
these conditions on $\varphi$ hold automatically
if $\cF$ is a stack satisfying~[1] and~[3].}
If $u \in U$ is a finite type point, then the following are equivalent:
\begin{enumerate}
\item $\varphi$ is formally smooth at $u$.
\item There is a smooth morphism of schemes $V \to U$
and a finite type point $v \in V$ mapping to $u$ such that the composite
$V \to U \to \cF$ is formally smooth at $v$.
\item For any smooth morphism of schemes $V \to U$,
and any finite type point $v \in V$ mapping to $u$,
the composite $V \to U \to \cF$ is formally smooth at $v$.
\end{enumerate}
\end{cor}
\begin{proof}
  If $V\to U$ is a smooth morphism, mapping the finite type point $v \in
  V$ to $u$, and $X\to\cF$ is a morphism from a finite type $S$-scheme $X$, then
  we consider the following commutative diagram.\[\xymatrix{V\times_\cF
    X\ar[d]\ar[r]&U\times_\cF X\ar[d]\ar[r]&X\ar[d]\\V\ar[r]&U\ar[r]&\cF}\] If
  (1) holds, then by the equivalence of (a) and (b) in part~(3) of Lemma~\ref{lem:versions of smoothness}, the projection
  $U\times_{\cF} X \to X$ contains the fibre over $u$ in its smooth locus.  The
  fibre over $v$ (with the respect to the projection $V \times_{\cF} X \to V$)
  is contained in the base-change to $V$ of this fibre over $u$, and since the
  morphism $V\times_{\cF} X \to U\times_{\cF} X$ is smooth, the smooth locus of
  the projection $V\times_{\cF} X \to X$ contains the fibre over $v$ in its
  smooth locus. Employing Lemma~\ref{lem:versions of smoothness} again, we see
  that~(1) implies~(3).

Considering the same diagram, and using the facts that smooth morphisms
are open, and that smoothness can be tested smooth-locally on
the source, we see in the same way that~(2) implies~(1).

Clearly~(3) implies~(2) (e.g.\ by taking $V = U$ and $v = u$), so we are done.
\end{proof}

\begin{cor}
\label{cor:formal smoothness at a point for stacks}
Suppose that $\cF$ is a stack 
over $S$ satisfying {\em [1]} and~{\em [3]}, 
that $\cU$ is an algebraic stack, locally of finite type over $S$,
and that $\cU \to \cF$
is a morphism of stacks over $S$.
If $u \in |\cU|$ is a finite type point of $\cU$,
then the following are equivalent:
\begin{enumerate}
\item
There exists a smooth morphism $V \to \cU$ whose source
is a scheme, and a finite type point $v \in V$ mapping to $u$,
such that the composite $V \to \cU \to \cF$ is formally smooth at $v$.
\item
For any smooth morphism $V \to \cU$ whose source
is a scheme, and any finite type point $v \in V$ mapping to $u$,
the composite $V \to \cU \to \cF$ is formally smooth at $v$.
\end{enumerate}
\end{cor}
\begin{proof}
Since $\cU$ is an algebraic stack, there does exist a smooth
surjection $V \to \cU$ whose source is a scheme, and the finite type points of $V$ are dense in the fibre over $u$.  Thus if (2)
holds, so does~(1).  Conversely, suppose that (1) holds,
for some choice of $V$ and $v$, and suppose that $V' \to \cU$
is a smooth morphism from a scheme to $\cU$,
and that $v' \in V'$ is a finite type point mapping to $u$.  We must show
that $V' \to \cF$ is formally smooth at $v'$.

Consider the fibre product $V\times_{\cU} V'$; this is an algebraic
space (since $\cU$ is an algebraic stack, and so satisfies~[3]),
and so we may find an \'etale (and in particular smooth)
surjection $W \to V\times_{\cU} V'$
whose source is a scheme.
Since $v$ and $v'$
both map to $u$, and since $W \to V\times_{\cF} V'$ is surjective,
we may find a finite type point $w \in W$ lying over both $v$ and $v'$.

Since $\cF$ is a stack satisfying~[1] and~[3], and since $V$ and $V'$ 
are each of finite type over $S$, the morphisms $V\to \cF$
and $V' \to \cF$ are both representable by algebraic spaces
and locally of finite type.  Thus,
since $W\to V$ and $W \to V'$ are both smooth morphisms,
we conclude from Corollary~\ref{cor:formal smoothness is smooth
local}, first that $W \to \cF$ is formally smooth at~$w$, 
and then that $V' \to \cF$ is formally smooth at~$v'$.
This completes the proof that (2) implies~(1).
\end{proof}

\begin{df}
\label{def:formal smoothness at a point for stacks}
If $\cU$ is an algebraic stack, locally of finite type over $S$,
and $\cF$ is a stack over $S$ satisfying~[1] and~[3],
then we say that a morphism $\varphi:\cU \to \cF$
of stacks over $S$ is formally smooth
at a finite type point $u \in |\cU|$ if the equivalent conditions
of Corollary~\ref{cor:formal smoothness at a point for stacks}
hold.
\end{df}

\begin{remark}
\label{rem:formal smoothness at a point for stacks}
Suppose we have a commutative diagram of morphisms of stacks over $S$
$$\xymatrix{\cU'' \ar^-{\varphi'}[r]\ar^-{\psi'}[d] & \cU' \ar^-{\psi}[d] \\
\cU \ar^-{\varphi}[r] & \cF}$$
in which $\cU$, $\cU'$, and $\cU''$ are algebraic stacks,
locally of finite type over $S$, and $\cF$ satisfies [1]
and~[3].  Suppose further that $\psi'$ and $\varphi'$ are
smooth.    Let $u'' \in |\cU''|$ be a finite type point,
with images $u \in |\cU|$ and $u' \in |\cU'|$.  
Then it follows directly from the definition (and Corollary~\ref{cor:formal smoothness at a point for stacks})
that $\varphi$ is formally smooth at $u$ if and only if
$\psi$ is formally smooth at $u'$.
(Both conditions hold if and only if the composite $\varphi\circ \psi'
= \psi\circ \varphi'$ is formally smooth at $u''$.)
\end{remark}

The following lemma extends those parts
of Lemma~\ref{lem:versions of smoothness} dealing with formal smoothness
at a point to the case of morphisms whose source is an algebraic stack.

\begin{lemma}
\label{lem:formal smoothness at a point for stacks}
Let $\cU$ be an algebraic stack, locally of finite type over $S$,
and let $\cF$ be a stack over $S$ satisfying {\em [1]} and~{\em [3]}.
Let $\varphi: \cU \to \cF$ be a morphism of stacks over $S$.
Let $u \in | \cU|$ be a point of $\cU$.
\begin{enumerate}
\item
The morphism $\varphi$ is formally smooth at $u$ if and only if,
for every morphism $\cX \to \cF$ whose source is an algebraic 
stack, the base-changed morphism $\cU \times_{\cF} \cX \to \cX$
contains the fibre over $u$ in its smooth locus.
\item If $\cF$ is also an algebraic stack,
then $\varphi$ is formally smooth at $u$ if and only 
if it is smooth in a neighbourhood of $u$.
\end{enumerate}
\end{lemma}
\begin{proof}
Given our assumptions on $\cF$, conditions (3)(a) and (3)(d)
of Lemma~\ref{lem:versions of smoothness} are equivalent.
The present lemma follows in a straightforward manner,
taking into account Definition~\ref{def:formal smoothness
at a point for stacks}, Lemma~\ref{lem:versions of smoothness}~(4)
(which we can apply, because~$\cF$ is locally of finite presentation
by Lemma~\ref{lem: limit preserving iff locally of finite presentation}), and the fact that smoothness
for morphisms from an algebraic stack to $\cF$ can
be checked smooth-locally on the source.\end{proof}

The preliminaries being dealt with,
we now state Axiom~[4], the condition of
{\em openness of versality}.  We assume that $S$ is locally Noetherian,
and that $\cF$ is a category fibred in groupoids over $S$.

\subsubsection*{Axiom~{\em [4]}}
If $U$ is a scheme locally of finite type over $S$, and $\varphi:U\to\cF$ is
versal at some finite type point $u\in U$, then there is an open neighbourhood
$U'\subseteq U$ of $u$ such that $\varphi$ is versal at every finite type point of $U'$.

\smallskip

\subsubsection{Alternatives to axiom~{\em [4]}} 
Following Artin \cite{MR0260746} we introduce a pair of axioms, labelled~[4a]
and~[4b], which are closely related to Axiom~[4]. These axioms are not needed
for our main results (and in particular for our applications to Galois
representations in Section~\ref{sec: finite flat}), but as we have modelled our
treatment of Artin's representability theorem on~\cite{MR0260746}, we have
followed Artin in introducing and discussing these variants on
Axiom~[4]. We will also make use of these axioms when discussing the
various examples in Section~\ref{sec: examples}.

In order to state~[4a], we first introduce some notation.
Suppose that $A$ is a DVR, with field of fractions $K$.
If $A'_K$ denotes an Artinian local thickening of $K$,
i.e.\ an Artinian local ring equipped with a surjection
$A'_K \to K$ (which then induces an isomorphism between
the residue field of $A'_K$ and $K$), then we let $A'$ denote
the preimage of $A$ under the surjection $A'_K \to K$;
it is a subring of $A'_K$, equipped with a surjection
$A' \to A,$ whose kernel is nilpotent.

We now state Axiom~[4a].  In the statement,
we suppose that $S$ is locally Noetherian,
and that $\cF$ is a category fibred in groupoids
over $S$.

\subsubsection*{Axiom~{\em [4a]}}
Let $A$ be a DVR over $\mathcal O_S$, whose residue field is of
finite type over~$\mathcal O_S$, and let $K$ denote the field of fractions
of~$A$.
Then for any morphisms $\Spec A\to\cF$, $\Spec A'_K \to \cF$,
with $A'_K$ being an 
Artinian local thickening of $K$,
for which the diagram
$$\xymatrix{
\Spec K \ar[d] \ar[r] & \Spec A \ar[d] \\
\Spec A'_K \ar[r] & \cF }
$$
commutes,
there exists a morphism $\Spec A' \to \cF$ making the enlarged
diagram
$$\xymatrix{
\Spec K \ar[d] \ar[r] & \Spec A \ar[d]\ar[dr] & \\
\Spec A'_K \ar[r] & \Spec A' \ar[r] & \cF }
$$
commute.
\smallskip

\begin{remark} 
Our formulation of Axiom~[4a] is slightly different to Artin's,
who only requires the condition of~[4a] to hold for $A$ that
are essentially of finite type over $\cO_S$.
To explain this, note that
the key role of~[4a] in the theory is to make Lemma~\ref{lem:Artin's
  3.10} below
true, and that,
in the proof of that lemma, we apply~[4a] after taking a certain
integral closure.   If $S$ is a {\em Nagata scheme} (see
e.g.~\cite[\href{http://stacks.math.columbia.edu/tag/033R}{Tag 033R}]{stacks-project})
then taking this integral closure keeps us in the world of essentially finite
type $\cO_S$-algebras, and so in the case of a Nagata base (which holds
in Artin's setting, since his base $S$ is assumed to be excellent)
the argument only requires Artin's more limited
form of~[4a].

We also observe that
if $A$ is essentially of finite type over $\cO_S$,
then $A'$ may be written as the inductive limit
$A' = \varinjlim_{i} A'_i,$ where $A_i'$ runs over the
essentially finite type $\cO_S$-subalgebras of $A'$ with the properties that
$A'_i \to A$ is surjective, and that the localisation of $A_i'$
at its generic point is equal to $A'_K$.  Thus, if $\cF$ also satisfies~[1],
then if we can find a morphism $\Spec A' \to \cF$ making the 
diagram of Axiom~[4a] commute,
we may similarly find a morphism $\Spec A'_i \to \cF$
making the diagram
$$\xymatrix{
\Spec K \ar[d] \ar[r] & \Spec A \ar[d]\ar[dr] & \\
\Spec A'_K \ar[r] & \Spec A'_i \ar[r] & \cF }
$$
commute, for some (sufficiently large) value of $i$.
In Artin's context, this allows an alternate phrasing of~[4a],
e.g.~in the form of Axiom~[4] of~\cite{MR0262237}.\end{remark}

Artin proves in~\cite[Thm.~3.7]{MR0260746}
that every locally separated algebraic space 
satisfies (his formulation of) Axiom~[4a].  In fact this is true
of arbitrary algebraic stacks (with the more expansive formulation
of the axiom given here).

\begin{lemma}
\label{lem:algebraic stacks satisfy 4a}
Every algebraic stack satisfies Axiom~{\em [4a]}.
\end{lemma}
\begin{proof}
Since $\Spec A'$ may be thought of as the pushout of $\Spec A$ and $\Spec A'_K$
over $\Spec K,$ this is a special case 
of~\cite[\href{http://stacks.math.columbia.edu/tag/07WN}{Tag 07WN}]{stacks-project}.
\end{proof}

We can now state Axiom~[4b].  We assume that $S$ is locally Noetherian,
and that $\cF$ is a category fibred in groupoids over $S$.

\subsubsection*{Axiom {\em [4b]}} 
If
$\varphi: U \to \cF$ is a morphism whose source is an $S$-scheme,
locally of finite type,
and $\varphi$ is smooth at a finite type point $u \in U$,
then $\varphi$ is smooth in a neighbourhood of $u$.

\smallskip

\begin{remark}
\label{rem:versions of smoothness}
If we assume that $\varphi: U \to \cF$ is representable by algebraic
spaces and locally of finite type,
then Lemma~\ref{lem:versions of smoothness}~(3)
shows that being {\em formally smooth} at the finite type
point $u$ is equivalent to requiring that for each morphism
$X \to \cF$ from a finite type $S$-scheme, there is an open 
set containing the fibre over $u$ in $U\times_{\cF} X$
at which the projection to $X$ is smooth.  
The morphism is {\em smooth at} $u$ if this open set can in 
fact be taken to be the preimage of a neighbourhood of $u \in U$.
Finally, Axiom~[4b] holds precisely when this neighbourhood of $u$ can
be chosen independently of $X$.
\end{remark}

\begin{remark}
In Corollary~\ref{cor:4a and b give 4} below, following Artin~\cite[Lem.~3.10]{MR0260746},
we show
that if $\cF$ is a stack satisfying~[1], [2](a), and~[3], and whose diagonal is furthermore quasi-compact, 
then Axioms~[4a] and~[4b] for $\cF$ together
imply Axiom~[4] for~$\cF$.
This in turn implies  
that in Artin's representability theorem (Theorem~\ref{thm: Artin 
representability}) we may replace Axiom~[4] by the combination
of Axioms~[4a] and~[4b], provided that we add to~[3]
the condition that the diagonal be quasi-compact (see Theorem~\ref{thm: Artin representability; variant}). 
\end{remark}

\begin{remark}
As Artin notes \cite[p.~39]{MR0260746},
the two axioms~[4a] and~[4b] are quite different in nature.  Axiom~[4a] is
finitary; for example, the fact that it holds for algebraic
stacks immediately implies that it also holds for Ind-algebraic stacks.
Axiom~[4b] is not finitary in nature,
and although it holds for algebraic stacks, it will typically not hold for
Ind-algebraic stacks.
(See Subsection~\ref{subsec:non-representable} below for a further
discussion of Ind-algebraic stacks and their comportment with regard
to Artin's axioms.)
\end{remark}

\subsection{Imposing Axiom [1] via an adjoint construction}\label{subsec:[1]-ification}We want to be able to apply our
results to stacks that are not necessarily limit preserving, and we
will need a way to pass from such a stack to one which is limit preserving
without losing too much information.

Fix a quasi-separated base-scheme $S$. Let $\AffS$ denote the category of affine $S$-schemes,
and let $\AfffpS$ denote the full subcategory of finitely presented
affine $S$-schemes.  The category $\AffS$ is closed under the formation
of filtering projective limits (see~\cite[\href{http://stacks.math.columbia.edu/tag/01YX}{Tag 01YX}]{stacks-project}). Let $\pro\AfffpS$ denote
the formal pro-category of $\AfffpS$.

\begin{thm}
  \label{thm: Thomason and Trobaugh, affine is pro-fp-affine}
The formation of projective limits
induces an equivalence
$$ \pro\AfffpS \iso \AffS.$$
\end{thm}
\begin{proof}
 This is a consequence of the theory of absolute Noetherian
 approximation developed in~\cite[Appendix C]{MR1106918}; for convenience, in
 this proof we will refer to the treatment of this material
 in~\cite[\href{http://stacks.math.columbia.edu/tag/01YT}{Tag
   01YT}]{stacks-project}, which develops the relative versions of the
 statements of~\cite[Appendix C]{MR1106918} that we need.

  Note firstly that if
 $X$ is an affine $S$-scheme, then $X$ is in particular quasi-compact and
 quasi-separated, so that by~\cite[\href{http://stacks.math.columbia.edu/tag/09MU}{Tag
   09MU}]{stacks-project}, we may write $X=\varprojlim X_i$ as the limit of a
 directed system of schemes $X_i$ of finite presentation over $S$, such that the
 transition morphisms are affine over
 $S$. By~\cite[\href{http://stacks.math.columbia.edu/tag/01Z6}{Tag
   01Z6}]{stacks-project}, we may assume that the $X_i$ are all affine. 
This proves that the purported equivalence is essentially surjective.

 Let $f:X\to Y$ be a morphism of $S$-schemes with $Y$ also affine. Then we can
 also write $Y=\varprojlim Y_j$ with $Y_j$ a directed system of affine schemes
 of finite presentation over $S$. Then we have
 \begin{align*}
   \Mor_{\pro\AfffpS}(\varprojlim_i X_i,\varprojlim_jX_j)&=\varprojlim_j\varinjlim_i\Mor_{\AffS}(X_i,Y_j)\\&=\varprojlim_j\Mor_{\AffS}(X,Y_j)\\&=\Mor_{\AffS}(X,Y),
 \end{align*}
where the first equality is the definition of a morphism in $\pro\AfffpS$, the
second is \cite[\href{http://stacks.math.columbia.edu/tag/01ZB}{Tag
   01ZB}]{stacks-project}, and the third is by the universal property
 corresponding to the statement that $Y=\varprojlim Y_j$. Thus the purported equivalence is fully faithful, as required.
\end{proof}

Now let $\cF \to \AfffpS$ be a category fibred in groupoids.
Passing to the corresponding pro-categories, we obtain a category
fibred in groupoids
$$\pro\cF \to \pro(\AfffpS) \iso \AffS.$$
On the other hand, given a category $\cF'$ fibred in groupoids over $\AffS$,
we may always restrict it to $\AfffpS$, to obtain a category
$\cF'{|_\AfffpS}$ 
fibred in groupoids
over $\AfffpS$.

Since $\AffS$ is closed under the formation of filtering projective limits,
we see that the same is true of any category fibred in groupoids $\cF'$ over $\AffS$,
and so for any such $\cF'$, evaluating projective limits induces a functor
\numequation
\label{eqn:projective limits}
\pro(\cF'{|_ \AfffpS}) \to \cF'
\end{equation}
over $\AffS$.
\begin{lemma}
\label{lem:adjunction equivalence}
If $\cF$ is a category fibred in groupoids
over $\AfffpS$, then the natural embedding
$\cF \to (\pro\cF)|_{ \AfffpS}$ is an equivalence.
\end{lemma}
\begin{proof}
This is formal: the fully faithful embedding $\AfffpS \to
\pro\AfffpS\iso \AffS$ identifies the essential image
of $\AfffpS$ with the subcategory 
of pro-systems that are isomorphic to a pro-system which is eventually constant.
Similarly, the essential image of $\cF$ in $(\pro\cF)$ consists 
of those pro-objects that are isomorphic to a pro-system which
is eventually constant, which by the previous remark are precisely the
pro-objects lying over elements of $\AfffpS$.
\end{proof}

\begin{lemma}
\label{lem:limit preserving and pro-categories}
If $\cF'$ is a category fibred in groupoids over $\AffS$,
then $\cF'$ is limit preserving if and only if the functor
{\em (\ref{eqn:projective limits})} is an equivalence.
\end{lemma}
\begin{proof}
This is almost formal.
(The non-formal ingredient is supplied
by~\cite[\href{http://stacks.math.columbia.edu/tag/01ZC}{Tag 01ZC}]{stacks-project}.)
\end{proof}

The previous two lemmas show that the formation of $\pro\cF$
from $\cF$ gives  an equivalence between the $2$-category
of categories fibred in groupoids over $\AfffpS$ and the full subcategory
of the $2$-category
of categories fibred in groupoids over $\AffS$ consisting of objects
satisfying~[1].

The following lemma establishes some additional properties of this construction.

\begin{lemma}
\label{lem:restriction to finitely presented algebras} 
\begin{enumerate}\item
If $\cF$ is a category fibred in groupoids over $\AfffpS$,
and $\cF'$ is a category fibred in groupoids over $\AffS$,
then there is an equivalence of categories 
$$\Mor_{\AfffpS}(\cF,\cF'{|_\AfffpS}) \iso
\Mor_{\AffS}(\pro\cF,\cF').$$
\item
If $\cF$ is a stack for the Zariski site  {\em (}resp.\ the \'etale site, resp.\ the fppf site{\em )}
on $\AfffpS$, then
$\pro\cF$ is a stack for the Zariski site {\em (}resp.\ the \'etale site, resp.\
the
fppf site{\em )} on $\AffS$.
\item
If $\cF'$ is a category fibred in groupoids over $\AffS$ whose diagonal and double diagonal
are both limit preserving on objects,
then the functor~{\em (\ref{eqn:projective limits})} is fully faithful.
\item 
If $\cF'$ is a category fibred in groupoids over $\AffS$ whose diagonal 
is limit preserving on objects and representable by algebraic spaces,
then the functor~{\em (\ref{eqn:projective limits})} is fully faithful,
and the diagonal of $\pro(\cF'{|_\AfffpS})$ is also
representable by algebraic spaces. 
\item
If $S$ is locally Noetherian, and if $\cF'$ is a category fibred in groupoids over $S$ which
admits versal rings at all finite type points, then $\pro(\cF'{|_\AfffpS})$ also admits versal rings at all finite type points.\end{enumerate}
\end{lemma}
\begin{proof}
(1) Passing to pro-categories,  
and then composing with the functor 
$$\pro(\cF'{|_ \AfffpS}) \to \cF'$$
of~(\ref{eqn:projective limits}),
induces a functor 
\begin{multline*}
\Mor_{\AfffpS}(\cF,\cF'{|_\AfffpS}) \to
\Mor_{\AffS}(\pro\cF,\pro(\cF'{|_\AfffpS}))
\\
\to \Mor_{\AffS}(\pro\cF,\cF'),
\end{multline*}
which is the required equivalence.
A quasi-inverse is given by the functor
\begin{multline*}
\Mor_{\AffS}(\pro\cF,\cF')
\to
\Mor_{\AfffpS}((\pro\cF){|_\AfffpS},\cF'{|_\AfffpS})
\\
\iso
\Mor_{\AfffpS}(\cF,\cF'{|_\AfffpS}),
\end{multline*}
obtained by first restricting to $\AfffpS$,
and then taking into account
the equivalence of Lemma~\ref{lem:adjunction equivalence}.

(2) This follows by a standard limiting argument, which we recall. We
need to show that all descent data is effective, and that the presheaves
$\Isom(x,y)$ are sheaves.
We being with the argument for descent data.

We must show that if $T$ is affine, and $T'\to T$ is a Zariski  (resp.\
\'etale, resp.\
\emph{fppf}) cover of $T$, equipped with a morphism $T'\to\pro\cF$ with descent
data, then there is a map $T\to\pro\cF$ inducing the given map $T'\to\pro\cF$.

By definition, the map $T'\to\pro\cF$ factors as $T'\to T''\to\pro\cF$, with
$T''$ of finite presentation over~$S$. By Lemma~\ref{lem:descending etale to
  finite presented affine}, we can find a commutative diagram
\[\xymatrix{T'\ar[d]\ar[r]&T'_j\ar[d]\ar[r]&T''\\ T\ar[r]&T_j&}\]in
which the square is Cartesian, and $T'_j$ is an affine scheme of finite presentation
over $S$. Furthermore (by the same lemma),
if $T'\to T$ is a Zariski (resp.\ \'etale, resp.\
\emph{fppf}) covering, then we may assume that the morphism $T'_j\to T_j$ is a
Zariski (resp.\ \'etale, resp.\ \emph{fppf}) covering.

Since $\cF$ is assumed to be a stack, it is enough to show that (after possibly
increasing $j$) the descent data for the morphism $T'\to \pro\cF$ arises as the
base change of descent data for the morphism $T'_j\to\cF$. This descent data is
given by an isomorphism between the two maps $T'\times_T T'\to\pro\cF$ given by
the two projections, which satisfies the cocycle condition on the triple
intersection.

Now, an isomorphism between the two maps $T'\times_T T'\to\pro\cF$ is equivalent
to factoring the induced map $T'\times_T T'\to\pro\cF\times_S\pro\cF$ through
$\Delta_{\pro\cF}$. So, we have a morphism $T_j'\times_{T_j}
T_j'\to\pro\cF\times_S\pro\cF$ which factors through the diagonal after pulling
back to $T$, and we want to show that it factors through the diagonal after
pulling back to some $T_{j'}$. This will follow immediately provided that
$\Delta_{\pro\cF}$ is limit preserving on objects. Similarly, to deal with the
cocycle condition, it is enough to show that the double diagonal
$\Delta_{\Delta_{\pro\cF}}$ is limit preserving on objects. Since $\pro\cF$ is
limit preserving by definition, the required limit preserving properties
of the diagonal and double diagonal follow from Lemma~\ref{lem:limit preserving vs. limit preserving on objects}.

We now turn to proving that if
$x,y: T \rightrightarrows \pro\cF$ are two morphisms,
then $\Isom(x,y)$ is a sheaf on the Zariski
(resp.\ \'etale, resp.\ \emph{fppf}) site of $T$. 
By definition, we may find a morphism $T \to T_j$, with $T_j$ affine of finite
presentation over $S$,
and morphisms $x_j,y_j: T_j \to \cF$ such that $x$ and $y$
are obtained as the pull-backs of $x_j$ and $y_j$. 
We may further find a projective system $\{T_{j'}\}$ of affine $S$-schemes
of finite presentation, having $T_j$ as final object, and an 
isomorphism $T \iso \varprojlim_{j'} T_{j'},$ inducing the given morphism
$T \to T_j$; we then write $x_{j'}$ and $y_{j'}$ for the composites 
$T_{j'} \to T_j \buildrel x_j,y_j \over \rightrightarrows \cF$.
We recall (see Lemma~\ref{lem:descending etale to
  finite presented affine}) that the
Zariski (resp.\ \'etale, resp.\ \emph{fppf})
 site of $T$ is then naturally identified with the projective limit
of the corresponding sites of the~$T_{j'}$; by~\cite[Thm.\
VI.8.2.3]{MR0354653}, the same is true of the corresponding topoi.
In particular, the various Isom presheaves $\Isom(x_j,y_j)$,
which by assumption are in fact sheaves on the 
Zariski (resp.\ \'etale, resp.\ \emph{fppf}) sites of the~$T_j$,
form a projective system whose projective limit can be identified 
with a sheaf on the 
Zariski (resp.\ \'etale, resp.\ \emph{fppf}) sites of the $T$.
Unwinding the definitions,
one furthermore finds that this projective limit sheaf
is naturally isomorphic to $\Isom(x,y)$.  Thus we find that $\Isom(x,y)$
is indeed a sheaf.

(3) To check that~(\ref{eqn:projective limits})
is fully faithful, it suffices to check that, for any
affine $S$-scheme $T$, if we write $T$ as a projective
limit $T = \varprojlim_{i} T_i$ of finitely presented affine $S$-schemes,
then the functor~(\ref{eqn:limit functor}) is fully
faithful. As was already noted in the proof of Lemma~\ref{lem:limit preserving
vs. limit preserving on objects}, this follows from the assumption
that the diagonal and double diagonal of $\cF'$ are limit preserving on objects.

(4) Proposition~\ref{prop: double diagonals of stacks satisfying [3]}
shows that our assumptions on $\cF'$ imply that its double diagonal
is also limit preserving on objects, and so it follows from part~(3)
that~(\ref{eqn:projective limits}) is fully faithful.
This in turn implies that the diagram
$$
\xymatrix{
\pro(\cF'{|_\AfffpS}) \ar^-{\Delta}[r] \ar[d] & \pro(\cF'{|_\AfffpS}) \times_S \pro(\cF'{|_\AfffpS})\ar[d] \\
\cF' \ar^-{\Delta}[r] & \cF' \times_S \cF'}
$$
is $2$-Cartesian,
and thus if the bottom arrow is representable by algebraic spaces,
the same is true of the top arrow.

(5) Since finite type Artinian $\mathcal O_S$-algebras
are objects of $\AfffpS$, we see that the functors
$\cF'$ and $\pro\cF'{|_\AfffpS}$ induce equivalent groupoids
when restricted to
such algebras. Thus the finite type points of $\cF'$ and $\pro\cF'{|_\AfffpS}$
are in natural bijection, in the strong sense that for each finite
type $\cO_S$-field~$k$ there is a natural
bijection between the morphisms $x:\Spec k\to\cF'$ and the morphisms
$x:\Spec k\to\pro\cF'{|_\AfffpS}$,  and a versal ring to $\cF'$ at
such a morphism is
also a versal ring to $\pro\cF'{|_\AfffpS}$ at the same morphism.
\end{proof}

Finally, we note the following basic result.
\begin{lem}
  \label{lem: full subcategory is compatible with pro}If $\cF$ is a
  category fibred in groupoids over $\AfffpS$, and $\cF'$ is a
  full subcategory fibred in groupoids, then~$\pro\cF'$ is a full
  subcategory fibred in groupoids of~$\pro\cF$.
\end{lem}
\begin{proof}
  This is immediate from the definitions.
\end{proof}
\subsection{Stacks satisfying [1] and~[3]}
\label{subsec: 1 and 3}
In this subsection we discuss some properties of stacks satisfying Artin's
axioms~[1] and~[3].

\begin{lemma}
\label{lem:finiteness}
If $\cF$ satisfies~{\em [1]} and~{\em [3]}, and $\cX$ is an algebraic
stack locally of finite presentation over $S$,
then any morphism of $S$-stacks $\cX \to \cF$ is locally of finite presentation
{\em (}in the sense of Definition~{\em \ref{def:properties defined by
base-change})}.
\end{lemma}
\begin{proof}
Lemma~\ref{lem:characterizing morphisms}
shows that we may verify this after composing the morphism
$\cX \to \cF$ with a smooth surjection $U \to \cX$ whose
source is a scheme, and thus may assume that $\cX$ is in fact an
$S$-scheme~$X$, locally of finite presentation.
It follows 
from~\cite[\href{http://stacks.math.columbia.edu/tag/06CX}{Tag 06CX}]{stacks-project}
that $X \to S$ is limit preserving on objects, so that $X$ is limit preserving
(as it is a stack in setoids).
It follows from Corollary~\ref{cor: morphism between limit preserving} that $X \to
\cF$ is limit preserving on objects.  Since this morphism is representable by
algebraic spaces (as $\cF$ satisfies~[3]),
applying~\cite[\href{http://stacks.math.columbia.edu/tag/06CX}{Tag
  06CX}]{stacks-project} again yields the lemma.
\end{proof}

\begin{lemma}
\label{lem:diagonal}
If $\cF$ satisfies~{\em [1]} and~{\em [3]},
then the diagonal $\Delta: \cF \to \cF \times_S \cF$ is locally of
finite presentation.
\end{lemma}
\begin{proof}

Since $\cF$
is limit preserving by assumption,
Lemma~\ref{lem:limit preserving vs. limit preserving on objects} shows that
$\Delta:\cF\to\cF\times_S\cF$ is limit preserving on objects.
Since $\Delta$ is representable by algebraic spaces by assumption, it is locally of finite
presentation by~\cite[\href{http://stacks.math.columbia.edu/tag/06CX}{Tag 06CX}]{stacks-project}.
\end{proof}

We now establish some simple results related to finiteness conditions
in the case when $S$ is locally Noetherian.  We first recall
that a morphism $\cX \to S$ from an algebraic stack $\cX$ to the
locally Noetherian scheme $S$ is locally of finite type if and only
if it is locally of finite presentation (since both conditions may be
verified after composing this morphism with a smooth surjection
from a scheme to $\cX$, which reduces us to the case when $\cX$
is itself a scheme, in which case the claim follows immediately
from the hypothesis that $S$ is locally Noetherian).

\begin{lemma}\label{lem:another finiteness lemma}
Suppose that $\cF$ satisfies~{\em [1]} and~{\em [3]},
that $S$ is locally Noetherian, 
and that $\cX \to \cF$ is a morphism whose source is an algebraic
stack.
Then the morphism $\cX \to \cF$ is locally of finite type if
and only if it is locally of finite presentation, 
and these conditions are in turn equivalent to the composite
$\cX \to \cF \to S$ being locally of finite type {\em (}or, equivalently,
locally of finite presentation, as noted above{\em )}.
\end{lemma}
\begin{proof}

If $\cX \to S$
is locally of finite type, and so locally of finite presentation,
then Lemma~\ref{lem:finiteness} shows that
$\cX \to \cF$ is locally of finite presentation. 
And certainly,
if $\cX \to \cF$ is locally of finite presentation
then it is locally of finite type.  

Next, we want to show that if $\cX \to \cF$ is locally 
of finite type, then it is locally of finite presentation.
To this end,
let $T$ be an affine $S$-scheme.
Since $\cF \to S$ is limit preserving, and so in particular limit
preserving on objects, we may factor any $S$-morphism
$T \to \cF$ as $T \to T' \to \cF$, where $T'$ is of finite presentation
over~$S$.  (It follows from~\cite[\href{http://stacks.math.columbia.edu/tag/09MV}{Tag
  09MV}]{stacks-project} that $T$ may be written as the limit of such~$T'$.)
The base-changed morphism $\cX_T \to T$ is obtained by pulling back the base-changed morphism
$\cX_{T'} \to~T'$.  This latter morphism is locally of finite type,
and its target
is of finite presentation over the locally Noetherian scheme $S$ (and hence locally
Noetherian itself),
and so it is in fact locally of finite presentation.
Thus the morphism $\cX_T \to T$ is also locally of finite 
presentation, and since $T$ and the morphism $T\to\cF$ were arbitrary,
we conclude that the morphism $\cX \to \cF$
is locally of finite presentation, as claimed.

It remains to show that if $\cX \to \cF$ is locally
of finite presentation, then $\cX \to S$ is locally of finite presentation.
Morally, we would like to prove this by arguing
that since $\cF$ satisfies~[1],
the morphism $\cF \to S$ is locally of finite presentation, and so conclude that
the composite $\cX \to S$ is locally of finite presentation.  
Unfortunately, while this is a valid argument if $\cF$
is an algebraic stack, it does not apply in the generality
we are considering here, where $\cF$ is assumed simply to
satisfy [1] and~[3].  Indeed, in this level of generality,
we haven't defined what it means for $\cF \to S$ to be
locally of finite presentation.

Thus we are forced to make a slightly more roundabout argument.
Since $\cX \to \cF$ is  locally
of finite presentation, it is
limit preserving on objects,
by Lemma~\ref{lem:lfp equals limit preserving on objects}.
The morphism $\cF \to S$ is also limit preserving on objects
(since $\cF$ satisfies~[1]), and hence the composite
$\cX \to S$ is limit preserving on
objects~\cite[\href{http://stacks.math.columbia.edu/tag/06CW}{Tag 06CW}]{stacks-project}.
It then follows
from Lemma~\ref{lem: limit preserving iff locally of finite presentation}
that $\cX \to S$ is locally of finite presentation (or, equivalently,
locally of finite type).
\end{proof}

Our next results are inspired by a lemma of Artin~\cite[Lem.~3.10]{MR0260746}.
Our first lemma isolates one of the steps in Artin's argument,
and generalises it to the stacky context.

\begin{lemma}
\label{lem:monomorphism}
Suppose that $S$ is locally Noetherian, that $T$ is a locally finite type $S$-scheme, that $\cF$ satisfies~{\em [1]} and {\em [3]},
that $t \in T$ is a finite type point,
and that $T \to \cF$ is an $S$-morphism which is formally smooth at
$t$. Then, replacing $T$ by an open neighbourhood of $t$ if necessary,
we may factor the morphism $T \to \cF$ as 
$$ T \to \cF' \to \cF,$$
where~$\cF'$ is an algebraic stack,
locally of finite presentation over $S$,  
the first arrow is a smooth surjection, and the second arrow is
locally of finite presentation, unramified,
representable by algebraic spaces,
and formally smooth at the image~$t'$ of $t$ in $|\cF'|$.
\end{lemma}
\begin{proof}
Write $R := T \times_{\cF} T$.  By the assumption that $\cF$ satisfies~[3], this is a locally finite type algebraic
space over $S$; indeed,
it is the base-change of the morphism $T\times_S T \to \cF\times_S \cF$
via the diagonal $\Delta:
\cF \to \cF\times_S\cF$, and Lemma~\ref{lem:diagonal} shows
that this latter morphism
is locally of finite presentation.
The projections $R \rightrightarrows T$ 
endow $R$ with the structure of a groupoid in algebraic
spaces over $T$.
Since the morphism $T\to\cF$ is locally of finite presentation by Lemma~\ref{lem:another finiteness lemma}~(1), each of these projections is also locally of finite presentation (or equivalently,
locally of finite type, since the base $S$ is locally Noetherian),
and the formal smoothness of $T\to \cF$ at $t$ implies that both projections
are smooth in a neighbourhood
of $(t,t)$ by~Lemma~\ref{lem:versions of smoothness}~(3).
Thus, applying Lemma~\ref{lem:groupoids} below, 
we see that, by shrinking $T$ around $t$ if necessary, we may find an open
subgroupoid $V \subseteq R$
that is actually a smooth groupoid over~$T$. We then define $\cF' := [T/V]$, and let $t'$  denote the image of $t$ in 
$|\cF'|$. 

Certainly, since $V \subseteq T\times_{\cF} T,$
the map $T \to \cF$ factors through $\cF'$. Since $T \to \cF'$ is smooth and $T$ is locally of finite type over
$S$, so is $\cF'$.
(The property of being locally of finite type is smooth local on the source;
see~\cite[\href{http://stacks.math.columbia.edu/tag/06FR}{Tag 06FR}]{stacks-project}).
The morphism $\cF' \to \cF$ is thus locally of finite presentation,
by Lemma~\ref{lem:another finiteness lemma}~(1).
Since $T \to \cF$ is formally smooth at $t$, and $T \to \cF'$ is smooth,
the morphism $\cF' \to \cF$ is also formally smooth at $t'$ {\em by definition}.
(See Definition~\ref{def:formal smoothness at a point for stacks}.)

It remains to show that
the morphism $\cF' \to \cF$ is unramified and representable by
algebraic spaces. We have already seen that it is locally of finite presentation,
and so, in particular, locally of finite type,
and (as explained in Example~\ref{ex:unramified morphisms})
we must show furthermore that the diagonal morphism
$\cF' \to \cF'\times_{\cF} \cF'$ is an open immersion. (Recall that since $\cF$ satisfies~[3] by assumption,
the fibre product $\cF'\times_{\cF} \cF'$ is an algebraic stack.)
Since the morphism $R = T\times_{\cF} T \to   \cF'\times_{\cF} \cF'$
is smooth and surjective (as $T\to \cF'$ is),
we can verify this after base-changing by this latter morphism.
Since $R=T\times_{\cF}T$, and $V=T\times_{\cF'}T$, one
verifies that the base-changed morphism $$\cF' \times_{\cF'\times_{\cF} \cF'} R \to R$$
is precisely the open immersion $V \to R.$ 
This establishes the claim.
\end{proof}

\begin{remark}
In the case when $\cF$ is simply a sheaf of sets (which is the
context of~\cite{MR0260746}), the algebraic stack $\cF'$ will
simply be an algebraic space, and if we replace it by an \'etale
cover by a scheme $X$, we obtain an unramified morphism $X \to \cF$
which is formally smooth at a point $x$ above $t$.  This is essentially
the conclusion of the second paragraph of the proof 
of~\cite[Lem.~3.10]{MR0260746}, and our argument is
an adaptation to the stacky context
of the argument given there.
\end{remark}

\begin{lemma}
\label{lem:groupoids}
Let $(U,R, s,t,c,e,i)$ be a groupoid in algebraic spaces locally
of finite type over a locally
Noetherian scheme $S$, with $U$ a scheme.
If $u \in U$ is a finite type point such that $s$ is smooth
at the point $e(u)$ {\em (}or equivalently, such that $t$ 
is smooth at the point $e(u)${\em )},
then there exists an open neighbourhood $U'$ of $u$ in $U$
such that, if $(R',s',t',c')$ denotes the restriction 
of the groupoid $(R,s,t,c)$ to $U'$,
there is an open subgroupoid $V \subseteq R'$ 
such that $s'_{| V}$ and $t'_{|V}$ are smooth,
i.e.\ such that $V$ is a smooth groupoid over $U'$.
\end{lemma}
\begin{proof}
Let $V'$ denote the open subspace of $R$ on which $s$ is smooth.  
Write $U' = e^{-1}(V');$ then $U'$ is
an open subscheme of $U$ containing $u$. 
If we replace $U$ by $U'$ and $R$ by $R' := R_{| U'}
= (s \times t)^{-1}(U' \times U')$,
so that $V'$ is replaced by $V' \cap R'$,
then we may, and do, assume that $e(U) \subseteq V'$.

Recall that we have the commutative diagram
\begin{equation*}
\xymatrix{
& U & \\
R \ar[d]_s \ar[ru]^t &
R \times_{s, U, t} R
\ar[l]^-{\text{pr}_0} \ar[d]^{\text{pr}_1} \ar[r]_-c &
R \ar[d]^s \ar[lu]_t \\
U & R \ar[l]_t \ar[r]^s & U
}
\end{equation*}
of~\cite[\href{http://stacks.math.columbia.edu/tag/043Z}{Tag
  043Z}]{stacks-project}, 
in which each square (including the top square) is Cartesian.
Pulling back the right-hand square of this diagram via the morphism
$V'\to R$, we form the Cartesian diagram
$$
\xymatrix{
R \times_{s, U, t} V'  
\ar[d]^{\text{pr}_1} \ar[r]_-c & R \ar[d]^s \\
V'  \ar[r]^s & U.
}
$$
The bottom arrow is smooth by the definition
of $V'$, and is furthermore surjective, since
we have put ourselves in a situation in which
$e(U) \subseteq V'$.  Thus the locus of smoothness
of the left-hand vertical arrow is precisely 
the preimage under $c$ of the locus of smoothness
of the right-hand vertical arrow; i.e.\ the 
locus of smoothness of $\text{pr}_1: R\times_{s,U, t} V'
\to V'$ is equal to $c^{-1}(V')$.  Now certainly
$\text{pr}_1:V'\times_{s,U,t} V' \to V'$ is smooth,
since it is a base-change of the smooth morphism
$s: V' \to U.$
Thus $V'\times_{s,U,t} V' \subseteq c^{-1}(V')$,
or equivalently, $c(V'\times_{s,U,t} V') \subseteq V'.$

Now define $V := V' \cap i(V')$.   Clearly $c(V \times_{s,U,t} V)
\subseteq V.$  Also $e(U) \subseteq V$, and $i(V) = V.$   
Taken together, these properties show that $V$ is an open
subgroupoid of~$R$.  Since $V \subseteq V',$ we see that
$s_{| V}$ is smooth.  Since $t = s i$ and $i(V) = V$,
we see that $t_{| V}$ is smooth as well.  Thus the lemma is proved.
\end{proof}

In the remainder of the section, we return to Axioms~[4a] and~[4b]; none of this
material is needed for our main theorems. The following lemma, and its corollary,
provide an analogue in the stacky context 
of~\cite[Lem.~3.10]{MR0260746} itself.
The argument is essentially the same as Artin's.

\begin{lemma}
\label{lem:Artin's 3.10}
Suppose that $S$ is locally Noetherian,
that $T$ is a locally finite type $S$-scheme,
and that $\cF$ satisfies~{\em [1]}, {\em [3]}, and~{\em [4a]}. 
If $\cX \to \cF$ is a morphism of finite type
from an algebraic stack to $\cF$ {\em (}in the sense
of Definition~{\em \ref{def:properties defined by base-change})},
and if $T \to \cF$ is an $S$-morphism which is formally
smooth at a finite type point $t \in T$,
then there exists a neighbourhood $T'$ of $t$ in $T$ such
that the base-changed morphism $T' \times_{\cF} \cX \to \cX$
is smooth.
\end{lemma}
\begin{proof}
Clearly we may replace $T$ by an affine open neighbourhood
of $t$, and thus suppose that $T$ is quasi-compact.
We may apply Lemma~\ref{lem:monomorphism} to the morphism
$T \to \cF$, and so, replacing $T$ by a neighbourhood of $t$ in $T$
if necessary, we factor $T \to \cF$ as $T \to \cF' \to \cF$
as in the statement of that lemma. 
We also choose a smooth surjection $U \to \cX$ whose source is a scheme.
We then consider the diagram 
$$\xymatrix{
T\times_{\cF} U \ar[r] \ar[d] & \cF'\times_{\cF} U \ar[r]\ar[d] & U \ar[d] \\
T\times_{\cF} \cX \ar[r] \ar[d] & \cF'\times_{\cF} \cX \ar[r] \ar[d]
& \cX \ar[d] \\
T \ar[r] & \cF' \ar[r] & \cF}
$$
Since $\cF$ satisfies~[3], the fibre product $T\times_{\cF} U$
is an algebraic space, while $\cF'\times_{\cF} U$,
$T\times_{\cF} \cX$, and $\cF'\times_{\cF} \cX$ are algebraic 
stacks.  Since the morphism $\cX \to \cF$ is quasi-compact,
and since $T$ (and hence also $\cF'$) is quasi-compact,
the fibre-products $T\times_{\cF} \cX$ and $\cF'\times_{\cF} \cX$
are furthermore quasi-compact.

Since the upper right horizontal arrow is an unramified
morphism, and so (by definition) has an \'etale diagonal,
it is in particular a DM morphism (i.e.\ has an unramified diagonal).
Since its target is a scheme, we thus find that
$\cF'\times_{\cF} U$ is in fact a DM stack. We may thus amplify
the preceding diagram to a commutative diagram
$$\xymatrix{
V \ar[r]\ar[d] & V' \ar[d] & \\
T\times_{\cF} U \ar[r] \ar[d] & \cF'\times_{\cF} U \ar[r]\ar[d] & U \ar[d] \\
T\times_{\cF} \cX \ar[r] \ar[d] & \cF'\times_{\cF} \cX \ar[r] \ar[d]
& \cX \ar[d] \\
T \ar[r] & \cF' \ar[r] & \cF}
$$
in which $V$ and $V'$ are quasi-compact schemes,
the vertical arrows with $V$ and $V'$ as their sources are \'etale,
and
the composites
$V \to T\times_{\cF} \cX$
and $V' \to \cF'\times_{\cF} \cX$,
as well as the four horizontal arrows on the left half of the diagram,
are smooth surjections.
(Note that the upper-most square in this diagram is {\em not} assumed
to be Cartesian.)

Let $C' \subseteq V$ denote the complement of
the smooth locus of the composite $V \to~U.$
It is a closed subset of $V$, and so its image $C$ in $T$
is a constructible subset of~$T$.  (The morphism $V \to T$
is a locally finite type morphism between Noetherian schemes,
and hence is finite type.   Thus it maps constructible sets
to constructible sets.)

Lemma~\ref{lem:another finiteness lemma}~(1) shows
that $U$ is locally of finite type over $S$, and so,
by Lemma~\ref{lem:versions of smoothness}~(3), the smooth
locus of the morphism $T\times_{\cF} U \to U$ contains
the fibre over $t$, and thus the same is true of the 
smooth locus of the morphism $V \to U$. 
In other words, the point $t$ does not lie in $C$.
To prove the lemma we must show that there is an open 
neighbourhood of $t$ disjoint from $C$; that is,
we must show that $t$ does not lie in the closure 
of $C$.

Suppose that $t$ does lie in the closure of $C$.
Then we claim that we may find a point $\tilde{t} \in C$
whose closure $Y$ contains $t$, and such that
$A := \mathcal O_{Y,t}$ is a one-dimensional domain. To see this, note firstly that $C$ is a finite union of irreducible
locally closed schemes, so $t$ is in the closure of one of these, and we may
replace $C$ by this irreducible component (with its induced reduced structure),
and thereby assume that $C$ is open in its closure. Replacing the closure of $C$
by an affine open subset containing $t$, and $C$ by a distinguished open subset
of this affine open, we put ourselves in the situation of having a Noetherian
domain $A$ with a closed point $t\in\Spec A$, and an element $f\in A$ such that
$\Spec A[1/f]$ is a proper non-empty subset. Replacing $A$ by its localisation at the maximal ideal
corresponding to $t$, we may assume that $A$ is a local domain, and that $f$ is
in the maximal ideal of $A$.

Since $A$ is local Noetherian, it is finite-dimensional, and we may choose a
prime $P$ not containing $f$ so that $\dim A/P$ is as small as
possible. Replacing $A$ by $A/P$, we have a Noetherian local domain $A$
containing a nonzero element $f\in A$, with the properties that $f$ is not a
unit, and $f$ is contained in every non-zero prime ideal of $A$. It remains to
check that $A$ is one-dimensional. To see this, note that the height one primes
of $A$ are precisely the isolated associated primes of $(f)$, so that $A$ has
only finitely many height one primes. Let $Q$ be any prime ideal of $A$. Take an
element $x\in Q$; then by Krull's principal ideal theorem, $(x)$ is contained in
some height one prime, so we see that $Q$ is contained in the union of the
height one primes. Since this is a finite union, we see by prime avoidance that
$Q$ is contained in some height one prime, as required.

Since $C'$ is closed in $V$, and so of finite type over~$T$,
we may find a point $\tilde{v} \in C'$ lying over $\tilde{t}$
such that residue field extension $\kappa(\tilde{t})\hookrightarrow
\kappa(\tilde{v})$ is finite.  We will obtain a contradiction
by showing that the morphism $V \to U$ is in fact formally smooth
at $\tilde{v}$ (and hence smooth at $\tilde{v}$, by Lemma~\ref{lem:versions
of smoothness}~(3)).

Consider a diagram 
$$\xymatrix{Z_0 \ar[d] \ar[r] & Z \ar[d] \\
V \ar[r]  & U}$$
as in the definition of formal smoothness at $\tv$ (Definition~\ref{def:versions of smoothness}~(2));
so $Z_0 \to Z$ is a closed immersion of Artinian local $\mathcal O_S$-algebras,
and the residue field $L$ at the closed point of $Z_0$ is a finite
extension of $\kappa(\tilde{v})$ (and thus also finite over
$K = \kappa(\tilde{t})$.)
We expand this diagram to the diagram
$$\xymatrix{\Spec L \ar[d] \ar[rrd] & & \\
Z_0 \ar[d]\ar[rd] \ar[rr] & & Z \ar[d] \\
V \ar[r]  & V'\ar[r] & U}$$
The morphism $V' \to U$ factors as $V' \to \cF'\times_{\cF} U \to U,$
and hence, as the composite of an unramified and an \'etale morphism,
is itself unramified.

We will show that the morphism $Z \to U$ lifts to a morphism
$Z \to V',$ compatible with the given morphism $\Spec L \to V'$.
Assuming this, we note that this morphism, when restricted
to $Z_0$, must coincide with the given morphism $Z_0 \to V'$;
this follows from the fact that $V' \to U$ is unramified, and thus
formally unramified.
Since $V \to V'$ is smooth, and hence formally smooth, we then see
that we may further lift the morphism $Z \to V'$ to a morphism
$Z \to V$ compatible with the given morphism $Z_0 \to V$.
This completes the proof that $V \to U$ is formally smooth at $\tilde{v}$,
and thus completes the proof of the lemma.

It remains to prove the existence of the lifting
$Z \to V'$.  Since $V' \to \cF'\times_{\cF} U$ is
\'etale, it suffices to obtain a lifting $Z \to \cF'\times_{\cF} U$. 
For this, it suffices in turn to obtain a morphism $Z \to \cF'$
lifting the composite $Z \to U \to \cF$,
and, for this, it suffices to obtain a morphism $Z \to T$
lifting $Z \to \cF$.

Write $Z = \Spec B'_L,$
and let $B$ denote the integral closure of $A$ in $L$.
Then by the Krull--Akizuki theorem (\cite[Thm.\ 11.7]{MR1011461} and its Corollary) $B$ is a semi-local Dedekind domain with field of fractions~$L$,
whose residue fields at its closed points are finite
extensions of the residue field of $A$ at its closed
point, i.e.\ of $\kappa(t)$. 
Let $B'$ be the inverse
image of $B\subseteq L$  in $B'_L$.
Since $\cF$ satisfies~[4a], we may extend the morphism
$\Spec B'_L \to \cF$ to a morphism $\Spec B'\to\cF,$
compatible with the given morphism
$$\Spec B \to \Spec A \to T \to \cF.$$ Lemma~\ref{lem:local lifting} below then shows that we may
lift the morphism $\Spec B' \to \cF$ to a morphism
$\Spec B' \to T$. 
Passing to the local ring at the generic point, this gives the required
morphism $Z \to T.$\end{proof}

\begin{lemma}
\label{lem:local lifting}
Suppose that $S$ is locally Noetherian,
that $\cF$ is a stack over $S$ satisfying {\em [1]} and~{\em [3]},
that $T$ is a scheme, locally of finite type over $S$,
and that $T \to \cF$ is a morphism over $S$
which is
formally smooth at a finite type point $t \in T$.
Suppose that $B$ is a Noetherian local $\mathcal O_S$-algebra
whose residue field is finite type over $\mathcal O_S$,
that $B'$ is an $\mathcal O_S$-algebra which
is a nilpotent thickening of $B$,
and that we have a commutative diagram of 
morphisms over $S$
$$\xymatrix{\Spec B \ar[r]\ar[d] & \Spec B' \ar[d] \\
T \ar[r] & \cF}
$$
such that the left-hand vertical arrow maps the closed point
of $\Spec B$ to the given point $t \in T$.
Then there is a lifting of the right-hand vertical arrow to a morphism
$\Spec B' \to T$.
\end{lemma}
\begin{proof}
Consider the projection $T\times_{\cF} \Spec B' \to \Spec B'.$
The commutative diagram in the statement of the lemma
gives rise to a morphism
\numequation
\label{eqn:lifting}
\Spec B\to T\times_{\cF} \Spec B'
\end{equation}
lifting the closed immersion $\Spec B \hookrightarrow \Spec B'$.
Since the closed point of $\Spec B$ maps to the point $t$ of $T$,
by assumption, it follows from Lemma~\ref{lem:versions of smoothness}~(3)
that the image of the closed point of $\Spec B$
under this morphism lies in the smooth locus of the projection,
and thus that the morphism~(\ref{eqn:lifting}) itself
factors through this smooth locus.   Since smooth morphisms
are, in particular, formally smooth, we may thus lift this morphism to
a section $\Spec B' \to  T\times_{\cF} \Spec B.$ Composing this section with the projection onto $T$
gives the morphism required by the statement of the lemma.
\end{proof}

\begin{cor}
\label{cor:Artin's 3.10}
Suppose that $S$ is locally Noetherian,
that $T$ is a locally finite type $S$-scheme,
that $\cF$ satisfies~{\em [1]}, {\em [3]}, and~{\em [4a]}, 
and furthermore that the diagonal of $\cF$ is quasi-compact.
If $T \to \cF$ is an $S$-morphism which is formally
smooth at a finite type point $t \in T$,
then this morphism is in fact smooth at $t$.
\end{cor}
\begin{proof}
Let $X \to \cF$ be a morphism whose source is a scheme of finite type over $S$,
and consider the usual ``graph diagram''
$$\xymatrix{X \ar^-{\Gamma}[r]\ar[d] & X \times_S \cF \ar[r] \ar[d] & \cF\\
\cF \ar^-{\Delta}[r] & \cF\times_S \cF }
$$
in which the square is Cartesian (by construction).
Lemma~\ref{lem:limit preserving vs. limit preserving on objects},
along with~\cite[\href{http://stacks.math.columbia.edu/tag/06CX}{Tag 06CX}]{stacks-project},
shows that the diagonal $\Delta:\cF \to \cF \times_S \cF$ is
locally of finite presentation, and it is quasi-compact 
by assumption.
Thus it is in particular of finite type, and so 
the same is true of the graph
$\Gamma: X \to X\times_S \cF$.  
The projection $X\times_S \cF \to \cF$ is also of finite type,
being the base-change of the finite type morphism of Noetherian schemes
 $X \to S$.
Thus the morphism $X \to \cF$ is of finite type,
and so Lemma~\ref{lem:Artin's 3.10} implies that we may find a neighbourhood
$U$ of $t$ in $T$ such that the base-changed morphism 
$U \times_{\cF} X \to X$ is smooth.  By definition, then,
the morphism $T \to \cF$ is smooth at $t$.
\end{proof}

\begin{cor}
\label{cor:4a and b give 4}
Suppose that $S$ is locally Noetherian, and that $\cF$ is a stack over
$S$ satisfying~{\em [1]}, {\em [2](a)},{\em [3]}, {\em [4a]}, and {\em [4b]},
whose diagonal is quasi-compact.  Then $\cF$ satisfies~{\em [4]}.
\end{cor}
\begin{proof}
Suppose that $\varphi: U \to \cF$ is a morphism from a scheme locally of
finite type over $S$ to $\cF$ which is versal at a  finite type
point $u \in U$.  Lemma~\ref{lem:versions of smoothness}~(2)
then shows that $\varphi$ is formally smooth at $u$. 
From Corollary~\ref{cor:Artin's 3.10} we conclude that $\varphi$
is in fact smooth at $u$, and Axiom~[4b] then implies that
$\varphi$ is smooth in a neighbourhood of~$u$.  Finally, 
Lemma~\ref{lem:versions of smoothness}~(1) shows that $\varphi$
is versal at each finite type point in this neighbourhood,
and thus $\cF$ satisfies Axiom~[4], as claimed.
\end{proof}

Our final result in this subsection strengthens the conclusion
of Lemma~\ref{lem:monomorphism}, in the presence of~[4a].

\begin{cor}
\label{cor:monomorphism}
Suppose, in the context of Lemma~{\em \ref{lem:monomorphism}},
that $\cF$ furthermore satisfies~{\em [4a]}, and
has quasi-compact diagonal.
Then, in addition to the conclusions of that lemma,
we may impose the condition that the morphism
$\cF' \to \cF$ be a monomorphism.
\end{cor}
\begin{proof}
By Corollary~\ref{cor:Artin's 3.10}, replacing $T$ by a neighbourhood
of $t$ if necessary, we may suppose that the projections
$R:= T\times_{\cF} T \rightrightarrows T$ are in fact smooth.  
Thus, in the proof of Lemma~\ref{lem:monomorphism},
we may take $V = R$. 
It is then easily verified
that $\cF' := [T/V] = [T/R] \to \cF$ is in fact a monomorphism.
Indeed, following the proof of Lemma~\ref{lem:monomorphism}, 
we deduce from the fact that $V = R$ that the diagonal
$\cF' \to \cF'\times_{\cF} \cF'$ is an isomorphism,
which is equivalent to $\cF' \to \cF$ being a monomorphism.
\end{proof}

\begin{remark}
Example~\ref{ex: line with infinitely many nodes} shows 
that we cannot remove the assumption that $\cF$ has quasi-compact diagonal 
from the preceding corollaries. \end{remark}

\subsection{Relationships between [2] and~[3]}
\label{subsec: 2 and 3}
Throughout this subsection we suppose that $\cF$ is a category fibred in
groupoids over the locally Noetherian base scheme
$S$.

\begin{lemma}
\label{lem:characterizing pro-representability}
Let 
$x:\Spec k \to \cF$ be a finite type point of $\cF$ {\em (}here $k$ is a finite type field over $\cO_S${\em )}.

(1)
If $\cF$ admits an effective Noetherian versal ring at $x$,
then,
for each object $A$ of~$\widehat{\mathcal C}_{\Lambda}$,
the functor $\cF_x(A) \to \widehat{\cF}_x(A)$ is essentially surjective.
Conversely, if this functor is essentially surjective 
for each object $A$ of $\widehat{\mathcal C}_{\Lambda}$,
and if $\cF$ admits a Noetherian versal ring at $x$,
then $\cF$ admits an effective Noetherian versal ring
at $x$.

(2)  If $\cF$ satisfies~{\em [3]}, then,
for each object $A$ of $\widehat{\mathcal C}_{\Lambda}$,
the functor $\cF_x(A) \to \widehat{\cF}_x(A)$ is fully faithful.
\end{lemma}
\begin{proof}
We begin by proving~(1).
Suppose first that $\cF$ admits an effective Noetherian versal ring  at $x$,
and let $A_x \in \widehat{\cC}_{\Lambda}$ be an effective versal ring
at the morphism $x$;
that is, we have an object $\eta$ of $\cF_x(A_x)$ whose image in
$\widehat{\cF}_x(A_x)$ is versal.
If $\xi: \Spf B \to \cF$ is any object of $\widehat{\cF}_x(B)$,
for some object $B$ of $\widehat{\cC}_{\Lambda}$,
then by definition of versality we may find a morphism $f:A_x \to B$
and an isomorphism $\xi \cong f^*\eta.$  Since $\eta$ is effective,
i.e.\ lies in the image of the functor
$\cF_x(A_x) \to \widehat{\cF}_x(A_x),$
we see that $\xi$ lies in the essential image of the functor
$\cF_x(B) \to \widehat{\cF}_x(B);$
i.e.\ this functor is essentially surjective.
On the other hand, if $\cF$ admits a Noetherian versal ring $A_x$ at~$x$, and if furthermore the functor $\cF_x(A_x)
\to \widehat{\cF}_x(A_x)$ is essentially surjective, then 
$\cF$ admits an effective Noetherian versal ring at $x$ by definition.

We turn to proving~(2). Suppose that $\eta$ and $\xi$ are two objects
of $\cF_x(A)$ (for some object $A$ of $\widehat{\cC}_{\Lambda}$).
Their product is a morphism $\eta\times\xi: \Spec A \to \cF\times_S \cF$.
The morphisms $\eta$ and $\xi$, when restricted to $\Spec k$
(the closed point of $\Spec A$), both induce the given morphism
$x: \Spec k \to \cF$,
and hence we have the outer square in a commutative diagram
$$\xymatrix{\Spec k \ar^-{x}[d] \ar[r] & \Spec A \ar^-{\eta\times\xi}[d]\ar@{-->}[dl] \\
\cF \ar^-{\Delta}[r] & \cF \times_S \cF}$$
The set of morphisms between $\eta$ and $\xi$ in the category $\cF_x(A)$
may be identified with the set of morphisms $\Spec A \to \cF$ which continue to
make the diagram commute.

If we let $\widehat{\eta}$ and $\widehat{\xi}$ denote the 
images of $\eta$ and $\xi$ in $\widehat{\cF}_x(A)$,
then the set of morphisms between $\widehat{\eta}$ and $\widehat{\xi}$
may similarly be identified with the set of morphisms
$\Spf A \to \cF$ for which the diagram
$$\xymatrix{\Spec k \ar^-{x}[d] \ar[r] & \Spf A \ar^-{\widehat{\eta}\times
\widehat{\xi}}[d]
\ar[dl] \\
\cF \ar^-{\Delta}[r] & \cF \times_S \cF}$$
commutes.

Since $\cF$ satisfies~[3], the fibre product $\cF\times_{\Delta,\cF\times_S \cF, \eta\times\xi} \Spec A$ is an algebraic space over $\Spec A$.  The claim
of full faithfulness in~(2) now follows from the following general fact~\cite[\href{http://stacks.math.columbia.edu/tag/0AQH}{Tag 0AQH}]{stacks-project}:
if $A$ is a complete local ring, if $T \to \Spec A$ is a morphism from
an algebraic space to $\Spec A$, and if $t: \Spec k \to T$ is a section
of this morphism over the closed point, then the restriction
map $\Mor_t(\Spec A,T) \to \Mor_t(\Spf A,T)$ from the set of sections 
over $\Spec A$ extending $t$ to the set of sections over $\Spf A$
extending $t$ is a bijection.\end{proof}

We now prove the lemma promised in Remark~\ref{rem:relationship of versal rings to pro-representability}.

\begin{lem}
\label{lem:versal rings give Art-triv if we know [3]}
Suppose that $\cF$ satisfies~{\em [3]}, and that the diagonal of $\cF$
is furthermore locally of finite type. 
If, for some finite type point $x$ of $\cF$,
there exists a versal morphism $U \to \widehat{\cF}_x$
with $U$ being {\em (}Noetherianly{\em )} pro-representable,
then $R := U\times_{\widehat{\cF}_x} U$
is also {\em (}Noetherianly{\em )} pro-representable, so that $\widehat{\cF}_x$
admits a presentation by a {\em (}Noetherianly{\em )} pro-representable smooth groupoid in functors.
\end{lem}
\begin{proof}Note firstly that~$(U,R)$ is a smooth groupoid in
  functors by Lemma~\ref{lem:relationship of versal rings to
    pro-representability}, so we only need to show that~$R$ is
  (Noetherianly) pro-representable. Suppose that $U$ is pro-represented by  $A= \varprojlim_{i \in I}
A_i$. 
If $B \in \cC_{\Lambda}$, then giving a morphism
$\Spec B \to  U \times_{\widehat{\cF}_x} U$
is equivalent to giving a pair of morphisms
$\Spec B \rightrightarrows U$ which induce the same morphism
to $\cF$.   By definition  these  morphisms factor through
$\Spec A_i$ for some $i$.    Thus $U\times_{\widehat{\cF}_x} U$
is the colimit over $I$
of the fibre products $\Spec A_i \times_{\widehat{\cF}_x} \Spec A_i$
(the fibre product denoting a fibre product
of categories cofibred in groupoids on $\cC_{\Lambda}$).
Since a colimit of pro-representable functors
is pro-representable (by the projective limit of the representing rings), it suffices to show that each of these
fibre products is pro-representable.

If we let $\xi: \Spec A_i \to \cF$ denote the composite
$\Spec A_i \to \Spf A \to \widehat{\cF}_x\hookrightarrow \cF,$
then $\Spec A_i \times_{\cF} \Spec A_i$ is an algebraic
space locally of
finite type over $S$ (because $\cF$ satisfies~[3],
and its diagonal is locally of finite type)
which represents the functor
$\Isom(\xi,\xi)$.
The fibre product $\Spec A_i \times_{\widehat{\cF}_x} \Spec A_i$
is the subfunctor of the restriction
of $\Isom(\xi,\xi)$ 
to $\cC_{\Lambda}$
consisting of isomorphisms which induce the identity from $x$
to itself.  The identity from $x$ to itself is a $k$-valued
point of $\Isom(\xi,\xi)$, and this subfunctor 
may equally well be described as the formal completion of
$\Isom(\xi,\xi)$ at this point. Since $\Isom(\xi,\xi)$ is represented by a
locally Noetherian algebraic space, this formal completion is represented by the formal spectrum
of the complete local ring of~$\Isom(\xi,\xi)$ at~$x$ in the
sense of Definition~\ref{df:complete local rings}
and thus is pro-representable.

It remains to show that if $A$ is Noetherian, then $U\times_{\widehat{\cF}_x} U$
is Noetherianly pro-representable. By~\cite[Prop.\ 5.1]{MR1603480}, we need to
check that if $k$ is the residue field at $x$, then the $k$-vector space of morphisms
$\Spec k[\epsilon]/\epsilon^2\to U\times_{\widehat{\cF}_x} U $ is
finite-dimensional. 

Now, any such morphism factors through a pair of morphisms
$\Spec A/\m_A^2\to U$.
Since $A$ is Noetherian, $A/\m_A^2$ is Artinian. Applying the argument of the previous
paragraph to $\Spec A/\m_A^2
\times_{\widehat{\cF}_x} \Spec A/\m_A^2$, we find that it is
pro-represented by a Noetherian complete local ring (the completion of
a locally Noetherian algebraic space at a closed point); any such ring admits only a finite-dimensional space of
maps to $k[\epsilon]/\epsilon^2$, so we are done.
\end{proof}\begin{cor}\label{cor: 3 and 1 and versal rings implies arttriv.}
  Suppose that $\cF$ satisfies~{\em [3]},  that its diagonal is furthermore locally of finite type,
and that $\cF$ admits versal
rings at all finite type points. Then $\cF$ is $\mathbf{Art}^{\mathbf{triv}}$-homogeneous.
\end{cor}
\begin{proof}This is immediate from
Lemmas~\ref{lem:versal rings give Art-triv if the relation is Art-triv}
and~\ref{lem:versal rings give Art-triv if we know [3]}.
  \end{proof}

The next result shows that
when $\cF$ satisfies~[2] and~[3],
we can strengthen condition~[2](a)
so as to allow $Y'$ and $Z$ to be complete
Noetherian local rings, rather than merely Artinian.

\begin{lemma}
\label{lem:complete homogeneity}
Suppose that~$\cF$ satisfies~{\em [3]} and~{\em [2](b)},
and that
$$\xymatrix{Y \ar[r] \ar[d] & Y' \ar[d] \\
Z \ar[r] & Z'}
$$
is a pushout diagram of $S$-schemes, with the horizontal arrows
being closed immersions, $Y$ being a finite type Artinian $\mathcal O_S$-scheme,
each of $Z$ and $Y'$ being the spectrum of a complete Noetherian local
$\mathcal O_S$-algebra whose residue field is finite type over $\mathcal O_S$,
and the left-hand vertical arrow being closed {\em (}i.e.\
corresponding to a local morphism of local $\mathcal O_S$-algebras{\em )}.
If either $\cF$ satisfies~{\em [2](a)}, or if~$\cF$ is an \'etale
stack  whose diagonal is locally of finite presentation and   the extension
of the residue field of $Y$ over the residue field of $Z$ is separable,
then the induced functor
$$\cF(Z') \to \cF(Y') \times_{\cF(Y)} \cF(Z)$$
is an equivalence of categories.
\end{lemma}
\begin{proof}Let $x$ be the underlying closed point of~$Y$. By
  Lemma~\ref{lem:fibre products of local rings}, $Z'$ is Noetherian (note that $Y\into Y'$
  is assumed to be a closed immersion, and the residue fields of~$Y,
  Y', Z$ are all finite type $\cO_S$-algebras, so that in particular
  the residue field of~$Y$ is a finite extension of that of~$Z$), so by Lemma~\ref{lem:characterizing pro-representability} and the
  assumption that~$\cF$ satisfies~[3] and~[2](b), we may replace~$\cF$
  by~$\widehat{\cF}_x$, and by the definition of~$\widehat{\cF}_x$, we can then
  reduce to the case that~$Y'$ and~$Z$ are Artinian. In the case that~$\cF$
  satisfies~[2](a) we are then done by the very definition of that condition,
and in the case that~$\cF$ is an \'etale stack whose diagonal is
locally of finite presentation and the extension of the residue field of~$Y$ over the residue field of~$Z$ is
  separable the result follows from Corollary~\ref{cor: 3 and 1 and
    versal rings implies arttriv.} and Lemma~\ref{lem:homogeneity conditions}.\end{proof}

\subsection{Artin's representability theorem}\label{subsec: Artin
  representability}
In this subsection, we continue to assume that $S$ is locally Noetherian.
In the statement of the next lemma, as well as for one direction of the main theorems,
we will suppose further that, for each
finite type point $s\in S$, the local ring $\cO_{S,s}$ is a $G$-ring. As
explained in Section~\ref{subsubsec:Grings}, this assumption is needed in order
to apply Artin approximation.

We begin with a lemma, which is essentially a rephrasing
of~\cite[\href{http://stacks.math.columbia.edu/tag/07XH}{Tag
  07XH}]{stacks-project}. 
It is the key application of Artin approximation which underlies
Artin's representability theorem.

\begin{lemma}
\label{lem:formally smooth charts}
Suppose that, for each
finite type point $s\in S$, the local ring $\cO_{S,s}$ is a $G$-ring.
If $\cF$ is a category in groupoids satisfying {\em [1]} and~{\em [2](a)}, 
if $k$ is a finite type $\cO_S$-field,
and if $x: \Spec k \to \cF$ is a morphism representing the finite type point $t \in | \cF|$, for which $\widehat{\cF}_x$
admits an effective Noetherian versal ring, 
then there exists a morphism $\varphi:U \to \cF$ whose source is 
a scheme locally of finite type over $S$, and a finite type point
$u \in U$ such that $\varphi(u) = t$,
and such that $\varphi$ is formally smooth at $u$.
\end{lemma}
\begin{proof}
It follows
from~\cite[\href{http://stacks.math.columbia.edu/tag/07XH}{Tag 07XH}]{stacks-project}
(i.e.\ Artin approximation),
together with our assumptions, that we may find a morphism $\varphi:U \to \cF$
with $U$ a scheme locally of finite type over $S$,
containing a point $u,$ such that $\phi(u) = t$,
and such that $\varphi$ is versal at $u$.
Lemma~\ref{lem:versions of smoothness}~(2), together with our assumption
of~[2](a), shows that in fact $\varphi$ is formally smooth at~$u$.
\end{proof}

\begin{lemma}
\label{lem:criterion for openness of versality}
Suppose that $\cF$ is a category fibred in groupoids which satisfies~{\em [3]},
and that there exists a morphism $\varphi:T \to \cF$ whose source
is a scheme, locally of finite type over~$S$,
which is smooth in a neighbourhood of a finite type point 
$t \in T$.  Then, if $\varphi':T' \to \cF$ is any morphism from
a scheme, locally of finite type over~$S$, to $\cF$,
which is formally smooth at a finite type point 
$t'\in T'$ whose image in $\cF$ coincides with the image of $t$,
there is a neighbourhood $V$ of $t'$ in $T'$ such that
the restriction of $\varphi'$ to $V$ is smooth.
\end{lemma}
\begin{proof}
Replacing $T$ by the hypothesised neighbourhood of $t$,
we may assume that $\varphi$ is smooth.  Note that since $\cF$
satisfies~[3], the morphism $\varphi$ is representable by algebraic
spaces, and so this is to be understood in the sense
of~\cite[\href{http://stacks.math.columbia.edu/tag/03YK}{Tag
03YK}]{stacks-project}, i.e.\ the base change of this morphism
over any morphism from a scheme to $\cF$ is smooth. In particular, the base-changed morphism of algebraic spaces
$T\times_{\cF} T' \to T'$ is smooth.  
Since the morphism $T' \to \cF$ is formally smooth at $t',$
the projection
$T\times_{\cF} T' \to T$
is smooth in a neighbourhood $U$ of the point $(t,t')$, 
by Lemma~\ref{lem:versions of smoothness}~(3). 
Now the the composite $U \to T \to \cF$ is the composite
of smooth morphisms, hence is smooth.  Rewriting this morphism
as $U \to T' \to \cF$, we see that this composite is also smooth.
If we let $V$ denote the image of $U$ in $T'$, then $V$ is an open
subset of $T'$ containing $t'$, and $U \to V$ is a smooth surjection.
Smoothness being a property that is smooth local on the source,
we see that $V \to \cF$ is a smooth morphism, as required.
\end{proof}

\begin{remark}
\label{rem:verifying axiom 4}
The preceding
lemma shows that in order to verify~[4] for a category fibred in groupoids 
satisfying~[2](a) and~[3], it suffices to find, for each finite type point
of $x \in |\cF|$, a morphism $\varphi:T \to \cF$ whose source is locally
of finite type over $S$, and contains a finite type point $t$
mapping to $x$, such that $\varphi$ restricts to a smooth morphism
on some neighbourhood of $t$.
(The role of~[2](a) is to ensure,
via Lemma~\ref{lem:versions of smoothness}~(2), 
that versality at a point coincides with formal smoothness at a point.)
\end{remark}

\begin{thm}
  \label{thm: Artin representability}
Suppose that $S$ is locally Noetherian.
Any algebraic stack, locally of finite presentation
over $S$, satisfies Axioms \emph{[1]}, \emph{[2]}, \emph{[3]} and \emph{[4]}. 
Conversely, suppose further that for each
finite type point $s\in S$, the local ring $\cO_{S,s}$ is a $G$-ring. 
Then if $\cF$ is an \'etale stack in groupoids over $S$ satisfying
\emph{[1]}, \emph{[2]}, \emph{[3]} and \emph{[4]}, then
  $\cF$ is an algebraic stack, locally of finite presentation
over $\cF$. 
\end{thm}
\begin{proof}
If $\cF$ is an algebraic stack, locally of finite presentation over $S$,
then it follows from
Lemma~\ref{lem: limit preserving iff locally of finite presentation}
that $\cF$ satisfies~[1], while Lemma~\ref{lem: algebraic stacks satisfy RS},
Lemma~\ref{lem:characterizing pro-representability},
\cite[\href{http://stacks.math.columbia.edu/tag/07WU}{Tag
  07WU}]{stacks-project},
\cite[\href{http://stacks.math.columbia.edu/tag/07WV}{Tag
  07WV}]{stacks-project},
\cite[\href{http://stacks.math.columbia.edu/tag/06IW}{Tag
  06IW}]{stacks-project},
\cite[\href{http://stacks.math.columbia.edu/tag/07X8}{Tag
  07X8}]{stacks-project}, and~
\cite[\href{http://stacks.math.columbia.edu/tag/07X1}{Tag 07X1}]{stacks-project} 
show that $\cF$ satisfies~[2]. (More precisely,~$\cF$ satisfies~[2](a)
by Lemma~\ref{lem: algebraic stacks satisfy
  RS}. By~\cite[\href{http://stacks.math.columbia.edu/tag/07X8}{Tag
  07X8}]{stacks-project} and Lemma~\ref{lem:characterizing
  pro-representability}, in order to prove that $\cF$ satisfies~[2](b), it suffices to show
that $\cF$ has Noetherian versal rings at each finite type
point. \cite[\href{http://stacks.math.columbia.edu/tag/06IW}{Tag
  06IW}]{stacks-project} gives a criterion for the existence of such
versal rings, which is satisfied by
\cite[\href{http://stacks.math.columbia.edu/tag/07WU}{Tag
  07WU}]{stacks-project},
\cite[\href{http://stacks.math.columbia.edu/tag/07WV}{Tag
  07WV}]{stacks-project}, and~
\cite[\href{http://stacks.math.columbia.edu/tag/07X1}{Tag
  07X1}]{stacks-project}.)  By definition~$\cF$ satisfies~[3]. Again
by definition, we may find a smooth surjection $U \to \cF$ with $U$ a
scheme.  If $x \in |\cF|$ is a point of finite type, then we may find
a finite type point $u \in U$ lying over $x$, and
Remark~\ref{rem:verifying axiom 4} then shows that $\cF$
satisfies~[4].

For the converse,
we follow the proof of~\cite[\href{http://stacks.math.columbia.edu/tag/07Y4}{Tag
    07Y4}]{stacks-project}. By definition, we need to show that $\cX$ admits a
  smooth surjection from a scheme. Taking the union of the morphisms 
obtained
  from Lemma~\ref{lem:formally smooth charts} as we run over all finite type
  points of $\cX$,  we obtain a smooth map $U\to\cX$ whose source is a scheme,
  whose image contains all finite type points of $\cX$. It remains to show that
  this is surjective. As in the proof of~\cite[\href{http://stacks.math.columbia.edu/tag/07Y4}{Tag
    07Y4}]{stacks-project}, this may be checked by pulling back to affine
  schemes of finite presentation over $S$, where it is immediate (as smooth maps
  are open, and the finite type points of a scheme are dense).
\end{proof}

\begin{thm}
  \label{thm: Artin representability; variant}
Any algebraic stack, locally of finite presentation
over $S$, satisfies
\emph{[1]}, \emph{[2]}, \emph{[3]}, \emph{[4a]}, and~{\em [4b]}.
Conversely, suppose further that for each
finite type point $s\in S$, the local ring $\cO_{S,s}$ is a $G$-ring. Then if $\cF$ is an \'etale stack in groupoids over $S$ satisfying
\emph{[1]}, \emph{[2]}, \emph{[3]}, \emph{[4a]}, and~{\em [4b]},
and if the diagonal of $\cF$ is furthermore quasi-compact,
then $\cF$ is an algebraic stack,
locally of finite presentation over $S$. \end{thm}
\begin{proof}For the first statement, in view of Theorem~\ref{thm: Artin
    representability}, we need to show that~[4a] and~[4b] are satisfied for
  algebraic stacks. Lemma~\ref{lem:algebraic stacks satisfy 4a} shows that~[4a]
  is satisfied, while Lemma~\ref{lem:versions of smoothness}~(4) shows
that~[4b] is satisfied.

To prove the second statement, note that
Corollary~\ref{cor:4a and b give 4} shows, under the given hypotheses
on $\cF$, that $\cF$ furthermore satisfies Axiom~[4].  
Theorem~\ref{thm: Artin representability} then shows that
$\cF$ is an algebraic stack, locally of finite presentation over $S$.
\end{proof}

Example~\ref{ex: line with infinitely many nodes} shows that the condition on the diagonal of $\cF$
is necessary in order to deduce that a stack
satisfying [1], [2], [3], [4a], and [4b] is algebraic.

\begin{remark}
In the proof of Artin representability, we don't require the full strength
of Axiom~[2](b); all we need is that~$\widehat{\cF}_x$ admits an
effective Noetherian versal
ring for at least one
morphism $x$ representing any given finite type point of $\cF$.
However, the following result, which describes the extent to which
this  hypothesis is independent of the choice
of representative of a finite type point, implies in particular
that, in the context of the preceding representability theorems,
this  hypothesis holds for at least one representative of
a given finite type point if and only if it holds for every such
representative.
\end{remark}

\begin{lemma}
\label{lem:independence of morphism for 2b}
  Suppose that $\cF$ satisfies~{\em [2](a)},
and suppose given $x: \Spec k \to \cF$, with $k$ a field 
which is a finite type $\cO_S$-algebra.
\begin{enumerate}
\item
If $\widehat{\cF}_x$
admits a Noetherian versal ring,
then for any other morphism $x': \Spec l\to
\cF$, with~ $l$ a field which is a finite type $\cO_S$-algebra, representing
the same point of~$\cF$, we have that $\widehat{\cF}_{x'}$ again
admits a Noetherian versal ring.
\item
If $\widehat{\cF}_x$
admits an  effective
Noetherian versal ring,
then for any finite extension $l$ of $k$, if we let $x'$ denote
the composite $\Spec l\to \Spec k \buildrel x \over \longrightarrow \cF$,
we have that $\widehat{\cF}_{x'}$ again
admits an effective Noetherian versal ring.
\item
Suppose that $\cF$ satisfies~{\em [3] (}in addition to {\em [2](a))},
and that, for some finite  extension $l$ of $k$, 
if we denote by $x'$ the composite $\Spec l \to \Spec k \buildrel x \over 
\longrightarrow \cF$, we have that $\widehat{\cF}_{x'}$ 
admits an effective
Noetherian versal ring. Assume also either that
$l/k$ is separable, or that $\cF$ is an fppf stack.
Then $\widehat{\cF}_x$ admits an effective
Noetherian versal ring.
\item
Suppose either that
{\em (a)} $\cF$ satisfies~{\em [3] (}in addition to {\em [2](a))},
and that either the residue field of the image of the composite
$\Spec k \to \cF \to S$ is perfect, or $\cF$ is an fppf stack; 
or {\em (b)} that 
$\cF$ satisfies~{\em [1]} and~{\em [3] (}in addition to {\em [2](a))},
and that the local rings of $S$ at finite type points are $G$-rings.
Suppose also that $\widehat{\cF}_x$
admits an  effective
Noetherian versal ring.
Then for any other morphism $x': \Spec l\to
\cF$, with $l$ a field which is a finite type $\cO_S$-algebra, representing
the same point of~$\cF$, we have that $\widehat{\cF}_{x'}$ again
admits an  effective
Noetherian versal ring.
\end{enumerate}
\end{lemma}
\begin{proof}
Given $l$ as in~(1), we may find a common finite extension $l'$ of $k$ and $l$.
After appropriate relabelling, then, to prove~(1)
we may assume that $l$ is a
finite extension of $k$ and that $x'$ is the composite $\Spec l \to \Spec k
\to \cF$, and we must show that $\widehat{\cF}_x$ admits a Noetherian
versal ring if and only
if $\widehat{\cF}_{x'}$ does.

We use the notation and results of Remark~\ref{rem: material on change of residue
  field}. Since~$\cF$ satisfies~[2](a), there is a natural equivalence of categories cofibred in groupoids 
$(\widehat{\cF}_{x})_{l/k}  \iso \widehat{\cF}_{x'}$.
We will deduce part~(1) of the lemma from
from~\cite[\href{http://stacks.math.columbia.edu/tag/06IW}{Tag 06IW}]{stacks-project}
(which gives a criterion for a category cofibred in groupoids over
$\cC_{\Lambda,k}$, or $\cC_{\Lambda,l}$, to admit a
Noetherian versal morphism), as follows. Since~$\cF$ satisfies~[2](a), it
follows from~\cite[\href{http://stacks.math.columbia.edu/tag/07WA}{Tag
  07WA},\href{http://stacks.math.columbia.edu/tag/06J7}{Tag 06J7}]{stacks-project} that both~$\widehat{\cF}_{x}$ and~
$(\widehat{\cF}_{x})_{l/k}$ satisfy all of the hypotheses
of~\cite[\href{http://stacks.math.columbia.edu/tag/06IW}{Tag
  06IW}]{stacks-project}, except possibly the condition of having a
finite-dimensional tangent
space. By~\cite[\href{http://stacks.math.columbia.edu/tag/07WB}{Tag
  07WB}]{stacks-project}, this condition holds for $\widehat{\cF}_{x}$
if and only if it holds for $(\widehat{\cF}_{x})_{l/k}$, so it
suffices to observe that the existence of a Noetherian versal ring
implies the finite-dimensionality of the tangent space, by~\cite[\href{http://stacks.math.columbia.edu/tag/06IU}{Tag 06IU}]{stacks-project}.

To prove part~(2), we again use the equivalence $(\widehat{\cF}_x)_{l/k}
\iso \widehat{\cF}_{x'}$.
Taking into account part~(1), it suffices to show that if $A_{x'}$ denotes
a versal ring at the morphism~ $x'$, then the versal morphism
$\Spf A_{x'} \to \widehat{\cF}_{x'}$ arises from a morphism
$\Spec A_{x'} \to \cF$.  The above-mentioned equivalence of categories
shows that this versal morphism arises as the composite 
of the natural morphism $\Spf A_{x'} \to \Spf (A_{x'} \times_l k)$
and a morphism $\Spf (A_{x'} \times_l k) \to \widehat{\cF}_x$.
Lemma~\ref{lem:characterizing pro-representability}~(1),
and our effectivity assumption regarding~$\widehat{\cF}_x$,
shows that this latter morphism is effective, i.e.\ is induced 
by a morphism
$\Spec (A_{x'} \times_l k)  \to \cF$. 
The composite $$\Spec A_{x'} \to \Spec (A_{x'} \times_l k) \to
\Spec A_x \to \cF$$ then induces the original versal morphism 
$\Spf A_{x'} \to \widehat{\cF}_{x'}$, and so~(2) is proved.

In order to prove~(3), we first fix a versal morphism
$\Spf A_x \to \cF$ at $x$. If $l/k$ is separable, then, by the topological invariance of
the \'etale site, we may find a finite \'etale local extension
$B$ of $A_x$ which induces the extension $l/k$ on residue
fields. Otherwise, let $\Lambda_k$, $\Lambda_l$ be Cohen rings for
$k$, $l$ respectively, and set $B:=A_x\widehat{\otimes}_{\Lambda_k}\Lambda_l$;
then $B/A_x$ is a finite extension by the topological version of
Nakayama's lemma, and it is a faithfully flat local extension of
complete local Noetherian rings by \cite[Lem.\ 0.19.7.1.2]{MR0173675}.  In either case, the composite $\Spf B \to \Spf A_x \to \cF$
may be regarded as a morphism $\Spf B \to \widehat{\cF}_{x'}$,
and so, by our effectivity assumption regarding $\widehat{\cF}_{x'}$,
together with Lemma~\ref{lem:characterizing pro-representability}~(1),
it is effective,
i.e.\ is induced by a morphism $\Spec B \to \cF$.

Consider the pull-backs
of this morphism along the two projections $\Spec B\otimes_{A_x} B
\rightrightarrows \Spec B.$  Since the morphism $\Spf B \to \cF$ factors
through $\Spf A_x$, we see that the two pull-backs become isomorphic 
over $\Spf (B\otimes_{A_x} B)$.  Since $\cF$ satisfies~[3],
Lemma~\ref{lem:characterizing pro-representability}~(2)
shows that the two pull-backs are themselves isomorphic.
We may check the cocycle condition in a similar way, and hence
obtain \'etale (or \emph{fppf}, in the case that $l/k$ is inseparable)
descent data on the morphism $\Spec B \to \cF$, which
(since $\cF$ is an \'etale stack, and an \emph{fppf} stack if $l/k$ is
inseparable) allows us to descend it to a morphism
$\Spec A_x \to \cF$, as required.

We turn to proving~(4).
Thus we suppose given
a morphism $x':\Spec l \to \cF$ representing the same point of $\cF$
that the given morphism $x: \Spec k \to \cF$ represents.  We are assuming
that $\widehat{\cF}_x$
admits an effective Noetherian versal ring, and we wish
to prove the corresponding fact for $\widehat{\cF}_{x'}$.  
If we let $l'$ be a common finite extension of $k$ and $l$,
and let $x''$ denote the composite $\Spec l' \to \Spec k\to \widehat{\cF}_x$,
then it follows from~(2) that $\widehat{\cF}_{x''}$ 
admits an effective Noetherian versal ring. 
So, relabelling $l'$ as $l$ and $l$ as $k$, we are reduced to the 
following problem: we are given the morphisms $x: \Spec k \to \cF$
and $x' : \Spec l \to \Spec k \buildrel x \over \longrightarrow \cF$,
for which
$\widehat{\cF}_{x'}$ 
admits an effective Noetherian versal ring, 
and we would like to conclude that 
$\widehat{\cF}_x$ also
admits an effective Noetherian versal ring.

If either the residue field of the image of 
$\Spec k$ in $S$ is perfect or $\cF$ is an \emph{fppf} stack, and if $\cF$ satisfies~[3],
then this follows from~part~(3).
Otherwise, we assume further that $\cF$ satisfies~[1],
and that the local rings of $S$ at finite type points are $G$-rings.
This allows us to apply Artin approximation, in the form of
Lemma~\ref{lem:formally smooth charts}, to the morphism $x'$, so as to
conclude that there exists a morphism
$\varphi: U \to \cF$ whose source is a scheme,
and a finite type point $u$ of $U$, lying over the image of $x'$ (which is
also the image of $x$), such that $\varphi$ is formally smooth at $u$.

At this stage $x'$ and $l$ have done their job in the argument,
and so we drop them from consideration; in fact,
we will recycle them as notation, in a manner which we now explain.
Pulling back $U$ over $x$, we obtain a $k$-scheme whose smooth locus
contains the fibre over $U$.  This fibre is non-empty,
by Lemma~\ref{lem:versions of smoothness}~(2),
and thus contains a point defined over a finite separable extension
$l$ of $k$.   If we let $x'$ denote the resulting
composite $\Spec l \to \Spec k \to \cF$,
then $\widehat{\cF}_{x'}$
admits an effective Noetherian versal ring (given by the complete
local ring of this fibre at this point).
Part~(3) now implies that the same is true for~$\widehat{\cF}_x$,
as required.
\end{proof}

\begin{rem}
  A version of Lemma~\ref{lem:independence of morphism for 2b} for the
  condition of admitting a presentation by a smooth Noetherianly
  pro-representable groupoid in functors, rather than admitting a
  Noetherian versal ring, can be proved in an almost identical fashion.
\end{rem}

\begin{remark}
Example~\ref{ex: conjugate lines meeting to infinite order}
shows that  parts~(3)
and~(4) of 
Lemma~\ref{lem:independence of morphism for 2b} are not true without the
assumption of~[3]. 
\end{remark}

\section{Scheme-theoretic images}
\label{sec:images}
Suppose that $\xi:\cX\to\cF$ is a proper morphism, where~$\cX$ is an
algebraic stack, and~$\cF$ is a stack whose diagonal is
representable by algebraic spaces and locally of finite
presentation. In Section~\ref{subsec:scheme images two} we define the
scheme-theoretic image~$\cZ$ of~$\xi$, which is initially a Zariski
substack of~$\cF$. Our main aim in this section is to prove
Theorem~\ref{thm:main-intro}, giving a criterion for~$\cZ$ to be an
algebraic stack, as well as to prove a number of related properties
of~$\cZ$.

\smallskip

{\em Interpreting Theorem}~\ref{thm:main-intro} {\em
as taking the quotient by a proper equivalence relation.}
Whether or not $\cZ$ satisfies~[2],
we can show (under mild hypotheses on~$\cF$) that the morphism $\xi:\cX \to \cF$ factors through 
a morphism $\overline{\xi}: \cX \to \cZ$ (see Lemma~\ref{lem:proper and surjective} below), and this morphism is ``scheme-theoretically
surjective''. If we define $\cR := \cX \times_{\cF} \cX$, then $\cR$ is an algebraic stack
(because $\cF$ satisfies~[3]) which defines a proper equivalence relation on $\cX$.
Thus, at least morally, we may regard $\cZ$ as the quotient of $\cX$ by
the equivalence relation $\cR$, and
Theorem~\ref{thm:main-intro} may be regarded as providing a context in
which the quotient $\cX/\cR$ may be defined as an algebraic stack.

Note that in general the quotient of an algebraic stack, or a scheme,
or even a variety, by a proper equivalence relation,
may not admit a reasonable interpretation as an object of algebraic
geometry. (See e.g.\ Examples~\ref{ex:sphere} and~\ref{ex:endos of line}.)
Our result shows
that when the desired quotient admits an interpretation as
the scheme-theoretic image of a morphism from $\cX$ to some stack $\cF$,
the quotient does indeed have a chance to be of an algebro-geometric nature.

One well-known theorem which concerns taking a quotient 
by a proper equivalence relation is Artin's result~\cite{MR0260747}
on the existence of contractions.
We close this discussion by briefly describing the relationship
between~\cite{MR0260747} and the present note.  

\smallskip

{\em The relationship with}~\cite{MR0260747}.
In~\cite{MR0260747}, Artin proves the existence of 
dilatations and contractions
of an algebraic space along a closed algebraic subspace,
given a formal model for the desired dilatation or contraction. (See
also Remark~\ref{rem: remark about Artin contractions} below.) 
In the case of contractions, his result can be placed
in the framework of Corollary~\ref{cor:representability-intro}.
Namely, taking $\cX$ in that theorem to be the algebraic
space on which one wishes to perform a contraction,
one can define a functor $\cF$ which is supposed to
represent the result of the contraction, and a proper morphism
$\xi: \cX \to \cF$.   One can then show that $\cF$ satisfies~[1],
[2], and~[3], and hence conclude
(via Corollary~\ref{cor:representability-intro})
that $\cF$
is representable by an algebraic space.   This is in fact
essentially how Artin proceeds, although there are some slight
differences between our approach and his:  Artin defines $\cF$
concretely  on the category of {\em Noetherian} $\cO_S$-algebra,
and then extends it to arbitrary $\cO_S$-algebras by taking limits.
Our approach would be to apply Artin's concrete definition only
to {\em finite type} $\cO_S$-algebras, and then to extend to
arbitrary $\cO_S$-algebras by taking limits.  The difference
between the two approaches is that in Artin's approach, the fact that
condition~[2] is satisfied by $\cF$ is rather automatic, whereas the proof
that $\cF$ satisfies~[1] becomes the crux of the argument.
In our approach, while the arguments remain the same, their interpretation
differs: essentially by definition, $\cF$ will satisfy~[1], while
the bulk of the argument can be seen as proving that it also satisfies~[2].
However, this attempt to link Artin's result to ours is a little
misleading, since in Artin's context, the verification that $\cF$
satisfies~[4] is straightforward.   
For some other $\cF$ (such as those that appear in the theory of 
moduli of Galois representations discussed in Section~\ref{sec: finite flat}),
the verification of~[4] seems to be less straightforward,
however,
and indeed the only approach we know is via the general arguments of
the present paper.

\subsection{Scheme-theoretic images (part one)}
\label{subsec:images part one}
We recall the following definition (see for
example~\cite[\href{http://stacks.math.columbia.edu/tag/01R8}{Tag
  01R8}]{stacks-project}; note that we do not include a quasi-separated
hypothesis, as it is not needed for the basic properties of scheme-theoretic
images that we use here).

\begin{defn}
\label{def:scheme-theoretic image and dominant for schemes}
Let $f: Y \to Z$ be a quasi-compact morphism of schemes.
The kernel of the natural morphism $\cO_Z \to f_*\cO_Y$
is a quasi-coherent ideal sheaf $\mathcal I$ on $Z$ (\cite[\href{http://stacks.math.columbia.edu/tag/01R8}{Tag
  01R8}]{stacks-project})
and we define 
the {\em scheme-theoretic image of $f$}
to be the closed subscheme $V(\mathcal I)$ of $Z$ cut out
by $\mathcal I$. 

We say that $f$ is {\em scheme-theoretically dominant}
if the induced morphism $\cO_Z \to f_* \cO_Y$ is injective; that is, if the
scheme-theoretic image of $f$ is $Z$.

We say that $f$ is {\em scheme-theoretically surjective} if it is
 scheme-theoretically dominant, and surjective on underlying
 topological spaces.
\end{defn}

\begin{remark}
\label{rem:scheme-theoretically dominant}
(1)  The morphism $f$ factors through $V(\mathcal I)$,
and $V(\mathcal I)$ is the minimal closed subscheme
of $Z$ through which $f$ factors. The induced morphism $f': Y \to V(\mathcal I)$
is scheme-theoretically dominant, and also
has dense image \cite[\href{http://stacks.math.columbia.edu/tag/01R8}{Tag
  01R8}]{stacks-project}.  

Thus, for a closed morphism, the notions of
scheme-theoretical dominance and of scheme-theoretical surjectivity
are equivalent.  

(2) The formation of scheme-theoretic images is compatible with arbitrary flat
base change~\cite[\href{http://stacks.math.columbia.edu/tag/081I}{Tag
  081I}]{stacks-project}.

(3) It follows easily from (2) that the formation of scheme-theoretic images is
{\em fpqc} local on the target, so that
the condition of being scheme-theoretically dominant, or surjective, may be checked
{\em fpqc} locally on the target.

(4) If $g:X\to Y$ is quasi-compact then the scheme-theoretic image of $fg$ is a
closed subscheme of the scheme-theoretic image of $f$. If $g$
is furthermore scheme-theoretically dominant
(e.g.\ {\em fpqc}),
then $\cO_Z \to f_* \cO_Y$ and $\cO_Z \to (fg)_* \cO_X$ have
the same kernel, and hence the scheme-theoretic images of $f$ and $fg$ coincide.

In particular, if $g: X \to Y$ is quasi-compact and scheme-theoretically dominant (resp.\
surjective),
then $f$ is scheme-theoretically dominant (resp.\ surjective)
if and only if the composite $fg$ is
scheme-theoretically dominant (resp.\ surjective).

(5)
An {\em fpqc} morphism is scheme-theoretically surjective.
\end{remark}

Points~(3), (4), and~(5) of the preceding remark
allow us to extend the notion of scheme-theoretically
dominant (resp.\ scheme-theoretically surjective) morphisms to the context
of morphisms of algebraic stacks in the following way.

Recall (\emph{cf.}\ Remark~\ref{rem: various properties of morphisms of algebraic stacks} and Lemma~\ref{lem:characterizing morphisms}) that an algebraic stack $\cY$ is quasi-compact if its underlying
topological space $|\cY|$ is quasi-compact,
or equivalently,
by~\cite[\href{http://stacks.math.columbia.edu/tag/04YC}{Tag 04YC}]{stacks-project}, if
there is a smooth surjection $U\to \cY$ with $U$ a quasi-compact scheme. A morphism of algebraic stacks $f:\cY\to\cZ$ is \emph{quasi-compact}
if for every morphism $\cV\to\cZ$ with $\cV$ a quasi-compact algebraic stack,
the fibre product $\cY\times_\cZ\cV$ is also quasi-compact \cite[\href{http://stacks.math.columbia.edu/tag/050U}{Tag 050U}]{stacks-project}.

\begin{defn}
  \label{defn: scheme-theoretic dominance for algebraic stacks}Let $f:\cY\to\cZ$ be
  a quasi-compact morphism of algebraic stacks. Let $V\to\cZ$ be a smooth
  surjection from a scheme, and let $V=\cup_i T_i$ be a cover of $V$ by 
  quasi-compact open subschemes. For each $T_i$, the fibre product $\cY\times_\cZ T_i$
  is quasi-compact, so admits a smooth surjection $U_i\to \cY\times_\cZ T_i$
  from a quasi-compact scheme. 

The composite morphism $U_i\to T_i$ is
  quasi-compact, and we say that $f:\cY\to\cZ$ is \emph{scheme-theoretically
dominant} if for all $i$, the morphism $U_i\to T_i$ is scheme-theoretically
dominant. We say that $f$ is \emph{scheme-theoretically surjective} if the
morphisms $U_i\to T_i$ are all scheme-theoretically
surjective.

It follows from Remark~\ref{rem:scheme-theoretically dominant}~(3)--(5) 
that this notion is well-defined,
independently of the choices of $V$ and the $T_i$ and $U_i$, and
that it agrees with Definition~\ref{def:scheme-theoretic image and dominant for
  schemes} if $\cY$ and $\cZ$ are schemes.
\end{defn}
Similarly, we can extend the definition of scheme-theoretic images to
quasi-compact morphisms of algebraic stacks.

\begin{defn}
  \label{defn: scheme theoretic image for algebraic stacks} Let $f:\cY\to\cZ$ be
  a quasi-compact morphism of algebraic stacks, and choose $V$, $T_i$ and $U_i$
  as in Definition~\ref{defn: scheme-theoretic dominance for algebraic stacks}.

  The composite morphism $U_i\to T_i$ is quasi-compact, hence admits a
  scheme-theoretic image $T'_i$.  The smooth equivalence relation $R_{T_i}:=
  T_i\times_{\cZ} T_i$ on $T_i$ restricts to a smooth equivalence relation $R_{T'_i}$ on
  $T'_i$, and the quotient stack $[T_i'/R_{T_i'}]$ is a closed substack of the
  quotient stack $[T_i/R_{T_i}]$, itself an open substack of $\cZ$.  The substacks
  $[T_i'/R_{T_i'}]$ glue together to form a closed substack of $\cZ$, which we
  define to be the {\em scheme-theoretic image} of $f$.

Again, it follows from Remark~\ref{rem:scheme-theoretically dominant}~(2)--(5) 
that this notion is well-defined,
independently of the choices of $V$ and the $T_i$ and $U_i$, and
that it agrees with Definition~\ref{def:scheme-theoretic image and dominant for
  schemes} if $\cY$ and $\cZ$ are schemes.
\end{defn}
Given this definition, the following remarks are an easy consequence of
Remark~\ref{rem:scheme-theoretically dominant}.
\begin{remark}
\label{rem:properties of scheme-theoretic images of algebraic stacks}(1)  The morphism $f$ factors through its scheme theoretic image $\cX$, and $\cX$ is the minimal closed substack
of $\cZ$ through which $f$ factors. The induced morphism $f': \cY \to \cX $
is scheme-theoretically dominant, and the induced map $|\cY|\to|\cX|$ has dense image. 

(2) The formation of scheme-theoretic images is compatible with arbitrary flat
base change.

(3) The formation of scheme-theoretic images is
{\em fpqc} local on the target, so that
the condition of being scheme-theoretically dominant, or surjective, may be checked
{\em fpqc} locally on the target.

(4) If $g:\cX\to \cY$ is quasi-compact then the scheme-theoretic image of $fg$ is a
closed substack of the scheme-theoretic image of $f$. If $g$ is also
scheme-theoretically dominant, then the scheme-theoretic images of $f$ and $fg$ coincide.

In particular, if $g: \cX \to \cY$ is quasi-compact and scheme-theoretically dominant (resp.\
surjective),
then $f$ is scheme-theoretically dominant (resp.\ surjective)
if and only if the composite $fg$ is
scheme-theoretically dominant (resp.\ surjective).

\begin{rem}
  \label{rem: we don't want to give the quasi-coherent sheaf definition of
    scheme-theoretic image for stacks}One can show that if
  $f:\cY\to\cZ$ is a quasi-compact and quasi-separated morphism of algebraic stacks, then the kernel
  of $\cO_{\cZ}\to f_*\cO_{\cY}$ is a quasi-coherent ideal sheaf, and cuts out
  the scheme-theoretic image of $f$ as defined above. (In fact, presumably
this is true even without the quasi-separatedness assumption,
although, as in the scheme-theoretic context, such a statement would
be slightly more delicate to prove, since then $f_*\cO_{\cY}$ 
need not be quasi-coherent.)  However, since we do not
  need to consider sheaves on stacks elsewhere in this paper, except very
  briefly in the proof of Lemma~\ref{lem:image conditions} below, we will not
  give the details here, but rather simply establish the few basic facts that we
  need as we use them.

  The one special case of this theory that we need is that if $f:\cX\to S$ is a
  quasi-compact algebraic stack over a Noetherian base scheme $S$, then $f$ is scheme-theoretically dominant if
  and only if $\cO_S\to f_*\cO_\cX$ is injective. To see this, let $g:U\to\cX$ be a
  smooth cover by a quasi-compact scheme; then by Definition~\ref{defn:
    scheme-theoretic dominance for algebraic stacks}, $f$ is
  scheme-theoretically dominant if and only if the composite $fg:U\to\cX\to S$
  is scheme-theoretically dominant, or equivalently if and only if
  $\cO_S\to(fg)_*\cO_U $ is injective. It remains to show that $f_*\cO_{\cX}\to
  (fg)_*\cO_U$ is injective. Since $f_*$ is left exact,  it is enough to show that $\cO_\cX\to g_*\cO_U$ is
  injective. By the definition of $\cO_\cX$, this can be checked after pulling back by a smooth cover of $\cX$
  by a scheme (for example, $U$ itself), 
so we reduce to
  the case of algebraic spaces, which is immediate from~\cite[\href{http://stacks.math.columbia.edu/tag/082Z}{Tag 082Z}]{stacks-project}.
\end{rem}

\end{remark}

We have the following simple lemma concerning scheme-theoretic dominance.

\begin{lemma}
\label{lem:isomorphism criterion}
Suppose that $f: \cY \to \cZ$ is a scheme-theoretically dominant quasi-compact morphism of algebraic stacks, and that $\cZ' \hookrightarrow \cZ$ is a closed
immersion for which the base-changed morphism $\cY' \to \cY$
is an isomorphism.  Then $\cZ' \hookrightarrow \cZ$ is itself an
isomorphism. 
\end{lemma}

\begin{proof}
The assumption that $\cY' \iso \cY$ can be rephrased
as saying that the morphism $\cY \to \cZ$ can be factored through the
closed substack $\cZ'$ of $\cZ$.  Since $\cY \to \cZ$ is scheme-theoretically
dominant, we see that necessarily $\cZ' = \cZ$. (This last statement is easily
reduced to the case of schemes, where it is immediate.) \end{proof}

\subsection{Scheme-theoretic images (part two)}
\label{subsec:scheme images two}
Our goal in this section is to generalise the construction of scheme-theoretic images
to certain morphisms of stacks whose domain is algebraic,
but whose target is of a possibly more general nature.
The general set-up, which will be in force throughout this section, will be as follows: we suppose given
a morphism $\xi:\cX \to \cF$ of stacks over a locally
Noetherian base-scheme~$S$, 
whose domain $\cX$ is assumed to be algebraic,
and whose target $\cF$ is assumed to have diagonal~$\Delta_{\cF}$
which is representable by algebraic spaces and locally of finite
presentation (so, in particular we are assuming throughout this
section that~$\cF$ satisfies~[3]).

We will furthermore typically assume
that either $\cF$ satisfies~[1], or else that $\cX$ is locally of finite
presentation. We
will sometimes reduce the latter situation to the former by using the results of
Section~\ref{subsec:[1]-ification}. (See Remark~\ref{rem: some generalities on reducing to the case that F
    satisfies [1]} below.) We will also frequently need to assume that $\xi$ is
  proper. 
  We work in maximal generality when we can, only introducing these hypotheses at
the points that they are needed (but see Remark~\ref{rem:laundry list}~(3)
below).

We begin with a lemma whose intent is to capture the properties
that will characterise when a morphism $\Spec A \to \cF$,
with $A$ a finite type  Artinian local $\mathcal O_S$-algebra,
factors through the scheme-theoretic image of $\xi$.

\begin{lemma}
\label{lem:image conditions}If $\Spec A$ is a finite type Artinian local $S$-scheme,
and  $\varphi:\Spec A \to \cF$ is a morphism over $S$, then
the following conditions are equivalent:
\begin{enumerate} 
\item There exists a complete Noetherian local $\mathcal O_S$-algebra~$B$,
and a factorisation of $\varphi$ into $S$-morphisms
$\Spec A \to \Spec B \to \cF$,
such that the base-changed morphism $\cX_B \to \Spec B$
is scheme-theoretically dominant.
\item There exists a complete Noetherian local $\mathcal O_S$-algebra~$B$,
and a factorisation of $\varphi$ into $S$-morphisms
$\Spec A \to \Spec B \to \cF$,
such that the base-changed morphism $\cX_B \to \Spec B$
is scheme-theoretically dominant,
and such that the morphism $\Spec A \to \Spec B$
is a closed immersion.
\end{enumerate}

If $\xi$ is furthermore proper, then these conditions are equivalent
to the following further two conditions.
\begin{enumerate}
\setcounter{enumi}{2}
\item There exists an Artinian local $\mathcal O_S$-algebra~$B$,
and a factorisation of $\varphi$ into $S$-morphisms
$\Spec A \to \Spec B \to \cF$,
such that the base-changed morphism $\cX_B \to \Spec B$
is scheme-theoretically dominant.
\item There exists an Artinian local $\mathcal O_S$-algebra~$B$,
and a factorisation of $\varphi$ into $S$-morphisms
$\Spec A \to \Spec B \to \cF$,
such that the base-changed morphism $\cX_B \to \Spec B$
is scheme-theoretically dominant,
and such that the morphism $\Spec A \to \Spec B$
is a closed immersion.
\end{enumerate}
\end{lemma}

\begin{remark}\label{rem:in Artinian case B is finite type}
We remark that if $A$ is a Artinian local $\mathcal O_S$-algebra,
which is furthermore of finite type (or, equivalently,
whose residue field is of finite type over $\cO_S$),
and if $B$ is a local $\mathcal O_S$-algebra for which
there exists a morphism of $\mathcal O_S$-algebras
$B \to A$, then the residue field of $B$ is necessarily
also of finite type over $\cO_S$.
Thus, in the context of the preceding lemma, the residue
field of the ring $B$ appearing in any of the conditions
of the lemma will necessarily be of finite type over $\cO_S$.
\end{remark}

The proof of Lemma~\ref{lem:image conditions} will  make use of the theorem on formal functions for algebraic stacks
in the form of the following theorem of Olsson. As in~\cite{MR2312554}, we work
with sheaves on the lisse-\'etale site.

\begin{thm}
 \label{thm: thm on formal functions for stacks} Let $A$ be a Noetherian adic
 ring, and let $I$ be an ideal of definition of $A$. Let $\cX$ be a proper algebraic
 stack over $\Spec A$. Then the functor sending a
 sheaf to its reductions is an equivalence of categories between the category
 of coherent sheaves $\cG$ on $\cX$ and the category of compatible systems of
 coherent sheaves $\cG_n$ on the reductions $\cX_n:=\cX\times_{\Spec A}\Spec(A/I^{n+1})$.

Furthermore, if $\cG$ is a coherent sheaf on $\cX$ with reductions $\cG_n$, then
the natural map \[H^0(\cX,\cG)\to\varprojlim_n H^0(\cX_n,\cG_n)\] is an
isomorphism of topological $A$-modules, where the left hand side has the
$I$-adic topology, and the right hand side has the inverse limit topology.
\end{thm}\begin{proof}A proper morphism of stacks is separated, so has
quasi-compact and quasi-separated diagonal, so that the more restrictive
definition of an algebraic stack in~\cite{MR2312554} is automatically
satisfied. The claimed result is then a special case of~\cite[Thm.\ 11.1]{MR2312554}.
 \end{proof}

\begin{proof}[Proof of Lemma~\ref{lem:image conditions}]
Evidently~(2) implies~(1).  Conversely, suppose~(1) holds.
Bearing in mind Remark~\ref{rem:in Artinian case B is finite type}, we
see that the residue field of~$A$ is finite over that of~$B$, so Lemma~\ref{lem: Cohen rings lemma} shows
that we may factor the morphism $\Spec A \to \Spec B$
as $\Spec A \to \Spec B' \to \Spec B$,
where $B'$ is again a complete Noetherian local ring,
such that the first morphism is a closed immersion,
and such that the second morphism is faithfully
flat.  Since scheme-theoretic dominance is preserved
under flat base-change, we see that $\cX_{B'} \to \Spec B'$
is scheme-theoretically dominant, and so replacing 
$B$ by $B'$, we find that~(2) is satisfied.

It is evident that~(4) implies~(3) implies~(1), with
no assumption on $\xi$.  Thus, to complete the proof
of the lemma, it suffices to show that~(2) implies~(4),
when $\xi$ is proper.  In fact, we will show that
if $\Spec A \to \Spec B$ is a morphism as in~(2),
then we may find an Artinian quotient $B'$ of $B$
such that this morphism factors through the closed
immersion 
$\Spec B' \to \Spec B$,
and such that the base-changed morphism $\cX_{B'} \to \Spec B'$
is scheme-theoretically dominant.
Certainly $\Spec A \to \Spec B$ factors through
the closed immersion $\Spec B' \to \Spec B$ for {\em some}
Artinian quotient of $B$, since $A$ itself is Artinian.
Thus it suffices to show that $B$ admits a cofinal collection
of Artinian quotients $B'$ for which $\cX_{B'} \to \Spec B'$
is scheme-theoretically dominant. This is the content of
Lemma~\ref{lem: formal functions implies formal image equals image
  formal} below.
\end{proof}

\begin{lem}
  \label{lem: formal functions implies formal image equals image formal}
  Assume that~$\xi$ is proper, and suppose that
  $\Spec B\to\cF$ is a morphism over~$S$, where~$B$ is a complete
  Noetherian local $\cO_S$-algebra. Then the base-changed morphism
  $\cX_B\to\Spec B$ is scheme-theoretically dominant if and only
  if~$B$ admits a cofinal collection of Artinian quotients~$B'$ for
  which $\cX_{B'}\to\Spec B'$ is scheme-theoretically dominant.
\end{lem}
\begin{proof}Suppose firstly that $B$ admits a cofinal collection of Artinian quotients~$B'$ for
  which $\cX_{B'}\to\Spec B'$ is scheme-theoretically dominant. Then
  the scheme-theoretic image of $\cX_B\to\Spec B$ is a closed
  subscheme of~$\Spec B$ containing each of the~$\Spec B'$, and since
  the quotients~$B'$ are cofinal, it must in fact be~$\Spec B$.

We now consider the converse. Since $\xi$ is proper, so is the base-changed morphism
$\cX_{B}\to \Spec B$. 
If $\cX_B\to\Spec B$ is scheme-theoretically dominant,
then the natural map $B\to H^0(\Spec
B,\xi_*\cO_{\cX_B})=H^0(\cX_B,\cO_{\cX_B})$ is injective by Remark~\ref{rem: we don't want to give the quasi-coherent sheaf definition of
    scheme-theoretic image for stacks}. 
Set $\cX_n=\cX_B \times_{\Spec B} \Spec (B/\mathfrak m_B^{n+1})$.
Noting that the structure sheaf $\cO_{\cX_B}$ is coherent,
Theorem~\ref{thm: thm on formal functions for stacks} applies to show that
we have an isomorphism of
topological $B$-modules \[H^0(\cX,\cO_X)\isoto \varprojlim_n
H^0(\cX_n,\cO_{\cX_n}).\] Let $I_n$ be the kernel of the composite
morphism \[B\to H^0(\cX,\cO_X)\to 
H^0(\cX_n,\cO_{\cX_n}),\] so that we have injections $B/I_n\to
H^0(\cX_n,\cO_{\cX_n}) $. The natural map 
$\cX_n \to \cX\times_{\Spec B}\Spec(B/I_n)$
is then an isomorphism by construction,
and so 
$\cX\times_{\Spec B}\Spec(B/I_n)\to\Spec B/I_n$ is scheme-theoretically
dominant (again by Remark~\ref{rem: we don't want to give the quasi-coherent sheaf
  definition of scheme-theoretic image for stacks}; since $B/I_n\to
H^0(\cX_n,\cO_{\cX_n}) $ is injective, the natural map $\cX_n \to
\Spec(B/I_n)$ is scheme-theoretically dominant). Thus the $B/I_n$ are the
sought-after cofinal collection of Artinian
quotients of $B$.\end{proof}

We now define, in our present context, a subgroupoid $\widehat{\cZ}$
of $\widehat{\cF}$. 

\begin{df}
If $\Spec A$ is a finite type Artinian local $S$-scheme,
then we let $\widehat{\cZ}(A)$ denote the full subgroupoid of $\cF(A)$
consisting of morphisms $\Spec A \to \cF$ that satisfy the equivalent
conditions (1) and (2) of Lemma~\ref{lem:image conditions}.
(If $\xi$ is proper, then we note that the objects of
$\widehat{\cZ}(A)$ can
equally well be characterised in terms of conditions (3) and (4)
of that lemma.)
\end{df}

We now define the scheme-theoretic image~$\cZ$ of $\xi$, using the
terminology and results of Section~\ref{subsec:[1]-ification} (note
that since~$S$ is assumed to be locally Noetherian, it is in particular
quasi-separated, so the results of Section~\ref{subsec:[1]-ification}
apply to categories fibred in groupoids over~$S$). As in Section~\ref{subsec:[1]-ification}, we let $\AffS$ denote the category of affine $S$-schemes,
and let $\AfffpS$ denote the full subcategory of finitely presented
affine $S$-schemes. 

To begin with, we will consider the restriction~$\cF_{| \AfffpS}$
to a category fibred in groupoids over $\AfffpS$,
and define a full subcategory of~$\cF_{| \AfffpS}$ which
(by an abuse of notation which we will justify below)
we denote $\cZ_{| \AfffpS}$.

\begin{df}
\label{def:stack-theoretic image}Let $\cZ_{| \AfffpS}$ be the full subcategory of $\cF|_{\AfffpS}$ defined as
follows: if $A$ is a finite type $\cO_S$-algebra,
then we let $\cZ_{| \AfffpS}(A)$ denote the full subgroupoid of $\cF(A)$ consisting
of points $\eta$ whose formal completion $\widehat{\eta}_t$,
at each finite type point $t$ of $\Spec A$,
factors through $\widehat{\cZ}$. 
\end{df}

\begin{lemma}
  \label{lem: Z is a Zariski stack first time around}$\cZ_{| \AfffpS}$ is a
  Zariski substack of~$\cF_{|\AfffpS}$. 
\end{lemma}
\begin{proof} 
To check that $\cZ_{| \AfffpS}$ is a full subcategory fibred
in groupoids of $\cF_{| \AfffpS},$
we need to check that if $T' \to T$ is a morphism in~$\AfffpS$,
and $T \to \cF$ satisfies the condition to lie in~$\cZ(T)$,
then  the same is true of
the pulled-back morphism $T' \to \cF$; but this is immediate from the
definition of~$\cZ$.
To see that it is actually a substack of~$\cF|_{\AfffpS}$, it is then
enough (by~\cite[\href{http://stacks.math.columbia.edu/tag/04TU}{Tag
  04TU}]{stacks-project}) to note that since the condition
of a morphism $T \to \cF$ factoring through $\cZ$ is checked
pointwise, it is in particular a Zariski-local condition.
\end{proof}

We remind the reader that
Theorem~\ref{thm: Thomason and Trobaugh, affine is pro-fp-affine}
gives an equivalence of
categories $\pro\AfffpS \iso \AffS$. We now give our definition of the
scheme-theoretic image of $\cX\to\cF$; see Remark~\ref{rem:laundry
  list}~(2) below for the relationship of this definition to the existing
one in the case that~$\cF$ is an algebraic stack.

\begin{df}
\label{df:scheme-theoretic images, finally}
We define $\cZ$ to be the pro-category $\pro\cZ_{| \AfffpS},$
which
(by Lemmas~\ref{lem: full subcategory is compatible with pro}
and~\ref{lem: Z is a Zariski stack first time around})
is a full subcategory in groupoids of $\pro\cF_{| \AfffpS}$.
Via the equivalence of
Theorem~\ref{thm: Thomason and Trobaugh, affine is pro-fp-affine},
we regard $\cZ$ as a category fibred in groupoids over $\AffS$.
We also note that Lemma~\ref{lem:adjunction equivalence} yields
an equivalence between the restriction of $\cZ$ to $\AfffpS$ 
and our given category $\cZ_{| \AfffpS}$, thus justifying the notation
for the latter.

Since $\Delta_{\cF}$ is assumed to be representable by algebraic
spaces and locally of finite presentation, it follows from
Lemma~\ref{lem:restriction to finitely presented algebras}~(4) that
$\pro(\cF|_{\AfffpS})$ is equivalent to a full subcategory of $\cF$,
so that we may regard $\cZ$ as a full subcategory of $\cF$
(the latter being thought of
as a category fibred in groupoids over the category of affine $S$-schemes).
We refer to $\stimage$ 
as the \emph{scheme-theoretic image} of the morphism $\xi: \cX \to \cF$.\end{df}

\begin{lemma}
  \label{lem: Z is a Zariski stack which satisfies [1]}$\cZ$ is a
  Zariski substack of~$\cF$ over~$\AffS$, and satisfies~{\em [1]}.
\end{lemma}
\begin{proof} This follows immediately
from  Lemmas~\ref{lem:restriction to finitely presented
    algebras}~(2),
 \ref{lem:limit preserving and pro-categories},
and~\ref{lem: Z is a Zariski stack first time around}.
\end{proof}

\begin{remark}
\label{rem:laundry list}

(1) Since~$\cZ$ is a Zariski stack over~$\AffS$, it naturally extends
to a Zariski stack over~$S$ (and is again a Zariski substack
of~$\cF$). We will see in 
Lemma~\ref{lem:we get a sheaf}
below that under some relatively mild additional assumptions,
$\cZ$ is in fact an \'etale substack of~$\cF$.

(2) In the case when $\cF$ is an algebraic stack which satisfies~[1], 
then its substack $\cZ$ coincides with the scheme-theoretic image of $\xi$
(in the sense of
Definition~\ref{defn: scheme theoretic image for algebraic stacks});
see Lemma~\ref{lem:algebraic case} below.

(3) If neither~$\cX$ nor $\cF$ satisfy~[1], then Definition~\ref{def:stack-theoretic image} is not a sensible one.  E.g. if $\cF$ is an algebraic
stack that does not satisfy~[1], and $\xi: \cF \to \cF$ is the identity
morphism, then evidently the scheme-theoretic image of $\xi$ (in the usual
sense, i.e. in the sense of
Definition~\ref{defn: scheme theoretic image for algebraic stacks}) is
just $\cF$, and so doesn't satisfy~[1] (by assumption) --- whereas we noted
above that the substack $\cZ$ of $\cF$ given by
Definition~\ref{def:stack-theoretic image} does satisfy~[1], by
stipulation.\end{remark}

Here and on several occasions below, we will have cause to consider 
how scheme-theoretic images interact with monomorphisms of stacks
$\cF' \hookrightarrow \cF$.  In light of this, it will be helpful to
recall that Lemma~\ref{lem:3.2 running assumptions descend to monos}
implies that if $\cF$ has a diagonal which is representable by locally
algebraic spaces and locally of finite presentation, 
then the same is true of any stack $\cF'$ admitting a monomorphism
$\cF' \hookrightarrow \cF$ (i.e.\ any substack~$\cF'$ of $\cF$).

\begin{lem}
\label{lem:images compatible with monos}Suppose that $\cX\to\cF$
factors as
$\cX\to\cF'\hookrightarrow\cF$, where $\cF'$ is a stack which satisfies~{\em[3]}, and
$\cF'\hookrightarrow\cF$ is a monomorphism. Then the scheme-theoretic
image of $\cX\to\cF'$ is contained in the scheme-theoretic image of $\cX\to\cF$.
\end{lem}
\begin{proof}A consideration of the definitions shows
that we need only check that if $\Spec
  B\to\cF'$ is such that $\cX\times_{\cF'}\Spec B\to\Spec B$ is
  scheme-theoretically dominant, then $\cX\times_{\cF}\Spec B\to\Spec
  B$ is scheme-theoretically dominant; but this is immediate, because $\cX\times_{\cF'}\Spec B=\cX\times_{\cF}\Spec B$.
\end{proof}

\begin{remark}
Throughout the
rest of this section we prove several variants and refinements of
Lemma~\ref{lem:images compatible with monos}. More specifically,
Lemma~\ref{lem:equality of images} gives a criterion for the
scheme-theoretic images of $\cX\to\cF'$ and $\cX\to\cF$ to actually be
equal, and Lemma~\ref{lem:xibar} gives a criterion for $\cX\to\cF$ to
factor through~$\cZ$. Lemmas~\ref{lem:minimality}
and~\ref{lem:minimality for finite type monomorphisms} relate the
existence of a factorisation $\cX\to\cF'\hookrightarrow\cF$ to the
property of~$\cF'$ containing~$\cZ$, and Proposition~\ref{prop: if there is a closed algebraic substack then Z is
    it} shows that~$\cZ$ can be characterised by a universal property
  if it is assumed to be algebraic.
\end{remark}

\begin{rem}\label{rem: some generalities on reducing to the case that F
    satisfies [1]}As we explained at the beginning of this section, our main results assume either
that $\cF$ satisfies~[1], or that $\cX$ satisfies~[1] (equivalently, $\cX$ is locally of finite
presentation over~$S$). In the proof of the main theorem, we will
argue by reducing the second case to the first case in the following
way:

We first note that it follows from
Lemmas~\ref{lem:adjunction equivalence},
\ref{lem:limit preserving and pro-categories},
and~\ref{lem:restriction to finitely presented algebras}~(2) and~(4)
that $\pro(\cF{|_\AfffpS})$ is a substack of $\cF$ which satisfies~[1] and~[3].
If we assume in addition that $\cF$ admits versal rings at all
finite type points,
then so does~$\pro(\cF{|_\AfffpS})$,
by~Lemma~\ref{lem:restriction to finitely presented algebras}~(5).

If $\cX$ satisfies~[1], then it follows 
from Lemma~\ref{lem:limit preserving and pro-categories}
that there is an equivalence $\pro\cX_{| \AfffpS} \iso \cX,$
and thus from 
Lemma~\ref{lem:restriction to finitely presented algebras}~(1)
that $\xi: \cX \to \cF$ can be factored through a morphism
$\cX \to
\pro\cF_{\AfffpS}.$
By construction the monomorphism $\cZ \hookrightarrow \cF$
also factors through $\pro\cF_{\AfffpS},$
and an examination of the definitions shows that
$\cZ$ is also the scheme-theoretic image (in the sense of
Definition~\ref{df:scheme-theoretic images, finally})
of the induced morphism $\cX\to\pro(\cF_{|_\AfffpS})$.
Taken together, the previous remarks will allow us to simply
  replace~$\cF$ by~$\pro(\cF_{|_\AfffpS})$, and thus assume that we
  are in the first case.  \end{rem}

\begin{lemma}
  \label{lem: scheme theoretic dominance for Artin proper}Assume that
  $\xi:\cX\to\cF$ is proper, and let $x:\Spec k\to\cF$ be a finite type
  point. Then $x$ is a point of $\cZ$ if and only if the fibre
  $\cX_x$ is non-empty.
\end{lemma}\begin{proof}If $\cX_x$ is non-empty, then $\cX_x\to\Spec k$ is scheme-theoretically surjective, and $x$
  is a point of $\cZ$ by definition. Conversely, if $x$ is a point of $\cZ$,
  then by Lemma~\ref{lem:image conditions} we can factor $x:\Spec k\to\cF$
  through a closed immersion $\Spec k\to\Spec B$ with $B$ an Artin local ring
  and $\cX_B\to\Spec B$ scheme-theoretically dominant. This implies that the
  fibre $\cX_x$ is non-empty, as required.
  \end{proof}

    \begin{defn}\label{defn:definition of scheme-theoretic images for versal rings}
      If $\Spf A_x\to\widehat{\cF}_x$ is a versal ring
      at $x$, then (by Remark~\ref{rem:scheme-theoretically
        dominant}~(4)) if we let $A_i$ run over the discrete Artinian
      quotients of $A_x$, the scheme-theoretic images $\Spec R_i$ of
      the morphisms
      $\cX_{A_i}\to\Spec A_i$ fit together to give a formal subscheme
      $\Spf R_x:=\Spf\varprojlim_iR_i$ of $\widehat{\cF}_x$, which we call the
      scheme-theoretic image of the base-changed morphism
      $\cX \times_{\cF} \Spf A_x \to \Spf A_x$.

      It follows from Remark~\ref{rem:pseudo-compact} that the natural morphism
$A_x\to R_x$ is surjective, and is a quotient map of topological rings.
Thus $\Spf R_x$ is even a {\em closed} formal subscheme of $\Spf A_x$.
    \end{defn}

   The following lemma shows that, when $\xi$ is proper, versal rings for $\cZ$ can
    be constructed from versal rings for $\cF$ by taking scheme-theoretic
    images.

  \begin{lem}
  \label{lem:versal ring for Z}Assume that $\xi:\cX\to\cF$ is
  proper. Let $x:\Spec k\to\cZ$ be a finite type
  point, which we also consider as a finite type point of $\cF$. Suppose that $\Spf
  A_x\to\widehat{\cF}_x$ is a versal ring at $x$, and let $\Spf R_x$
  be the scheme-theoretic image of $\cX_{\Spf A_x}\to\Spf A_x$. Then the
  morphism $\Spf R_x\to\widehat{\cF}_x$ factors through a versal morphism $\Spf R_x\to\widehat{\cZ}_x$.
\end{lem}
\begin{proof}
We claim that a morphism $\Spec A\to\Spf A_x$, with $A$ an object of
$\cC_{\Lambda}$, 
factors through $\Spf R_x$ if and only if the composite $\Spec
A\to\Spf A_x\to\widehat{\cF}_x$ factors through~$\widehat{\cZ}_x$. In
the notation of Definition~\ref{defn:definition of scheme-theoretic images for versal rings}, if
$\Spec A\to\widehat{\cF}_x$ factors through~$\Spf R_x$, then it in fact factors
through $\Spec R_i$ for some $i$, and hence 
through $\widehat{\cZ}_x$ by definition (as the morphism $\cX_{R_i}
  \to\Spec R_i$ is scheme-theoretically dominant).

Conversely, if the composite $\Spec
A\to\Spf A_x\to\widehat{\cF}_x$ factors through~$\widehat{\cZ}_x$, then by
Lemma~\ref{lem:image conditions}, we have a factorisation $\Spec A\to\Spec B\to\cF$, where $B$ is
an object of~$\cC_{\Lambda}$, the morphism $\Spec
A\to\Spec B$ is a closed immersion, and $\cX_B\to \Spec B$ is
scheme-theoretically dominant.
 By the versality
of $\Spf A_x\to\widehat{\cF}_x$,
we may lift the morphism $\Spec B\to\widehat{\cF}_x$ to a morphism $\Spec B\to\Spf A_x$, which furthermore we may factor as $\Spec B \to \Spec A_i \to \Spf A_x$, for some value of $i$.   Since $\cX_B \to \Spec B$ is scheme-theoretically dominant, the morphism $\Spec B \to \Spec A_i$ then factors through
$\Spec R_i$,
and thus through $\Spf R_x$, as claimed.

In particular, we see that the composite $\Spf R_x\to\Spf A_x\to\widehat{\cF}_x$
factors through a morphism $\Spf R_x\to\widehat{\cZ}_x$. It remains to check
that this morphism is versal. This is formal. Suppose we are given a commutative
diagram \[\xymatrix{\Spec A_0\ar[r]\ar[d]&\Spf R_x\ar[r]\ar[d]&\Spf A_x\ar[d]\\
  \Spec A\ar[r]&\widehat{\cZ}_x\ar[r]&\widehat{\cF}_x }\]where the left hand
vertical arrow is a closed immersion, and $A_0, A$ are objects
of~$\cC_\Lambda$. By the versality of $\Spf A_x\to\widehat{\cF}_x$, we may lift the
composite $\Spec A\to\widehat{\cF}_x$ to a morphism $\Spec A\to\Spf A_x$. Since
the composite $\Spec A\to\Spf A_x\to\widehat{\cF}_x$ factors through
$\widehat{\cZ}_x$, the morphism $\Spec A\to \Spf A_x$ then factors through $\Spf
R_x$, as required.
 \end{proof}

If $\cZ$ were to behave like the scheme-theoretic image of a morphism
of algebraic stacks, i.e.\ like a closed substack,
then we would expect to be able to test the property
of a morphism to $\cF$ factoring through $\cZ$ by precomposing with
a scheme-theoretically dominant morphism.  
We cannot prove this general statement at this point of the development,
but we begin our discussion of this general problem
by establishing some special cases.

\begin{lemma}
\label{lem:scheme-theoretically dominant --- Artinian case}
Let $ Y \to  Z \to \cF$ be a composite of morphisms over $S$, with
$ Y$ and $ Z$ being local Artinian schemes of finite type over $S$,
for which the morphism $ Y \to  Z$ is scheme-theoretically
surjective, and such that the composite $ Y \to \cF$ factors
through $\cZ$.  Suppose that:

(a) either $\xi$ is proper and $\cF$
is~$\mathbf{Art}^{\mathbf{triv}}$-homogeneous, 
or $\cF$ satisfies~{\em [2](b)}; and

(b) either $\cF$
is $\mathbf{Art}^{\mathbf{fin}}$-homogeneous,
i.e.\ $\cF$ satisfies~{\em (RS)},
or the residue field extension
in $Y \to Z$ is separable. 
\newline
Then the morphism $ Z \to \cF$ also factors through~$\cZ$.
\end{lemma}
\begin{proof}
Set $Y = \Spec A$ and $Z = \Spec A'$;
the scheme-theoretically surjective
morphism $Y \to Z$ then corresponds to an injective  morphism
$A' \to A$ of $\cO_S$-algebras.
By assumption we may find a surjection of $\cO_S$-algebras
$B \to A$, whose source is a complete Noetherian local
$\cO_S$-algebra, such that the given morphism $\Spec A \to \cF$
factors through a morphism $\Spec B \to \cF$
for which the induced morphism $\cX_{B} \to \Spec B$
is scheme-theoretically surjective.

Let $B':= A'\times_{A} B.$ 
We have a commutative diagram\bigskip
$$
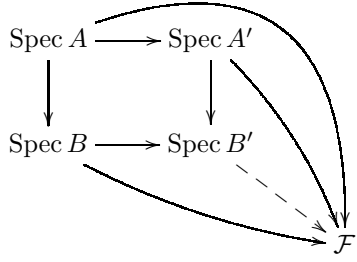
\xymatrix{\Spec A \ar[r]\ar@/^4.5pc/[rrdd] \ar[d] & \Spec A' \ar[d]\ar@/^/[rdd] &  \\
\Spec B \ar[r]\ar@/_/[rrd] & \Spec B'\ar@{-->}[dr] & \\
& &  \cF}
$$
Suppose firstly that $\xi$ is proper and $\cF$ is
$\mathbf{Art}^{\mathbf{triv}}$-homogeneous. By Lemma~\ref{lem:image conditions},
we may assume that $B$ is Artinian.
We claim that we may fill in this diagram with a morphism
$\Spec B' \to \cF$; this follows since we are assuming that either $\cF$ is
$\mathbf{Art}^{\mathbf{fin}}$-homogeneous, or that
the residue field extension is separable, in which case this follows from
Lemma~\ref{lem:homogeneity conditions}.
 Otherwise, if $\xi$ is not proper, then 
by assumption $\cF$ satisfies~[2](b).
Lemma~\ref{lem:complete homogeneity} and our hypotheses then
show that
we may once again fill in this diagram with a morphism
$\Spec B' \to \cF.$

It remains to show that $\cX_{B'} \to \Spec B'$
is scheme-theoretically dominant.  But this follows
directly from the fact that each of $\cX_B\to \Spec B$
and $\Spec B \to \Spec B'$
are scheme-theoretically dominant (the latter because
$B'\to B$ is injective),
and a consideration of the commutative diagram
$$\xymatrix{\cX_B \ar[r] \ar[d] & \cX_{B'} \ar[d] \\
\Spec B \ar[r] & \Spec B'}
$$
(This is a particular case of Remark~\ref{rem:properties of scheme-theoretic images of algebraic stacks}~(4).)
\end{proof}

\begin{lemma}\label{lem:scheme-theoretically dominant --- smooth case}Suppose that either
$\xi$ is proper and $\cF$ is~$\mathbf{Art}^{\mathbf{triv}}$-homogeneous, or
that $\cF$ satisfies~{\em [2](b)}. If $\cU \to \cY$ is a smooth surjective
morphism of algebraic stacks
over~$S$, and if $\cY \to \cF$ is a morphism of stacks over $S$
such that the composite $\cU \to \cY \to \cF$ factors through
$\cZ$, then the morphism $\cY \to \cF$ itself factors through~$\cZ$.
\end{lemma}
\begin{proof}
Let $T \to \cY$ be a morphism,
with $T$ a scheme over $S$.  We must show that the composite
$T \to \cY \to \cF$ factors through $\cZ$.
By assumption the morphism $\cU\times_{\cY} T \to T \to \cY \to \cF$
(which admits the alternate factorisation
$\cU\times_{\cY} T \to \cU \to \cY \to \cF$) factors through $\cZ$.
Since $\cU\times_{\cY} T$ is an algebraic stack,
it admits a smooth surjection from a scheme $U$.

The composite $U \to T$ is a smooth morphism,
and so 
we may
apply~\cite[\href{http://stacks.math.columbia.edu/tag/055V}{Tag 055V}]{stacks-project} to find a morphism $V \to U$ for which
the composite $V\to U \to T$ is \'etale and surjective.  
Replacing $U$ by $V$,
we are reduced to showing that if $U \to T$ is a surjective \'etale 
morphism of $S$-schemes,
and the composite $U \to T \to \cF$ factors through $\cZ$,
then the morphism $T \to \cF$ factors through $\cZ$.
(This is essentially the \'etale sheaf property of $\cZ$;
see Lemma~\ref{lem:we get a sheaf} and its proof below.)

Since $\cZ$ is a Zariski substack of $\cF$, we may check this statement Zariski locally on $T$, and hence assume that $T$ is affine.
Since $T$ is then quasi-compact, we may find a quasi-compact
open subset of $U$ which surjects onto $T$, and then replacing
$U$ by the disjoint union of the members of a finite affine
cover of this quasi-compact open subset, we may further assume
that $U$ is affine.
Since $\cZ$ satisfies~[1], we may further reduce to the case
that $T$ and $U$ are locally of finite presentation over $S$,
and hence of finite type over $S$ (since $S$ is locally
Noetherian). (More precisely: by Lemma~\ref{lem:descending etale to
  finite presented affine}, we may write $T=\varprojlim_i T_i$,
$U=\varprojlim_i U_i$, with $U_i\to T_i$ an \'etale surjection of
affine schemes locally of finite presentation over~$S$, such that
$U\to T$ is the pull-back of each $U_i\to T_i$. Since $\cZ$
satisfies~[1], the morphism $U\to\cZ$ factors through $U_i$ for
some~$i$, and thus so does the composite $U\to \cZ\to\cF$. Arguing as
in the proof of Lemma~\ref{lem:restriction to finitely presented
  algebras}~(2), the \'etale descent data for the morphism $U\to\cF$
arises as the base change of \'etale descent data for the morphism
$U_{i'}\to\cF$ for some $i'\ge i$, and since~$\cF$
is an \'etale stack, it follows that $T\to\cF$  factors through~$T_{i'}$,
as required.)

In this case we are reduced to checking that 
if $u: \Spec l \to  U$ is a finite type point
lying over the finite type point $t: \Spec k \to T$, then
the morphism $\Spf \widehat{\cO}_{T,t} \to \widehat{\cF}_{t}$
factors through $\widehat{\cZ}_{t}$
provided that the morphism $\Spf \widehat{\cO}_{U,u} \to \widehat{\cF}_{u}$
factors through~$\widehat{\cZ}_u$. 
That this is true follows from
Lemma~\ref{lem:scheme-theoretically dominant --- Artinian case}
and the fact that the local \'etale (and hence faithfully flat)
morphism $\widehat{\cO}_{T,t} \to \widehat{\cO}_{U,u}$
can be written as the inverse limit of faithfully flat 
(and hence injective) morphisms of Artinian local rings with separable residue field 
extensions (to which the cited lemma applies).
\end{proof}

The following lemma gives a useful criterion for when a morphism
$\cT \to \cF$ factors through $\stimage$.

\begin{lemma}
\label{lem:scheme-theoretic images}
Assume  that either
$\xi$ is proper and $\cF$ is~$\mathbf{Art}^{\mathbf{triv}}$-homogeneous, or
that $\cF$ satisfies~{\em [2](b)}. Let $\cT$ be an algebraic stack over $S$, 
and let $\eta:\cT \to \cF$ be a morphism for which
the base-changed morphism $\cX_\cT \to \cT$ is scheme-theoretically
dominant. Assume further that  either $\cF$ satisfies~{\em [1]}, or that $\cT$
is locally of
finite type over~$S$. Then~$\eta$ factors through $\stimage$. \end{lemma}
\begin{proof}
Since $\cT$ is an algebraic stack, we may find a smooth surjection $U \to \cT$
with $U$ a scheme.  Since the formation of scheme-theoretic images
commutes with flat base-change, we see
that $\cX_U \to U$ must be scheme-theoretically dominant. On the other hand, Lemma~\ref{lem:scheme-theoretically dominant --- smooth case} 
shows that $\cT \to \cF$ factors through $\cZ$ if and only if $U \to \cT \to \cF$
does.  Thus we reduce to the case that $\cT$ is a scheme~$T$.

Since $\stimage$ is a Zariski substack of $\cF$ by Lemma~\ref{lem: Z is a Zariski stack which satisfies [1]},
we can check the assertion of
the lemma Zariski locally on $T$, and hence we may assume that $T = \Spec A$ is
affine. If $\cT$ is locally of finite type over~$S$, then we may assume that $T$
is of finite type over~$S$. In the case that we
are assuming instead that $\cF$ satisfies~[1], write $A = \varinjlim_i A_i$ as the inductive limit of its
finitely generated $\cO_S$-subalgebras. Then for some $i$ there exists $\eta_i: \Spec A_i \to \cF$ such that
$\eta$ is obtained as the composite of $\eta_i$ with $\Spec A \to \Spec A_i$;
Remark~\ref{rem:properties of scheme-theoretic images of algebraic stacks} (4)
then shows that $\cX_{A_i}
\to \Spec A_i$ is scheme-theoretically dominant.  
Thus we may replace $A$ by $A_i$,
and hence again assume that $T$ is finite type over $S$.

Let $t$ be a finite type point of $T$. Since $\cX_T \to T$ is
scheme-theoretically dominant, the same is true of the (flat) base-change
$\cX_{\widehat{\cO}_{T,t}} \to \Spec \widehat{\cO}_{T,t}.$
By definition, then, the composite morphism
$\Spf \widehat{\cO}_{T,t} \to T \to \cF$ factors
through~$\widehat{\cZ}$. 
Since this holds for all finite type points~$t$,  the morphism $T \to
\cF$ factors through $\cZ$ by definition.
\end{proof}

Our next lemma spells out the basic properties satisfied
by $\cZ$.

\begin{lemma}\label{lem:we get a sheaf} Suppose either
that $\xi$ is proper and that $\cF$ is~$\mathbf{Art}^{\mathbf{triv}}$-homogeneous, or that
$\cF$ satisfies~{\em [2](b)}. Assume also either that $\cF$ satisfies~{\em [1]}, or that
$\cX$ is locally of finite presentation over $S$. Then the scheme-theoretic image $\stimage$ forms a substack of $\cF$,
and satisfies Axioms~{\em [1]} and~{\em [3]}. 
If $\cF$ satisfies {\em [2](a)}, then so does $\cZ$.
If $\xi$ is proper and $\cF$ satisfies {\em [1]} and~{\em [2](b)},
then $\cZ$ satisfies~{\em [2](b)}.
\end{lemma}
\begin{proof}
As already noted,
by its very definition,
we see that $\cZ$  satisfies~[1].
It then follows by standard limit arguments 
that in order to show that $\cZ$ is a stack on the big \'etale site of $S$,
it suffices to do so after restricting $\cZ$ to the category 
of finite type $\cO_S$-algebras; more precisely, this follows easily
from Lemma~\ref{lem:restriction to finitely presented algebras}~(3) and~\cite[\href{http://stacks.math.columbia.edu/tag/021E}{Tag
  021E}]{stacks-project}.

Thus,
since $\cF$ is a stack for the \'etale topology,
and since $\cZ$ is defined to be a full subcategory of $\cF$,
in order to show that $\cZ$ is a stack on the big \'etale site of $S$, it suffices to verify that
the property of a morphism $\eta: T \to \cF$ factoring through $\cZ$
(for an $S$-scheme $T$) is \'etale local on $T$.
This follows from Lemma~\ref{lem:scheme-theoretically dominant --- smooth case}.

Since $\cZ$ (thought of as a category fibred in groupoids) is
a full subcategory of~$\cF$, we see that the diagonal $\cZ \to \cZ \times_S \cZ$
is the base-change under the morphism $\cZ \times_S \cZ \to \cF \times_S \cF$
of the diagonal $\cF \to \cF \times_S \cF$.  Thus, since $\cF$ satisfies~[3],
the same is true of $\cZ$.

Suppose now that $\cF$ satisfies~[2](a).
To verify~[2](a) for $\cZ$,
suppose that we are given a pushout diagram
$$\xymatrix{Y \ar[d]\ar[r] & Y' \ar[d] \\ Z \ar[r] & Z' }$$
of finite type Artinian $S$-schemes, 
with the horizontal arrows being closed immersions.
By assumption, the induced functor 
$$\cF(Z') \to \cF(Y') \times_{\cF(Y)} \cF(Z)$$
induces an equivalence of categories. Since $\cZ$ is a substack of $\cF$,
we see that the induced functor
$$\cZ(Z') \to \cZ(Y') \times_{\cZ(Y)} \cZ(Z)$$
is fully faithful, and, furthermore, that in 
order to verify that it induces an equivalence
of categories, it suffices to show that an element of $\cF(Z')$
whose image under pull-back in $\cF(Y') \times \cF(Z)$ lies in $\cZ(Y') \times \cZ(Z)$
itself necessarily lies in $\cZ(Z')$.

Suppose given such an element of $\cF(Z')$.
By hypothesis, we may find $\mathcal O_S$-algebras $\tY'$ and $\tZ$,
each the spectrum of a complete Noetherian local $\mathcal O_S$-algebra
with finite type residue field, and closed immersions
$Y' \hookrightarrow \tY'$ and $Z \hookrightarrow \tZ$,
and morphisms $\tY' \to \cF$ and $\tZ \to \cF$ inducing the given
morphisms $Y' \to \cF$ and $Z \to \cF$, such that the
pulled-back morphisms $\cX_{\tY'} \to \tY'$ and $\cX_{\tZ} \to \tZ$
are scheme-theoretically dominant. If $\xi$ is proper, then we can and do assume
that $\tY'$ and $\tZ$ are Artinian, and then $\tY' \coprod \tZ \to \cF$ factors through
$\tZ' := \tY' \coprod_Y \tZ$ by our assumption that $\cF$
satisfies~[2](a). Otherwise, $\cF$ satisfies~[2](b) by assumption, and
Lemma~\ref{lem:complete homogeneity} shows that the induced morphism 
$\tY' \coprod \tZ \to \cF$ factors through
$\tZ' := \tY' \coprod_Y \tZ.$   Since the tautological morphism $\tY' \coprod \tZ
\to \tZ'$ is scheme-theoretically dominant,
the pulled-back morphism $\cX_{\tZ'} 
\to \tZ'$
is scheme-theoretically dominant. 
Since the given morphism $Z' := Y' \coprod_Y Z\to \cF$ factors
through the natural morphism $Z' \to \tZ',$
by definition the morphism $Z' \to \cF$ does indeed factor through $\cZ$.

Finally, suppose that $\xi$ is proper,
and that $\cF$ satisfies [1] and~[2](b), and let $\Spec A_x\to
\cF$ be an effective versal morphism at a finite type point $x$ of $\cZ$
(regarded as a finite type point of $\cF$). Let $\Spec R_x$ be the
scheme-theoretic image  of $\cX_{\Spec A_x}\to \Spec A_x$; then the morphism
$\Spec R_x\to\cF$ factors through $\cZ$,
by Lemma~\ref{lem:scheme-theoretic images}.
Lemma~\ref{lem: formal functions implies formal image equals image formal},
applied with~$B=R_x$,
then shows that the scheme-theoretic image of
$\cX_{\Spf A_x} \to \Spf A_x$ is equal to $\Spf R_x$,
and it now follows from Lemma~\ref{lem:versal ring for Z}
that the induced morphism $\Spf R_x \to \cZ$ is versal at~$x$.
\end{proof}

The next lemma gives a refinement of Lemma~\ref{lem:images compatible with monos}.

\begin{lem}
\label{lem:equality of images}
Suppose that $\xi$ admits a factorisation $\cX \buildrel \xi' \over
\longrightarrow \cF' \hookrightarrow \cF$,
with $\cF'$ also satisfying~{\em [3]}, and 
with the second arrow being a monomorphism.
Suppose furthermore that the monomorphism $\cZ \hookrightarrow \cF$ factors
through $\cF'$, and that either 
$\xi$ is proper and $\cF$ is~$\mathbf{Art}^{\mathbf{triv}}$-homogeneous,
or that $\cF$ satisfies~{\em [1]} and~{\em [2](b)}. Then the scheme theoretic image of $\xi'$ is equal to $\cZ$.
\end{lem}
\begin{proof}
If we let $\cZ'$ denote the scheme-theoretic image of $\xi'$,
then Lemma~\ref{lem:images compatible with monos} shows that
$\cZ' \hookrightarrow \cZ$. We must prove the reverse inclusion;
that is, we must show that for any morphism $T \to \cZ$, 
with $T$ an affine scheme over $S$,
the composite $T \to \cZ \hookrightarrow \cF$  factors
through $\cZ'$.  By definition the morphism $T \to \cZ$
factors as $T \to T' \to \cZ$, with $T'$ of finite type
over $S$, and so it is no loss of generality to assume
from the beginning that $T$ is finite type over $S$.
By assumption the composite $T \to \cZ \hookrightarrow \cF$ factors through
$\cF'$, and now by definition of $\cZ'$, 
we see that it suffices to check that this composite factors
through $\cZ'$ under the additional hypothesis that $T = \Spec A$
for some finite type Artinian $\mathcal O_S$-algebra $A$.

Since the morphism $T \to \cF$ factors through $\cZ$,
by definition there exists a morphism $\Spec A \to \Spec B \to \cF$
as in Lemma~\ref{lem:image conditions}~(1) such that $\cX_B \to \Spec B$
is scheme-theoretically dominant. In the case that~$\cF$ does not
satisfy~[1], we are assuming that~$\xi$ is proper, so by
Lemma~\ref{lem:image conditions} we can assume that~$B$ is Artinian,
which by Remark~\ref{rem:in Artinian case B is finite type} implies in
particular that in this case, we can assume that $\Spec B$ is of finite
type over~$S$.  Lemma~\ref{lem:scheme-theoretic images}
then implies that $\Spec B \to \cF$ factors through $\cZ$,
and so in particular through $\cF'$.  It now follows from the definition
that $T \to \cF'$ factors through $\cZ'$, as required.
\end{proof}

\begin{lem}
\label{lem:xibar}
Assume either that $\cF$ satisfies~{\em [1]} and~{\em [2](b)}, or that
$\cF$ is~$\mathbf{Art}^{\mathbf{triv}}$-homogeneous,  $\xi$ is
proper, and either $\cF$ satisfies~{\em [1]}, or $\cX$ is locally of finite presentation over $S$. Then the morphism $\xi: \cX \to \cF$ factors
through a morphism $\xibar:\cX\to\stimage$, and the scheme-theoretic image of
$\xibar$ is just~$\stimage$. \end{lem}
\begin{proof}Since
the diagonal map $\cX \to \cX\times_{\cF} \cX$ gives a section to the projection
$\cX\times_{\cF} \cX \to \cX$, it is immediate from
Lemma~\ref{lem:scheme-theoretic images} that $\xi$ factors through~$\cZ$. That
the scheme-theoretic image of~$\xibar$ is~$\cZ$ is immediate from Lemma~\ref{lem:equality of images}.
\end{proof}

\begin{lemma}
\label{lem:proper and surjective}  Assume either that $\cF$ satisfies~{\em [1]}, or that
$\cX$ is locally of finite presentation over $S$. Suppose further that
$\xi$ is proper,
and that
$\cF$ is $\mathbf{Art}^{\mathbf{triv}}$-homogeneous.
Then the morphism $\xibar$ is proper and surjective.
\end{lemma}

\begin{remark}
Note that Lemma~\ref{lem:we get a sheaf} shows that $\cZ$ is a stack satisfying~[3],
and so it makes sense
(following Definition~\ref{def:properties defined by base-change})
to assert that $\xibar$ is proper and
surjective.
\end{remark}

\begin{proof}[Proof of Lemma~\ref{lem:proper and surjective}.]
Following Definition~\ref{def:properties defined by base-change},
we have to show that if $\cY$ is an algebraic
stack, and if
$\mathcal Y \to \mathcal Z$ is a morphism of stacks,
then the base-changed morphism of algebraic stacks $\mathcal Y\times_{\mathcal Z} \mathcal X \to \mathcal Y$
is proper and surjective.  Since by definition $\mathcal Z \to \mathcal F$ is 
a monomorphism, the fibre product over $\mathcal Z$ 
is isomorphic to the fibre product $\mathcal Y \times_{\mathcal F} \mathcal X,$
taken over $\mathcal F$,
and so the properness of $\xi$ immediately implies the properness of $\overline{\xi}$. 

To verify the surjectivity, we note that surjectivity can be checked 
after pulling back by a surjective morphism.  Replacing $\mathcal Y$
by a cover of $\mathcal Y$ by a scheme, we may assume that $\mathcal Y$
is in fact a scheme.  Since we may also check surjectivity locally
on the target, we may in fact assume that $\mathcal Y$ is an affine 
scheme $T$.  Since $\cZ$ satisfies ~[1], we may factor the morphism
$T \to \cZ$ through an affine $S$-scheme of finite type, and
hence (since surjectivity is preserved under base-change)
may further assume that $T$ is of finite type over $S$.

Now $T\times_{\cZ} \mathcal X \to T$ is a proper morphism of algebraic stacks,
and so to check that it is surjective, it suffices to check
that its image contains each finite type point of~$T$. Thus it suffices
to show that if $t: \Spec k \to T$ is any finite type point, then the fibre of $\mathcal X$ over $t$ is non-empty. 
Since $\xi$ is proper, it follows from
Lemma~\ref{lem: scheme theoretic dominance for Artin proper}
that the fact that the composite $t:\Spec k\to T\to\cZ\to\cF$ factors through $\cZ$
implies
that the fibre $\mathcal X_t \neq
\emptyset.$ 
\end{proof}

\begin{lemma}
\label{lem:minimality}
If $\cF'$ is a closed Zariski substack of~$\cF$,
and if the morphism $\xi: \cX \to \cF$ factors through~$\cF'$,
then $\cF'$ contains $\cZ$.

Conversely, assume either that $\cF$ satisfies~{\em [1]} and~{\em [2](b)}, or that
 $\xi$ is proper,
$\cF$ is~$\mathbf{Art}^{\mathbf{triv}}$-homogeneous, and either $\cX$
is locally of finite presentation over $S$, or $\cF$ satisfies~{\em [1]}. Then if $\cF'$ is a closed substack of $\cF$ which contains~$\cZ$,
then the morphism $\xi: \cX \to \cF$ factors through $\cF'$.
\end{lemma}
\begin{proof}
Suppose that $\cF'$ is a closed Zariski substack of $\cF$,
and that $\xi$ factors through $\cF'$; we will show that
$\cF'$ contains $\cZ$. Since $\cZ$ is (by definition) a Zariski substack of $\cF$, it suffices to show that if $T \to \cZ$ is a morphism with $T$
an affine $S$-scheme, then the composite $T \to \cZ \to \cF$
factors through $\cF'$.  By definition, given such a morphism
$T \to \cZ$, we may find a finite type affine $S$-scheme $T'$
and a factorisation $T \to T' \to \cZ$; thus we may and do
assume that $T$ is of finite type over $S$.

In order to show that $T \to \cF$ factors through $\cF'$,
it suffices to show that the base-changed closed immersion
$T \times_{\cF} \cF' \hookrightarrow T$ is an isomorphism.
Since $T$ is of finite type over $S$, for this it suffices
to show that if $t \in T$ is a finite type point,
then for any $n \geq 0,$
the induced morphism $\Spec \widehat{\cO}_{T,t}/\mathfrak m_t^n \to T$
factors through $T \times_{\cF} \cF'$, or equivalently,
that the composite $\Spec \widehat{\cO}_{T,t}/\mathfrak m_t^n \to T \to \cF$
factors through $\cF'$.
Thus we are reduced to the case when $T = \Spec A$,
with $A$ an Artinian $\cO_S$-algebra of finite type.

By the definition of $\cZ$, we may find a complete Noetherian local $\cO_S$-algebra
$B$, a morphism $\Spec A \to \Spec B$ over $S$, and a morphism $\Spec B \to \cF$ inducing the given morphism $\Spec A \to \cF$, such that the base-changed
morphism $\cX_B \to \Spec B$ is scheme-theoretically dominant.
Since $\xi$ factors through $\cF'$ by assumption,
the morphism $\cX_B \to \Spec B$ factors through the closed immersion
$\Spec B \times_{\cF} \cF' \hookrightarrow \Spec B$.  Since the former morphism
is also scheme-theoretically dominant, we see that this closed immersion
is necessarily an isomorphism, and hence that the morphism $\Spec B \to \cF$ 
factors through $\cF'$.  Thus the morphism $\Spec A \to \cF$
also factors through $\cF'$, and we are done.

For the converse, we note that, under the additional hypotheses, the morphism $\xi$ factors through a morphism
$\xibar: \cX \to \cZ$, by Lemma~\ref{lem:xibar}; so if $\cF'$ contains $\cZ$,
we see that $\xi$ factors through $\cF'$, as required. 
\end{proof}
Under the assumption that $\xi$ is proper, we may strengthen
Lemma~\ref{lem:minimality} as follows.

\begin{lem}
  \label{lem:minimality for finite type monomorphisms}Assume either that $\cF$ satisfies~{\em [1]}, or that
$\cX$ is locally of finite presentation over $S$. Assume also that 
$\xi$ is proper and that $\cF$ is~$\mathbf{Art}^{\mathbf{triv}}$-homogeneous.
If $\cF'$ is a substack of~$\cF$, and if the monomorphism $\cF'\into\cF$ is
representable by algebraic spaces and of finite type, then the morphism $\xi: \cX \to \cF$ factors through $\cF'$
if and only if $\cF'$ contains $\cZ$.
\end{lem}
\begin{proof}The ``if'' direction was proved in
  Lemma~\ref{lem:minimality}. For the ``only if'' direction,
set $\cZ':=\cZ\times_\cF\cF'$. We begin by showing that the finite type
  monomorphism $\cZ'\to\cZ$ is a closed immersion; it suffices to show that it
  is proper, and since monomorphisms are automatically separated, and since
it is finite type by assumption, it in fact
  suffices to show that it is universally closed. 

Firstly, note that since $\cX\to\cZ$ is proper and surjective by Lemma~\ref{lem:proper and
  surjective}, the composite $|\cX|\to|\cZ'|\to|\cZ|$ is surjective, while
$|\cZ'|\to|\cZ|$ is injective. Thus $|\cX|\to|\cZ'|$ is surjective. Similarly,
$|\cX|\to|\cZ|$ is closed, and a trivial topological argument then shows that
$|\cZ'|\to|\cZ|$ is closed. Since both properness and surjectivity are preserved
by arbitrary base changes, we conclude that $\cZ'\to\cZ$ is universally closed,
as claimed, and thus a closed immersion.

 By Lemma~\ref{lem:proper and surjective}, the morphism
$\xibar:\cX\to\cZ$ is proper and surjective; but it also factors
through the closed substack~$\cZ'$ of~$\cZ$, so we must
have~$\cZ'=\cZ$. Thus  $\cF'$ contains~$\cZ$, as required.
\end{proof}

\begin{remark}
If $\cF$ is an algebraic stack that is locally of finite presentation
over~$S$, then in particular it satisfies Axioms~[1], [2], and [3]
(by Theorem~\ref{thm: Artin representability}), 
and so all of the previous results apply.
\end{remark}
\begin{rem}\label{rem: representability of monomorphism}
  Note that at this point in the development of the theory, we do not know that
  the monomorphism $\cZ\to\cF$ is representable by algebraic spaces, which means
  for example that it is hard to use Lemma~\ref{lem:minimality for finite type
    monomorphisms} to uniquely characterise~$\cZ$ by a universal
  property. However, once we have shown that~$\cZ$ is an algebraic stack, it is
  immediate from the assumption that~$\cF$ satisfies~[3] that the monomorphism
  $\cZ\to\cF$ \emph{is} representable by algebraic spaces, and 
  Proposition~\ref{prop:if Z to F is representable, then it is a
    closed immersion} below can then be applied to deduce that~$\cZ$ is in fact a closed substack
  of~$\cF$.\end{rem}
\begin{lemma}
\label{lem:algebraic case}
If $\cF$ is an algebraic stack, locally of finite presentation over~$S$, then $\cZ$
coincides with the scheme-theoretic image of $\xi$ {\em (}in the sense of Definition~{\em \ref{defn: scheme theoretic image for algebraic stacks})}.
\end{lemma}
\begin{proof}
Let $\cY \hookrightarrow \cF$ denote the scheme-theoretic image
of $\xi$.  Since $\cY$ is a closed substack of $\cF$,
and since the morphism $\xi$ factors through $\cY$, 
Lemma~\ref{lem:minimality} shows that $\cY$ contains $\cZ$.
On the other hand, since $\cY \to \cF$ is a monomorphism and $\xi$
factors through $\cY$, we see that the
projection $\cX_{\cY} := \cX \times_{\cF} \cY \to \cY$
is naturally identified with the canonical morphism $\cX \to \cY$,
and so in particular is scheme-theoretically dominant.  
Lemma~\ref{lem:scheme-theoretic images} then implies that $\cY$ is contained
in $\cZ$.
\end{proof}

\begin{remark}
As noted in Remark~\ref{rem:laundry list}~(3), the conclusion
of this lemma needn't hold if $\cF$ doesn't satisfy~[1].
\end{remark}
The following result will be useful in Section~\ref{subsec:relatively representable}.
\begin{prop}
  \label{prop: if there is a closed algebraic substack then Z is
    it}Suppose
  that there is a closed substack~$\cY$ of~$\cF$ such that~$\cY$ is an
  algebraic stack, locally of finite presentation over~$S$, and that the
  morphism $\cX\to\cF$ factors as a composite $\cX\to\cY\to\cF$.
  Then $\cZ$ coincides with the scheme-theoretic image of $\cX$ in $\cY$
{\em (}in the sense of
Definition~{\em \ref{defn: scheme theoretic image for algebraic stacks})}. 
\end{prop}
\begin{proof}
Replacing $\cY$ by the scheme-theoretic image of $\cX \to \cY$,
we may assume that it coincides with this scheme-theoretic image; we
must then show that $\cY$ equals $\cZ$.
  By Lemma~\ref{lem:minimality}, we see that $\cY$ contains~$\cZ$. To show the
  reverse inclusion, we note first that Lemma~\ref{lem:algebraic case}
  implies that $\cY$ is the scheme-theoretic image of the morphism
$\cX \to \cY$;   the desired inclusion now follows by Lemma~\ref{lem:images compatible with monos}.
\end{proof}

The following lemma will be used to help prove openness of versality for~$\cZ$.
\begin{lemma}
\label{lem:dominance over a neighbourhood} Assume that $\cF$ satisfies~{\em
  [1]}, that $\xi$ is proper, and that $\cF$ is~$\mathbf{Art}^{\mathbf{triv}}$-homogeneous.
If $\cY$ is an algebraic stack, locally of finite type over $S$,
and if $\cY \to \cZ$ is a morphism over $S$ which is formally
smooth at a finite type point $y \in |\cY|$, then there 
exists an open neighbourhood $\cU$ of $y$ in $\cY$ such
that the induced morphism $\cX_{\cU} \to \cU$ is
scheme-theoretically dominant {\em (}where $\cX_{\cU}$
denotes the base-change of $\xibar: \cX \to \cZ$
via the composite $\cU \to \cY \to \cZ$; or equivalently,
the base-change of $\xi: \cX \to \cF$ via the composite
$\cU \to \cY \to \cZ \to \cF${\em )}.
\end{lemma}
\begin{proof}
By definition (see Definition~\ref{def:formal smoothness at a point for stacks}) we may find
a smooth morphism $V \to \cY$ from a scheme $V$ to $\cY$,
and a finite type point $v \in V$ mapping to $y$,
such that the composite $V \to \cY \to \cZ$
is formally smooth at $v$.  Write $B := \widehat{\cO}_{V,v},$
let $\cX_B \to \Spec B$ denote the 
base-change of $\xibar: \cX \to \cZ$ via the composite
$\Spec B \to V \to \cY \to \cZ$ (or, equivalently,
the base-change of $\xi: \cX \to \cF$ via the
composite $\Spec B \to V \to \cY \to \cZ \to \cF$),
and let $\Spec B/I$ denote the scheme-theoretic image of this
base-change.

Let $A$ be any Artinian quotient of $B$.  We will show that 
the surjection $B \to A$ necessarily contains $I$ in its kernel.
This will show that $I = 0$ (since $A$ was arbitrary) and
hence that $\cX_B \to \Spec B$ is scheme-theoretically dominant. 
From this one sees that there is a neighbourhood $U$ of $v$ in $V$
such that $\cX_U \to U$ is scheme-theoretically dominant (here
$\cX_U$ has the evident meaning).  Letting $\cU$ denote
the image of $U$ in $\cY$ (an open substack of $\cY$), we
see that $\cX_{\cU} \to \cU$ is scheme-theoretically dominant,
as required.   (Scheme-theoretic dominance can be checked
{\em fpqc} locally, and in particular smooth locally, on the target.)

The composite $\Spec A \to \Spec B \to V \to \cY \to \cZ \to \cF$
factors through $\cZ$ by construction, and so by definition we
may find a complete Noetherian local $\cO_S$-algebra $C$ with
finite type residue field, a closed immersion of $\cO_S$-schemes
$\Spec A \to \Spec C$, and a morphism $\Spec C \to \cF$ inducing
the above morphism $\Spec A \to \cF$, such that the base-changed
morphism $\cX_C \to \Spec C$ is scheme-theoretically surjective.
The morphism $\Spec C \to \cF$ factors through $\cZ$,
by Lemma~\ref{lem:scheme-theoretic images}.  
Since $V \to \cZ$ is formally smooth at $v$, we may
lift the morphism $\Spec C \to \cZ$ to a morphism $\Spec C \to \Spec B$,
extending the composite morphism $\Spec A \to \Spec B$. 
Since $\cX_C \to \Spec C$ is surjective,
we conclude that the morphism $\Spec A \to \Spec B$ factors
through $\Spec B/I$, as required.
\end{proof}

\begin{prop}
  \label{prop:if Z to F is representable, then it is a closed immersion}Assume either that $\cF$ satisfies~{\em [1]}, or that
$\cX$ is locally of finite presentation over $S$. Assume also that $\xi$ is proper,
that
$\cF$ is $\mathbf{Art}^{\mathbf{triv}}$-homogeneous,
and that the
monomorphism $\cZ\to\cF$ is representable by algebraic spaces.
Then~$\cZ$ is a closed substack of~$\cF$.
\end{prop}
\begin{proof}
We must show (under either set of finiteness assumptions
on $\cX$ and $\cF$) that
$\cZ \to \cF$ is a closed immersion.  Since it is a monomorphism,
it suffices to show that it is proper.  As monomorphisms are separated,
and as it is locally of finitely presentation  (this follows from Lemma~\ref{lem:morphism between limit
  preserving}, applied to the composite $\cZ\into\cF\to S$, together
with Lemma~\ref{lem:limit preserving vs. limit preserving on objects} and~\cite[\href{http://stacks.math.columbia.edu/tag/06CX}{Tag
  06CX}]{stacks-project}) 
it suffices to show that it is universally closed and
quasi-compact. These are properties that (by definition) can be checked after pulling back 
by a morphism $W\to\cF$ from a scheme.

Pulling back the morphisms $\cX \to \cZ \to \cF$ via such a morphism,
we obtain morphisms $\cX_{W} \to \cZ_{W} \to W$
with the composite being 
proper, and the first being surjective (by Lemma~\ref{lem:proper
and surjective}).   We must show that $\cZ_{W} \to W$ is closed,
and that, if $W$ is quasi-compact, the same is true of $\cZ_{W}$.
These are properties that can be checked on the underlying topological spaces.
So we consider the  continuous morphisms $|\cX_{W}| \to |\cZ_{W}| \to |W|$,
with the first being surjective (by~\cite[\href{http://stacks.math.columbia.edu/tag/04XI}{Tag 04XI}]{stacks-project}) and the
composite being closed.  Furthermore, if $|W|$ is quasi-compact,
the same is true of $|\cX_{W}|$. It follows immediately that
the second arrow is closed, and that $|\cZ_{W}|$ is quasi-compact
if $|W|$ is.  This completes the proof of the proposition. \end{proof}

We will now prove  Theorem~\ref{thm:main-intro}, as well as a variant where we
assume that $\cF$ satisfies~[1], but make no finiteness assumption on~$\cX$. (In
fact, as remarked above, we will reduce Theorem~\ref{thm:main-intro} to this case.)
\begin{theorem}
\label{thm:main-intro assuming [1]} Let $S$ be  a locally Noetherian
scheme, all of whose local rings $\mathcal O_{S,s}$ at finite type
points
$s \in S$ are $G$-rings. 
Suppose that $\xi: \cX \to \cF$ is a proper morphism, where $\cX$
is an algebraic stack and $\cF$ is a
stack  over $S$
satisfying~{\em [3]}, and that $\cF$ admits versal rings at all
finite type points.  Assume also either that $\cF$ satisfies~{\em [1]}, or that
$\cX$ is locally of finite presentation over $S$. Let $\cZ$ denote the scheme-theoretic image of $\xi$ as in Definition~{\em \ref{def:stack-theoretic image}}, and
suppose that $\cZ$ satisfies~{\em [2]}.  

Then $\cZ$ is an algebraic stack, locally of finite presentation over $S$; 
the inclusion $\cZ \hookrightarrow \cF$
is a closed immersion; and the morphism $\xi$ factors 
through a proper, scheme-theoretically surjective morphism $\overline{\xi}: \cX
\to \cZ$. 
Furthermore, if $\cF'$ is a substack of $\cF$ for which
the monomorphism $\cF' \hookrightarrow \cF$ is of finite type
{\em (}e.g.\ a closed substack{\em )} with the property that
$\xi$ factors through $\cF'$, then $\cF'$ contains~$\cZ$.
\end{theorem}
\begin{proof} Note firstly that by Lemma~\ref{lem:diagonal}, the diagonal of~$\cF$
is locally of finite presentation, so it follows from Corollary~\ref{cor: 3 and 1 and versal rings
    implies arttriv.} that~$\cF$ is $\mathbf{Art}^{\mathbf{triv}}$-homogeneous. Taking into account Lemmas~\ref{lem:xibar},  \ref{lem:proper and surjective}
and~\ref{lem:minimality for finite type monomorphisms}, Remark~\ref{rem: representability of monomorphism}, and Proposition~\ref{prop:if Z to F is representable, then it is a closed immersion},
we see that we only need to prove that $\cZ$ is an
  algebraic stack, locally of finite presentation over $S$. By Theorem~\ref{thm: Artin representability}, it is enough to
  show that it is an \'etale stack over $S$, and satisfies Axioms~[1], [2], [3] and~[4]. It follows from
  Lemma~\ref{lem:we get a sheaf} that $\cZ$ is an \'etale stack and satisfies Axioms~[1] and~[3], and we
  are assuming that it satisfies Axiom~[2]. It remains to show that $\cZ$ also
  satisfies Axiom~[4], i.e.\ openness of versality.  We will assume
  from now on that we are in the case that $\cF$ satisfies~[1], and explain at the end how to reduce the other case to this
  situation.

Consider an $S$-morphism $T \to \cZ$,
where $T$ is a locally finite type $S$-scheme,
and let $t \in T$ be a finite type point at which this
morphism is versal. 
In fact this morphism is then formally smooth at $t$, by Lemma~\ref{lem:versions of smoothness}~(2). Lemma~\ref{lem:monomorphism} allows us, shrinking $T$ if necessary, to factor the morphism $T \to \cZ$ as
$$T \to \cZ' \to \cZ,$$
where $\cZ'$ is an algebraic stack, locally of finite presentation over $S$,
the first arrow is a smooth surjection,
and the second arrow is locally of finite presentation, unramified,
and representable by algebraic spaces, and formally smooth at the image $t' \in |\cZ'|$ of $t$.

Consider the following diagram, both squares in which are defined to be Cartesian:
\numequation
\label{eqn:first Cartesian diagram}
\xymatrix{T\times_{\cZ'}\cX'\ar[r]\ar[d]& \cX' \ar[r] \ar[d] & \cX \ar[d] \\ T \ar[r] &  \cZ' \ar[r] & \cZ }
\end{equation}
The vertical arrows are proper, and hence closed,
the first of the horizontal arrows are smooth surjections,
and the second of the horizontal arrows are unramified,
locally of finite presentation, and representable by algebraic spaces. By
Lemma~\ref{lem:formal smoothness at a point for stacks}~(1),
there is an open substack $\cU$ of $\cX'$, containing the fibre over~$t'$,
such that the induced morphism $\cU \to \cX$ is smooth.

Since $\cX' \to \cZ'$ is closed,
the complement $\cU'$ of the image in~$\cZ'$ of the complement 
in~$\cX'$ of $\cU$ is an open substack of $\cZ'$, containing $t'$, 
whose preimage in $\cX'$ maps smoothly to $\cX$.
Thus, replacing $\cZ'$ by $\cU'$, and $T$ by its preimage,
we may further assume that $\cX' \to \cX$ is smooth,
and hence (being also unramified) \'etale.

Let $t_1 \in T$ be another finite type point.  We will
show that the morphism $T \to \cZ$ is formally smooth (and so also versal)
at $t_1$.  This will establish~[4] for $\cZ$.
 Let $t_1'$ denote the image of $t_1$ in $\cZ'$.  Since $T\to \cZ'$
 is smooth, it suffices to show that $\cZ' \to \cZ$ is formally smooth
 at $t_1'$.

Let $\overline{t}_1$ denote the image of $t_1$ in $\cZ$. 
By Lemma~\ref{lem:formally smooth charts},
we may find an $S$-scheme~$T_1$, locally of finite type,
a  morphism $T_1 \to \cZ$, and a point $\tilde{t}_1 \in T_1$
mapping to $\overline{t}_1$, such that the morphism $T_1 \to \cZ$ is formally smooth
at $\tilde{t}_1$.
We apply Lemma~\ref{lem:monomorphism} once more,
to obtain a factorisation $T_1 \to \cZ'' \to \cZ$ (possibly
after shrinking $T_1$ around~$\tilde{t}_1$), where $\cZ''$ is an
algebraic stack, locally of finite presentation over $S$,
the first arrow is a smooth surjection, and the second arrow is 
locally of finite presentation, unramified, representable by algebraic
spaces, and formally smooth at the image $t''\in | \cZ''|$ of $\tilde{t}_1$.
If we let $\cX''$ denote the base-change of $\xibar:\cX\to \cZ$ over $\cZ''$,
then Lemma~\ref{lem:dominance over a neighbourhood} shows
that, shrinking $\cZ''$ around $t''$ if necessary,
we may assume that the morphism $\cX'' \to \cZ''$ is scheme-theoretically
dominant, and thus scheme-theoretically surjective (being 
a base-change of the morphism $\xibar$, which is surjective
by Lemma~\ref{lem:proper and surjective}).

Now consider the Cartesian square
$$\xymatrix{\cZ'\times_{\cZ} \cZ'' \ar[r] \ar[d] & \cZ''\ar[d] \\
\cZ' \ar[r] & \cZ }
$$
Applying Lemma~\ref{lem:formal smoothness at a point for stacks},
we find that the morphism $\cZ'\times_{\cF} \cZ'' \to \cZ'$ contains
the fibre over $t''$ in its smooth locus.  Since $t'_1$ and $t''$ map
to the same point of $|\cF|$, this fibre contains a point $t'''$ lying over $t'_1$ and $t''$,
and so we may find an open substack $\cZ''' \subseteq \cZ'\times_{\cF} \cZ''$,
such that  $t''' \in |\cZ'''|$,
and such that the morphism $\cZ''' \to \cZ'$ is smooth, and hence
(being also unramified) \'etale.

In summary, we have a diagram 
$$\xymatrix{\cZ''' \ar[r]\ar[d] & \cZ'' \ar[d] \\
\cZ' \ar[r] & \cZ}
$$
in which the horizontal arrows are locally of finite presentation
and unramified,
the left-hand vertical arrow is \'etale,
the right-hand vertical arrow is formally smooth at $t''$,
and for which there is a point $t''' \in |\cZ'''|$ lying
over $t_1' \in |\cZ'|$ and $t'' \in |\cZ''|$.
Our goal is to show that the lower horizontal arrow is
formally smooth at the point~$t_1'$.  
For this, it suffices (by Remark~\ref{rem:formal smoothness at a point
for stacks}) to show that the upper horizontal arrow
is smooth (or equivalently, \'etale, since it is unramified).

Thus we turn to showing that $\cZ'''\to\cZ''$ is \'etale. Form the Cartesian diagram of algebraic stacks
$$\xymatrix{\cX''' \ar[r] \ar[d] & \cX'' \ar[d] \\ \cZ''' \ar[r] 
& \cZ''}$$
Regarding the top arrow as being the composition of an open immersion with the
base-change of the top
arrow of~(\ref{eqn:first Cartesian diagram}), 
we see that it is \'etale.  Recall
that
the right-hand vertical arrow $\cX''\to\cZ''$ is scheme-theoretically surjective.
Both vertical arrows are proper, thus universally closed,
and the bottom arrow is locally of finite presentation and unramified. Lemma~\ref{lem:unramified with etale base-change} below thus implies that the bottom
arrow is \'etale, and so Axiom~[4] for $\cZ$ is proved.

We have therefore shown that if $\cF$ satisfies~[1], then $\cZ$ is an algebraic
stack, locally of finite presentation over~$S$. Assume now that we are
in the case that $\cX$ is locally of finite presentation over
$S$. Recall that by Lemma~\ref{lem:diagonal}, the diagonal of~$\cF$
is locally of finite presentation. By Remark~\ref{rem: some
  generalities on reducing to the case that F satisfies [1]}, $\xi$ factors
as \[\cX\stackrel{\xi'}{\to}\pro(\cF{|_\AfffpS})\into\cF,\] and
$\pro(\cF{|_\AfffpS})$ is an \'etale stack which satisfies [1] and [3] and admits versal
rings at all finite type points; and the
scheme-theoretic image of $\xi'$ is $\cZ$. Regarding $\xi'$ as the pull-back of
$\xi$ along the embedding
$\pro(\cF{|_\AfffpS})\into\cF$, we see that $\xi'$ is proper, and therefore satisfies the
assumptions of the Theorem; so it follows from the case already proved that
$\cZ$ is an algebraic stack, locally of finite presentation over~$S$. \end{proof}

\medskip

Although we have stated and
proved the following lemma in what seems to be its natural level of generality,
the only application of it that we make is in the case when the morphism
$\cZ' \to \cZ$ is in fact representable by algebraic spaces.

\begin{lemma}
\label{lem:unramified with etale base-change}
Let $\cY \to \cZ$ 
be a quasi-compact morphism of algebraic stacks that 
is scheme-theoretically surjective and universally closed,
and let $\cZ'\to \cZ$ be a morphism of algebraic stacks that is
locally of finite presentation and unramified. 
If the base-changed morphism $\cY'\to \cY$ is \'etale, then $\cZ'\to \cZ$ is also \'etale.
\end{lemma}  
\begin{proof}
Since $\cZ'\to \cZ$ is unramified, its diagonal is \'etale
and in particular unramified,
so $\cZ'\to \cZ$ is a DM morphism in the sense
of~\cite[\href{http://stacks.math.columbia.edu/tag/04YW}{Tag
  04YW}]{stacks-project}. Choose a scheme mapping $Z$ via a smooth surjection to $\cZ$, and pull
everything back over~$Z$;
in this way we reduce to
the case $\cZ=Z$. As explained in~\cite[\href{http://stacks.math.columbia.edu/tag/04YW}{Tag
  04YW}]{stacks-project}, it now follows
from~\cite[\href{http://stacks.math.columbia.edu/tag/06N3}{Tag
  06N3}]{stacks-project} that there is a scheme $Z'$ and a surjective \'etale
morphism $Z'\to\cZ'$.\footnote{As noted above,
the only application of this lemma that we will make is in the case
when $\cZ' \to \cZ$ is representable by algebraic spaces.  In this case,
the pull-back of the algebraic stack $\cZ'$ over the scheme $Z$ is
already an algebraic space, and so by definition admits an \'etale
surjection from a scheme.  Thus, in this case, we can avoid appealing
to the theory of DM morphisms and their relationship to Deligne--Mumford
stacks.}

 Replacing $\cZ'$ by $Z'$, we have reduced to the case that
the morphism $\cZ'\to\cZ$ is the morphism of schemes $Z'\to Z$.

We must show that $Z'\to Z$ is \'etale at each point $z'\in Z'$. Fix such a
point, with image $z\in Z$. Then we may
replace $Z$ by its base change to the strict henselisation
$\Spec\cO_{Z,z}^{\textrm{sh}}$, and thus assume that $Z$ is a local strictly
henselian scheme with closed point $z$. Then the \'etale local structure of
unramified morphisms~\cite[\href{http://stacks.math.columbia.edu/tag/04HH}{Tag
  04HH}]{stacks-project} shows that after shrinking $Z'$ around $z'$, we can
arrange that $Z'\to Z$ is a closed immersion, so that in particular $Z'$ is also
local with closed point $z'=z$. Then $Z'\to Z$ is \'etale if and only if $Z'=Z$.

Since $Z'\to Z$ is a closed immersion, so is $\cY'\to\cY$. Since $\cY'\to\cY$ is
\'etale by assumption, it is an \'etale monomorphism, and therefore it is an
open (as well as closed) immersion. The complement in $\cY$ of the image of
$\cY'\to\cY$ is therefore closed, and has closed image in $Z$, as $\cY\to Z$ is
closed by assumption; but this image has empty special fibre, and so this
complement is empty, and $\cY'=\cY$.  It now follows from Lemma~\ref{lem:isomorphism
 criterion} that $Z'=Z$, as required.
\end{proof}

\begin{remark}
\label{rem:counterexample}
In the statement of Lemma~\ref{lem:unramified with etale base-change}, it would not suffice to assume
that $\cY\to \cZ$ is merely scheme-theoretically surjective. 
For example, if $\cZ = Z$ is taken to be the scheme given as the union of
two lines $Z_1$ and $Z_2$ crossing at the point $y$ in the plane,
if $\cY = Y$ is taken to be the scheme obtained as the disjoint union $Z_1 \,\coprod \,  Z_2 \! \setminus \! \{y\},$
if $Y \to Z$ is taken to be the obvious morphism (which is scheme-theoretically
surjective), 
and if $\cZ' = Z'$ is taken to be $Z_1$,
then the closed immersion $Z'\to Z$ is unramified,
and the base-changed morphism $Y' \to Y$ is an open immersion,
and so in particular \'etale.  
On the other hand, the closed immersion $Z' \to Z$ is certainly not \'etale.
\end{remark}

\subsection{Base change}\label{subsec:relatively representable}In
this section we study the behaviour of our scheme-theoretic images under base
change by a morphism which is representable by algebraic spaces. We will only
need to consider cases where we know (in applications, as a
consequence of Theorem~\ref{thm:main-intro}) that the scheme-theoretic image
of the morphism being base-changed is an algebraic stack, and we have
therefore restricted to this case, and have not investigated the
compatibility with base change in more general situations.

  We first of all consider the question of base change of the target stack~$\cF$.

\begin{prop}
  \label{prop: representable by algebraic spaces base change - new
    version}Suppose that we have a commutative
diagram of \'etale stacks in groupoids over a locally Noetherian base scheme~$S$ \[\xymatrix{\cX'\ar[r]\ar[d]&\cX\ar[d]\\ \cF'\ar[r]&\cF}\]{\em (}which is not assumed to be Cartesian{\em )}, in which $\cX,\cX'$ are
algebraic stacks, $\cF'\to\cF$ is representable by algebraic spaces
and locally of finite presentation,
and $\Delta_{\cF}$ is representable by algebraic spaces and locally of
finite presentation. Write $\cZ$ and $\cZ'$ for the scheme-theoretic
images of $\cX\to\cF$ and $\cX'\to\cF'$ respectively,
and assume further that $\cZ$ is an
  algebraic stack, that~$\cZ$ is a closed substack of~$\cF$, and that
  the morphism $\cX\to\cF$ factors through a scheme-theoretically dominant morphism $\cX\to\cZ$.

Then $\cZ'$ is an algebraic stack, locally of finite presentation over~$S$. In
fact, the morphism $\cX'\to\cF'$ factors through the algebraic stack
$\cZ\times_{\cF}\cF'$
{\em (}which is in turn a closed substack of $\cF'${\em )}, and $\cZ'$
is the scheme-theoretic image of the induced morphism of algebraic stacks
$\cX'\to\cZ\times_{\cF}\cF'$.
\end{prop}
\begin{proof}Since $\cX\to\cF$ factors through~$\cZ$ by assumption,
  the composite $\cX'\to\cX\to\cF$ factors through~$\cZ$, and so
  $\cX'\to\cF'$ factors through the algebraic stack
  $\cZ\times_{\cF}\cF'$, which is a closed substack of~$\cF'$. By Proposition~\ref{prop: limit preserving and base
    change}, $\Delta_{\cF'}$ is representable by algebraic spaces and
  locally of finite presentation. The result then follows from Proposition~\ref{prop:
    if there is a closed algebraic substack then Z is it}.
\end{proof}

\begin{cor}
  \label{cor: base change with main theorem incorporated, for use in
    CEGS}Suppose that we have a commutative
diagram of \'etale stacks in groupoids over a locally Noetherian base scheme~$S$ \[\xymatrix{\cX'\ar[r]\ar[d]&\cX\ar[d]\\ \cF'\ar[r]&\cF}\]{\em (}which is not assumed to be Cartesian{\em )}, in which $\cX,\cX'$ are
algebraic stacks, $\cF'\to\cF$ is representable by algebraic spaces
and locally of finite presentation,
and $\Delta_{\cF}$ is representable by algebraic spaces and locally of
finite presentation. 

Suppose that $\cX \to \cF$ is  proper, and that 
 $\cX$
is locally of finite presentation over~$S$. Write $\cZ$, $\cZ'$
 for the scheme-theoretic images of $\cX\to\cF$,
$\cX'\to\cF'$ respectively
Suppose that $\cZ$ satisfies~{\em [2]}, and  that $\cF$ admits
versal rings at all finite type points. Then $\cZ$ and $\cZ'$ are both algebraic stacks, locally of finite presentation over~$S$. In
fact, the morphism $\cX'\to\cF'$ factors through the algebraic stack
$\cZ\times_{\cF}\cF'$
{\em (}which is in turn a closed substack of $\cF'${\em )}, and $\cZ'$
is the scheme-theoretic image of the induced morphism of algebraic stacks
$\cX'\to\cZ\times_{\cF}\cF'$.

If $x$ is a finite type point of~$\cZ'$, let $R_x$ be a
versal ring for the corresponding finite type point of $\cF$, and write
$\Spf S_x$ for the complete local ring at~$x$ of $\Spf
R_x\times_{\cF}\cF'$, in the sense of
Definition~\emph{\ref{df:Ind complete
  local rings}} below. Then a 
versal ring for~$\cZ'$ at~$x$ is given by the scheme-theoretic image of the
induced morphism $\cX'_{\Spf S_x}\to\Spf S_x$.\end{cor}
\begin{proof}
Everything except for the claim about versal rings is immediate from
Theorem~\ref{thm:main-intro} and Proposition~\ref{prop: representable by algebraic spaces base change - new
    version}. The description of the versal rings follows from
  Lemmas~\ref{lem:versal rings and base change} and~\ref{lem:versal ring for Z}.
\end{proof}

We now consider base changes $S'\to S$. It is of course unreasonable to expect that the
formation of scheme-theoretic images is compatible with such base changes unless
we make a flatness hypothesis; but under this condition, we are able
to prove the following general result.

\begin{prop}
  \label{prop: base change of S, new version}Suppose that $\cX\to\cF$
  is a morphism of \'etale stacks over the locally Noetherian base
  scheme~$S$, where $\cX$ is algebraic and~$\Delta_{\cF}$ is
  representable by algebraic spaces.
  Suppose that the scheme-theoretic image~$\cZ$ of
  $\cX\to\cF$ is an algebraic stack,  that $\cZ\to\cF$ is a closed
  immersion, and that $\cX\to\cF$ factors through a
  scheme-theoretically dominant
  morphism $\cX\to\cZ$.

Let $S'\to S$ be a flat morphism of locally Noetherian schemes, and
write $\cX_{S'}$, $\cF_{S'}$, and $\cZ_{S'}$ for the base changes of
$\cX$, $\cF$, and
$\cZ$ to~$S'$. Write $\cZ'$ for the scheme-theoretic image of
$\cX_{S'}\to\cF_{S'}$. 

Then $\cZ'=\cZ_{S'}$; so $\cZ'$ is an algebraic stack, 
$\cZ'\to\cF_{S'}$ is a closed immersion, and $\cX_{S'}\to\cF_{S'}$
factors through a scheme-theoretically dominant
  morphism $\cX\to\cZ'$.\end{prop}
\begin{proof}
 Note firstly that $\cZ_{S'}$ is an algebraic stack, that
 $\cZ_{S'}\to\cF_{S'}$ is a closed immersion, and that
 $\cX_{S'}\to\cF_{S'}$ factors through $\cZ_{S'}$. Since $S'\to S$ is
 flat, we see that $\cX_{S'}\to\cZ_{S'}$ is scheme-theoretically dominant by
 Remark~\ref{rem:properties of scheme-theoretic images of algebraic
   stacks}~(2). Since the formation of the diagonal is compatible with
 base change, $\Delta_{\cF_{S'}}$ is representable by algebraic spaces
 and locally of finite presentation. The result follows immediately from
 Proposition~\ref{prop: if there is a closed algebraic substack then Z
   is it}. 
\end{proof}

\section{Examples}\label{sec: examples}

\subsection{Quotients of varieties by proper equivalence relations}
\label{subsec:proper}
We indicate some simple examples of quotients of varieties (over $\mathbb C$,
so that we may form the quotients in the sense of topological spaces)
by proper equivalence relations which are not algebraic objects.

\begin{example}
\label{ex:sphere}
Consider the equivalence relation on $\mathbb P^2$ which contracts a line to a
point.   In the topological category, the quotient of $\mathbb P^2(\mathbb C)$
by this  equivalence relation is~$S^4$ (a $4$-sphere), which is not a K\"ahler manifold (since it has vanishing~$H^2$).  Thinking more algebraically, one can show that there is no proper morphism of algebraic spaces $\mathbb P^2 \to X$
for which $\mathbb P^2 \times_X \mathbb P^2$ coincides with
this equivalence relation;
indeed, if there were such a morphism,
then the theorem on formal functions would show that the complete local ring of $X$ at the point obtained
as the image of the contracted $\mathbb P^1$
is equal to the ring of global sections of 
the structure sheaf of the formal completion of $\mathbb P^2$ along
$\mathbb P^1$; but since the conormal bundle of $\mathbb P^1$ in $\mathbb P^2$ equals $\mathcal O(-1)$, this ring of global sections is just equal to
$\mathbb C$,
and so cannot arise as the complete local ring at a point of a two-dimensional
algebraic space.  Thus there isn't any reasonable way to take the quotient of
$\mathbb P^2$ by this equivalence relation in the world of algebraic spaces
over $\mathbb C$.
\end{example}

\begin{example}
\label{ex:endos of line}
Consider the space $X$ of endomorphisms of $\mathbb P^1$ of degree $\leq 1$.
Such an endomorphism is described by a linear fractional transformation
$x \mapsto \dfrac{ax+b}{cx+d}$ for which the matrix $\begin{pmatrix}a & b \\ c
& d\end{pmatrix}$ is non-zero, and so $X$ is a quotient of $\mathbb P^3$
(thought of as the projectivisation of the vector space $M_2(\mathbb C)$
of $2\times 2$ matrices).  However, the locus of singular matrices
(i.e.\ the subvariety cut out by the vanishing of the determinant),
which gives rise to the constant endomorphisms, is two-dimensional
(indeed a quadric surface in $\mathbb P^3$),
while the space of constant endomorphisms 
is obviously just equal to $\mathbb P^1$.  Thus $X$ can be thought of
as a quotient of $\mathbb P^3$ in which a quadric is contracted to
a projective line (via one of the projections $\mathbb P^1 \times \mathbb P^1
\to \mathbb P^1$).  Similarly to the preceding example,
the space $X$ can't be realised in the world of algebraic spaces.
(The first author would like to thank V.~Drinfeld for explaining this
example to him.)
\end{example}

\begin{remark}\label{rem: remark about Artin contractions}
It follows from \cite{MR0260747} that in contexts like
those considered in the preceding examples, in which we wish
to contract a closed subvariety of a proper variety in some manner,
{\em if} we have
a formal model for neighbourhood of a contraction, then we {\em can} perform 
the contraction in the category of algebraic spaces.    The computation
with complete local rings on the (hypothetical) quotient in 
Example~\ref{ex:sphere} can be thought of as 
showing that, for this example, such a formal model doesn't exist.

A key intuition behind the main theorem of the present paper is that,
as one of the hypotheses of the theorem,
we do assume that we have well-behaved complete local rings at
each finite type point of the quotient we are trying to analyse.
\end{remark}

\subsection{Ind-algebraic stacks} \label{subsec:non-representable}

Ind-algebraic stacks give simple examples of stacks which satisfy some but not all
of Artin's axioms, and which admit surjective morphisms from algebraic
stacks. We refer the reader to~\cite{Emertonformalstacks} for the
details of the theory of Ind-algebraic stacks, contenting ourselves
here with recalling some general facts, and giving various examples (several of which are taken from 
\cite{MR0262237}) which illustrate the roles of the various axioms in the
results of Section~\ref{sec: artin's axioms and stacks in groupoids}.
Example~\ref{ex: the zipper in a surface}
also illustrates the necessity of the properness
hypothesis in the statement of Theorem~\ref{thm:main-intro},
while Example~\ref{ex: infinitely many embedded points}
illustrates the necessity of the hypothesis of scheme-theoretic dominance
(rather than mere surjectivity)
in the statement of Corollary~\ref{cor:representability-intro}.

\begin{df}
\label{def:ind stacks}
We say that a stack $\cX$ over $S$ is an {\em Ind-algebraic stack}
if we may find
a directed system of algebraic stacks $\{\cX_i\}_{i \in I}$ over $S$ and an
isomorphism $\varinjlim_{i \in I} \cX_i\iso \cX$, the inductive limit
being computed in the $2$-category of stacks. 
If $\cX$ is a stack in setoids (i.e.\ is equivalent to a sheaf of sets),
then we say that $\cX$ is an {\em Ind-algebraic space},
resp.\ an {\em Ind-scheme},
if the $\cX_{i}$ can in fact be taken to be algebraic spaces,
resp.\ schemes.
\end{df}
\begin{rem}
  More properly, the inductive limit is a 2-colimit, but as
  in~\cite{Emertonformalstacks}, we find it more suggestive to use the usual notation
  for direct limits.\end{rem}
\begin{remark}
\label{rem:don't require closed immersions}
Our definitions of Ind-algebraic stacks, Ind-algebraic spaces, and Ind-schemes
are broader than usual.  Indeed,
at least in the case of Ind-schemes, it is conventional to require the 
transition morphisms $\cX_i \to \cX_{i'}$ to be closed immersions. 
We have adopted a laxer definition simply because it provides a convenient
framework in which to discuss many of the examples of
Subsection~\ref{subsec:examples} below.
\end{remark}

\begin{remark}\label{rem:Ind-stackification}
Given a directed system $\{\cX_i\}_{i \in I}$ as in Definition~\ref{def:ind stacks},
then for any $S$-scheme $T$
there is morphism
of groupoids
\numequation
\label{eqn:morphism of inductive limits}
\varinjlim_{i \in I} \cX_i(T) \to \bigl(\varinjlim_{i \in I} \cX_i\bigr)(T).
\end{equation}
Although this is not an equivalence in general, it {\em is} an equivalence
if $T$ is quasi-compact and quasi-separated 
(e.g.\ if $T$ is an affine $S$-scheme). 
We record this fact in the following lemma (which is presumably
well-known);
see e.g.~\cite[\href{http://stacks.math.columbia.edu/tag/0738}{Tag 0738}]{stacks-project}  
for the analogous statement in the context of sheaves. 

Also, although in the above definition all stacks involved are (following
our conventions) understood to be stacks for the \'etale site,
in the statement of the lemma we consider stacks for other topologies
as well.  
\end{remark}

\begin{lemma}
\label{lem:inductive limits of stacks}
Let $\{\cX_i\}_{i \in I}$ be an inductive system of Zariski,
\'etale, fppf, or fpqc stacks,
and consider the inductive limit $\varinjlim_{i \in } \cX_i,$
computed as a stack for the topology under consideration.
If $T$ is an quasi-compact $S$-scheme, then the morphism~{\em
(\ref{eqn:morphism
of inductive limits})} is faithful.
If, in addition, either  $T$ is quasi-separated, or the transition maps in the inductive system are monomorphisms,
then it is in fact an equivalence of groupoids.
\end{lemma}
\begin{proof}
We write $\cX := \varinjlim_{i \in I} \cX_i$ (the inductive limit
being taken as stacks).
If $\xi_i,\eta_i$ are a pair of objects in $\cX_i(T)$, inducing objects
$\xi_{i'}, \eta_{i'}$ in $\cX_{i'}(T)$ for each $i' \geq i,$
and objects $\xi,\eta$ in $\cX(T)$,
then the set of morphisms between $\xi_{i'}$ and $\eta_{i'}$ 
is equal to the space of global sections of the sheaf
$\cX_{i'} \times_{\Delta, \cX_{i'}\times_S \cX_{i'}, \xi_{i'} \times \eta_{i'}} T$,
while the set of morphisms between $\xi$ and $\eta$ is equal
to the space of global sections of
$\cX \times_{\Delta,\cX\times_S \cX,\xi\times \eta} T,$
which is also equal to the set of global sections of the inductive limit of
sheaves
$\varinjlim_{i' \geq i} 
\cX_{i'} \times_{\Delta, \cX_{i'}\times_S \cX_{i'}, \xi_{i'} \times \eta_{i'}} T.$

If $T$ is quasi-compact and quasi-separated,
then the inductive limit of global sections maps isomorphically
to the global sections of the inductive limit
\cite[\href{http://stacks.math.columbia.edu/tag/0738}{Tag 0738}]{stacks-project},
and so we find that the morphism~(\ref{eqn:morphism of inductive limits})
is fully faithful for such $T$.
The same reference shows that if $T$ is merely assumed
to be quasi-compact, then the morphism on global sections is injective,
and hence that the morphism~(\ref{eqn:morphism of inductive limits})
is faithful.
Finally, if the transition morphisms in the inductive system of sheaves
under consideration are injective,
then this reference again shows that the map on global
sections is an isomorphism provided merely that $T$ is quasi-compact.
If the transition morphisms $\cX_{i'} \to \cX_{i''}$ are monomorphisms,
then the corresponding transition morphisms 
$\cX_{i'} \times_{\Delta, \cX_{i'}\times_S \cX_{i'}, \xi_{i'} \times \eta_{i'}} T
\to
\cX_{i''} \times_{\Delta, \cX_{i''}\times_S \cX_{i''}, \xi_{i''} \times \eta_{i''}} T$
are indeed injective,
and hence in this case the morphism~(\ref{eqn:morphism of inductive limits})
is again fully faithful.

We turn to considering the essential surjectivity of~(\ref{eqn:morphism
of inductive limits}).   If $T$ is quasi-compact, 
then for any object
of $\cX(T)$, we may find a cover $T'\to T$
(in the appropriate topology: Zariski,
\'etale, \emph{fppf}, or \emph{fpqc}, as the case may be)
such that $T'$ is again quasi-compact, and a morphism
$T' \to \cX_i$ for some $i$ such that the diagram
$$\xymatrix{T' \ar[r] \ar[d] & \cX_i \ar[d] \\
T \ar[r] & \cX }
$$
commutes.
We obtain an induced morphism $T'\times_T T' \to \cX_i\times_{\cX} \cX_i$.
The target of this morphism may be written
as $\varinjlim_{i'\geq i} \cX_i\times_{\cX_{i'}}\cX_i.$
If $T$ is also quasi-separated, so that $T'\times_T T'$ is again quasi-compact,
then by what we have already observed,
this morphism arises from a morphism
$T'\times_T T' \to \cX_i\times_{\cX_{i'}}\cX_i$ for some sufficiently large
value of $i'$.  Consequently, we obtain a factorisation of the
induced morphism $T'\times_T T' \to \cX_{i'}\times_S \cX_{i'}$ through
the diagonal.   When we pass from $\cX_{i'}$ to $\cX$,
we obtain the descent data to $T$ of the composite $T'\to T\to \cX$,
and so the faithfulness result already proved shows
that, if we enlarge $i'$ sufficiently, this factorisation
provides descent data to $T$ for the morphism $T'\to \cX_{i'}$.  
Since $\cX_{i'}$ is a stack, we obtain a morphism $T \to \cX_{i'}$
inducing the original morphism $T \to \cX$.

If the transition morphisms in the inductive system $\{\cX_i\}$ are
monomorphisms, then the various fibre products $\cX_i\times_{\cX_{i'}}\cX_i$
are all isomorphic to $\cX_i$ (via the diagonal), and hence the natural map
$\cX_i \to \varinjlim_{i'\geq i} \cX_i\times_{\cX_{i'}} \cX_i$ is an
isomorphism.  Thus in this case, the morphism $T'\times_T T' \to
\cX_i\times_{\cX} \cX_i$ necessarily
arises from a morphism $T'\times_T T' \to \cX_i$,
without any quasi-compactness assumption on $T'\times_T T'$.
Thus, in this case, we obtain the essential surjectivity 
of~(\ref{eqn:morphism of inductive limits}) under the assumption
that $T$ is merely quasi-compact.
\end{proof}

\begin{remark}
\label{rem:Ind alg stacks can be computed Zariskiwise}
The preceding lemma (in particular, the fact that the
isomorphism~(\ref{eqn:morphism of inductive limits}) provides
an explicit description of $\bigl(\varinjlim_{i \in I}\cX_i\bigr)(T)$
which makes no reference to the site over which the inductive limit
is computed, when $T$ is an affine $S$-scheme), 
shows that to form the stack
$\varinjlim_{i\in I} \cX_i$, for an inductive systems of stacks in
the \'etale, {\em fppf}, or {\em fpqc} topologies, it in fact suffices
to form the corresponding inductive limit as a stack for the Zariski topology.
\end{remark}

\begin{remark}
\label{rem:inductive limits evaluated on stacks}
If $\{\cX_i\}_{i\in I}$ is an inductive system of stacks
(for the \'etale topology), 
and if $\cY$ is a category fibred in groupoids over $S$,
then, analogously to~(\ref{eqn:morphism of inductive limits}),
we have a natural morphism of groupoids
$$\varinjlim_{i \in I} \Mor(\cY,\cX_i)
\to \Mor(\cY,\varinjlim_{i \in I} \cX_i).$$
In particular if $\cY$ is an algebraic stack, then the obvious
extension of~ Lemma~\ref{lem:inductive limits of stacks} holds:
namely, if $\cY$ is quasi-compact then this morphism is faithful,
and if furthermore either $\cY$ is quasi-separated,
or the transition morphisms in the inductive system
are monomorphisms, then this morphism is an equivalence. 
(This can be checked by choosing a smooth surjection $T \to \cY$ from
an $S$-scheme to $\cY$, whose domain can be taken to be affine if $\cY$
is quasi-compact, and using the fact that morphisms from $\cY$
to any stack on the \'etale site can be identified with morphisms
from $T$ to the same stack with appropriate descent data; we leave
the details to the reader.)\end{remark}

\begin{remark}
\label{rem:axioms for Ind-algebraic stacks}
We briefly discuss the comportment of Ind-algebraic stacks with regard
to Artin's axioms. Suppose that
$\cX := \varinjlim_{i \in I } \cX_i$ is an Ind-algebraic stack.

(1)
If each $\cX_i$ satisfies~[1], then it follows from Remark~\ref{rem:Ind-stackification}
that $\cX$ also satisfies~[1].   

(2) It follows from Remark~\ref{rem:Ind-stackification} that 
$\cX$ satisfies~[2](a), since this is true of each~$\cX_i$.
The question of whether or not a particular Ind-stack satisfies~[2](b) is 
more involved, and we discuss it in some detail for Ind-algebraic spaces
in~\ref{subsubsec:versal rings} below.  

(3)  Typically, 
 $\cX$ need not satisfy~[3].  If $T$ is a quasi-compact and quasi-separated $S$-scheme,
then any morphism $T \to \cX \times_S \cX$ factors through a morphism
$T \to \cX_i \times_S \cX_i$
for some $i \in I$, and thus $\cX \times_{\cX \times_S \cX} T$ is
the base-change with respect to this morphism
of
$$\cX \times_{\cX \times_S \cX} (\cX_i \times_S \cX_i)
= \cX_i \times_{\cX} \cX_i = \varinjlim_{i' \geq i} \cX_i \times_{\cX_{i'}} 
\cX_i
$$
(equipped with its natural morphism to $\cX_i \times_S \cX_i$);
consequently the diagonal of $\cX$ is representable by Ind-algebraic
spaces. If each of the transition morphisms $\cX_i \to \cX_{i'}$ is a monomorphism,
then so is the morphism $\cX_i \to \cX$,
and we find that $\cX_i \iso \cX_i \times_{\cX} \cX_i$.  Thus in this case
we find that $\cX$ does satisfy~[3].  One can give examples 
in which the transition morphisms are {\em not} monomorphisms,
but for which nevertheless the diagonal of $\cX$ is representable by algebraic
spaces; see e.g.\
Examples~\ref{ex: line with infinitely many nodes} and~\ref{ex: the zipper
in a surface} below.  However, certainly in general the diagonal of 
$\cX$ is not representable by algebraic spaces, even for examples which seem
tame in many respects; see e.g.\ Example~\ref{ex: line with infinitely
many cusps}.

(4) Typically an Ind-algebraic stack will not satisfy~[4].  However,
since algebraic stacks satisfy~[4a]
(Lemma~\ref{lem:algebraic stacks satisfy 4a}), we see that
Ind-algebraic stacks will satisfy~[4a].  Typically they will not satisfy~[4b],
though.  In some cases, however, they may satisfy both~[4a] and~[4b],
but not~[4]; Example~\ref{ex: line with infinitely many nodes} illustrates
this possibility (and demonstrates that the quasi-separatedness hypotheses 
in Corollary~\ref{cor:4a and b give 4} and
Theorem~\ref{thm: Artin representability; variant}
are necessary).
\end{remark}

\subsubsection{Versal rings for Ind-algebraic spaces}
\label{subsubsec:versal rings}
If $x$ is a finite type point of the Ind-algebraic stack $\cX
:= \varinjlim_{i \in I} \cX_i$,
then $x$ arises from a finite type point $x_i$ of $\cX_i$ for some $i \in I$,
and if we write $x_{i'}$ to denote the corresponding finite type point
of~ $\cX_{i'}$, for each $i' \geq i$,
then one can attempt to construct a versal deformation ring at $x$
by taking a limit over the versal rings of each $x_{i'}$.
This limit process is complicated in general by the non-canonical
nature of the versal deformation ring of a finite type point in a stack,
and so we will restrict our precise discussion to the case of Ind-algebraic
spaces, for which we can find canonically defined minimal versal rings.

To be precise, suppose that $S$ is locally Noetherian.
Note that an algebraic space $X$ over $S$ satisfies~[1] if and only
if it is locally of finite presentation over $S$ (Lemma~\ref{lem:
limit preserving iff locally of finite presentation}),
which in turn holds if and only if it is locally of finite type over $S$
\cite[\href{http://stacks.math.columbia.edu/tag/06G4}{Tag 06G4}]{stacks-project}.  

\begin{df}
\label{df:Ind-lft alg space}
We say that a sheaf of sets $X$ on the \'etale site of $S$ 
is an {\em Ind-locally finite type algebraic space over} $S$ 
if there is an isomorphism
$X \cong \varinjlim_{i \in I}X_i$, where $\{X_i\}_{i \in I}$
is a directed system of algebraic
spaces, locally of finite type over $S$. 
\end{df}

In this case, we will be able to show that
each finite type point $x$ of $X$ admits a representative which
is a monomorphism, unique up to unique isomorphism, and 
for such a monomorphism $\Spec k \to X$, we will construct a canonical
minimal versal ring to $X$ at $x$.  
Despite this, it doesn't follow that $X$ satisfies~[2](b).  Loosely speaking,
there are two possible obstructions to~[2](b) holding: Firstly, it may
be that the versal ring at a point $x$ is not Noetherian. This happens 
when the dimension of the various $\cX_i$ at the points $x_i$ increases without
bound, so that the Ind-algebraic space $X$ is infinite dimensional; a typical
example is given by infinite-dimensional affine space (Example~\ref{ex:
infinite dimensional affine space}).
Secondly, even if the versal ring $x$ is Noetherian, it may not be effective.
This happens when the infinitesimal germs of $X_i$ at the points $x_i$
collectively fill out a space of higher dimension than each of the individual
germs of $X_i$ individually does; typical
examples are given by the union of infinitely many lines passing
through a fixed point (Example~\ref{ex: infinitely many lines through
a point}), or a curve with a cusp of infinite order (Example~\ref{ex: infinite
order cusp}).
Roughly speaking, $X$ will satisfy~[2](b) when the infinitesimal
germs of the $X_i$ at the points $x_i$ eventually stabilise, or equivalently,
when the projective limit defining the versal deformation ring at $x$ 
eventually stabilises.
(This is not quite true; for example, in the inductive system $X_i$,
one could alternately add additional components, and then contract them to
a point.  It will be true in all the examples we give below for which~[2](b)
is satisfied.)

We now present the details of the preceding claims.

\begin{lemma}
\label{lem:monomorphisms at points of Ind-spaces}
Any finite type point of an Ind-locally finite type
algebraic space  $X = \varinjlim_{i \in I} X_i$
 over the locally Noetherian scheme $S$ admits a representative
$x:\Spec k \to X$ which a monomorphism.  This representative is unique up
to unique isomorphism, the field $k$ is a finite type $\cO_S$-field,
and any other representative $\Spec K \to X$ of the given point
factors through 
the morphism $x$ in a unique fashion.
Furthermore, if $i \in I$ is sufficiently large,
then $x:\Spec k \to X$ factors in a unique manner as a composite
$\Spec k \stackrel{x_i} \longrightarrow X_i \to X,$ 
and the morphism $x_i:\Spec k \to X_i$ is again a monomorphism.
\end{lemma}
\begin{proof}
From the definitions, one sees that $|X| = \varinjlim_{i \in I} |X_i|.$
Thus the given point of $X$ arises from a point of $X_i$ for some $i\in I$.
Let $x_i: \Spec k \to X_i$ be the monomorphic representative of this
point whose existence is given by Lemma~\ref{lem:monomorphisms at points}.
For each $i' \geq i,$ let $x_{i'}: \Spec k_{i'} \to X_{i'}$ be
the monomorphic representative of the image of this point in $X_{i'}$.
Let $k'$ denote the residue field of the image of this point in $S$.
Then we have natural containments $k \supseteq k_{i'} \supseteq k'$,
for each $i' \geq i$.   Thus, since $k$ is finite over $k'$,
we see, replacing $i$ by some sufficiently large $i' \geq i$
if necessary, that we may assume that $k_{i'} = k$ for all $i' \geq i,$
and thus conclude that $x_{i'}$ is simply the composite $x_i: \Spec k \to 
X_i \to X_{i'}$, for all $i' \geq i$.
If we define $x$ to be the composite $x:\Spec k \to X_i \to X$,
then, since $x_{i'}$ is a monomorphism for each $i' \geq i$,
one easily verifies that $x$ is a monomorphism.

The remaining claims of the lemma are proved identically to the analogous
claims of Lemma~\ref{lem:monomorphisms at points}.
\end{proof}

If $x$ is a finite type point of the Ind-locally finite type algebraic 
space  $X = \varinjlim_{i \in I} X_i$ over $S$, and if (by abuse of notation) we also write
$x: \Spec k \to X$ to denote the monomorphic representative of $x$ provided
by the preceding lemma, and (following the lemma) write $x_i: \Spec k \to
X_i$ to denote the induced monomorphisms (for sufficiently large $i$),
then we may consider the complete
local rings $\widehat{\cO}_{X_i,x_i}$ (for sufficiently large
$i$), in the sense of Definition~\ref{df:complete local rings}. If $i' \geq i$  (both sufficiently large),
then we obtain a canonical local morphism of complete
local $\cO_S$-algebras
$\widehat{\cO}_{X_{i'},x_{i'}} \to \widehat{\cO}_{X_i,x_i}.$

\begin{df}
\label{df:Ind complete local rings}In the above context, we write
$$\widehat{\cO}_{X,x} := \varprojlim_{i\in I} \widehat{\cO}_{X_i,x_i},$$
considered as a pro-Artinian ring, 
and we refer to $\widehat{\cO}_{X,x}$ as the complete local ring
to $X$ at $x$. By the following lemma it is canonically defined, independent of the
description of~$X$ as an Ind-locally finite type algebraic space.
\end{df}

\begin{lemma}\label{lem: complete local ring for ind algebraic space}
If $X$ is an Ind-locally finite type algebraic space over the locally
Noetherian scheme $S$, and if $x: \Spec k \to X$ is the
monomorphic representative of a finite type point of $X$, which
{\em (}by abuse of notation{\em )} we also denote by $x$,
then $\widehat{\cO}_{X,x}$ is a versal ring for $X$ at
$x$. Furthermore, the morphism $\Spf\widehat{\cO}_{X,x}\to X$ is a formal
monomorphism, and therefore the ring $\widehat{\cO}_{X,x}$, equipped
with this morphism, is unique up to unique isomorphism.
\end{lemma}
\begin{proof}
Any commutative diagram
$$\xymatrix{\Spec A \ar[r] \ar[d] & \Spec B \ar[d] \\
\Spf \widehat{\mathcal O}_{X,x} \ar[r] & \widehat{X}_x}$$
in which $A$ and $B$ are finite type Artinian local $\cO_S$-algebras,
is induced by a commutative diagram
$$\xymatrix{\Spec A \ar[r] \ar[d] & \Spec B \ar[d] \\
\Spf \widehat{\mathcal O}_{X_i,x_i} \ar[r] & \widehat{X}_{i,x_i}}$$
for some sufficiently large value of $i$.  Since the lower horizontal
arrow of this diagram is versal, we may lift the right-hand vertical arrow
to a morphism $\Spec B \to \Spf \widehat{\cO}_{X_i,x_i}$, 
and hence obtain a lifting of the right-hand vertical arrow of the
original diagram
to a morphism $\Spec B \to \Spf \widehat{\cO}_{X,x}$. 
This establishes the versality.

The property of being a formal monomorphism follows easily from the
fact that (by Proposition~\ref{prop:complete local rings}) the
morphisms $\Spf \widehat{\cO}_{X_i,x_i}\to X_i$ are formal
monomorphisms. That the data of $\widehat{\cO}_{X,x}$ together with   the morphism
$\Spf\widehat{\cO}_{X,x}\to X$ is unique up to unique isomorphism is
then immediate from Lemma~\ref{lem: formal versal monomorphism is
  unique up to unique isomorphism.}.
\end{proof}

\subsubsection{Ind-algebraic stacks as scheme-theoretic images}
If $\cX := \varinjlim_{i \in I} \cX_i$ is an Ind-algebraic stack,
then there is an evident morphism of stacks $\coprod_{i \in I} \cX_i 
\to \cX$, whose source is an algebraic stack (being the disjoint union
of a collection of algebraic stacks), 
and which can be verified to be
representable by algebraic spaces precisely when $\cX$ satisfies~[3].
Thus Ind-stacks give us examples of maps from algebraic stacks to
not-necessarily-algebraic stacks, to which we can try to apply the
machinery of Section~\ref{sec:images} (although the fact that
the morphism $\coprod_{i \in I} \cX_i 
\to \cX$ is not quasi-compact in general
is an impediment to such applications).

As one illustration of this,
we note that Example~\ref{ex: the zipper in a surface} below shows that our
main results will not extend in any direct way to morphisms of finite
type which are not assumed to be proper.

\subsection{Illustrative examples}
\label{subsec:examples}

Here we present various illustrative examples, and explain 
how they relate to the general theory.  Several of them are due originally
to Artin~\cite[\S 5]{MR0262237}.

\bigskip
\begin{example}\label{ex: squiggly example with thickened origin}
\cite[Ex.\ 5.11]{MR0262237}: $X$ is the sheaf of sets obtained by taking the union of the two schemes $\Spec k[x,y][1/x]$ and $\Spec
  k[x,y]/(y)$ in the $(x,y)$-plane. (More precisely, $X$ is the sheaf
obtained as the pushout of these schemes along their common intersection.)
The sheaf $X$ satisfies~[1], [2], [3], and~[4b], but  doesn't satisfy~[4a].
(And hence is not an algebraic space, and so doesn't satisfy~[4].)
 
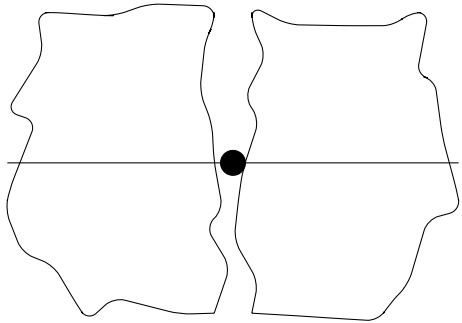
\[ \begin{tikzpicture}[scale=0.5,decoration={random steps,segment
      length=0.5cm,amplitude=2.5mm}]
    \draw (-6,0) -- (6,0); \path (0,0) node [shape=circle,draw=black,fill=black]
    {}; \draw [decorate, rounded corners=1.5mm] (0.5,-4) -- (0.5,4) -- (5,4) --
    (6,-1) -- (4,-4) -- (0.5,-4); \draw [decorate, rounded corners=1.5mm]
    (-0.5,-4) -- (-0.5,4) -- (-5,4) -- (-6,-1) -- (-4,-4) -- (-0.5,-4);
  \end{tikzpicture}\]
\end{example}

We now give a series of examples of various kinds of Ind-algebraic stacks
over a field $k$.
All the stacks we consider will be the inductive limit of stacks
satisfying~[1], and hence will satisfy~[1], as well as~[2](a) and~[4a].
In fact, all the examples other than Example~\ref{ex: limit of classifying
spaces} will be Ind-locally finite algebraic spaces (in fact,
even Ind-locally finite type schemes)  over $k$,
and so will admit versal rings at finite type points by Lemma~\ref{lem: complete local ring for ind algebraic space}.
(As discussed above, though, this doesn't necessarily imply that they
satisfy~[2](b).)

\begin{example}
\label{ex: limit of classifying spaces}
Let $\cX_i$ be a directed system of algebraic stacks,
locally of finite presentation over a locally Noetherian scheme $S$,
for which the transition morphisms are smooth,
and consider the Ind-algebraic stack $\cX := \varinjlim \cX_i$.
Since each~$\cX_i$ satisfies~[1], so does~$\cX$,
and being an Ind-algebraic stack,
it also satisfies~[2](a) and~[4a].
Since the transition morphisms $\cX_i \to \cX_{i'}$
(for $i' \geq i$) are smooth, and since each of the algebraic stacks
$\cX_i$ satisfies~[4], 
we find that $\cX$ also satisfies~[4] (and hence also~[4b]).
If $x:\Spec k \to \cX$ is a finite type point, then $x$ 
factors as $\Spec k \buildrel x_i \over \longrightarrow \cX_i \to \cX$
for some $i$, and one easily verifies (again using the
smoothness of the transition morphisms) that if $\Spf R \to \cX_i$
is a versal ring to $x_i$, then the composite $\Spf R \to \cX_i \to \cX$
is a versal ring to $\cX$.  Thus $\cX$ satisfies~[2](b).  In conclusion,
such an Ind-algebraic stack necessarily satisfies each of
our axioms except possibly~[3].

One example of such an Ind-algebraic stack is obtained by
taking $\{G_i\}$ to be a directed system of smooth
algebraic groups over a field $k$
and setting $\cX_i := [\, \cdot/G_i \,],$
so that
$\cX := \varinjlim\, [\, \cdot/G_i \,]$. (Here ``$\cdot$'' stands for ``the point'', i.e.\ $\Spec k$.)
If we set $G:= \varinjlim_{i \in I} G_i,$ then $\cX$ may be regarded
as the classifying stack $[\, \cdot/G \, ].$
The fibre product $\cdot \times_{\cX} \cdot$ is then naturally identified
with $G$,
and so if the Ind-algebraic group
$G$ is not a scheme (as it typically will not be),
then $\cX$ does not satisfy~[3].
\end{example}

\begin{example}
\label{ex: infinite dimensional affine space}
We take $X$ to be ``infinite dimensional affine space''.
More formally, we write $X := 
\varinjlim \mathbb A^n$,
with the transition maps being the evident closed immersions:
$$\mathbb A^n \cong \mathbb A^n \times\{0\} \hookrightarrow \mathbb A^n
\times \mathbb A^1 \cong \mathbb A^{n+1}.$$
Since the transition maps are closed immersions, the Ind-scheme
$X$ satisfies [3] (and is quasi-separated). 
The complete local rings at points are non-Noetherian (they are power series
rings in countably many variables), and so $X$ does not satisfy~[2](b).
It does satisfy~[4] (and so also~[4b]) vacuously: because the complete
local rings are so large,
one easily verifies that
a morphism from a finite type $k$-scheme to $X$
is never versal at a finite type point.
\end{example}

\begin{example}
  \label{ex: line with infinitely many nodes}We take $X$ 
to be the line with infinitely many
  nodes~\cite[Ex.\ 5.8]{MR0262237} (considered as an Ind-algebraic
space in the evident way).
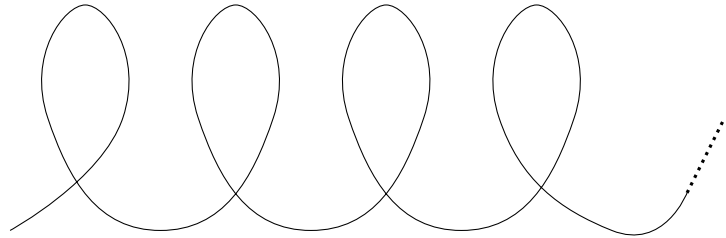
\[\begin{tikzpicture}[node distance = 2cm, auto, scale=0.5, transform shape]
\draw  plot [smooth, tension=1] coordinates { (0,0) (3,3) (2,6) (1,3) (4,0)
  (7,3) (6,6) (5,3) (8,0) (11,3) (10,6) (9,3) (12,0) (15,3) (14,6) (13,3)
  (16,0) (18,1)};
\draw [dotted, very thick] (18,1) -- (19,3);
\end{tikzpicture}\]
The Ind-scheme $X$ satisfies~[2](b), and,
although the transition maps are not closed immersions, it also
satisfies~[3]; however, it is not quasi-separated, since the diagonal morphism
$X \to X \times X$ is not quasi-compact.  Concretely, if $\mathbb A^1 \to X$
is the obvious morphism, namely the one that identifies countably many pairs
of points 
on $\mathbb A^1$ to nodes, then $\mathbb A^1\times_X \mathbb A^1$ is
the union of the diagonal copy of $\mathbb A^1$, and countably many discrete
points (encoding the countably many identifications that are made to create
the nodes of $X$).

The Ind-scheme $X$ doesn't satisfy~[4],
but it does satisfy~[4b], vacuously. 
(One can check that a morphism from a finite type $k$-scheme
to $X$ cannot be smooth at any point.) This example
shows that the quasi-separatedness hypotheses are necessary
in Corollaries~\ref{cor:4a and b give 4}
and~\ref{cor:monomorphism} and 
Theorem~\ref{thm: Artin representability; variant}.
 \footnote{Artin states that this
example does not satisfy~[4b], but seems to be in error on this point.}
\end{example}
\begin{example}
  \label{ex: line with infinitely many cusps}We take $X$ to 
be the line with infinitely many
  cusps.
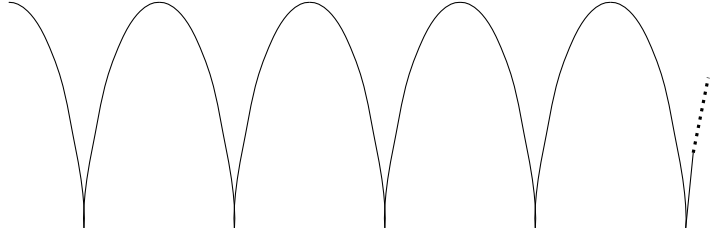
\[\begin{tikzpicture}[node distance = 2cm, auto, scale=0.5, transform shape]
\draw  plot [smooth, tension=1] coordinates { (0,6) (1,5) (1.8, 2) (2,0)};
\draw  plot [smooth, tension=1] coordinates { (2,0) (2.2,2) (3,5) (4,6)};
\draw  plot [smooth, tension=1] coordinates { (4,6) (5,5) (5.8,2) (6,0)};
\draw  plot [smooth, tension=1] coordinates { (6,0) (6.2,2) (7,5) (8,6)};
\draw  plot [smooth, tension=1] coordinates { (8,6) (9,5) (9.8,2) (10,0)};
\draw  plot [smooth, tension=1] coordinates { (10,0) (10.2,2) (11,5) (12,6)};
\draw  plot [smooth, tension=1] coordinates { (12,6) (13,5) (13.8,2) (14,0)};
\draw  plot [smooth, tension=1] coordinates { (14,0) (14.2,2) (15,5) (16,6)};
\draw  plot [smooth, tension=1] coordinates { (16,6) (17,5) (17.8,2) (18,0)};
\draw  plot [smooth, tension=1] coordinates { (18,0) (18.2,2) };
\draw [dotted, very thick] (18.2,2) -- (18.6,4);
\end{tikzpicture}\]
The Ind-scheme $X$ satisfies~[2](b).   As in the previous example,
the transition morphisms are not closed immersions, and in this case
$X$ does not satisfy~[3]; indeed, if $\mathbb A^1 \to X$ is the natural
morphism contracting a countable set of points to the cusps of $X$,
then $\mathbb A^1 \times_X \mathbb A^1$ is the Ind-scheme obtained by
adding non-reduced structure to the diagonal copy of $\mathbb A^1$
at each of the points that is contracted to a cusp.

The Ind-scheme $X$ does not satisfy~[4], but just as in the previous 
example, it does satisfy~[4b] vacuously. \end{example}

\begin{example}\label{ex: line with infinitely many lines crossing it}
We take $X$ to be the line
  with infinitely many lines crossing it \cite[Ex.\ 5.10]{MR0262237}.
  \[\begin{tikzpicture}[node distance = 2cm, auto, scale=1.5, transform shape]
\draw (0,0) -- (7,0);
\draw (1,-1) -- (1,1);
\draw (2,-1) -- (2,1);
\draw (3,-1) -- (3,1);
\draw (4,-1) -- (4,1);
\draw (5,-1) -- (5,1);
\draw [dotted, very thick] (6,0.2) -- (8,0.8);
\end{tikzpicture}\]
This example satisfies~[2](b) and~[3] (and it is quasi-separated, since the
transition maps are closed immersions), but doesn't satisfy~[4b]
(and hence doesn't satisfy~[4]).
\end{example}

\begin{example}
\label{ex: infinitely many embedded points}
As a variation on the preceding example, we consider the Ind-scheme $X$
given by adding infinitely many embedded points to a line.

\[\begin{tikzpicture}[node distance = 2cm, auto, scale=1.5, transform shape]
\draw (0,0) -- (7,0);
\draw[black,fill=black] (1,0) circle (.7ex);
\draw[black,fill=black] (2,0) circle (.7ex);
\draw[black,fill=black] (3,0) circle (.7ex);
\draw[black,fill=black] (4,0) circle (.7ex);
\draw[black,fill=black] (5,0) circle (.7ex);
\draw[black,fill=black] (6,0) circle (.7ex);
\draw [dotted, very thick] (7,0) -- (8.2,0);
\end{tikzpicture}\]

As with the preceding example, this example satisfies [2](b) and [3],
is quasi-separated, and doesn't satisfy [4b] (and hence doesn't
satisfy~[4]). 

One point of interest related to this example
is that 
the natural morphism $\mathbb A^1 \to X$ is a closed immersion
(in particular it is both quasi-compact and proper),
and is surjective; however it is not scheme-theoretically dominant.
Since $X$ is not an algebraic space, this shows the importance
of scheme-theoretic dominance (rather than mere surjectivity) as a
hypothesis in Corollary~\ref{cor:representability-intro}.
\end{example}

\begin{example}
  \label{ex: infinitely many lines through a point}
We take $X$ to be the union of infinitely many
  lines through the origin in the plane~\cite[Ex.\ 5.9]{MR0262237}.
  \[\begin{tikzpicture}[node distance = 2cm, auto, scale=2.5, transform shape]
\draw (-1,0) -- (1,0);
\draw (0,-1) -- (0,1);
\draw (-0.707,-0.707) -- (0.707,0.707);
\draw (-0.5,-0.866) -- (0.5,0.866);
\draw (-0.2, -0.98) -- (0.2,0.98);
\draw [line width=2pt, line cap=round, dash pattern=on 0pt
off2\pgflinewidth,out=0,in=0, looseness=2] (0,0.8) -- (-0.1,0.78) --
(-0.2,0.7) -- (-0.3,0.6);
\draw [line width=2pt, line cap=round, dash pattern=on 0pt
off2\pgflinewidth,out=0,in=0, looseness=2] (0,-0.8) -- (0.1,-0.78) --
(0.2,-0.7) -- (0.3,-0.6);
\end{tikzpicture}\]
This example satisfies~[3], but not~[2](b):
the complete local ring at the origin is equal to $k[[x,y]]$,
and the corresponding morphism $\Spf k[[x,y]] \to X$ is not effective.
It satisfies~[4] (and hence also~[4b]) vacuously: the 
non-effectivity of the complete local ring at the origin shows
that one cannot find a morphism from a finite type $k$-scheme to $X$
which is versal at a point lying over the origin in $X$.
\end{example}

\begin{example}
  \label{ex: infinite order cusp}Let $X_n$ be the plane curve cut out by the equation $y^{2^n} = x^{2^n +1},$
define the morphism $X_n \to X_{n+1}$ via $(x,y) \mapsto (x^2,xy)$,
and let $X := \varinjlim X_n$; so $X$ is a line with a cusp of infinite order.
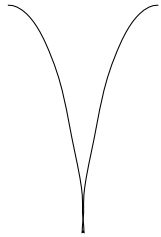
\[\begin{tikzpicture}[node distance = 2cm, auto, scale=0.5, transform shape]
\draw  plot [smooth, tension=1] coordinates { (0,6) (1,5) (1.8, 2)
  (2,0.3) (2,0)};
 \draw  plot [smooth, tension=1] coordinates { (2,0) (2,0.3) (2.2,2) (3,5) (4,6)};
\end{tikzpicture}\]
The Ind-scheme $X$ does not satisfy~[2](b): the complete local ring at the cusp
is equal to $k[[x,y]]$, and the morphism $\Spf k[[x,y]] \to X$ is not
effective.  Also, $X$ does not satisfy~[3]: if $\mathbb A^1  = X_0 \to X$ is
the natural morphism, then $\mathbb A_1 \times_X \mathbb A^1$ is a formal 
scheme, which is an infinite-order
thickening up of the diagonal copy of $\mathbb A^1$
at the origin.
\end{example}

\begin{example}\label{ex: conjugate lines meeting to infinite order}
In~\cite[Ex.\ 5.3]{MR0262237}, 
Artin gives the example of two lines meeting to infinite order
as an Ind-scheme satisfying~[2](b) and [4], but not [3].  
We give a variant of Artin's example here, which illustrates the necessity
of~[3] in Lemma~\ref{lem:independence of morphism for 2b}~(3).

We take $k = \mathbb R$, and we define $X_n := \Spec \mathbb R[x,y]/
(y^2 + x^{2n}),$ with the transition morphism
$X_n \to X_{n+1}$ given by $(x,y) \mapsto (x,xy),$   
and set $X := \varinjlim_{n} X_n.$

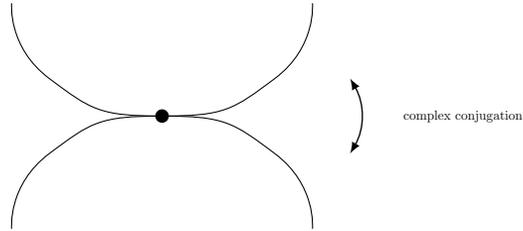
\[\begin{tikzpicture}[node distance = 2cm, auto, scale=0.5, transform shape]
\draw  plot [smooth, tension=1] coordinates { (0,3) (1,1)  (4,0) (7,1)
  (8,3)};
\draw  plot [smooth, tension=1] coordinates { (0,-3) (1,-1)  (4,0) (7,-1) (8,-3)};
 \path (4,0) node [shape=circle,draw=black,fill=black]
    {};
\draw [latex-latex, black, line width=0.5pt, bend left]  (9,1) to
(9,-1);
\node at (12,0) {complex conjugation};
\end{tikzpicture}\]

As with Artin's example, the Ind-scheme $X$ satisfies~[4].  
The complete local ring of $X$ at the origin is equal to $\mathbb R[[x]],$
but the natural morphism $\Spf \mathbb R[[x]] \to X$ is not effective;
thus $X$ does not satisfy~[2](b).
On the other hand, if we consider the map $\Spec \mathbb C \to X$
induced by the origin, then we obtain a versal morphism
$\Spf \mathbb C[[x]] \to X$ which is effective,
although the resulting morphism $\Spec \mathbb C[[x]] \to X$
is not unique; we can map $\Spec \mathbb C[[x]]$ along either
of the branches through the origin.
\end{example}

\begin{example}\label{ex: the zipper} We take $X$ to be the indicated
Ind-scheme (``the zipper'').
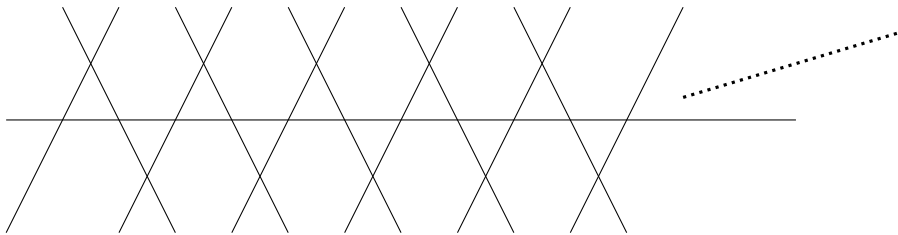
   \[\begin{tikzpicture}[node distance = 2cm, auto, scale=1.5, transform shape]
\draw (0,0) -- (7,0);
\draw (0,-1) -- (1,1);
\draw (1.5,-1) -- (0.5,1);
\draw (1,-1) -- (2,1);
\draw (2.5,-1) -- (1.5,1);
\draw (2,-1) -- (3,1);
\draw (3.5,-1) -- (2.5,1);
\draw (3,-1) -- (4,1);
\draw (4.5,-1) -- (3.5,1);
\draw (4,-1) -- (5,1);
\draw (5.5,-1) -- (4.5,1);
\draw (5,-1) -- (6,1);
\draw [dotted, very thick] (6,0.2) -- (8,0.8);
\end{tikzpicture}\]
It has the same formal properties as the Ind-scheme of
Example~\ref{ex: line with infinitely many lines crossing it},
namely it satisfies~[2](b) and~[3], but not [4] or~[4b].
\end{example}

\begin{example}\label{ex: the zipper in a surface}
We will give an example of an Ind-algebraic {\em surface} which contains
the zipper of the preceding example as a closed sub-Ind-scheme.
   
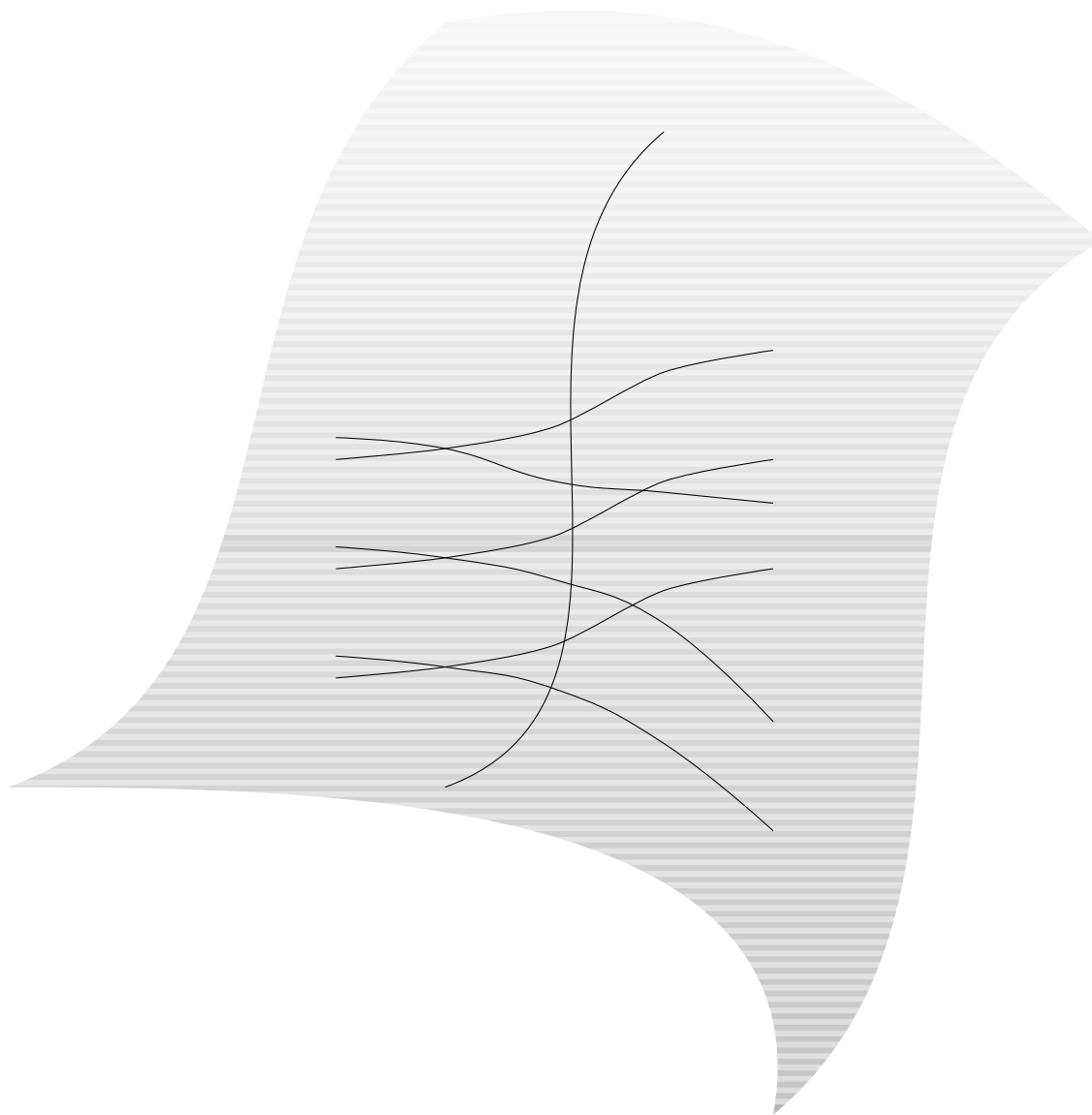
\[\begin{tikzpicture}[node distance = 2cm, auto, scale=1.5, transform shape]

\shade[top color=gray!10,bottom color=gray!90,opacity=.30] 
  (-1,-4) to[out=20,in=220] (3,3)  to[out=10,in=140] (9,1)
 to[out=210,in=40] (6,-7) to[out=80,in=0] (-1,-4);
\path[draw] (3,-4) to[out=20,in=220] (5,2.0);
\draw plot [smooth, tension=0.8] coordinates {(2,-1.8) (3,-1.9) (4,-2.1)  (5,-2.5) (6,-3.4) };
\draw plot [smooth] coordinates {(2,-3) (3,-2.9) (4,-2.7)  (5,-2.2) (6,-2) };

\draw plot [smooth, tension=0.8] coordinates {(2,-0.8) (3,-0.9) (4,-1.2)  (5,-1.3) (6,-1.4) };
\draw plot [smooth] coordinates {(2,-2) (3,-1.9) (4,-1.7)  (5,-1.2) (6,-1) };

\draw plot [smooth, tension=0.8] coordinates {(2,-2.8) (3,-2.9) (4,-3.1)  (5,-3.6) (6,-4.4) };
\draw plot [smooth] coordinates {(2,-1) (3,-0.9) (4,-0.7)  (5,-0.2) (6,0) };
\end{tikzpicture}\]
We begin by setting $\overline{X}_0 := \mathbb A^2$; we also choose
a closed point $P_0 \in \overline{X}_0$.  We let 
$\overline{X}_1$ be the blow-up of $\overline{X}_0$ at $P_0$; it contains 
an exceptional divisor $E_1$, and we choose a closed point $P_1 \in E_1$.
We proceed to construct surfaces $\overline{X}_n$ inductively: 
each $\overline{X}_n$ is a smooth surface obtained by blowing
up $\overline{X}_{n-1}$ at a point $P_{n-1}$.  The surface $\overline{X}_n$
contains an exceptional divisor~$E_n$,
as well as the strict transform of the exceptional divisor $E_{n-1}$
on $\overline{X}_{n-1}$.  We choose a closed point $P_n \in E_n$
which does not lie in the strict transform of $E_{n-1}$,
and then define $\overline{X}_{n+1}$ to be the blow-up of $\overline{X}_n$
at $P_n$.

There are natural open immersions $\overline{X}_n\setminus P_n \subseteq
\overline{X}_{n+1} \setminus P_{n+1}$; the Ind-scheme $\overline{X}$
obtained by taking
the inductive limit of these open immersions is in fact a scheme, and is the
standard example of a locally finite type smooth irreducible surface which
is not of finite type.  Inside $\overline{X}$ we have the union
$E: = \cup_{n \geq 1} E_n\setminus \{P_n\}$, which is an infinite chain of $\mathbb
P^1$'s.   We now form an Ind-scheme $X$ by identifying a countable collection
of points on $E_1 \setminus P_1$ with a point on each $E_n \setminus \{P_n\}$
(for $n \geq 2$); the image of $E$ in $X$ is then a copy of the zipper
of Example~\ref{ex: the zipper}.

The Ind-surface $X$ satisfies~[2](b) and~[3], and is quasi-separated.  
It does not satisfy~[4] or~[4b].

We note that each of the composite morphisms $\overline{X}_n \setminus \{P_n\}
\to \overline{X} \to  X$
is of finite type, so that our formalism of scheme-theoretic images
applies to it.  Each of these morphisms
is in fact scheme-theoretically dominant, in the sense
that its scheme-theoretic image is all of $X$.  (Morally, the Ind-surface
$X$ is irreducible.)   Since $X$ is not an algebraic space, this
example shows that Theorem~\ref{thm:main-intro}
and Corollary~\ref{cor:representability-intro}
don't extend in any direct way to the case of morphisms of finite type
that are not proper.
\end{example}
\section{Moduli of finite height $\varphi$-modules and Galois representations}\label{sec: finite
  flat}
 In this section we will
combine Theorem~\ref{thm:main-intro} with the results of~\cite{KisinModularity,MR2562795}
to construct moduli stacks of finite height and finite flat
representations of the absolute Galois groups of $p$-adic
fields.

We begin by proving some foundational results about $\varphi$-modules of finite
height and \'etale $\varphi$-modules.  We then introduce various moduli stacks of finite height $\varphi$-modules and \'etale
 $\varphi$-modules closely related to those considered
by Pappas and Rapoport in~\cite{MR2562795}, to which we apply the machinery of the earlier parts of the paper.

\subsection{Projective modules over power series and Laurent series
  rings}\label{sec: modules over power Laurent series rings}We begin with a  discussion of some foundational results concerning finitely generated
modules over the power series ring $A[[u]]$ and the Laurent series
ring $A((u))$, where $A$ is an arbitrary (not necessarily Noetherian)
commutative ring. In particular, we recall some deep results of
Drinfeld~\cite{MR2181808} on the \emph{fpqc}-local nature of
the projectivity of such modules.
We are grateful to Drinfeld for sharing with us some of
his unpublished notes on the subject; several arguments in this
section are essentially drawn from these notes.

\subsubsection{Projective and locally free modules}
The basic objects we are interested in are finitely generated projective
modules over the rings $A[[u]]$ and $A((u))$,
where $A$ is some given ground ring.
Of course, a finitely generated projective module over $A[[u]]$ 
(resp.\ over $A((u))$) is Zariski locally free on $A[[u]]$ (resp.\ $A((u))$),
and so we can speak of a finitely generated 
projective $A[[u]]$- or $A((u))$-module being
{\em of rank $d$}: this just means that it is Zariski locally free of rank $d$
on $\Spec A[[u]]$ (resp.\ $\Spec A((u))$).
However, our point of view is that we want to regard the ``base'' of our modules
as being $\Spec A$, and so we will be interested in understanding the
behaviour of $A[[u]]$- or $A((u))$-modules locally on $\Spec A$.

This prompts the following definition, which is used in~\cite{MR2562795}.

\begin{defn}\label{defn: locally free for modules over power Laurent series}
  If $\gM$ is a finitely generated $A[[u]]$-module,
  then we say that that $\gM$ is \emph{fpqc locally free of rank $d$ over $A$}
	  if there exists a faithfully flat $A$-algebra $A'$ such that
	  $\gM\otimes_{A[[u]]} A'[[u]]$ is free of rank $d$ over $A'[[u]]$.
	  Similarly, if $M$ is a 
  finitely generated $A((u))$-module,
  then we say that that $M$ is \emph{fpqc locally free of rank $d$ over $A$}
	  if there exists a faithfully flat $A$-algebra $A'$ such that
	  $\gM\otimes_{A((u))} A'((u))$ is free of rank $d$ over $A'[[u]]$.
	  
	  We use analogous terminology for other topologies besides
	  the {\em fpqc} topology.  E.g.\ if $A'$ can be chosen
	  so that $\Spec A' \to \Spec A$ is an {\em fppf}, {\em \'etale},
	  {\em Nisnevich}, or {\em Zariski} cover (the last notion
	  being understood in the sense that $\Spec A' \to \Spec A$
	  should be surjective, and locally on the source an open immersion),
	  then we say that $\gM$ or $M$ (as the case may be) is
	  {\em fppf, \'etale, Nisnevich, or Zariski locally free
		  of rank $d$ over $A$}.
\end{defn}

\begin{rem}
  \label{rem: not sure if projective implies locally free}One of the main objects of our discussion in this section is to understand, as best
we can, the relationship between the various local freeness properties
just defined, and the property of being finitely generated and projective
(over $A[[u]]$ or $A((u))$).

As we note in Lemma~\ref{lem:power series for fppf extensions} below,
if $A$ is Noetherian and
$A \to B$ is an \emph{fppf} morphism, then the induced morphisms
$A[[u]] \to B[[u]]$ and $A((u)) \to B((u))$ are themselves
faithfully flat.  Thus an $A[[u]]$- or $A((u))$- module
which is \emph{fppf} locally free of finite rank
in the sense of Definition~\ref{defn: 
locally free for modules over power Laurent series} is in fact
finitely generated and projective (since being finitely generated
and projective is a property of a module that can be
checked \emph{fpqc} locally).

In the case that $A$ is not Noetherian, or that the morphism
$A \to B$ is merely faithfully flat, but not of finite presentation, we aren't able to gain
the same level of control over either of
the morphisms $A[[u]] \to B[[u]]$ or $A((u)) \to B((u))$,
and so the precise relationship between projectivity and the notions of
  local freeness
  introduced in Definition~\ref{defn: locally free for modules over
    power Laurent series} is not completely clear to us. 

The most general statement that we were able to prove in
the context of $A[[u]]]$-modules is given
in Proposition~\ref{prop:improving fpqc for A[[u]]} below,
in which we show that
  if an $A[[u]]$-module is finitely generated and projective, then it
  is Zariski locally free (and thus also \emph{fpqc} locally free);
  and that if an $A[[u]]$-module is \emph{fpqc} locally free,
  $u$-torsion free and $u$-adically complete and separated, then it is
  finitely generated and projective.

    In the context of $A((u))$-modules, the relationship between these
    notions is less clear. In Lemma~\ref{lem:locally free over u^n} we show
    that a finitely generated projective $A((u))$-module is Nisnevich
    locally free as an $A((u^n))$-module for all $n$ sufficiently
    large, but in Example~\ref{example: not etale locally free} we
    give an example of a finitely generated
    projective $A((u))$-module of rank one which
    is not \'etale locally free as an $A((u))$-module (although it is
    Zariski locally free as an $A((u^n))$-module for all~$n\ge 2$).
\end{rem}

\begin{rem}
	The second main object of our discussion is to describe the extent to which
	various notions of projectivity/local
	freeness for $A[[u]]$- and $A((u))$-modules are genuinely local
	notions, and (closely related) the descent properties
	of these notions.

	It is obvious that the various local freeness notions presented
	in Definition~\ref{defn: locally free for modules over power Laurent series}
	are local in the relevant topology.  Below, we will recall
	Drinfeld's result that finitely generated projective $A((u))$-modules
	satisfy descent in the {\em fpqc} topology.
\end{rem}

Before turning to these main points of our discussion,
we note the following result.

\begin{lem}
  \label{lem: components of A((u))}Let $A$ be a commutative ring.  Each of
  the 
  natural maps $\Spec A[[u]]\to \Spec A$ and $\Spec A((u))\to\Spec A$ induces
  an isomorphism on the Boolean algebras of simultaneously open and closed subsets
  of its source and target.
 \end{lem}\begin{proof} We need to show that the injections $A\into A[[u]]\into
  A((u))$ induce bijections on the sets of idempotents. Writing
 $\Idem(R)$ for the set of idempotents in the ring $R$, we note
 that the isomorphism
  $A[[u]] \iso \varprojlim A[u]/(u^n)$ induces a bijection
  $\Idem(A[[u]]) \iso \varprojlim \Idem(A[u]/u^n).$
  Since idempotents lift uniquely through nilpotent ideals,
  each of the transition morphisms in this latter projective
  system is a bijection, and hence we find that the morphism
  $\Idem(A[[u]]) \to \Idem(A)$ (induced by the map $A[[u]] \to A[u]/(u) = A$)
  is a bijection.  Its inverse is then given by the map
  $\Idem(A) \to \Idem(A[[u]])$ induced by the inclusion $A \hookrightarrow
  A[[u]]$; in particular, this map is also a bijection.

  Now 
  let $e\in A((u))$ be an idempotent,
  and write
  $e= \sum_{n = -\infty}^{\infty} a_n u^n$.  We claim that $a_i$
  is nilpotent for each $i < 0$.  
  Indeed, if $A$ is reduced, then it is immediate from the equation $e^2 = e$
  that $a_i = 0$ when $i < 0$.  The claim in general
  then follows by considering
  the image of $e$ in $A_{\red}((u))$.
Let $I$ denote the ideal of $A$ generated by the $a_i$ for $i <0 .$
  By the claim we have just proved, together with the fact that $a_i = 0$
  for all but finitely many $i<0$, we find that $I$ is finitely
  generated by nilpotent elements, and so is a nilpotent ideal
  of $A$.  Thus the kernel of the morphism $A((u))\to (A/I)((u))$
  is also nilpotent, and so we obtain a commutative square 
  $$\xymatrix{\Idem(A) \ar[r]\ar[d] & \Idem\bigl( A((u)) \bigr)\ar[d] \\
	  \Idem(A/I) \ar[r] & \Idem\bigl( (A/I)((u))\bigr)}$$
  in which the vertical arrows are bijections and the horizontal
  arrows are injections.   By the definition of $I$, and the
  discussion above, we see that
  the image of $e$ in $(A/I)((u))$
  in fact lies in $(A/I)[[u]],$ and so, by what we have already proved
  in fact lies in $A/I$.  A consideration of the preceding commutative 
  square then shows that $e \in A$,  
  as required.
  \end{proof}

As a consequence of the preceding result,
we have the following reassuring statement,
which shows that in those contexts in which we have multiple ways to
define the locally free rank of an $A[[u]]$ or $A((u))$-module,
the definitions coincide.

\begin{lem}
  \label{lem: about locally free implies rank}If $M$ is a finitely
  generated projective $A[[u]]$- {\em (}resp.\ $A((u))$-{\em )}module,
	  and it is also
  fpqc locally free of rank~$d$, then it is of rank~$d$
  as a projective module over $A[[u]]$ {\em (}resp.\ $A((u))${\em )}.\end{lem}
\begin{proof}
The rank of a finitely generated
projective module over a ring is locally constant, and so we
may find a partition of $\Spec A[[u]]$ (resp.\ $\Spec A((u))$)
into disjoint open subsets on each of which the rank of $M$ is constant.
  Lemma~\ref{lem: components of A((u))}
  shows that this partition of $\Spec A[[u]]$ is induced by a corresponding
  partition of $\Spec A$.     Working separately
  over each of these open and closed subsets of $A$, we may assume
  that $M$ is a projective $A[[u]]$- (resp.\ $A((u))$-)module
  of some fixed rank~$n$.   

  By assumption there is a faithfully flat morphism $A \to B$
  such that $B[[u]]\otimes_{A[[u]]}M$ (resp.\ $B((u))\otimes_{A((u))} M$)
  is free of rank $d$ over $B[[u]]$ (resp.\ $B((u))$). 
  This module is also finitely generated projective of rank $n$,
  and thus we deduce that $n = d,$ as required.
\end{proof}

We now develop those aspects of our discussion that don't require
Drinfeld's theory of Tate modules.

\begin{lemma}
\label{lem:power series for fppf extensions}
If $A$ is Noetherian and $A \to B$ is an fppf morphism,
then each of the induced morphisms $A[[u]]\to B[[u]]$ and $A((u))\to B((u))$
is faithfully flat.
\end{lemma}
\begin{proof}
Note that $B[[u]]$ 
is flat over $B\otimes_A A[[u]]$ (being the $u$-adic completion
of the latter ring, which is finitely presented over the Noetherian
ring $A[[u]]$, and thus is itself Noetherian), which is in turn flat
over $A[[u]]$.   The maximal ideals of $A[[u]]$ are all of the form
$(\mathfrak m, u)$, where $\mathfrak m$ is a maximal ideal of $A$.
Given such a maximal ideal, since $\Spec B \to \Spec A$ is surjective,
we may find a prime ideal $\mathfrak p$ of $B$ which maps to $\mathfrak m$,
and then $(\mathfrak p, u)$ is a prime ideal of $B[[u]]$ which maps
to the maximal ideal $(\mathfrak m, u)$ of $A[[u]]$.  

Thus $A[[u]] \to B[[u]]$ is a flat morphism for which the induced map
$\Spec B[[u]] \to \Spec A[[u]]$ contains all maximal ideals in its image.  
Since flat morphisms satisfy going-down, this morphism is in fact surjective,
and thus $A[[u]] \to B[[u]]$ is faithfully flat.   
The morphism $A((u))\to B((u))$ is obtained from this one via extending
scalars from $A[[u]]$ to $A((u))$, and so is also faithfully flat.
\end{proof}

As already discussed in
Remark~\ref{rem: not sure if projective implies locally free},
for $A[[u]]$-modules, the conditions of being finitely generated
projective and of being \emph{fpqc} locally free of finite rank are closely
related, the precise nature of this relationship being the subject of
the following results. The first part of the
next proposition is closely related to~\cite[Prop.\ 7.4.2]{MR2713872}.

\begin{prop}
  \label{prop: criterion for projectivity over A[[u]]}Let $\gM$ be an $A[[u]]$-module. Then the following conditions are equivalent:
  \begin{enumerate}
  \item $\gM$ is a finitely generated projective $A[[u]]$-module.
  \item $\gM$ is $u$-torsion free and $u$-adically complete and
    separated, and $\gM/u\gM$ is a finitely generated projective $A$-module.
  \end{enumerate}
Moreover if these conditions hold then there is an isomorphism of
$A[[u]]$-modules $(\gM/u\gM)\otimes_A A[[u]]\isoto \gM$, which reduces
to the identity modulo~$u$. In particular if furthermore~$\gM/u\gM$ is
a free $A$-module, then~$\gM$ is a free $A[[u]]$-module.
\end{prop}
\begin{proof}If $\gM$ is projective, then it is a direct summand of a
  finite free $A[[u]]$-module, and is therefore $u$-torsion free and
  $u$-adically complete and separated; and certainly $\gM/u\gM$ is a
  projective $A$-module. For the reverse implication, note that by~\cite[Prop.\
0.7.2.10(ii)]{MR3075000}, we need only show that for each integer
$n\ge 1$, $\gM/u^n\gM$ is a projective $A[u]/u^n$-module.

To show this, note that firstly that since $\gM$ is $u$-torsion free,
for each $m,n\ge 1$ we have a short exact sequence of $A$-modules
\[0\to\gM/u^n\gM\stackrel{u^m}{\to}\gM/u^{m+n}\gM\to\gM/u^m\gM\to 0.\]
By the equivalence of conditions (1) and (4) of~\cite[Thm.\ 22.3]{MR1011461}
 (the local flatness criterion), we see that $\gM/u^n\gM$ is a flat
 $A[u]/u^n$-module for each $n$. It follows from the same 
short exact sequence and induction on $n$ that $\gM/u^n\gM$ is an
$A$-module of finite presentation; so
by~\cite[\href{http://stacks.math.columbia.edu/tag/0561}{Tag
  0561}]{stacks-project}, it is an $A[[u]]/u^n$-module of finite
presentation, and is therefore projective, as required.

Finally, if these conditions hold, then since $\gM/u\gM$ is projective,
we can choose an $A$-linear section to the natural surjection $\gM\to\gM/u\gM$. We
therefore have a morphism of projective $A[[u]]$-modules
$(\gM/u\gM)\otimes_A A[[u]]\to \gM$, which reduces to the identity
modulo~$u$. Write~$\gN:=(\gM/u\gM)\otimes_A A[[u]]$. Then the morphism
$\gN\to\gM$ is surjective by the topological version of Nakayama's
lemma, so we need only prove that it is injective. For this, note that
if $K$ is the kernel of the morphism, then since~$\gM$ is projective,
the short exact sequence \[0\to K\to\gN\to\gM\to 0 \]splits, so that $K$ is a finitely generated
projective $A[[u]]$-module, and 
$K/uK=0$. Since~$K$ if finitely generated projective, it is in
 particular $u$-adically separated, so we have~$K=0$, as required.
\end{proof}

\begin{prop}
\label{prop:improving fpqc for A[[u]]}\leavevmode
\begin{enumerate}
		\item If $\gM$ is a finitely generated projective $A[[u]]$-module, then $\gM$ is Zariski
  locally on $\Spec A$ free of finite rank as an $A[[u]]$-module.
\item If  $\gM$ is $u$-torsion free and is $u$-adically
  complete and separated, and $\gM$ is fpqc locally on $\Spec A$ free of finite rank as an
  $A[[u]]$-module, then $\gM$ is a finitely generated projective $A[[u]]$-module.
\item If $A$ is Noetherian, then $\gM$ is a finitely generated projective $A[[u]]$-module
  if and only if $\gM$ is fppf locally on $\Spec A$ free of finite
  rank as an $A[[u]]$-module.
\end{enumerate}
\end{prop}
\begin{proof}We begin with (1). If $\gM$ is a finitely generated
  projective $\gS_A$-module, then by Proposition~\ref{prop: criterion
    for projectivity over A[[u]]}, $\gM/u\gM$ is a finitely generated
  projective $A$-module. It is therefore Zariski locally free. It
  follows that it is enough to show that if
  $(\gM/u\gM)\otimes_{A}B$ is a free $B$-module of finite
  rank, then $\gM\otimes_{A[[u]]}B[[u]]$ is a free
  $B[[u]]$-module of finite rank. But
  $\gM\otimes_{A[[u]]}B[[u]]$ is a finitely generated projective
  $B[[u]]$-module, so it follows from Proposition~\ref{prop:
    criterion for projectivity over A[[u]]} that
  $\gM\otimes_{A[[u]]}B[[u]]$ is isomorphic to
  $((\gM/u\gM)\otimes_{A}B)\otimes_{B}B[[u]]$, and is
  therefore free of finite rank.

  For (2), if $\gM$ is \emph{fpqc} locally free, then $\gM/u\gM$
  is
  \emph{fpqc} locally free of finite rank, and is in particular
  \emph{fpqc} locally finitely generated
  projective. By~\cite[\href{http://stacks.math.columbia.edu/tag/058S}{Tag
    058S}]{stacks-project}, the property of being finitely generated and projective
is \emph{fpqc} local, so we see that~ $\gM/u\gM$ is finitely generated
  projective. Thus $\gM$ is finitely generated  projective,
  by Proposition~\ref{prop: criterion for projectivity over A[[u]]}.

For (3), one implication is immediate from~(1) (which furthermore allows
us to strengthen \emph{fppf} locally free to Zariski locally free).
The converse follows from Lemma~\ref{lem:power series for fppf extensions},
again using~\cite[\href{http://stacks.math.columbia.edu/tag/058S}{Tag
    058S}]{stacks-project}.
\end{proof}

\subsubsection{Tate modules}\label{subsec: Tate modules}As discussed in Remark~\ref{rem: not sure if projective implies locally free},
for general faithfully flat morphisms $A \to B$ (i.e.\ outside the
context of Lemma~\ref{lem:power series for fppf extensions}), we are
not able to gain much direct control over the induced morphisms
$A[[u]] \to B[[u]]$ or $A((u)) \to B((u))$, and so we are not able to
apply standard descent results in the context of Definition~\ref{defn:
  locally free for modules over power Laurent series}.  However,
in~\cite{MR2181808}, Drinfeld is able to establish descent results for
these morphisms, provided that one restricts attention to modules that
are finitely generated and projective. Drinfeld's basic descent result is stated in the language of
{\em Tate modules}, 
and we begin by recalling the definition of this notion from~\cite[\S3]{MR2181808}, as well as some related definitions.

\begin{df}
\label{def:Tate modules}
Let $A$ be a commutative ring.  (In fact,
\cite{MR2181808} defines Tate modules over not necessarily commutative rings, but we
will only need the commutative case.) An \emph{elementary Tate
  $A$-module} is a topological $A$-module which is isomorphic to
$P\oplus Q^*$, where $P,Q$ are discrete projective $A$-modules,
and $Q^* := \Hom_A(Q,A)$ equipped with its natural projective limit
topology (where we write~$Q^*$ as the projective limit of
$(Q')^*$, where~$Q'$ is a finite direct summand of~$Q$, and give
$(Q')^*$ the discrete topology). A \emph{Tate $A$-module} is a direct summand of an elementary Tate
$A$-module.

A morphism of Tate modules is a continuous morphism of the underlying $A$-modules.
\end{df}

\begin{df}
\label{def:lattice}
A submodule $L$ of a Tate module $M$ is a \emph{lattice} if it is
open, and if furthermore for every open submodule $U\subseteq L$, the
quotient $L/U$ is a finitely generated $A$-module. We say that $L$ is
\emph{coprojective} if $M/L$ is a projective $A$-module (equivalently,
a flat $A$-module; see~\cite[Rem.\ 3.2.3(i)]{MR2181808}). A Tate
module contains a coprojective lattice if and only if it is elementary
\cite[Rem.\ 3.2.3(ii)]{MR2181808}.
\end{df}

The most important example of these definitions
for our purposes is the following.

\begin{example}
We endow $A[[u]]$ with its $u$-adic topology,
and endow $A((u))$ with the unique topology in which $A[[u]]$ (equipped
with its $u$-adic topology) is embedded as an open subgroup.
Equivalently,
we write
$A((u)) = \varinjlim_n \dfrac{1}{u^n} A[[u]]$, 
and endow $A((u))$ with its inductive limit topology, 
where each term in the inductive limit is endowed with its $u$-adic 
topology.
Or, again equivalently, we write $\displaystyle
A((u)) = A[[u]]\oplus \dfrac{1}{u}A[1/u],$
in which the first factor is endowed with its $u$-adic topology,
and the second factor is discrete.

  Both $A[[u]]$ and $A((u))$ are then elementary Tate $A$-modules, and
  $A[[u]]$ is a coprojective lattice in $A((u))$. Furthermore,
  by~\cite[Ex.\ 3.2.2]{MR2181808}, any finitely generated projective
  $A((u))$-module has a natural topology, making it a Tate $A$-module.
  (Indeed, we may write such a module~$M$ as a direct summand of
  $A((u))^n$ for some $n \geq 1$, where this latter module is
  endowed with its product topology.)
 In fact, we have the the
  following theorem~\cite[Thm.\ 3.10]{MR2181808}.
\end{example}

\begin{thm}\label{thm: Drinfeld characterisation of projective A((u)) modules}
There is a natural bijection between finitely generated projective $A((u))$-modules, and pairs $(M,T)$
consisting of a Tate $A$-module $M$ and a topologically nilpotent
automorphism $T:M\to M$, by giving $M$ the $A((u))$-module structure
determined by
$um:=T(m)$. {\em (}Here $T$ is topologically nilpotent if and only if for
each pair of lattices $L,L'\subseteq M$, we have $T^nL\subseteq L'$ for
all sufficiently large~$n$.{\em )}
\end{thm}

\subsubsection{Descent}\label{subsec: descent Drinfeld results}
Drinfeld's fundamental descent result is the following theorem
~\cite[Thm.\ 3.3]{MR2181808}, which shows
that the notion of a Tate $A$-module is
\emph{fpqc}-local on $\Spec A$.

\begin{thm}
\label{thm:fpqc descent for Tate modules}
If $A'$ is a faithfully flat
$A$-algebra, then the functor $M \mapsto A'\cotimes_A M$ induces
an equivalence between the category of Tate $A$-modules
and the category of Tate $A'$-modules with descent data to $A$.
Furthermore, the identification of $\Hom$-spaces given by
this equivalence respects the natural topologies.  
\end{thm}

The statement concerning topologies on $\Hom$-spaces is left implicit
in~\cite{MR2181808}, but is easily checked.\footnote{Any Tate module is a direct summand of one of the form $P \oplus Q^*$, where $P$ and $Q$ are both
{\em free} $A$-modules.  Since the formation of Homs is compatible with finite
direct sums, this reduces us to verifying the claim in the case
of $\Hom(M,N)$ where $M$ and $N$ are $A$-Tate modules which are either
free or dual to free, in which case it is straightforward.}

\begin{remark}
\label{rem:Drinfeld vs. Raynaud--Gruson}
It is a theorem of Raynaud--Gruson~\cite[Ex.\ 3.1.4, Seconde partie]{MR0308104} 
that the property of a module being projective can be checked \emph{fpqc} locally.  This implies that, for any faithfully flat morphism $A \to A'$,
the functor $M \mapsto A'\otimes_A M$ induces an equivalence
between the category of projective $A$-modules
and the category of projective $A'$-modules equipped with descent data to $A$. 
Since projective modules are particular examples of Tate modules
(they are precisely the discrete Tate modules), 
Drinfeld's Theorem~\ref{thm:fpqc descent for Tate modules} incorporates
this descent result of Raynaud--Gruson as a special case.
\end{remark}

The relationship between Tate modules and finitely generated projective
$A((u))$-modules given by
Theorem~\ref{thm: Drinfeld characterisation of projective A((u)) modules}
then implies that
the notion of a finitely generated projective $A((u))$-module is also
\emph{fpqc} local on $\Spec A$; this is~\cite[Thm.\ 3.11]{MR2181808},
which we now recall.   In addition, we prove some slight variants 
of this result that we will need below.

\begin{thm}
  \label{thm: fpqc locality of projective and locally free}
  The following notions are local for the fpqc topology on $\Spec A$.
  \begin{enumerate}
  \item A finitely generated projective $A((u))$-module.
\item A projective $A((u))$-module of rank~$d$.
  \item A finitely generated projective $A((u))$-module which is fpqc
    locally free of rank~$d$.
  \item A finitely generated projective $A[[u]]$-module.
\item A projective $A[[u]]$-module of rank~$d$.
  \item A finitely generated projective $A[[u]]$-module which is fpqc
    locally free of rank~$d$.
  \end{enumerate}

\end{thm}

\begin{rem}
  \label{rem: what fpqc locality means for projective etc}More
  precisely, saying that the notion of a finitely generated projective
  $A((u))$-module is local for the \emph{fpqc} topology on $\Spec A$ means
  the following (and the  meanings of the other statements in Theorem~\ref{thm: fpqc
    locality of projective and locally free} are entirely
  analogous):

If $A'$ is any faithfully flat $A$-algebra, set
  $A'':=A'\otimes_A A'$. Then the category of finitely generated
  projective $A((u))$-modules is canonically equivalent to the
  category of finitely generated projective $A'((u))$-modules $M'$
  which are equipped
  with an isomorphism \[M' \otimes_{A'((u)),a\mapsto 1\otimes a}A''((u))
  \isoto M' \otimes_{A'((u)),a\mapsto a\otimes 1}A''((u)) \] which satisfies the usual cocycle condition.
\end{rem}

\begin{proof}[Proof of Theorem~\ref{thm: fpqc locality of projective
    and locally free}]Since the notion of being \emph{fpqc} locally
  free of rank~$d$ is \emph{fpqc} local by definition, and the rank
  of a finitely generated projective module can be computed \emph{fpqc}
  locally by Lemma~\ref{lem: about locally free implies rank}, it suffices to
  prove statements (1) and (4). The first of these is~\cite[Thm.\
  3.11]{MR2181808}.
As noted above, it follows from Theorem~\ref{thm:fpqc descent for Tate
modules} together with
Theorem~\ref{thm: Drinfeld characterisation of projective A((u)) modules}.
(The fact that the property of an automorphism
being topologically nilpotent satisfies descent follows from
the compatibility with topologies on $\Hom$-spaces stated 
in Theorem~\ref{thm:fpqc descent for Tate modules}.)

For~(4), let $A'$ be a faithfully flat
  $A$-algebra, and let $L'$ be a projective $A'[[u]]$-module equipped with
  descent data. Then $M':=L'\otimes_{A'[[u]]}A'((u))$ is a projective
  $A'((u))$-module equipped with descent data, and $L'$ is a
  coprojective lattice in $M'$ by Lemma~\ref{lem: coprojective versus
    projective} below.
The short exact sequence of Tate $A'$-modules 
$$0 \to L' \to M' \to M'/L' \to 0$$
is split, and admits descent data to $A$.  Thus, by Theorem~\ref{thm:fpqc
descent for Tate modules}, we may descend this to a (split) short exact
sequence of Tate $A$-modules.  
By~(1) (or, perhaps better, by its proof),
the endomorphism $u$ of $M'$  descends to an endomorphism of~$M$,
which equips $M$ with the structure of a finitely generated and
projective  $A((u))$-module.
Since $u$ preserves the submodule $L'$ of $M'$, 
it preserves the descended submodule $L$ of $M$.

Since $M'/L'$ is discrete (or, equivalently, a projective
$A'$-module), the same is true of $M/L$, and thus we
see that $L$ is open in $M$, and coprojective
(\emph{cf.}~Remark~\ref{rem:Drinfeld vs. Raynaud--Gruson}).
In fact, $L$ is a lattice in $M$.  Indeed, since $u$ is a topologically
nilpotent automorphism of $M$, the submodules $u^n L$ ($n \geq 0$) 
form a neighbourhood basis of zero in $L$.  
Since
$$A'\otimes_A (L/u^n L) = (A'\cotimes_A L)/ u^n(A'\cotimes L)
= L'/u^n L'$$ is finitely generated, and since
  the property of being a finitely generated $A$-module is local for
  the \emph{fpqc} topology on $\Spec A$, 
we find that $L/u^n L$ is finitely generated over~$A$. 
To complete the proof, we note
that $L$ is a projective $A[[u]]$-module by another application
  of Lemma~\ref{lem: coprojective versus projective}, as required. \end{proof}
We learned the following lemma from Drinfeld.\begin{lem}
  \label{lem: coprojective versus projective}Let $M$ be a finitely
  generated projective $A((u))$-module, and let $L$ be an
  $A[[u]]$-submodule of $M$. Then the following are equivalent:
  \begin{enumerate}
  \item $L$ is a finitely generated projective $A[[u]]$-module with
    $A((u))L=M$.
  \item $L$ is a coprojective lattice in $M$.
  \end{enumerate}

\end{lem}
\begin{proof}
  If (1) holds, then $L$ is certainly open in $M$, and each $L/u^nL$
  is a finitely generated $A$-module, so $L$ is a lattice in $M$. Since
  $M/L\cong L\otimes_{A[[u]]}\bigl(
  A((u))/A[[u]]\bigr)$, and $A((u))/A[[u]]$ is a free $A$-module,
  $M/L$ is a projective $A$-module, so that $L$ is a coprojective
  lattice, as required.

Conversely, suppose that $L$ is a coprojective lattice. Since $L$ is a
lattice, we certainly have $A((u))L=M$. The short exact sequence of
$A$-modules \[0\to L/uL\to M/L\stackrel{u}{\to} M/L\to 0\] splits
(because $M/L$ is projective), so that $L/uL$ is a direct summand of
the projective $A$-module $M/L$, and is thus itself projective. Since
$L$ is a lattice, $L/uL$ is finitely generated. It follows from
Proposition~\ref{prop: criterion for projectivity over A[[u]]} that
$L$ is a finitely generated projective $A[[u]]$-module, as required.
\end{proof}

\begin{remark}
\label{rem:what we don't know}
As we proved in Proposition~\ref{prop:improving fpqc for A[[u]]},
a finitely generated and projective $A[[u]]$-module is 
\emph{fpqc} (indeed, even Zariski)
locally free of finite rank; thus the only difference
between the situations of parts~(4), (5) and~(6) of Theorem~\ref{thm:
fpqc locality of projective and locally free} is that, in parts~(5) and~$(6)$,
the locally
free rank of the $A[[u]]$-module in question is prescribed.

\end{remark}

\begin{lemma} 
	\label{lem:free over u^n} Let $M$ be a finitely
	generated  $A((u))$-module, which is projective as an
        $A((u^n))$-module, for some $n$ that is invertible in $A$.  Then $M$ is
	projective over $A((u))$.
\end{lemma}
\begin{proof}
	 Identify $A((u))$ with $A((u^n))[X]/(X^n - u^n),$
	 so that we may regard $M$ as a module over $A((u^n))[X]/(X^n - u^n)$
	 which is projective as an $A((u^n))$-module.
	 Now consider the base-change $M' := A((u))\otimes_{A((u^n))} M$;
	 this is a module over $A((u))[X]/(X^n - u^n)$ which is projective
	 as an $A((u))$-module. 
	 Since $n$ and $u$ are both invertible in $A((u))$,
	 the quotient $A((u)) := A((u))[X]/(X-u)$ of $A((u))[X]/(X^n - u^n)$
	 is a direct summand of $A((u))[X]/(X^n - u^n)$.
	 Thus $M \cong M'/(X-u)M'$ is a direct summand of $M'$,
	 and hence is projective as an $A((u))$-module.
\end{proof}

The following result relates the property of an $A((u))$-module
being finitely generated and projective to the property of it
being locally free, in the sense of Definition~\ref{defn: locally free for modules over power Laurent series}.

\begin{lemma}
  \label{lem:locally free over u^n}Let~$M$ be a finitely generated
  projective $A((u))$-module. Then there exists an~$n_0\ge 1$ such that
  for all $n\ge n_0$, $M$ is Nisnevich locally free as an $A((u^n))$-module.
\end{lemma}
\begin{proof}  By~ \cite[Thm.\
   3.4]{MR2181808}, we may make a Nisnevich localisation so that
   $M$ is an elementary Tate $A$-module, i.e.\ contains a coprojective lattice $L$.

 Since multiplication
 by $u$ is topologically nilpotent, there is an integer $n_0\ge 0$ such
 that $u^nL\subseteq L$ for all $n\ge n_0$; thus, for each such value of $n$,
 we see that $L$ is naturally an $A[[u^n]]$-module,  
 and that the natural morphism $A((u^n))\otimes_{A[[u^n]]} L \to M$
is an isomorphism. By Lemma~\ref{lem: coprojective versus projective},
$L$ is a finitely generated projective $A[[u^n]]$-module. By
Proposition~\ref{prop:improving fpqc for A[[u]]}~(1), after making a
further Zariski localisation, we may suppose that $L$ is free of
finite rank as an $A[[u^n]]$-module, so that~$M$  is free of finite rank
as an~$A((u^n))$-module, as required. 
  \end{proof}

  The following example shows that, in the context of the preceding 
  lemma, we can't necessarily take $n_0 = 1$.

\begin{example}\label{example: not etale locally free}
Let $A = k[x,y]/(y^2-x^3),$
	for some field $k$, and let $I \subseteq A((u))$ denote the ideal
	generated by $(u^2 - x, u^3 - y)$. 
One can check that $I$ is freely generated over $A((u^2))$ 
        by $u^2 -x$ and $u^3 - y$, and if the characteristic of $k$
is different from~$2$, $I$ is projective over $A((u))$
	by Lemma~\ref{lem:free over u^n}.

	Alternatively, and more conceptually, one can deduce
	this projectivity (with no assumption on the characteristic of~$k$)
by noting that $(u^2,u^3)$ is a smooth point
	(over the complete non-archimedean field $k((u))$)
	of the rigid analytic curve $y^2 = x^3$ lying in the closed
	polydisk $|x|,|y| \leq 1$  over $k((u))$,
	and that $I$ is ideal sheaf of $(u^2,u^3)$ in the Tate algebra
	$A((u))$ of the curve.  

	One can check that the $A[[u^2]]$-submodule $L$ of $I$
	generated by $u^2 -x$ and $u^3 - y$ is an $A[[u^2,u^3]]$-submodule
	of $I$, and hence is closed under multiplication by $u^n$
	for any $n \geq 2$. Since~$I$ is freely generated
        over~$A((u^2))$ by $u^2 -x$ and $u^3 - y$, we see that $I/L$
        is a free~$A$-module, so that $L$ is a coprojective lattice in~$I$. The proof of Lemma~\ref{lem:locally free
		over u^n} shows that $L$ is then Zariski locally
	free over $A[[u^n]]$ for any $n \geq 2$, and thus that $I$    
	is Zariski locally free over $A((u^n))$, for any $n \geq 2$.

	We claim that $I$ is {\em not} an \'etale locally free~$A((u))$-module.
To see this, it suffices to show that if $R$ denotes
	the strict Henselisation of $A$ at the maximal ideal $(x,y)$,
	then $I\otimes_{A((u))} R((u))$ is not free over $R((u))$,
	i.e.\ that the ideal $(u^2-x,u^3 - y)$ is not a principal
	ideal in $R((u))$.   We prove this in the following lemma.

\end{example}

	\begin{lemma}\label{lem: not etale locally principal} If $R$ denotes the strict Henselisation
		of $k[x,y]/(y^2 - x^3)$ at the maximal ideal $(x,y)$, 
		then the ideal $(u^2-x,u^3-y)$ of $R((u))$ is not
		principal.
	\end{lemma}
	\begin{proof}
		Let $f \in R[[u]]$ be non-zero, with non-zero constant
		term $f_0 \in R$.  (Any non-zero element of $R((u))$ may be be
		multiplied by some power of $u$ so as to satisfy
		this condition, and hence any principal ideal of $R((u))$
		has a principal generator satisfying this condition.)
		We claim, then, that $fR((u)) \cap R[[u]] = fR[[u]].$

		Indeed, if $g \in fR((u)) \cap R[[u]],$ then 
		$g = fh$ for some $h \in R((u))$.  Choose $N$ minimally
		so that $h u^N \in R[[u]]$, and suppose that $N >0$.
		If we reduce the equation
		$gu^N = f hu^N $
		modulo $u$,
		we find that
		$$0 = f_0 \times \text{ the constant term of } h u^N $$
		(an equation in $R$),
		and so (since $R$ is a domain, as the cuspidal cubic
		$y^2 = x^3$ is geometrically unibranch at its singular
		point $(0,0)$) we find that the constant term of $h u^N $
		is zero.  This contradicts the minimality of $N$,
		and thus shows that $N = 0$, so that in fact
		$h \in R[[u]]$
		and $g \in f R[[u]]$, as claimed.
		Thus,
                if $fR((u))$ is any principal ideal of $R((u))$, 
		with $f$ chosen as above, then
		$$R[[u]]/(fR((u))\cap R[[u]] + uR[[u]]) =  R[[u]]/(f,u)R[[u]] 
		= R/f_0 R.$$ Since~$R$ (which is a one-dimensional
                local ring) is not regular, the quotient
                $R/f_0 R$ is necessarily of dimension $> 1$ over $k$.

	        On the other hand,
		we find that $$R[[u]]/\bigl( (u^2 - x, u^3 - y)\cap R[[u]]
		+ uR[[u]]) = R/(x,y) = k.$$ 
		Taking into account the result of the preceding paragraph,
		this shows that indeed~$(u^2 - x, u^3 - y)$
		is not a principal ideal in $R((u))$.
		\end{proof}

\subsection{Modules of finite height and \'etale $\varphi$-modules}\label{sec: Kisin modules and phi modules}
In this section we discuss finite height $\varphi$-modules and \'etale $\varphi$-modules. With an eye to future applications (for example, the case of
Lubin--Tate $(\varphi,\Gamma)$-modules), we work in a more general
context than that usually considered.

\subsubsection{Definitions}
Fix a finite extension $k/\Fp$, and write $\gS:=W(k)[[u]]$. Let~$q$ be
some power of~$p$, and let~$\varphi$ be a ring endomorphism of~$\gS$
which is congruent to the~$q$-power Frobenius endomorphism
modulo~$p$.

\begin{lem}\label{lem: phi is what you think on Witt vectors}$\varphi$ induces the usual $q$-power Frobenius
on~$W(k)$.
\end{lem}
\begin{proof}
  It is enough to note that~$W(k)$ is generated as a~$\Zp$-algebra by
  a primitive~$(\#k-1)$st root of unity, and that the $(\#k-1)$st roots of
  unity are distinct modulo~$p$.
\end{proof}
\begin{lem}
  \label{lem: phi(u) divisible by u W(k) version}
  For each $M,a\ge 1$, $\varphi(u^{M+a-1})\in (u^{Mq},p^a)$, and $u^{(M+a-1)q}\in(\varphi(u^M),p^a)$. In particular, $\varphi$ is
  continuous with respect to the $(p,u)$-adic topology. \end{lem}
\begin{proof}
  Write $\varphi(u)=u^{q}+p Y$. Then
  $\varphi(u^{M+a-1})=(u^q+pY)^{M+a-1}$, and
  $u^{(M+a-1)q}=(\varphi(u)-pY)^{M+a-1}$, and the result follows from
  the binomial theorem.
\end{proof}

We fix a finite extension~$E/\Qp$ with ring of integers~$\cO$ and
uniformiser~$\varpi$. 
If $A$ is an~$\cO/\varpi^a$-algebra for some $a\ge 1$, we write
$\gS_A:=(W(k)\otimes_{\Zp}A)[[u]]$; we equip~$\gS_A$ with its
$u$-adic topology. Let $\OEA$ equal
$\gS_A[1/u]$.

\begin{lem}
  \label{lem: extending phi to A}$\varphi$ admits a unique continuous $A$-linear
  extension to~$\gS_A$, which in turn admits a unique continuous $A$-linear extension to~$\OEA$.
\end{lem}
\begin{proof}
  This follows immediately from Lemma~\ref{lem: phi(u) divisible by u
    W(k) version}.
\end{proof}

\begin{lemma}\label{lem: phi faithfully flat} $\varphi$ is faithfully
  flat and finite on~$\gS_A$ and~$\OEA$.
\end{lemma}\begin{proof} It is enough to show this for~$\gS_A$, as it then follows
  for the localisation~$\OEA$. It then suffices to show that~$\gS_A$ is a
  finite free $\varphi(\gS_A)$-module of rank~$q$, with a basis given
  by~$u^i, 0\le i\le q-1$.

  If $A$ is an $\Fp$-algebra this is clear, and it follows from
  Nakayama's lemma that in general $\gS_A$ is generated as an
  $\varphi(\gS_A)$-module by ~$u^i, 0\le i\le q-1$. It remains to
  check that the~$u^i$ are linearly independent. To this end, suppose
  that we have a relation $\sum_{i=0}^{q-1}u^i\varphi(a_i)=0$, where the
  $a_i$ are all contained in~$ p^n\gS_A$ for some~$n\ge 0$. We will
  show that in fact the~$a_i$ are all contained in~$ p^{n+1}\gS_A$;
  since ~$p^n\gS_A=0$ for~$n$ sufficiently large, it follows that
  $a_i=0$ for all~$i$, as required.

  Write $\varphi_0:\gS_A\to\gS_A$ for the particular lift of the $q$-power
  Frobenius on $\gS_{A/pA}$ defined via $\varphi_0(u)=u^q$. Then for any $x\in p^n\gS_A$, we
  have $\varphi(x)-\varphi_0(x)\in p^{n+1}\gS_A$, so our assumptions
  imply that $\sum_{i=0}^{q-1}u^i\varphi_0(a_i)\in p^{n+1}\gS_A$.  Writing the~$a_i$
  out as power series in~$u$, and equating coefficients, we see
  immediately that $a_i\in p^{n+1}\gS_A$, as required.
\end{proof} 

We fix a polynomial $F \in (W(k)\otimes_{\Zp}\cO)[u]$ that is congruent to a positive power of $u$ modulo $\varpi$. The following elementary
lemma will be useful below.

\begin{lem}\label{lem: F versus u mod p^a}For all integers $a,h\ge 1$ there is an
  integer~$n(a,h)\ge 1$ depending only on~$a$, $h$ and~$F$ such that if $A$
  is an $\cO/\varpi^a$-algebra,
  then $u^{n(a,h)}$ is divisible by~$F^h$ in~$\gS_A$, and $F^{n(a,h)}$ is divisible by~$u^{h}$.
\end{lem}
\begin{proof}Write $F=u^n-\varpi X$ for some~$n\ge 1$ and
  $X\in\gS_A$. By the binomial theorem, $F^{a+h-1}$ is divisible by
  $u^{nh}$ in~$\gS_A$, and thus by~$u^h$; similarly, writing $u^n=F+\varpi X$,
  $u^{n(a+h-1)}$ is divisible by~$F^h$. Putting these together, we see
  that we can take~$n(a,h):=n(a+h-1)$.
\end{proof}
We will also use the following result.

\begin{lem}
  \label{lem: phi(u) divisible by u} If $A$ is a $\cO/\varpi^a$-algebra,
  then for each $M\ge 1$, $u^{(M+a-1)q}$  is divisible
  by~$\varphi(u^M)$, and $\varphi(u^{(M+a-1)})$ is divisible by~$u^{Mq}$.
\end{lem}
\begin{proof}This is immediate from Lemma~\ref{lem: phi(u) divisible by u W(k) version}.\end{proof}

If~ $\gM$ is an $\gS_A$-module (resp.\ $M$ is an $\OEA$-module) then we
write~$\varphi^*\gM$ for~$\gM\otimes_{\gS_A,\varphi}\gS_A$ (resp.\ $\varphi^*M$
for~$M\otimes_{\OEA,\varphi}\OEA$). Since~$\varphi$ is faithfully
flat, the
functors $\gM\mapsto\varphi^*\gM$, $M\mapsto\varphi^*M$ are exact. 
\begin{cor}
  \label{cor: exponent of phi}If $A$ is a $\cO/\varpi^a$-algebra,  $\gM$ is an $\gS_A$-module, and
  $u^M\gM\ne 0$, then $u^{(M-a+1)q}\varphi^*\gM\ne 0$.
\end{cor}
\begin{proof}If $u^{(M-a+1)q}\varphi^*\gM=0$, then by Lemma~\ref{lem:
    phi(u) divisible by u} we have
  $\varphi^*(u^M\gM)=\varphi(u^M)\varphi^*\gM=0$. Since~$\varphi$ is faithfully
  flat, this implies that $u^M\gM=0$, a contradiction.  
\end{proof}

The following lemma is a straightforward generalisation (with a very
similar proof) of~\cite[Prop.\ 2.2]{MR2562795} to our setting. Let~$R$
be an $\cO/\varpi^a$-algebra, and let~$u\in R$ be a nonzerodivisor, such
that $R$ is $u$-adically complete and separated. (For example, we
could take~$R=\gS_A$ for some $\cO/\varpi^a$-algebra~$A$.) For  $n\ge
0$, write \[U_n=1+u^nM_d(R),\] \[V_n=\{A\in\GL_d(R[1/u])|A,A^{-1}\in u^{-n}M_d(R)\}.\]
\begin{lem}
    \label{lem: phi conjugacy approximation} Suppose that $n>(2m+(a-1)q)/(q-1)$.
\begin{enumerate}
\item For each $g\in U_n$, $A\in V_m$, there is a unique $h\in U_n$
  such that $g^{-1}A\varphi(g)=h^{-1}A$.
\item For each $h\in U_n$, $A\in V_m$ there is a unique $g\in U_n$
  such that $g^{-1}A\varphi(g)=h^{-1}A$.
\end{enumerate}

  \end{lem}
  \begin{proof} We follow the proof of~\cite[Prop.\
    2.2]{MR2562795}. For the first part, note that we can solve for $h^{-1}$,
    namely $h^{-1} = g^{-1}A\varphi(g) A^{-1},$
    so that the uniqueness of~$h$ is clear,
    and what we must show is that $h \in U_n$,
    or equivalently, that $h^{-1} \in U_n$.
    We can write $g^{-1}=I+u^nX$
    with $X\in M_d(R)$, and by Lemma~\ref{lem: phi(u) divisible by
      u} we can write $\varphi(g)=I+u^{(n-a+1)q}Y$ with $Y\in
    M_d(R)$. Then
    $g^{-1}A\varphi(g)A^{-1}=(I+u^nX)(I+u^{(n-a+1)q}AYA^{-1})$. Since~$A\in
    V_m$ we have $AYA^{-1}\in u^{-2m}M_d(R)$, so that
    $g^{-1}A\varphi(g)A^{-1}\in U_n$, as required.

For the second part we begin by showing uniqueness of~$g$, for which
it is enough to show that if $g^{-1}A\varphi(g)=A$ then~$g=1$. Write
$g=I+X$, so that the preceding relation between $A$ and $g$
may be rewritten as $X=A\varphi(X)A^{-1}$.
It is enough to
check that we have~$X\in u^sM_d(R)$ for all $s\ge n$. We prove this by
induction on~$s$, the case~$s=n$ being by hypothesis. If~$X\in
u^sM_d(R)$, then as above we have~$A\varphi(X)A^{-1}\in
u^{(s-a+1)q-2m}M_d(R)$, and since~$s\ge n$ we have $(s-a+1)q-2m>s$, as
required.

Finally we must show existence of~$g$. For this, let $A'=h^{-1}A$, and
set $A_0=A, h_0=h$. We inductively define sequences $(h_i)$, $(A_i)$
by setting $A_i=h_{i-1}^{-1}A_{i-1}\varphi(h_{i-1})$,
$h_i=(A')^{-1}A_i$. These equalities imply that
$h_i=A'\varphi(h_{i-1})(A')^{-1}$, so since $A'\in V_m$ and $h_0\in
U_n$, an easy induction as above shows that~$h_i\in U_{n+i}$ for
all~$i$. If we now set $g_i=h_0h_1\cdots h_i$, then $g_i$ tends to
some limit $g\in U_n$. Since for all~$i$ we have
$g_i^{-1}A\varphi(g_{i-1})=h^{-1}A$, in the limit we have
$g^{-1}A\varphi(g)=h^{-1}A$, as required.
\end{proof}

\begin{defn}\label{defn: phi module of finite E-height}Let $h$ be a non-negative
integer, and let~$A$ be an $\cO/\varpi^a$-algebra. A \emph{$\varphi$-module of
  height~$F$ with $A$-coefficients} is  a pair $(\gM,\varphi_M)$
consisting of a finitely generated $u$-torsion free $\gS_A$-module~$\gM$, and a
$\varphi$-semilinear map $\varphi_\gM:\gM\to~\gM$, with the further
properties that if we
write \[\Phi_{\gM}:=\varphi_{\gM}\otimes 1:\varphi^*\gM\to\gM,\] then $\Phi_{\gM}$ is injective, and the cokernel of $\Phi_{\gM}$ is killed
by $F$. A \emph{$\varphi$-module of finite height with
  $A$-coefficients}, or a \emph{finite height $\varphi$-module with
  $A$-coefficients}, is a $\varphi$-module with $A$-coefficients which
is of height~$F$ for some~$F$.

A
  morphism of finite height $\varphi$-modules is a morphism of the underlying
  $\gS_A$-modules which commutes with the morphisms~$\Phi_\gM$.

  We say that a finite height $\varphi$-module  is \emph{projective  of
    rank~$d$} if it is a finitely generated projective~$\gS_A$-module
  of constant rank~$d$.
\end{defn}

\begin{remark}\label{rem:Kisin module motivations}
We will primarily be interested in finite height $\varphi$-modules (resp.\ \'etale $\varphi$-modules) that are furthermore projective over $\gS_A$ (resp.\ over $\gS_A[1/u])$,
to which
the base-change and descent results of Drinfeld \cite{MR2181808}
recalled in Subsection~\ref{subsec: descent Drinfeld results}
are applicable.
However, we sometimes need to make
constructions that take us outside the category of projective
modules, and in particular we need to consider finite height $\varphi$-modules which are
not projective, but become projective after inverting~$u$.\end{remark}

\begin{defn}\label{defn: etale phi module} Let~$A$ be a $\cO/\varpi^a$-algebra.
  An \emph{\'etale $\varphi$-module with $A$-coefficients} is a pair
  $(M,\varphi_M)$ consisting of a finitely generated $\OEA$-module~$M$,
  and a $\varphi$-semilinear map $\varphi_M:M\to M$ which induces an
  isomorphism of $\OEA$-modules $\Phi_M:=\varphi_M\otimes 1:\varphi^*M\to M$.

  A morphism of \'etale $\varphi$-modules is a morphism
  of the underlying $\OEA$-modules which commutes with the
  morphisms~$\Phi_M$.

We say that~$M$ is projective (resp.\ free) of rank~$d$ if it is a
finitely generated projective (resp.\ free) $\OEA$-module of constant
rank~$d$. If~$\tau$ is any topology on the category of
$\cO/\varpi^a$-modules lying between the Zariski topology and the
\emph{fpqc} topology, then we say that $M$ is $\tau$-locally free of
rank~$d$ if it is projective of rank~$d$, and if $\tau$-locally on
$\Spec A$, it is free of rank~$d$.

\end{defn}

\begin{remark}
\label{rem: inverting u}If $(\gM,\varphi)$ is a $\varphi$-module of height~$F$ with
$A$-coefficients, then by Lemma~\ref{lem: F versus u mod p^a}, $(\gM[1/u],\varphi)$ is an \'etale $\varphi$-module with
$A$-coefficients.
\end{remark}

We will sometimes prove results about projective \'etale
$\varphi$-modules by reducing to the free case, using the following
lemma.

\begin{lemma}\label{lem: projective phi module is summand of
    free}If~$M$ is a  projective \'etale $\varphi$-module
  with $A$-coefficients, then~$M$ is a direct summand of a free \'etale $\varphi$-module
  with $A$-coefficients. 
  \end{lemma}
\begin{proof}
  	Since $M$ is projective, we may find another
	finitely generated
	projective $\OEA$-module $P$ such $M \oplus P \iso F,$ for some 
	finite rank free module $F$.  Then $$\varphi^*M \oplus \varphi^*P 
	\iso \varphi^*F,$$ and since $M$ is an \'etale $\varphi$-module,
	we have $\varphi^*M \cong M,$ while since $F$ is free, we have
	$\varphi^*F \cong F$.
	Thus $$M\oplus P \cong  M \oplus \varphi^*P,$$
	and so taking the direct sum with another copy of $P,$ we find that
	$$F \oplus P \cong F \oplus \varphi^*P \cong \varphi^*(F \oplus P).$$
	In other words, the finitely generated 
	projective module $Q:= F\oplus P$ admits the structure of an
	\'etale $\varphi$-module, and $$Q \oplus M\cong F \oplus F$$
	is free of finite rank.
\end{proof}

\begin{lem}
  \label{lem: existence of a lattice in a projective}
  Let $A$ be a Noetherian
   $\cO/\varpi^a$-algebra, and let~$M$ be an
  \'etale $\varphi$-module with $A$-coefficients. Then
  for some~$F$, there is a $\varphi$-module~$\gM$ of height~$F$ with
  $A$-coefficients such that $\gM[1/u]=M$. If~$M$ is
  furthermore a free $\OEA$-module, then we may choose~$\gM$ to be a
  free $\gS_A$-module.
\end{lem}
\begin{rem}
  \label{rem: lattice need not be projective}Note that in the
  case when~$M$ is projective but not necessarily free, we do not claim
  that the $\varphi$-module~$\gM$ in Lemma~\ref{lem: existence of a
    lattice in a projective} can be chosen to be projective.
\end{rem}
\begin{proof}[Proof of Lemma~\ref{lem: existence of a lattice in a projective}]
  By definition~$M$ is finitely generated as an $\OEA$-module, so we may
choose a generating set, and let~$\gM$ be the $\gS_A$-span of this
generating set; if~$M$ is free then we may and do also choose~$\gM$ to be free. It follows easily from Lemma~\ref{lem: phi(u)
  divisible by u} that if we scale $\gM$ by a large enough power of~$u$, we
may assume that~$\gM$ is $\varphi$-stable. Since~$M$ is
$u$-torsion free, so is~$\gM$, and since~$\Phi_M$ is an isomorphism,
$\Phi_{\gM}$ is injective and 
the cokernel of~$\Phi_\gM$ is killed by some power of~$u$. Thus
$\gM$ is a finite height $\varphi$-module, as required.
\end{proof}

\begin{rem}
  \label{rem: Galois context}If~$q=p$, then for certain choices of
  $\varphi$ and~$F$, the theories of finite height $\varphi$-modules
  and \'etale $\varphi$-modules with Artinian coefficients admit
  interpretations in terms of Galois representations; more precisely,
  in terms of representations of the absolute Galois groups of certain
  perfectoid fields. We refer to~\cite[\S 2]{emertongeepicture} for a more
  thorough discussion of this; in Section~\ref{subsubsec: Galois reps}
  below we explain the connection in a particular case, that of
  \emph{Breuil--Kisin modules}.
\end{rem}

 \subsection{Lifting rings
 }\label{sec: effectivity of projectivity}We now prove the existence of universal
 lifting rings for \'etale $\varphi$-modules, as well as variants for
 finite height $\varphi$-modules. These will be used below to show that our moduli
 stacks admit versal rings,  and to verify the effectivity hypothesis in our application of
 Theorem~\ref{thm:main-intro}.

We remark that if we were in one of the settings mentioned
in Remark~\ref{rem: Galois context}, then we could use the equivalence
of categories between \'etale $\varphi$-modules and Galois representations
with Artinian coefficients to study the versal rings in which we are interested 
using standard techniques from the formal deformation theory of Galois
representations; in particular, in the case of Breuil--Kisin modules,
we could use the results of~\cite{MR2782840}. Even in the general framework that we have adopted, it seems plausible
that we could use the
Artin--Schreier theory constructions that underlie that equivalence 
to replace the formal deformation theory of \'etale $\varphi$-modules
by the formal deformation theory of some more finitistic objects. 
However, we have found it more direct, and interesting in its own right,
to argue with the formal deformation theory of \'etale $\varphi$-modules.
We caution the reader that this leads us, in what follows,
to consider some rather large pro-Artinian rings!

\subsubsection{Lifting
   rings}The main results of this subsection are Proposition~\ref{prop: existence of infinite height deformation ring} and
Theorem~\ref{thm:
   existence of height F deformation ring}.

We
 begin by studying the formal deformation theory of \'etale
 $\varphi$-modules. Let~$E/\Qp$ be a finite extension with ring of
 integers~$\cO$, uniformiser~$\varpi$ and residue field~$\F$, and let $M$ be an \'etale $\varphi$-module
 with $\F$-coefficients which is free of rank~$d$. Fix an integer~$a\ge 1$, and (following the notation of~\S~\ref{subsec:axiom 2}) write~$\cC_{\cO/\varpi^a}$ for
 the category of Artinian local $\cO/\varpi^a$-algebras for which the
 structure map induces an isomorphism on residue fields. 

 Fix a choice of (ordered) $\cO_{\cE,\F}$-basis of~$M$, or equivalently,
an identification of $\cO_{\cE,\F}$-modules
\numequation
\label{eqn:fixed iso}
M \iso \cO_{\cE,\F}^d.
\end{equation}

\begin{df}
\label{def:etale liftings}
 A \emph{lifting} of~$M$
 to an object~$R$ of~$\cC_{\cO/\varpi^a}$ is a triple consisting of an \'etale
 $\varphi$-module~$M_R$ which is free of rank~$d$, a choice of (ordered)
 $\cO_{\cE,R}$-basis of~$M_R$, and an isomorphism $M_R\otimes_R
 \F\cong M$ of \'etale $\varphi$-modules which takes the chosen basis
 of~$M_R$ to the fixed basis of~$M$.  Equivalently, a lifting $M$ 
 consists of an \'etale 
 $\varphi$-module~$M_R$ endowed with an isomorphism of $\cO_{\cE,R}$-modules
 \numequation
 \label{eqn:lifted iso}
 M_R \iso \cO_{\cE,\F}^d
 \end{equation}
 such that
 the \'etale $\varphi$-module structure on $\cO_{\cE,\F}^d$
 which is obtained by reducing~\eqref{eqn:lifted iso}
 modulo $\mathfrak m_R$ coincides with 
 the \'etale $\varphi$-module structure on $\cO_{\cE,\F}^d$
 induced by the isomorphism~\eqref{eqn:fixed iso}.

 By regarding the isomorphisms~\eqref{eqn:fixed iso}
 and~\eqref{eqn:lifted iso} as identifications, we see that 
 the liftings of $M$ admit yet another equivalent description:
 namely, the identification~\eqref{eqn:fixed iso} allows us
 to regard $M$ as simply a choice of matrix $\Phi \in \GL_d(\cO_{\cE,\F})$,
 which describes the \'etale $\varphi$-module structure
 $$\cO_{\cE,\F}^d = \varphi^*\cO_{\cE,\F}^d = \varphi^*M \to M = \cO_{\cE,\F}^d,$$
 and the identification~\eqref{eqn:lifted iso} allows us to regard the lifting $M_R$
 as a matrix $\Phi_R \in \GL_d(\cO_{\cE,R})$ lifting $\Phi$; this is the matrix
 describing the \'etale $\varphi$-module structure
 $$\cO_{\cE,R}^d = \varphi^*\cO_{\cE,R}^d = \varphi^*M_R \to M_R = \cO_{\cE,R}^d.$$

We denote by $D^{\square}:\cC_{\cO/\varpi^a}\to\Sets$ 
 the functor taking~$R$ to the set of isomorphism classes of liftings
 of~$M$ to~$R$.  (The functorial structure on $D^{\square}$ is the
 obvious one induced by extension of scalars.)
\end{df}

\begin{remark}
\label{rem:action}
There is a natural action of the subgroup
$1 + \mathfrak m_R M_d(\cO_{\cE,R})$
of $\GL_d(\cO_{\cE,R})$
on $D^{\square}(R)$, given by change of basis.   In terms of the description
of liftings via a choice of isomorphism~\eqref{eqn:lifted iso}, 
this is given by composing the isomorphism~\eqref{eqn:lifted iso} with
the automorphism of $\cO_{\cE,R}^d$ induced by a matrix 
$g \in 1 + \mathfrak m_R M_d(\cO_{\cE,R})$.
In terms of the description of liftings in terms of a matrix~$\Phi_R$,
this is given by the twisted conjugation action $\Phi_R \mapsto g \Phi_R \varphi(g)^{-1}.$
\end{remark}

We begin our study of liftings by establishing the pro-representability
of~$D^{\square}$.

\begin{prop}
  \label{prop: existence of infinite height deformation ring}The functor
  $D^{\square}$ is pro-representable by an object $R^{\square}$ of
  $\pro\cC_{\cO/\varpi^a}$.
\end{prop}
\begin{proof}
  We follow Dickinson's appendix to~\cite{MR1860043}, which uses Grothendieck's
  representability theorem to prove the existence of universal deformation
  rings. By~\cite[Cor.\ to Prop.\ 3.1, \S A]{MR1603480}, it is enough to show that
$D^{\square}$
is left exact.
Left exactness is equivalent to preserving fibre products and terminal
objects. There is a unique terminal object of $\cC_{\cO/\varpi^a},$
namely $\F$, and since
$D^{\square}(\F)$ consists of the trivial lifting, it is also a terminal
object. It remains to check that $D^{\square}$ preserves fibre products in
$\cC_{\cO/\varpi^a}$; but this is obvious if we think of liftings as 
being a choice of $\Phi_R$ lifting~$\Phi$.
\end{proof}

\begin{remark}
\label{rem:big ring}
If $M_R$ is a rank $d$ free
 \'etale $\varphi$-module over the Artinian ring $R$ lifting $M$,
then Lemma~\ref{lem: existence of a lattice in a projective}
shows that we may write $M_R = \gM_R[1/u]$ for a free finite height $\varphi$-module
$\gM_R$ of some height $F$ (depending on $M_R$).
The unicity statement of Lemma~\ref{lem: phi conjugacy
  approximation}~(2) (with~$h=1$) then shows that if $N$ is sufficiently large (depending on $F$),
the elements of $1 + u^N \mathfrak m_R M_d(\gS_R)$ act freely
(in the sense of the action described in Remark~\ref{rem:action})
on the element
of $D^{\square}(R)$ represented by $M$.

Thus, for a non-trivial thickening $R$ of $k$,
if the set $D^{\square}(R)$ is non-empty
(which e.g.\ necessarily will be the
case if $R$ is a $k$-algebra, since in that case we can simply base-change
$M$ from $k$ to $R$),
then it is quite enormous!  In particular, the pro-representing ring $R^{\square}$
constructed in Proposition~\ref{prop: existence of infinite height deformation ring}
does not admit a countable basis of neighbourhoods of zero.
\end{remark}

We now consider liftings of bounded height. If~$M$ is an \'etale
$\varphi$-module with $A$-coefficients, then we say
  that~$M$ admits a model of height~$F$ if there is a (not necessarily
  projective) $\varphi$-module $\gM$ of height~$F$ with $A$-coefficients such that
  $\gM[1/u]=M$.

\begin{df}
\label{def:Kisin liftings} Let~$D^{\square}_F$ be the subfunctor of~$D^{\square}$ whose elements are the liftings
of~$M$ which admit a model of height~$F$. 
\end{df}

\begin{prop}
  \label{prop: existence of finite height deformation ring}The functor
  $D^{\square}_F$ is pro-representable by a quotient $R_F^{\square}$ of
  $R^{\square}$.
\end{prop}
\begin{proof}Let~$A$ be an Artinian quotient of~$R^{\square}$. 
It follows from 
Lemma~\ref{lem: sums and subquotients of height h}
below
that there is a unique maximal quotient~$A^F$ of~$A$ for which the
corresponding \'etale $\varphi$-module admits a model of
height~$F$. We then take~$R_F^{\square}:=\varprojlim_AA^F$.
\end{proof}

\begin{lemma}\label{lem: sums and subquotients of height h}
  Let~$A$ be Noetherian. The direct sum of two \'etale
  $\varphi$-modules with $A$-coefficients which have models
  of height~$F$ also has a model of height~$F$. Any subquotient of an \'etale
  $\varphi$-module with $A$-coefficients with a model of height~$F$ 
  also has a model of height~$F$.
\end{lemma}
\begin{proof}
  The statement about direct sums is trivial. For subquotients, let
  $0\to M'\to M\to M''\to 0$ be a short exact sequence of \'etale
  $\varphi$-modules with $A$-coefficients, and suppose that $\gM$ is a
  finite height $\varphi$-module of height~$F$ with
  $\gM[1/u]=M$. Let~$\gM',\gM''$ be respectively the kernel and image
  of the induced map $\gM\to M''$; then it is easy to check
  that~$\gM'$, $\gM''$ are $\varphi$-modules of height~$F$, and that
  $\gM'[1/u]=M'$, $\gM''[1/u]=M''$.\end{proof}

\begin{remark}
\label{rem:free action}
As noted in Remark~\ref{rem:big ring},
there is an $N$ (depending on $F$, and which we fix once
and for all) such that for each object $R$ of $\cC_{\cO/\varpi^a}$,
the elements of
$1 + u^N \gm_R M_d(\gS_R)$
act freely
on $D^{\square}_F(R)$.
This again implies that the rings $R_F^{\square}$ are huge.

However, since $F$ is now fixed, we may also fix the integer $N$,
and then systematically quotient out by the free action of
$1 + u^N \gm_R M_d(\gS_R)$.  This leads to the following definition.
\end{remark}

\begin{df}
\label{def:Kisin deformations}
We let $D_F :\cC_{\cO/\varpi^a}\to\Sets$ 
denote the functor defined by
$D_F(R) := D^{\square}_F(R)/\bigl(1 + u^N \gm_R M_d(\gS_R)\bigr)$.
\end{df}

In the rest of this section we develop the theory needed to prove
Theorem~\ref{thm: existence of height F deformation ring}, which shows
that~$D_F$ is pro-representable by a Noetherian ring.

\begin{lem}\label{lem: bounding Kisin modules inside each other}Let~$A$ be an
  Artinian $\cO/\varpi^a$-algebra. Suppose that 
  $M$ is an \'etale
  $\varphi$-module with $A$-coefficients, and
  that~$\gM$, $\gM'$ are two models of~$M$ of height~$F$. If~$j$ is
  minimal such that $u^j\gM\subseteq
  \gM'$, and~$k$ is such that $u^k\varphi^*\gM\subseteq\varphi^*\gM'$, then
  $k> (j-a)q$.
\end{lem}
\begin{proof}By assumption we have $u^{j-1}(\gM/\gM\cap\gM')\ne 0$, so
  by Corollary~\ref{cor: exponent of phi} we
  have~$u^{(j-a)q}\varphi^*(\gM/\gM\cap\gM')\ne 0$. It is therefore
  enough to check that
  $\varphi^*\gM\cap\varphi\gM'=\varphi^*(\gM\cap\gM')$; but this
  follows from the flatness of~$\varphi$ (applied to the map
  $\gM\oplus\gM'\to M$, $(x,y)\mapsto x-y$).
\end{proof}

The following lemma is a generalisation of~\cite[Cor.\
3.2.6]{MR2543474}, and the proof is similar.
\begin{lem}
  \label{lem: maximal Kisin modules exist}Let~$A$ be an Artinian
  $\cO/\varpi^a$-algebra. If $M$ is an \'etale $\varphi$-module with
  $A$-coefficients, and~$M$ has a model of height~$F$, then it has
  both a minimal model~$\gM_{\min}$ of height~$F$ and a maximal
  model~$\gM_{\max}$ of height~$F$. Furthermore, $\gM_{\max}/\gM_{\min}$ is an
  $A$-module of finite length.\end{lem}
\begin{proof}Let $\gM$ be a model of~$M$ of height~$F$. We claim that
  there is an~$i\ge 0$ such that any other model~$\gM'$ of height~$F$
  satisfies $u^i\gM\subseteq\gM'\subseteq u^{-i}\gM$. 

  To see this, we follow the proof of~\cite[Prop.\
  2.1.7]{KisinModularity}, and choose~$r$ minimal such
  that~$u^r\gM\subseteq\Phi_{\gM}(\varphi^*\gM)\subseteq u^{-r}\gM$ (note
  that~$r$ exists by Lemma~\ref{lem: F versus u mod p^a}), and
  choose~$j$ minimal such that $u^j\gM\subseteq
  \gM'$, and~$l$ minimal such that $\gM'\subseteq u^{-l}\gM$. We must
  show that $j,l$ are bounded independently of~$\gM'$. It follows from Corollary~\ref{cor: exponent of phi}
  that if ~$u^k\varphi^*\gM\subseteq \varphi^*\gM'$ then $k>
  (j-a)q$. By Lemma~\ref{lem: F versus u mod p^a}
  there is a constant~$n(a)$ such that $F$  divides $u^{n(a)}$
  in~$\gS_A$, so that \[\Phi_{\gM}(\varphi^*\gM)\subseteq
    u^{-r}\gM\subseteq u^{-j-r}\gM'\subseteq
    u^{-n(a)-j-r}\Phi_{\gM'}(\varphi^*\gM').  \] It follows from Lemma~\ref{lem: bounding Kisin modules inside each other} that
  $n(a)+j+r>(j-a)q$, so that~$j$ is bounded independently of~$\gM'$.

Similarly, if we choose~$l$ minimal such that $\gM'\subseteq
u^{-l}\gM$, then \[\Phi_{\gM'}(\varphi^*\gM')\subseteq \gM'\subseteq
  u^{-l}\gM\subseteq u^{-r-l}\Phi_{\gM}(\varphi^*\gM), \]so that by Lemma~\ref{lem:
  bounding Kisin modules inside each other} again we have $r+l>(l-a)q$, and~$l$ is
also bounded independently of~$\gM'$, as required.

Since $u^{-i}\gM/u^{i}\gM$ has finite length as an $A$-module, it 
follows that the category of models of~$M$ of height~$F$ is Artinian
and Noetherian.
To see that maximal and minimal models of height~$F$ exist, it is now
enough to note that if $\gM$, $\gM'$ are models of~$M$ of height~$F$,
then so are~$\max(\gM,\gM'):=\gM+\gM'$ and
$\min(\gM,\gM'):=\gM\cap\gM'$. \end{proof}

\begin{thm}
  \label{thm: existence of height F deformation ring}The
  functor~$D_F$ is pro-representable by a Noetherian
ring $R_F$.
  Furthermore, we may choose a natural transformation $D_F \to D_F^{\square}$
which is a section to the natural transformation $D_F^{\square} \to D_F$;
and thus for any~$R$, we obtain a natural isomorphism
$$D_F(R) \times \bigl( 1 + u^N \mathfrak m_R M_d(\gS_R) \bigr)
\iso
D_F^{\square} (R).$$
\end{thm}
\begin{proof}
Define the group-valued functor $\cG$ on $\cO/\varpi^a$-algebras
via $\cG(A) := 1  + u^N M_d(\gS_A)$
(the group structure being given via multiplication of matrices).
Then $\cG$ is in fact an affine group  scheme, equal to 
$$
\Spec 
(\cO/\varpi^a)[ \{x_{i,j,k,n}\}_{1\leq i \leq f, 1\leq j,k \leq d, N \leq n < \infty}],
$$
where we have written $f := [k:\F_p]$
and have chosen a basis for $W(k)$  over~$\Z_p$,
so that $\gS_A \iso A[[u]]^f$ and $1 + u^N M_d(\gS_A) \iso
\bigl(1 + u^N M_d(A[[u]])\bigr)^f$;
the variables $x_{i,j,k,n}$ then parameterize the coefficient  of $u^n$
in the $(j,k)$ entries of the $i$th summand.

Now let $S$ denote the pro-Artinian completion of 
$$
(\cO/\varpi^a)[ \{x_{i,j,k,n}\}_{1\leq i \leq f, 1\leq j,k \leq d, N \leq n < \infty}]
$$
at the closed point defined by setting $\varpi$,  as well
as  all the $x_{i,j,k,n}$, equal to $0$.
(This point corresponds to the identity $1 \in \M_d(\gS_{\F})$.)
If we set $G := \Spf S$,
then for any object $R$ of $\cC_{\cO/\varpi^a},$
we have
$G(R) = 1 + u^N\mathfrak m_R M_d(\gS_R).$
The free action of  $G(R)$ on $D_F^{\square}(R)$, for each $R$,
induces an equivalence relation
$$G\times D_F^{\square} \rightrightarrows D_F^{\square},$$
the quotient of $D_F^{\square}$ by which is of course just~$D_F$.
The product $G\times D_F^{\square}$ is pro-representable
by 
$R_F^{\square} \cotimes_{\cO/\varpi^a} S$. 
In order to proceed, we need to show that this completed tensor product 
is {\em topologically flat} over $R_F^{\square}$, in the sense of~\cite{MR0213358}.
For this, we have have to give a  more explict description of~$S$, and so we make
a small digression in order to do so.

Let $V$ denote the free $\cO/\varpi^a$-module spanned by
the variables
$$\{x_{i,j,k,n}\}_{1 \leq i \leq f,  1 \leq j,k \leq d, N \leq n < \infty},$$
endowed with its discrete topology.
We let $U \mapsto U^*$ denote Pontrjagin duality on topological $\cO/\varpi^a$-modules.
If $U$ is discrete, then $U^*$  is profinite, and conversely.  In particular,
if $U$ is a discrete  finite $\cO/\varpi^a$-module, then $U^*$ is again
both discrete and finite.
Continuing to assume that $U$ is discrete  and finite,
and writing $V^*_{\delta}$ to denote $V^*$ endowed with its discrete topology,
we find that
\begin{multline*}
\Hom_{\cO/\varpi^a}(V,U)  
= \Hom_{\cO/\varpi^a,\cont}(U^*,V^*)
\\ = \Hom_{\cO/\varpi^a}(U^*,V^*)
= \Hom_{\cO/\varpi^a,\cont}\bigl( (V^*_{\delta})^*, U\bigr).
\end{multline*}
Thus the canonical embedding $V\hookrightarrow (V^*_{\delta})^*$  identifies 
$(V^*_{\delta})^*$  with the pro-Artinian (equivalently, pro-finite) completion
of $V$. 
Thus $S = \cO/\varpi^a[[(V^*_{\delta})^*]]$,
in the notation of~\cite[\S 1.2.5]{MR0213358}.
In particular, $S$ {\em is} topologically flat over $\cO$
and 
$R_F^{\square}\cotimes_{\cO/\varpi^a} S = R_F^{\square}[[(V^*)_{\delta}]]$
is topologically flat  over $R_F^{\square}$.


It follows from~\cite[Thm.\ 1.4]{MR0213358} that the kernel $R_F$ of the corresponding pair of morphisms
$$R_F^{\square} \rightrightarrows
R_F^{\square}\cotimes_{\cO/\varpi^a} S = R_F^{\square}[[(V^*)_{\delta}]].$$
pro-represents the quotient~$D_F$.  
The same result also shows that $R_F^{\square}$ is topologically
flat over $R_F$.  

Since $R_F^{\square}\cotimes_{\cO/\varpi^a} S$ is
a power series ring 
over $R_F^{\square}$,
the morphism $G\times D_F^{\square} \to D_F^{\square}$ satisfies an 
infinitesimal lifting criterion of the type considered in Definition~\ref{def:versal ring}
above.   Thus this morphism is {\em versal}, in the terminology employed in
that definition, or {\em smooth}, in the terminology
  of~\cite[\href{http://stacks.math.columbia.edu/tag/06HG}{Tag
    06HG}]{stacks-project}.
It is then no doubt a matter of general principles that
$D_F^{\square} \to D_F$ also satisfies this infinitesimal lifting property.
In our particular case we can see this directly: since any
surjection $R\to R'$ manifestly induces a surjection $G(R) \to G(R')$,
one immediately confirms that $D_F^{\square} \to D_F$ satisfies the infinitesimal
lifting property. An evident generalization 
  of~\cite[\href{http://stacks.math.columbia.edu/tag/06HL}{Tag
    06HL}]{stacks-project} to our pro-Artinian setting, a detailed
  proof of which may be be found at~\cite{169051},
then shows that $R_F^{\square} \cong R_F[[\{x_i\}_{i \in I}]],$ for some index set~$I$.
The morphism $R_F \to R_F^{\square}$ obtained by mapping each $x_i$ to $0$
then determines a functorial section $D_F\to D_F^{\square}$ with 
the property stated in the theorem.

It remains to show that~$R_{F}$ is
  Noetherian. By~\cite[Prop.\ 5.1, \S A]{MR1603480}, this is
  equivalent to showing that $D_F(\F[\epsilon])$ is a
  finite-dimensional $\F$-vector space. Elements of $D_F(\F[\varepsilon])$ determine self-extensions $0\to
M\to M'\to M\to 0$ with the property that~$M'$ admits a model of
height~$F$. Let~$\gM'$ be such a model; then the image of~$\gM'$
in~$M$ is a model of~$M$ of height~$F$, and thus contains~$\gM_{\min}$,
the minimal model of height~$F$ (whose existence is guaranteed by
Lemma~\ref{lem: maximal Kisin modules exist}). We can replace~$\gM'$
by the preimage of~$\gM_{\min}$, and accordingly we can assume that
the image of~$\gM'$ in~$M$ is~$\gM_{\min}$.

Similarly, the kernel of~$\gM'\to M$ is a model of~$M$, and is
therefore contained in~$\gM_{\max}$; replacing~$\gM'$ by its sum
with~$\gM_{\max}$, we may assume that in fact~$\gM'$ is an extension of
$\gM_{\min}$ by~$\gM_{\max}$. (Having made this replacement, $\gM'$ may
only be of height~$F^2$, rather than of height~$F$, but this does not
matter for our argument.)
It suffices to show that the $\F$-vector space of such
extensions, considered up to the equivalence relation induced by
$(1+u^N\varepsilon M_d(\gS_k))$-conjugacy, is finite-dimensional.

After choosing bases for~$\gM_{\max}$
and~$\gM_{\min}$, the possible matrices for~$\varphi_{\gM'}$ are of
the form\[
  \begin{pmatrix}
    A_{\max}&B\\0&A_{\min}
  \end{pmatrix}
\]where~$A_{\min}$ and~$A_{\max}$ are respectively the matrices
of~$\varphi_{\gM_{\min}}$ and~$\varphi_{\gM_{\max}}$. Conjugating by
matrices of the form~$
\begin{pmatrix}
  1& u^N\varepsilon X\\0&1
\end{pmatrix}
$, we see that we are free to replace~$B$
by~$B+u^NXA_{\min}-A_{\max}\varphi(u^NX)$. It therefore suffices to
show that for some sufficiently large~$M$, for every~$Y\in M_d(\gS_k)$
we can write
\[u^MY=u^NXA_{\min}-A_{\max}\varphi(u^NX)\] for some~$X\in M_d(\gS_k)$.

Since~$\gM_{\min}$ has height~$F$, it follows from Lemma~\ref{lem: F
  versus u mod p^a} that for some~$t\ge 0$ and some~$Z_{\min}\in  M_d(\gS_k)$
we can write $A_{\min}Z_{\min}=u^t\Id_d$. It therefore suffices to
show that we can always solve the
equation \numequation\label{eqn: stupid phi
  equation}u^{M-N}Y=u^tX-u^{(q-1)N}A_{\max}\varphi(X)Z_{\min}\end{equation}
(as a solution 
to this equation with~$Y$ replaced by~$Y Z_{\min}$ provides a solution
to our original equation).

For any~$V\in u^t M_d(\gS_k)$,
write~$\delta(V):=u^{(q-1)N-t}A_{\max}\varphi(V)Z_{\min}\in
M_d(\gS_k)$. Note that if~$V\in u^sM_d(\gS_k)$ for some~$s\ge t+1$,
then $\delta(V)\in u^{qs+(q-1)N-t}M_d(\gS_k)$, and in 
particular~$\delta(V)\in u^{s+1}M_d(\gS_k)$. It follows that the
sum~$W:=V+\delta(V)+\delta^2(V)+\dots \in u^sM_d(\gS_k)$ converges,
and $W-\delta(W)=V$. Therefore, if we set~$M:=N+2t+1$,
take~$V=u^{M-N-t}Y$, and write~$W=X$, we have the required solution
to~\eqref{eqn: stupid phi equation}.
  \end{proof}

If $E'/E$ is a finite extension with ring of integers~$\cO'$, then
we have an obvious map from $\varphi$-modules (and finite height $\varphi$-modules) with
respect to~$\cO$ from those with respect to~$\cO'$, given by tensoring
with~$\cO'$ over~$\cO$. We end this section with the following
reassuring results.

\begin{lem}\label{lem:model of height F and extension of scalars}
  If $M$ is an \'etale $\varphi$-module with $A$-coefficients, then
  $M$ admits a model of height~$F$ if and only
  if~$M\otimes_{\cO}\cO'$ admits a model of height~$F$.
\end{lem}
\begin{proof}
  Since the inclusion $\cO \hookrightarrow \cO'$ is split as an inclusion
  of $\cO$-modules (e.g.\ because~$\cO'$ is a faithfully flat and finite  $\cO$-algebra),
  this is immediate from
  Lemma~\ref{lem: sums and subquotients of height h}.
\end{proof}

If $\F'$ is the residue field of~$\cO'$, then we have corresponding
universal lifting rings~$R_{\cO'}^\square$, $R_{\cO',F}^{\square}$ and~$R_{\cO',F}$
for liftings to $\cO'/\varpi^a$-algebras of $M_{\F}\otimes_{\F}\F'$.

\begin{cor}
  \label{cor: change of O and deformation ring}We have natural
  isomorphisms of $\cO'/\varpi^a$-algebras $R_{\cO'}^{\square}\cong
  R^{\square}\otimes_{\cO}\cO'$, $R_{\cO',F}^{\square}\cong
  R_{F}^{\square}\otimes_{\cO}\cO'$, and $R_{\cO',F}\cong
  R_{F}\otimes_{\cO}\cO'$
\end{cor}
\begin{proof}
This follows easily from Lemma~\ref{lem:change of versal ring under
  separable extension}. Alternatively, we can argue slightly more
explicitly (but essentially equivalently) as follows. We give the
argument for~$R^{\square}$, the other cases being essentially
identical.

Tensoring over~$\cO'$ with~$\cO$ and considering the universal
  property gives a natural map
  $R_{\cO'}^{\square}\to R^{\square}\otimes_{\cO}\cO'$. Let
  $R_{\cO',\F}^{\square}$ be the subring of~$R_{\cO'}^{\square}$
  consisting of elements whose reductions modulo the maximal ideal lie
  in~$\F$. Considering the matrix of the universal \'etale
  $\varphi$-module over~$R_{\cO'}^{\square}$ with respect to a basis
  which lifts (the extension of scalars of) a basis of~$M_\F$, we see
  that this universal \'etale $\varphi$-module is defined
  over~$R_{\cO',\F}^{\square}$. It follows from the universal property
  that there is a natural map
  $R^{\square}\to R_{\cO',\F}^{\square}$.

  Considering the composites
  $R_{\cO'}^{\square}\to R^{\square}\otimes_{\cO}\cO'\to
  R_{\cO',\F}^{\square}\otimes_{\cO}\cO'\to R_{\cO'}^{\square}$ and $R^{\square}\to
  R_{\cO',\F}^{\square}\to R_{\cO'}^{\square}\to
  R^{\square}\otimes_{\cO}\cO'$, we easily obtain the first
  claim. The second claim then follows from Lemma~\ref{lem:model of height F and extension of scalars}.
\end{proof}

\subsection{Moduli of finite height $\varphi$-modules and of \'etale $\varphi$-modules}
\label{subsec:moduli}
In this subsection we will define the moduli stacks that we are interested
in, and prove our key results regarding them. 
We begin by establishing some terminology, which will be important for
all that follows.

We fix an integer $a \geq 1$, and proceed to define various
categories fibred in groupoids (which will in fact be stacks,
although some only in the Zariski topology) over $\cO/\varpi^a$.

\begin{df}
If $a,d \geq 1$ are positive integers,
then for any $\cO/\varpi^a$-algebra $A$,
we define $\cR^a_d(A)$ to be the groupoid of
 \'etale $\varphi$-modules with $A$-coefficients which are projective
of rank~$d$ over $\OEA$.
If $A \to B$ is a morphism of $\cO/\varpi^a$-algebras,
and if $M$ is an object of $\cR^a_d(A)$,
then the pull-back of $M$ to $\cR^a(B)$ is defined to be the tensor product
$\OEB\otimes_{\OEA} M$.
\end{df}

The resulting category fibred in groupoids $\cR^a_d$ is in fact
a stack in groupoids in the \emph{fpqc} topology over $\cO/\varpi^a$,
as follows from the results of~\cite{MR2181808},
and more specifically from Theorem~\ref{thm: fpqc locality of projective and locally free} above. 

\begin{df}
	If $a,d \geq 1$ are positive integers, and if $F \in (W(k)\otimes_{\Zp}\cO)[u]$ is a polynomial
	that is congruent to a positive power of $u$ modulo $\varpi$,
        then for any $\cO/\varpi^a$-algebra~$A$,
	we define $\cC^a_{d,F}(A)$ to be the groupoid of
	$\varphi$-modules of height $F$ with $A$-coefficients which are
	projective of rank $d$ over $\gS_A$.
        If $A \to B$ is a morphism of $\cO/\varpi^a$-algebras,
        and if $\gM$ is an object of $\cC^a_{d,F}(A)$,
        then the pull-back of $\gM$ to $\cC^a_{d,F}(B)$
	is defined to be the tensor product $\gS_B\otimes_{\gS_A} \gM$.
\end{df}

Again, it follows from Theorem~\ref{thm: fpqc locality of projective and locally free} that
the resulting category fibred in groupoids $\cR^a_d$ is in fact
a stack in groupoids in the \emph{fpqc} topology over $\cO/\varpi^a$. There is an obvious morphism
$\cC^a_{d,F}\to\cR^a_d$,
defined by sending $(\gM,\varphi)$ to $(\gM[1/u],\varphi)$.

\begin{remark}
	Our notation, and the entire set-up that we have just
	introduced, is very much inspired by 
the work of Pappas and Rapoport~\cite{MR2562795}. 
Indeed, in the case when~$q=p$, $\varphi(u)=u^p$, and $F\in W(k)[u]$ is an
Eisenstein polynomial,  our stack $\cC^a_{d,F}$ coincides with the stack~$\cC^a_{h,W(k)[u]/F}$
defined in~\cite[\S 3.b]{MR2562795}. However,
our stack $\cR^a_d$ is subtly different from the stack 
denoted in the same manner in \cite{MR2562795}.  In that 
reference, the \'etale $\varphi$-modules under consideration
are not required to be projective, but are required
to be {\em fpqc} locally free.  However, it seems to us 
that it is necessary to impose this projectivity 
in order to obtain a stack, while the local freeness
hypothesis seems unnatural from the point of view 
of our intended applications (which is that $\cR^a_d$ should provide
models for moduli stacks of local Galois representations ---
and a direct summand of a family of representations
should again form such a family); also, at a technical
level, the effectivity result 
of Theorem~\ref{thm: hypotheses of the main thm hold in our
  setting}~(5) below
depends on working with projective \'etale $\varphi$-modules
that are not necessarily locally free over the coefficient ring.
\end{remark}

In spite of the difference between our definition 
of $\cR^a_d$ and that
of~\cite{MR2562795}, 
we nevertheless
rely on many of the arguments of that reference.  In order
to make the connection between our set-up and that 
of~\cite{MR2562795}, 
it is helpful to introduce the following auxiliary objects,
in which we require projectivity of our \'etale $\varphi$-modules,
but also impose freeness conditions,
as in~\cite{MR2562795}. 

\begin{df}
	We define $\cR^a_{d,\free}$ to be the full subgroupoid
	of $\cR^a_d$ classifying free \'etale $\varphi$-modules
	of rank $d$.

	If $\tau$ is any topology on the category of
	$\cO/\varpi^a$-modules lying between the Zariski topology and 
	the \emph{fpqc} topology,
	then we define $\cR^a_{d,\tau-\free}$ to be the
	full subgroupoid of $\cR^a_d$ classifying projective
	\'etale $\varphi$-modules of rank $d$ that
	are furthermore $\tau$-locally free over the
	ring of coefficients.
	\end{df}

Taking into account the fact that
$\cR^a_d$ is an \emph{fpqc} stack,
the category fibred in groupoids $\cR^a_{d,\tau-\free}$ 
is evidently a stack in the topology $\tau$.
Indeed, one immediately checks that it is the $\tau$-stackification
of $\cR^a_{d,\free}$.

\begin{remark}
As already noted, the category fibred
in groupoids $\cR^a_{d,\free}$,
as well as the stacks $\cR^a_{d,\tau-\free}$, will
play a purely auxiliary role. 
Furthermore, we need only make one choice of topology
$\tau$ and work with that particular choice throughout;
e.g.\ we could simply take $\tau$ to be the Zariski topology.
\end{remark}

From now on, for the duration of the paper,
we fix a choice of~$a\ge 1$, and omit it from the notation.

\begin{lemma}
	The morphism $\cC_{d,F} \to \cR_d$ factors through
	$\cR_{d,\tau-\free}$ {\em (}for any choice of $\tau${\em )}.
\end{lemma}
\begin{proof}
	This follows from Proposition~\ref{prop:improving fpqc for
          A[[u]]} ~(1), which shows that a projective finite height $\varphi$-module is actually Zariski locally free.\end{proof}

\begin{prop}
	\label{prop:freeness stratification}
	If $\Spec A \to \cR_d$ is a morphism with $A$ a Noetherian
	$\cO/\varpi^a$-algebra, then there exists a scheme-theoretically
	surjective morphism $\Spec B \to \Spec A$ such that the
	composite morphism $\Spec B \to \cR_d$ factors through
	$\cR_{d,\free}.$
\end{prop}
\begin{proof}
	Let $M$ be the \'etale $\varphi$-module over $A$ classified
	by the given $A$-valued point of $\cR_d$,
	and let $M_{\red}$ denote the base-change of $M$ over $A_{\red}$.
	By Lemma~\ref{lem: existence of a lattice in a projective}, we may find a (not necessarily projective) finite height $\varphi$-module $\gM
        \subset M_{\red}$ with~$\gM[1/u]=M_{\red}$. The quotient
	$\gM/u\gM$ is then a coherent sheaf on $\Spec A_{\red}$, 
	and so we may find a dense open subset $U 
	\subset \Spec
	A_{\red}$ such that $\gM/u\gM$ restricts to a free
	sheaf over $U$.  Thus, by Proposition~\ref{prop: criterion for
          projectivity over A[[u]]}, $\gM$ restricts to a free finite height $\varphi$-module over $U$, and so $M_{\red}$ restricts to a free \'etale $\varphi$-module
       over $U$ (necessarily of rank~$d$).

       We may regard $U$ as an open subset of $\Spec A$,
       and without loss of generality we may in fact
       assume that $U = \Spec A_f$ for some $f \in A$.
       Because $A_f$ is Noetherian, the nilradical of $A_f$ is nilpotent,
       and so the kernel of the morphism $A_f((u)) \to (A_f)_{\red}((u))$
       is also nilpotent.
       Thus the restriction of $M$ to $U$
       is a projective $A_f((u))$-module which becomes free modulo
       a nilpotent ideal.
       By a standard Nakayama-type argument,
       we see that this restriction
       itself is a free \'etale
       $\varphi$-module
      over $U$. 

       We now note that we may choose a closed subscheme
       $Z \hookrightarrow \Spec A$ whose underlying closed subset is equal
       to $\Spec A \setminus U,$ and for which the obvious morphism
       $U \coprod Z \to \Spec A$ is scheme-theoretically dominant,
       in addition to being surjective.
       (Since $A$ is Noetherian, the kernel $A[f^{\infty}]$ of 
       $A \to A_f$ is equal to $A[f^n]$ for some $n\geq 1$,
       and we may take $Z = \Spec A/f^n$.)
       The proposition now follows by an evident Noetherian induction.
       \end{proof}

\begin{prop}\label{prop: hom isom}Let $A$ be an $\cO/\varpi^a$-algebra, and let $M,N$ be projective
  \'etale $\varphi$-modules of finite rank with~$A$-coefficients. Then
  the functors on $A$-algebras taking~$B$ to $\Hom(M_B,N_B)$ and
  $\Isom(M_B,N_B)$ are both represented by affine schemes of finite
  presentation over~$A$.
  \end{prop}
  \begin{proof}By Lemma ~\ref{lem: projective phi module is summand of
    free}, there are projective
    \'etale $\varphi$-modules $P,Q$ of finite rank such that the
    \'etale $\varphi$-modules $F:=M\oplus P$ and $G:=N\oplus Q$ are
    both free of finite rank. We now follow the proof of~\cite[Cor.\
    2.6(b)]{MR2562795}.  Choosing bases of~$F,G$ as~$\OEA$-modules, an
    element of $\Hom(F_B,G_B)$ is given by a matrix~$g$ with
    coefficients in~$\OEB$. If the matrices of~$\varphi_M$,
    $\varphi_N$ with respect to the chosen bases are respectively~$X$,
    $Y$ then the condition that~$g$ respects $\varphi$ is that
    $\varphi(g)=Y^{-1}gX$. 

    Choose an integer $n\ge 0$ such that $X, X^{-1}, Y, Y^{-1}$ all
    have entries with poles of degree at most~$n$, and let~$s\ge 0$ be
    minimal such that~$g$ has poles of degree at
    most~$s$. By Corollary~\ref{cor: exponent of phi}, $\varphi(g)$
    has poles of degree greater
    than~$(s-a)q$. Since~$\varphi(g)=Y^{-1}gX$, we see that we must
    have $(s-a)q<2n+s$, whence $s<(2n+aq)/(q-1)$.

    Writing the matrix~$g$ as $\sum_{i=-s}^\infty g_iu^i$,
    $g_i\in M_d(W(k)\otimes_{\Zp}B)$, the equation
    $g=Y\varphi(g)X^{-1}$ and~ Lemma~\ref{lem: phi conjugacy approximation}
    show that the~$g_i$ for~$i\le (2n+(a-1)q)/(q-1)$ determine all of
    the~$g_i$. It follows that~$\Hom(F_B,G_B)$ is represented by an
    affine scheme of finite presentation over~$A$. 

Let~$e\in\End(F)$,
    $f\in \End(G)$ be the idempotents corresponding to~$M$, $N$
    respectively. Since $\Hom(M_B,N_B)\subset\Hom(F_B,G_B)$ is given
    by those~$g$ satisfying $g(1-e)=0$ and $(1-f)g=0$, we see that it
    is represented by a closed subscheme of the scheme
    representing~$\Hom(F_B,G_B)$, and is therefore of finite
    presentation (for example by  Lemma~\ref{lem:another finiteness lemma}).
    Finally, the result
    for~$\Isom(M_B,N_B)$ follows by regarding it as the subfunctor of
    pairs $(\alpha,\beta)\in\Hom(M_B,N_B)\times\Hom(N_B,M_B)$
    satisfying~$\alpha\beta=\Id_{N_B}$, $\beta\alpha=\Id_{M_B}$.
  \end{proof}

The following theorem generalises some of the main results 
of~\cite{MR2562795} to our setting. The proofs are almost identical,
and we content ourselves with explaining the changes that need to be
made to the arguments of~\cite{MR2562795}, rather than writing them
out in full.

\begin{thm}\label{main result of PR}

(1) The stack
$\cC_{d,F}$ is an algebraic stack of finite presentation over
  $\Spec\cO/\varpi^a$, with affine diagonal. 

(2) The morphism
  $\cC_{d,F}\to\cR_{d,fpqc-\free}$ is representable by
  algebraic spaces, proper, and of finite presentation.

(3) The diagonal morphism $\Delta:\cR_{d,fpqc-\free}\to\cR_{d,fpqc-\free}\times_{\cO/\varpi^a}\cR_{d,fpqc-\free}$ is
representable by algebraic spaces, affine, and of finite
presentation.\end{thm}
  \begin{proof}
In the case that~$\cO=\Zp$, $q=p$, and
    $\varphi(u)=u^p$, it follows from the main results
    of~\cite{MR2562795} that~$\cC_{d,F}$
    is an algebraic stack of finite type over $\Spec \cO/\varpi^a$,
    and that~(2) holds. (Strictly speaking, \cite{MR2562795} assume that~$F$
    is an Eisenstein polynomial, but their arguments go through
    unchanged with our slightly more general choice of~$F$.) In the
    case of general~$\cO$, $q$
    and~$\varphi$, the arguments go over essentially unchanged
    provided that one replaces the use of~\cite[Prop.\ 2.2]{MR2562795}
    with an appeal to
    Lemma~\ref{lem: phi conjugacy approximation}, and that in~\cite[\S
    3]{MR2562795} one replaces~$eah$ by the quantity~$n(a,h)$
    appearing in Lemma~\ref{lem: F versus u mod p^a}.

    Part~(3) is immediate from Proposition~\ref{prop: hom isom}. To prove the remaining claims of~(1), we have to show that
    $\cC_{d,F}$ is in fact of finite presentation over $\cO/\varpi^a$,
    with affine diagonal.  These facts are certainly implicit in 
    the arguments of~\cite{MR2562795}, but for the reader's convenience,
    we explain how they follow formally from the results already
    established.
    Since $\cO/\varpi^a$ is Noetherian, and since
    we know already that $\cC_{d,F}$ is finite type over $\cO/\varpi^a$,
    the diagonal morphism $\cC_{d,F} \to \cC_{d,F}\times_{\cO/\pi^a} \cC_{d,F}$
    is automatically quasi-separated (being a representable morphism
    between finite type algebraic stacks over $\cO/\varpi^a$), and so
    to show that $\cC_{d,F}$ is of finite presentation over $\cO/\varpi^a$,
    it suffices to show that this diagonal morphism is quasi-compact.
    Since affine morphisms are quasi-compact, this will follow once
    we show that $\cC_{d,F}$ has affine diagonal.  For this, we 
    factor the diagonal of $\cC_{d,F}$ as 
    $$\cC_{d,F} \to \cC_{d,F}\times_{\cR_{d,fpqc-\free}}\cC_{d,F} \to \cC_{d,F}\times_{\cO/\varpi^a}
    \cC_{d,F}.$$
    The first of these morphisms is a closed immersion, since $\cC_{d,F} 
    \to \cR_{d,fpqc-\free}$ is representable and proper (by (2)), while the second
    is affine, being a base-change of the diagonal
    morphism of $\cR_{d,fpqc-\free} \to \cR_{d,fpqc-\free} \times_{\cO/\varpi^a} \cR_{d,fpqc-\free}$ (which
    is affine, by (3)).
    Their composite is thus an affine morphism, as claimed.
  \end{proof}
  In Theorem~\ref{thm: morphism from C to Rproj} we prove the analogue
of Theorem~\ref{main result of PR} for~$\cR_d$. In order to do so we
need the following Lemma.
  \begin{lem}
    \label{lem: representability of automorphisms of lattices}Let $A$
    be an $\cO/\varpi^a$-algebra, and let~$M$ be an \'etale
    $\varphi$-module with $A$-coefficients, which is free of rank~$d$ as an
    $\OEA$-module. Let~$T$ be an automorphism of~$M$. Then the functor
    on~$A$-algebras taking~$B$ to the set of $T$-invariant projective $\varphi$-modules  $\gM_B\subset M_B$ of rank~$d$ and height~$F$ is representable by a
    projective $A$-scheme.\end{lem}
  \begin{proof}Let~$\Gr$ denote the affine Grassmannian classifying
    projective $\gS_A$-lattices in~$M$; this is an Ind-projective
    $A$-scheme. We begin by showing that the subfunctor of $\Gr\times\Gr$ given by
    \[\{(\gM,\gN):  \gM\subseteq\gN\}\] is a closed Ind-subscheme of
    $\Gr\times\Gr$. To
    see this, we have to show that for any $A$-algebra $B$, and any
    pair of lattices $\gM_B,\gN_B\in\Gr(B)$, the locus
    in $\Spec B$ over which $\gM_B\subseteq\gN_B$ is
    closed. Equivalently, we need to show that the locus over which
    the morphism $\gM_B\to M_B/\gN_B$ vanishes is closed. To see this,
    note that we may factor 
    this map as $\gM_B\to P\to Q\to M_B/\gN_B$,
    where $P$ is a finite rank (as a $B$-module) projective quotient (hence direct summand)
    of~$\gM_B$,
    and $Q$ is a finite rank projective
    direct summand of $M_B/\gN_B$. Since the maps $\gM_B\to P$ and $Q\to M_B/\gN_B$
    are both split, we are in fact considering the locus over which
    the morphism $P\to Q$ vanishes, and this is obviously closed, as
    it is given by the vanishing of matrix entries.

The endomorphism~$T$ induces an automorphism $T_*$ of $\Gr$ (taking
$\gM$ to $T(\gM$)), and we let $\Gamma_T:=T_*\times\id:\Gr\to\Gr\times\Gr$ be the
graph of~$T$. Pulling back the closed locus considered above
by~$\Gamma_T$, we see that there is a closed Ind-subscheme $\Gr^T$
of~$\Gr$, classifying the lattices~$\gM$ with $T(\gM)\subset\gM$.

The $\varphi$-module $M$ corresponds to a morphism $\Spec
    A\to\cR_{d,fpqc-\free}$, and for each $F$, the fibre product $\Gr^{F}:=\Spec A\times_{\cR_{d,fpqc-\free}}\cC_{d,F}$ is a closed subscheme of~$\Gr$
    (it is a scheme by Theorem~\ref{main result of PR}~(2)). Then the
    intersection of $\Gr^{F}$ and $\Gr^T$ is closed in~$\Gr^{F}$, and is therefore projective over~$\Spec A$, as required.    
  \end{proof}

We let~$\End(\cR_{n,\free})$ be the category fibred in groupoids
over~$\cO/\varpi^a$ with $\End(\cR_{n,\free})(A)=\{(M,f)\}$ where
$M\in\cR_{n,\free}(A)$ and $f\in\End(M)$. It contains a subcategory
fibred in groupoids~$\Proj_{n,d}$, classifying those pairs $(M,f)$ for
which $\Im f$ is projective of rank~$d$. There are natural morphisms
$\Proj_{n,d}\subseteq\End(\cR_{n,\free})\to\cR_{n,\free}$ and $\Proj_{n,d}\to\cR_d$, which
respectively take~$(M,f)$ to~$M$ and to~$\Im f$. These morphisms fit
into the following commutative diagram.\[ 
\begin{tikzcd}
 \cC_{n,F}\arrow{dd}& &  & \Proj_{n,d} \arrow{ldd}
\arrow{rdd} \arrow[hookrightarrow]{r}
& \End(\cR_{n,\free}) \arrow [bend left]{lldd}\\ \\
\cR_n & & \cR_{n,\free} \arrow[left hook->]{ll}  & &\cR_d\\
\end{tikzcd}
\]

  \begin{thm}
    \label{thm: morphism from C to Rproj}(1) The morphism
  $\cC_{d,F}\to\cR_d$ is representable by
  algebraic spaces, proper, and of finite presentation.

(2) The diagonal morphism $\Delta:\cR_d\to\cR_d\times_{\cO/\varpi^a}\cR_d$ is
representable by algebraic spaces, affine, and of finite
presentation.

(3) $\cR_d$ satisfies~\emph{[1]}.
  \end{thm}
  \begin{proof}We begin with~(1). Let~$B$ be an $A$-algebra, and let
    $\Spec B\to\cR_d$ be a morphism, corresponding to a projective
    \'etale $\varphi$-module $M_B$ of rank~$d$. We need to show that $\Spec
    B\times_{\cR_d}\cC_{d,F}$ is representable by a proper
    algebraic space over~$\Spec B$ of finite presentation. This can be
    checked \'etale locally, so in particular
    by Lemma~\ref{lem:locally free over u^n} we can assume that $M_B$ is free over
    $(W(k)\otimes_{\cO/\varpi^a}B)((u^n))$ for some~$n$. The claim follows from Lemma~\ref{lem:
      representability of automorphisms of lattices}, applied with~$u$
    replaced by~$u^n$, and~$T$ being given by multiplication by~$u$.  

Part~(2) is immediate from Proposition~\ref{prop: hom isom}. For~(3),
by Proposition~\ref{prop: double diagonals of stacks satisfying [3]}
and part~(2), it is enough to show that $\cR_d\to\Spec\cO/\varpi^a$ is
limit preserving on objects. To this end,
suppose that we have a morphism $T\to\cR_d$, where~$T=\varprojlim T_i$
is a limit of affine schemes. By Lemma~\ref{lem: projective phi module is summand of
    free}
we can lift the morphism $T\to\cR_d$ to a morphism $T\to\Proj_{n,d}$
for some~$n$. The composite morphism $T\to\Proj_{n,d}\to\cR_{n,\free}\to\cR_n$
lifts to a morphism $T\to\cC_{n,F}$ for some~$F$ (because every free
\'etale $\varphi$-module contains a free finite height $\varphi$-module,
by Lemma~\ref{lem: existence of a lattice in a projective})
and
since~$\cC_{n,F}$ is locally of finite presentation, this morphism
factors through~$T_i$ for some~$i$.

Consequently, the composite $T\to\Proj_{n,d}\to\cR_{n,\free}$ factors
through~$T_i$, and it suffices to prove that the morphism
$\Proj_{n,d}\to\cR_{n,\free}$ is limit preserving on objects. This
follows from Proposition~\ref{prop: hom isom}, which shows that the
morphism~$\Proj_{n,d}\to\cR_{n,\free}$ is representable by schemes
of finite presentation (note that the condition that an endomorphism
of a free module be idempotent is a closed condition, and is therefore
of finite presentation by Lemma~\ref{lem:another finiteness lemma}).\end{proof}

By Theorem~\ref{thm: morphism from C to Rproj}, the running assumptions of
Section~\ref{subsec:scheme images two} apply to the morphism
$\cC_{d,F}\to\cR_{d}$; so we may use
Definition~\ref{def:stack-theoretic image} to define the scheme-theoretic image of
$\cC_{d,F}\to\cR_{d}$, which we denote by
$\cR_{d,F}$. The main result of this section is Theorem~\ref{thm: hypotheses of the main thm hold in our
    setting} below, showing that $\cR_{d,F}$ is an algebraic
  stack. Before proving it, we study the versal rings of~$\cR_d$ and~$\cR_{d,F}$.

Let~$\F'/\F$ be a finite extension, and let~ $M_{\F'}$ be an \'etale
$\varphi$-module with $\F'$-coefficients, corresponding to a finite
type point $x:\Spec\F'\to\cR_d$.  Write $\cO'$ for
 the ring of integers in the compositum of~$E$
 and~$W(\F')[1/p]$, so that~$\cO'$ has residue field~$\F'$. 

As in~Section~\ref{sec: effectivity of projectivity}, we  fix a choice of (ordered)
 $\cO_{\cE,\F'}$-basis of~$M_{\F'}$, and we let $D^\square:\cC_{\cO'/\varpi^a}\to\Sets$ be
 the functor taking~$R$ to the set of isomorphism classes of liftings
 of~$M_{\F'}$ to~$R$. By Remark~\ref{rem:action}, the group functor $H$ defined via $H(R) := R^{\times} +
\mathfrak m_R M_d(\cO_{\cE,R})$ acts on $D^{\square}$ via change of basis.
 By Proposition~\ref{prop: existence of infinite height deformation
   ring}, $D^\square$ is pro-representable by an object~ $R^{\square}$
 of~$\pro\cC_{\cO'/\varpi^a}$.   

\begin{lemma}
\label{lem:R versal}
The natural morphism $D^{\square} \to \widehat{\cR}_{d,x}$,
defined by mapping any lift $M_R$ of $M_{\F'}$ over some test
object~$R$ to the underlying \'etale $\varphi$-module
{\em (}i.e.\ forgetting the choice of basis of $M_R$,
as well as the chosen isomorphism between $M_R/\mathfrak m_R$ and~$M_{\F'}${\em )},
is versal, and is also $H$-equivariant, for the change-of-basis action
of $H$ on~$D^{\square}$ and for the trivial action of $H$ on~$\widehat{\cR}_{d,x}$.
\end{lemma}
\begin{proof}
The claimed equivariance is clear, since the morphism is defined
in part by forgetting the chosen bases.
To see the claimed versality, it suffices to show that if $M_{A}$ is an \'etale
$\varphi$-module with $A$-coefficients, where $A$ is a finite Artinian
$\cO'/\varpi^a$-algebra, and if $M_B$ is an \'etale $\varphi$-module with
$B$-coefficients, with $B$ a finite Artinian $\cO'/\varpi^a$-algebra admitting a
surjection onto $A$, such that the base change $(M_B)_A$ of $M_B$ to
$A$ is isomorphic to $M_A$, and we have compatible choices of basis
for~$M_A$ and~$M_B$, then we may find an
$M'_B$ (together with such a basis) which lifts $M_A$, and is isomorphic to $M_B$.  The existence
of such a lift amounts to showing that we can lift the choice of
basis, and this is clear from the surjectivity of~$\GL_d(\OEB)\to\GL_d(\OEA)$. \end{proof}

The preceding lemma shows
in particular that~$\cR_d$
admits versal rings at all finite type points.

Let $\cC_{\Spf R^{\square}}$ denote the pull-back of $\cC_{d,F} \to \cR_d$
along the versal morphism $D^{\square} = \Spf R^{\square} \to \cR_d$,
and let~$R^{\square,\cC}$ be the scheme-theoretic image (in the sense of
Definition~\ref{defn:definition of scheme-theoretic images for versal
  rings}) of the morphism $\cC_{\Spf R^\square}\to \Spf R^\square$. By
Lemmas~\ref{lem:versal ring for Z} and~\ref{lem:R versal}, there is a versal morphism $\Spf
R^{\square,\cC}\to\widehat{\cR}_{d,F,x}$.

\begin{df}
\label{def:framed C-defs}
We let $D^{\square,\cC}$ denote the subfunctor of $D^{\square}$ represented
by $\Spf R^{\square,\cC}$. 
\end{df}

\begin{remark}
\label{rem:H equivariance}
The equivariance statement of Lemma~\ref{lem:R versal}
implies that the $H$-action on $D^{\square}$ restricts to an $H$-action
on~$D^{\square,\cC}$,
and the the morphism
$D^{\square,\cC} := \Spf R^{\square,\cC} \to
\widehat{\cR}_{d,F,x}$
is $H$-equivariant, with respect to the induced $H$-action
on~$D^{\square,\cC}$,
and the trivial $H$-action on~$\widehat{\cR}_{d,F,x}$.
\end{remark}

Our goal will be to show that an appropriately chosen subgroup functor
of $H$ acts freely on $D^{\square,\cC}$, and that the corresponding 
quotient $D^{\cC}$ of $D^{\square,\cC}$ also admits a versal morphism
to $\widehat{\cR}_{d,F,x}$, and is {\em Noetherianly} pro-representable.
To this end, we will relate $D^{\square,\cC}$ to the subfunctor $D^{\square}_F$
of $D^{\square}$, where,
as in Section~\ref{sec: effectivity of projectivity}, we let~$D_F^\square$ be the subfunctor of~$D^{\square}$ consisting of
those~$M$ which have a model of height~$F$. In fact, we will show that $D^{\square,\cC}$ is a subfunctor of~$D^{\square}_F$;
equivalently, we will show that $R^{\square,\cC}$, which is
{\em a priori} a quotient of $R^{\square}$,
is actually a quotient
of $R^{\square}_F$
(the quotient of $R^{\square}$ that pro-represents~$D^{\square}_F$,
whose existence is proved in
Proposition~\ref{prop: existence of finite height deformation ring}).
  We begin with a useful general criterion for such
  a factorisation of the map~$R^{\square}\to R^{\square,\cC}$ to exist.

We momentarily place ourselves in the general deformation-theoretic context
of Subsection~\ref{subsec:axiom 2}; that is,
we fix a Noetherian ring $\Lambda$, and a finite ring map $\Lambda\to k$,
whose target is a field.
We will also allow ourselves to use the language of formal algebraic spaces 
from~\cite[\href{http://stacks.math.columbia.edu/tag/0AHW}{Tag
  0AHW}]{stacks-project}, but not in a serious way.
(If $R$ is a local ring with maximal ideal $\mathfrak m$,
then $\Spf R$ is simply the Ind-scheme $\varinjlim_i \Spec R/\mathfrak m^i$,
and giving a finite type morphism of formal algebraic spaces
$X \to \Spf R$ amounts to giving a collection of compatible finite type
morphisms of algebraic spaces $X_i \to \Spec R/\mathfrak m^i$.)

\begin{lem}\label{lem: criterion to contain scheme theoretic image
    versal ring}
	Let $R \to S$ be a continuous surjection of objects in
	$\pro\cC_{\Lambda}$, 
	let $X\to \Spf R$ be a finite type morphism
	of formal algebraic spaces,
     and make the following assumption: 
	if $A$ is any finite-type Artinian local $R$-algebra
	for which the canonical morphism $R\to A$ factors through 
	a discrete quotient of $R$,
	and for which the canonical morphism
	$X\times_{\Spf R} \Spec A \to \Spec A$
	admits a section,
	then the canonical morphism $R \to A$ furthermore factors
        through $S$.

        Then if $A$ is any discrete Artinian quotient of $R$ for which
	the base-changed morphism
	$X\times_{\Spf R} \Spec A \to \Spec A$
	is scheme-theoretically dominant,
	the surjection $R\to A$ factors through $S$.
\end{lem}
\begin{proof}
  The desired conclusion is equivalent to the claim 
  that the closed immersion
  \numequation
  \label{eqn:closed imm}
  \Spf S\times_{\Spf R} \Spec A \to \Spec A
  \end{equation}
  is an isomorphism.
  By Yoneda's lemma it is enough to
  show that the following condition (*) holds whenever $B$ is a finite
  type $A$-algebra: 

(*) Any morphism $\Spec B\to \Spec A$ which can be factored through
$X\times_{\Spf R} \Spec A$ necessarily factors through
the closed immersion~(\ref{eqn:closed imm}).

(Indeed, if (*) holds, then, since the morphism
$X\times_{\Spf R} \Spec A \to \Spec A$
is of finite type, by assumption, we see that it
factors through~(\ref{eqn:closed imm}).
On the other hand, this morphism is scheme-theoretically
dominant, by assumption; thus~(\ref{eqn:closed imm}) is an isomorphism,
as required.)

  By considering the product of the localisations of a finite type
  $A$-algebra $B$ at all its maximal ideals, this will follow if we
  prove~(*) when $B$ is the localisation of a finite type $A$-algebra
  at one of its maximal ideals.  Since such a localisation is
  Noetherian, this in turn will follow if we prove~(*) when $B$ is the
  completion of a finite type $A$-algebra at one of its maximal
  ideals.  Considering the reduction of such a completion modulo the
  various powers of its maximal ideal, we then reduce further to
  proving~(*) in the the case when $B$ is a finite type Artinian local
  $A$-algebra. But in this case, condition~(*) holds by assumption.
  \end{proof}

\begin{prop}\label{prop: versal ring factors through finite height}
$D^{\square,\cC}$ is a subfunctor of $D^{\square}_F$.
\end{prop}
\begin{proof} The claim of the proposition amounts to showing that the surjection
$R^{\square} \to R^{\square,\cC}$ factors through $R^{\square}_F$.
      By Lemma~\ref{lem: formal functions implies formal image equals image formal}, we may write $R^{\square,\cC}$ as the inverse limit
      of Artinian quotients $A$, for each of which the base-changed
      morphism $\cC_{d,F,a} \to \Spec A$ is scheme-theoretically dominant.
      It suffices to show that each of the composite surjections
      $R^{\square} \to R^{\square,\cC} \to A$ factors through $R^{\square}_F$.
      This will follow 
      from Lemma~\ref{lem: criterion to contain scheme theoretic image
        versal ring}, taking~$R=R^{\square}$, ~$S=R^{\square}_F$,
      and~$X=\cC_{d,F,\Spf R^{\square}}$, provided we show that the hypotheses of that
lemma hold.

To this end, let~$A$ be a finite type Artinian local $R^{\square}$-algebra for which
  the canonical morphism $\cC_{d,F,\Spec A}\to\Spec A$ admits a section,
  and let~$M_A$ denote the \'etale $\varphi$-module corresponding
  to the induced morphism $\Spec A\to\Spf R^{\square}$.
  The existence of
  the section to $\cC_{d,F,\Spec A}$
  is, by definition, equivalent to the existence
  of a projective $\varphi$-module~$\gM_A$ of height~$F$ such that $\gM_A[1/u]=M_A$.
  Again by definition, we have
  $\gM_A\otimes_A\kappa(A)=M_{\F'}\otimes_{\F'}\kappa(A)$. By Corollary~\ref{cor: change of O and deformation ring}, the functor $D'_F$ of liftings
  of~$M_{\F'}\otimes_{\F'}\kappa(A)$ which have a model of type~$F$ is
  pro-represented by~$R^{\square}_F\otimes_{W(\F')}W(\kappa(A))$, so in
  particular the existence of~$\gM_A$ implies that the morphism 
  $\Spec A\to\Spf R^{\square}$ factors through~$\Spf R^{\square}_F$, as required.
\end{proof}

It follows from Proposition~\ref{prop: versal ring factors through finite height},
together with Remark~\ref{rem:free action},
that the action of $H(R) := R^\times + \mathfrak m_R M_d(\cO_{\cE,R})$
on $D^{\square,\cC}(R)$ (for any test object $R$),
whose existence was noted in Remark~\ref{rem:H equivariance},
restricts to a free action of $G(R) := 1 + u^N \mathfrak m_R M_d(\gS_R).$
We then make the following definition (in analogy to
Definition~\ref{def:Kisin deformations}).

\begin{df}
\label{def:unframed C-defs}
We let $D^\cC :\cC_{\cO/\varpi^a}\to\Sets$ 
denote the functor defined by
$D^{\cC}(R) := D^{\square,\cC}(R)/G(R)$. \end{df}

By construction there is a Cartesian square
$$\xymatrix{ D^{\square,\cC} \ar[r]\ar[d] & D^{\square}_F \ar[d] \\
D^\cC \ar[r] & D_F}$$
and the section $D_F \to D_F^{\square}$ of
Theorem~\ref{thm: existence of height F deformation ring}
then restricts to a section $D^{\cC} \to D^{\square,\cC}$.
An argument almost identical to that used in the proof of
Theorem~\ref{thm: existence of height F deformation ring}
shows that $D^{\cC}(R)$ is pro-representable by some $R^{\cC}$;
since $D^{\cC}$ is a subfunctor of $D_F$, 
this pro-representing object is a quotient of the ring $R_F$
that pro-represents $D_F$.   As the latter ring is Noetherian
(by Theorem~\ref{thm: existence of height F deformation ring}),
so is~$R^{\cC}$.


With these various definitions and observations in place,
we are now ready to prove our main theorem.\begin{thm}
  \label{thm: hypotheses of the main thm hold in our
    setting}The hypotheses of Theorem~{\em \ref{thm:main-intro}} hold for the morphism
  $\cC_{d,F}\to\cR_{d}$. That is:
  \begin{enumerate}
  \item $\cC_{d,F}$ is an algebraic stack, locally of finite presentation
    over~$\Spec \cO/\varpi^a$.
  \item $\cR_{d}$ satisfies~\emph{[3]}, and its diagonal is locally of
    finite presentation.
  \item $\cC_{d,F}\to\cR_{d}$ is a proper morphism.
  \item $\cR_{d}$ admits versal rings at all finite type points.
  \item $\cR_{d,F}$ satisfies~\emph{[2]}.
  \end{enumerate}
Accordingly,  $\cR_{d,F}$ is
an algebraic stack of finite presentation over $\Spec\cO/\varpi^a$, and the
morphism $\cC_{d,F}\to\cR_{d}$ factors as
$\cC_{d,F} \to \cR_{d,F} \to \cR_d$, with the first morphism being a 
proper surjection, and the second a closed immersion.
\end{thm}
\begin{proof}Points (1), (2), (3) and~(4) follow from Theorems~\ref{main
    result of PR} and~\ref{thm: morphism from C to Rproj} together
  with 
  Lemma~\ref{lem:R versal},
  so it only remains to
  check~(5). For this, it follows from Lemma~\ref{lem: RS is Art-fin
    homogeneous}, Corollary~\ref{cor: 3 and 1 and versal rings implies
    arttriv.}, and Lemma~\ref{lem:we get a sheaf} that we need only
  check that that~$\cR_{d,F}$ admits effective Noetherian versal rings
  at all finite type points.

The $H$-equivariance that was commented upon in Remark~\ref{rem:H equivariance}
implies that the versal morphism $D^{\square,\cC} \to \widehat{\cR}_{d,F,x}$
factors through the quotient $D^\cC$ of $D^{\square,\cC}$.  
The induced morphism $D^\cC \to \widehat{\cR}_{d,F,x}$ is again versal
(as one immediately checks, using the chosen section $D^{\cC} \to
D^{\square,\cC}$).
As we have already observed, the functor $D^{\cC}$ is pro-representable by a
Noetherian ring~$R^{\cC}$.

To complete the verification of~(5), 
we need to check that the morphism
$\Spf R^{\cC}  = D^{\cC} \to\widehat{\cR}_{d,F,x}$ is effective.
To this end, note that by Theorem~\ref{thm: effectivity of height F
  deformation ring} below, the morphism ~$\Spf R^{\cC}\to\widehat{\cR}_{d,x}$ is
induced by a morphism~$\Spec R^{\cC}\to\cR_{d}$.
It remains to check that this
morphism $\Spec R^{\cC}\to\cR_d$ factors through~$\cR_{d,F}$. By Lemma~\ref{lem:scheme-theoretic images}, it suffices to show that
the morphism $\cC_{d,F,R^{\cC}}\to\Spec R^{\cC}$ is
scheme-theoretically dominant. This is precisely the statement of Lemma~\ref{lem:dominance over R}. 

The theorem now follows from Theorem~{\ref{thm:main-intro}},
except that we have only proved that~$\cR_{d,F}$ is locally of finite
presentation. In order to show that it is of finite presentation
over~$\Spec\cO/\varpi^a$, we must show that it is quasi-compact and
quasi-separated. Since~$\cC_{d,F}$ is quasi-compact and the
map~$\cC_{d,F}\to\cR_{d,F}$ is surjective, it follows
from~\cite[\href{http://stacks.math.columbia.edu/tag/050X}{Tag
  050X}]{stacks-project} that~$\cR_{d,F}$ is
quasi-compact. Since~$\cR_{d}$ has affine diagonal,
	by Theorem~\ref{main result of PR}, and since $\cR_{d,F}$
	is a closed substack of $\cR_d$,
	the diagonal of~$\cR_{d,F}$ is also affine.
Thus $\cR_{d,F}$ is quasi-separated, as required.
\end{proof}

We now describe $\cR_{d}$ as an Ind-stack. Note that the inductive
limit in the statement of the following theorem could equivalently be
computed with respect to any cofinal system of~$F$s (see~\cite[\S2]{Emertonformalstacks}).

\begin{thm}\label{thm: R is an ind algebraic stack and satisfies 1}
If $F|F'$, then the natural morphism $\cR_{d,F} \to \cR_{d,F'}$ is a closed
immersion.  Furthermore, the natural morphism
$\varinjlim_F \cR_{d,F}\to \cR_{d}$
{\em (}where, following Remark~{\em \ref{rem:Ind alg stacks can be computed
Zariskiwise}}, the inductive limit is computed as a stack on any of
the Zariski, \'etale, fppf, or fpqc sites{\em )}
is an isomorphism of stacks.
In particular, the stack $\cR_d$ is an Ind-algebraic stack which
satisfies~{\em[1]}.\end{thm}\begin{proof}
	We first show that each of the morphisms
    $\cR_{d,F}\into\cR_{d,F'}$ is a closed immersion; indeed, this follows
    from the fact that in the chain of
    monomorphisms \[\cR_{d,F}\into\cR_{d,F'}\into\cR_d\]
    both the composite and the second morphism are closed immersions.

    It remains to be shown that the natural morphism
    $\varinjlim_F \cR_{d,F} \to \cR_d$ is an isomorphism.  Since
    $\cR_d$ satisfies~[1], it suffices to show that if
    $\Spec A \to \cR_d$ is a morphism with $A$ a Noetherian
    $\cO/\varpi^a$-algebra, then this morphism factors through the
    closed substack $\cR_{d,F}$ for some~$F$. Equivalently, we must
    show that for some~$F$, the closed embedding
    $\cR_{d,F}\times_{\cR_d}\Spec A\into\Spec A$ is an isomorphism. It
    therefore suffices to show that we may find a
    morphism $\Spec B \to \Spec A$ which is scheme-theoretically dominant (equivalently,
so that the corresponding morphism $A \to B$ is injective) such that
the induced morphism $\Spec B \to \cR_d$ factors through~$\cR_{d,F}$. By Proposition~\ref{prop:freeness
	stratification}, we can find a scheme-theoretically dominant
      morphism $\Spec B\to \Spec A$ such that $\Spec B\to \cR_d$
      factors through $\cR_{d,\free}$. Since a free \'etale
      $\varphi$-module contains a free finite height $\varphi$-module,
by Lemma~\ref{lem: existence of a lattice in a projective},
      we see that the morphism $\Spec B\to\cR_{d,\free}$ factors through
      $\cC_{d,F}$ for some~$F$, and thus through~$\cR_{d,F}$, as required.
\end{proof}

\begin{rem}
  \label{rem: R is not an algebraic stack}$\cR_d$ is presumably not an algebraic
  stack. Indeed, since it is Ind-algebraic and satisfies~[1], if it
  were algebraic, it would be locally finite
  dimensional. Since~$\cR_d$ is the inductive limit of its closed
  substacks~$\cR_{d,F}$, it would follow that for each finite type
  point~$x$ of~$\cR_{d}$, there would be a uniform bound on the
  dimension of~$\cR_{d,F}$ at~$x$, independently of~$F$. However, this
  dimension can be computed in terms of the versal rings at~$x$, and  it is presumably straightforward  to use the arguments
  of~\cite{MR2782840} to compute the dimensions of the rings~$R_F$
  and~$R_\cC$ and thereby obtain a contradiction (for example, in the
  case considered in Section~\ref{subsubsec: Galois reps} below, the
  results of~\cite{MR2782840} directly imply that the algebraic stacks
  ~$\cR_{d,E^h}$ are equidimensional, with dimension growing linearly in~$h$).\end{rem}

\subsubsection{An alternative approach}As we now explain, by slightly altering the definitions of $\cC_{d,F}$ and $\cR_d$, 
we could avoid appealing to the descent results 
of~\cite{MR2181808}, without substantially altering our conclusions.

Namely, setting $S := \Spec \cO/\varpi^a$,
we define $\widetilde{\cC}_{d,F} := \pro\bigl((\cC_{d,F})_{| \AfffpS}\bigr)$
and $\widetilde{\cR}_d := \pro\bigl((\cR_{d})_{| \AfffpS}\bigr)$.
Without appealing to the results of~\cite{MR2181808},
we know that these are categories fibred in groupoids.
Proposition~\ref{prop:improving fpqc for A[[u]]} shows that 
$(\cC_{d,F})_{| \AfffpS}$
is in fact an {\em fppf} stack.
Using faithfully flat descent results from rigid analytic geometry,
one can similarly show that
$(\cR_{d})_{| \AfffpS}$
is an {\em fppf} stack.
Furthermore, the arguments of~\cite{MR2562795}, as adapted 
and modified in the present paper, show that 
$(\cC_{d,F})_{| \AfffpS}$ is furthermore
represented by an algebraic stack of finite type
over~$S$.  
Lemma~\ref{lem:limit preserving and pro-categories} then
implies that this same algebraic stack represents $\widetilde{\cC}_{d,F}$,
while Lemma~\ref{lem:restriction to finitely presented algebras}~(2) 
implies that $\widetilde{\cR}_d$ is an \emph{fppf} stack,
which satisfies axiom~[1] by 
Lemma~\ref{lem:limit preserving and pro-categories}.
The arguments of~\cite{MR2562795} are easily adapted
to prove the analogue of Theorem~\ref{main result of PR}
for $\widetilde{\cC}_{d,F}$ and $\widetilde{\cR}_d$.  
Furthermore, the proofs in the subsequent section immediately
adapt to establish the analogue of  Theorem~\ref{thm: R is an ind algebraic stack and satisfies 1}.

Appealing to the results  
of~\cite{MR2181808}, as we do, we in fact prove that $\widetilde{\cC}_{d,F}
= \cC_{d,F}$ and that $\widetilde{\cR}_d = \cR_d.$  However, the primary appeal
of this approach is aesthetic: it allows us to give straightforward and
natural definitions of the stacks that we will study.  In practice, and
in applications, it seems that little would be lost by adopting the 
slightly weaker and more circumlocutious approach described
here.

\subsubsection{Galois representations}\label{subsubsec: Galois reps}
Let~$K/\Qp$ be a finite extension with
residue field~$k$.
We now specialise to the case that~$q=p$, $\varphi(u)=u^p$, and $\cO=\Zp$.
Let~$E$ be the minimal polynomial over~$W(k)$ of a fixed
uniformiser~$\pi$ of~$K$; then we refer to a
$\varphi$-module of height at most~$E^h$ as a \emph{Breuil--Kisin
  module of height at most h}. Fix a uniformiser $\pi$ of $K$, and elements
$\pi_n\in\overline{K}$, $n\ge 0$, such that $\pi_{n+1}^p=\pi_n$ and $\pi_0=\pi$.
Set $K_\infty=\cup_{n\ge 0}K(\pi_n)$, and
$G_{K_\infty}:=\Gal(\overline{K}/K_\infty)$.

The connections between Breuil--Kisin  modules, \'etale $\varphi$-modules,
and Galois representations are as follows. Let $A$ be a
finite\footnote{Recall that an Artinian $\Z_p$-algebra that is finitely
generated as a $\Z_p$-module is necessarily finite as a set, so that the word
``finite'' here can be interpreted either in the commutative algebra sense,
or literally.}  Artinian
$\Zp$-algebra. By \cite[Lem.\
1.2.7]{KisinModularity} (which is based on the results
of~\cite{MR1293971}, and makes no use of the running hypothesis
in~\cite{KisinModularity} that $p\ne 2$),
there is an equivalence of abelian categories between the category of continuous
representations of $G_{K_\infty}$ on finite $A$-modules, and
the category of \'etale $\varphi$-modules with $A$-coefficients. Furthermore, if we write
$T(M)$ for the $A$-module with $G_{K_\infty}$-action corresponding to the
\'etale $\varphi$-module $M$, then $M$ is a free $\gS_A[1/u]$-module of rank $d$ if
and only if $T(M)$ is a free $A$-module of rank $d$.

If $M$ is an \'etale $\varphi$-module, and $V=T(M)$ is the corresponding representation of $G_{K_\infty}$, then we say that $M$
has height at most~$h$ if and only if it there is a Breuil--Kisin
module $\gM$ of height at most~$h$
 with $\gM[1/u]=M$, and we say that
$V$ has height at most~$h$ if and only if $M$ has
height at most~$h$. 

Suppose that $p\cdot A=0$. We say that a continuous representation of $G_K$ on a finite free $A$-module is
\emph{flat} if it arises as the generic fibre of a finite flat group scheme over
$\cO_K$ with an action of $A$. It follows from ~\cite[Thm.\ 1.1.3, Lem.\ 1.2.5]{KisinModularity} (together with the results of~\cite{MR2923180}
in the case that~$p=2$) that restriction to $G_{K_\infty}$ induces an
equivalence of categories between the category of flat representations of $G_K$
on finite free $A$-modules and the category of representations of
$G_{K_\infty}$ of height at most~$1$ on finite free $A$-modules.

The above discussion 
shows that in the preceding context,
we may (somewhat informally) regard $\cR^a$ as a moduli space
of $d$-dimensional continuous representations of $G_{K_{\infty}}$ over $\Z/p^a \Z$.
Furthermore,
if $A$ is a reduced finite-$\Z/p^a$ algebra,
and thus a product of finite fields,
then any Breuil--Kisin module $\gM$ over $A$ is necessarily free.  
Thus an $A$-valued point of $\cR^a$ corresponds to a Galois representation
of height $h$ if and only if it factors through
$\cR^a_{d,E^h}$.
Indeed, we have the following result.

\begin{thm}\label{thm: existence of moduli space of finite flat
    representations}The $\Fpbar$-points of
$\cR^a_{d,E^h}$ naturally biject with isomorphism classes of Galois
representations $G_{K_\infty}\to\GL_d(\Fpbar)$ of height at most~$h$. \end{thm}\begin{proof}
Since $\cR_{d,E^h}$ is a finite
type stack over $\Z/p^a\Z$, any $\Fpbar$-point comes from an $\F$-point for
some finite extension $\F/\Fp$, and by the definition of $\cR_d$ (and the
correspondence between \'etale $\varphi$-modules and continuous $G_{K_\infty}$-representations
explained above), we see that we
need to prove that a morphism $\Spec\F\to\cR_d$ factors through~$\cR_{d,E^h}$ if and
only if the corresponding \'etale $\varphi$-module has height at most~$h$ (possibly
after making a finite extension of scalars).

 By the
definition of $\cC_{d,E^h}$, this latter condition is equivalent to the assertion that the morphism $\Spec\F\to\cR_d$
factors through
the morphism $\cC_{d,E^h}\to\cR_d$, while by
Lemma~\ref{lem: scheme theoretic dominance for Artin proper}, the former condition is equivalent to the assertion
that the fibre $\Spec \F \times_{\cR_d} \cC_{d,E^h}$ be non-empty.  Since
this fibre product is a finite type $\F$-algebraic space (by Theorem~\ref{thm: morphism from C to Rproj}~(1)),
if it is non-empty it contains
a point defined over a finite extension of $\F$.  
Thus these conditions are indeed equivalent, if we allow ourselves to replace
$\F$ by an appropriate finite extension.
\end{proof}

\begin{cor}
  \label{cor:the finite flat version of the main theorem}There is an algebraic
  stack of finite type over $\Spec\Fp$, whose $\Fpbar$-points naturally
  biject with isomorphism classes of finite flat Galois representations $G_K\to\GL_d(\Fpbar)$.
\end{cor}
\begin{proof}
 In view of the equivalence
  between finite flat representations of $G_K$ and representations of
  $G_{K_\infty}$ of height at most~$1$ explained above, this follows immediately from Theorem~\ref{thm: existence of moduli space of
    finite flat representations} (applied in the case
  $a=h=1$).\end{proof}

\subsection{Effectivity of projectivity}\label{subsec: proving effective
  versality}In this final subsection we prove Theorem~\ref{thm:
  effectivity of height F deformation ring}, which was used in the
proof of Theorem~\ref{thm: hypotheses of the main thm hold in our
  setting} in order to show that $\cR_{d,F}$ admits \emph{effective}
Noetherian versal rings at all finite type points. The basic idea of
the proof is to show that the universal family of finite height
$\varphi$-modules on~$\cC_{d,F}$ can be pushed forward to give a
finite height $\varphi$-module on $\Spec R^{\cC}$, and the
corresponding \'etale $\varphi$-module on $\Spec R^{\cC}$ gives the
required promotion of the morphism $\Spf R^{\cC}\to\cR_{d,F}$ to a
morphism $\Spec R^{\cC}\to\cR_{d,F}$. It is, however, not obvious (at
least to us) that this construction gives a \emph{projective} \'etale
$\varphi$-module on $\Spec R^{\cC}$, and that the resulting
morphism $\Spec R^{\cC}\to\cR_{d,F}$ agrees with the given morphism
$\Spf R^{\cC}\to\cR_{d,F}$, and we have to do some work to establish
both of these claims. 

\subsubsection{Interpreting scheme-valued points of $\cC_{d,F}$}
\label{subsubsec:scheme-valued points}
By definition, if $A$ is an $\cO/\varpi^a$-algebra,
then giving an $A$-valued point of $\cC_{d,F}$ is the  same
as giving a rank $d$ projective $\gS_A$-module $\gM$ equipped
with a $\varphi$-module structure of height~$F$.
Suppose more generally that  $T$  is a scheme over $\cO/\varpi^a$,
and let $\widehat{Z}$ denote  the $u$-adic completion of $T\times_{\Spec\Z_p} \Spec\gS,$
i.e.\ the formal scheme obtained by completing this fibre product along its closed
subscheme cut  out  by $u = 0$.
Lemma~\ref{lem: extending phi to A} shows that the endomorphism~ $\varphi$ of $\gS$ induces
a morphism $\varphi_{\widehat{Z}}: \widehat{Z} \to \widehat{Z}$ of formal schemes
over~$T$.

Giving  a  $T$-valued point of $\cC_{d,F}$
is then the same as giving a rank $d$ locally free coherent sheaf $\widehat{\gM}$
on~$\widehat{Z}$,
equipped with an injective morphism $\Phi_{\widehat{\gM}}: \varphi_{\widehat{Z}}^* \widehat{\gM} \to \widehat{\gM}$
whose cokernel is killed by~$F$.

In our discussion below, we will be in the following situation: the scheme $T$
will be projective over a Noetherian affine scheme $\Spec A$.   The morphism $T \to \Spec A$
then induces a morphism $Z: = T\times_{\Spec A} \Spec \gS_A \to  \Spec \gS_A.$
The $u$-adic completion of the source is then equal to the formal scheme
$\widehat{Z}$ introduced above  (although  the map
$Z := T\times_{\Spec A} \Spec \gS_A \to  T \times_{\Spec \Z_p} \Spec  \gS$
need not be an isomorphism in general).
In this situation, the Grothendieck existence theorem shows
that $\widehat{\gM}$ arises as the $u$-adic completion
of a rank~$d$ locally free sheaf $\gM$ on~$Z$.
If $\pi: Z \to \Spec \gS_A$ denotes the projection,
then $\pi_* \gM$ is a coherent sheaf on $\Spec \gS_A$,
which is to say, a finite type $\gS_A$-module.  Since multiplication
by $u$ on $\gM$ is injective (as  $\gM$ is locally free), the
same is true of multiplication by $u$  on $\pi_* \gM$   (since $\pi_*$ is left  exact).

The endomorphism $\varphi_{\widehat{Z}}$ of $\widehat{Z}$ arises as
the $u$-adic completion of an endomorphism $\varphi_{Z}$ of~$Z$; and
in fact $\varphi_Z$  is just the base-change to $Z$ of
the endomorphism $\varphi$ of~$\gS_A$.  This latter endomorphism is finite and
flat by Lemma~\ref{lem: phi faithfully flat}, and thus so is $\varphi_Z$.  Consequently, if $\pi: Z \to \Spec \gS_A$
denotes the projection, we find that there  is a natural isomorphism
$\varphi^* \pi_* \gM \iso \pi_* \varphi_Z^* \gM.$ 

Now the  morphism $\Phi_{\widehat{\gM}}: \varphi_{\widehat{Z}}^* \widehat{\gM} \to \widehat{\gM}$ algebraizes to a morphism
$\Phi_{\gM}: \varphi_Z^* \gM \to \gM,$
and  since the  cokernel of $\Phi_{\widehat{\gM}}$   is supported (set-theoretically)
on the  locus where $u =  0$ (since it is annihilated by~$F$), it coincides
with the cokernel of~$\Phi_{\gM}$.   
Thus this latter  cokernel is also killed by~$F$.

Since $\pi_*$ is left-exact, we may push forward the
short exact sequence
$$0 \to \varphi_Z^* \gM \buildrel \Phi_{\gM} \over
\longrightarrow  \gM \to \gM/\varphi_Z^*\gM \to 0$$
to obtain a left exact sequence
\[0\to \pi_*\varphi_Z^*\gM \to\pi_*\gM\to
\pi_*(\gM/\varphi_{Z}^*\gM),  \]
which, using the previously noted isomorphism, we may rewrite 
as a left exact  sequence
$$0 \to \varphi^* \pi_*\gM \to \pi_*\gM \to \pi_*(\gM/\varphi_Z^*\gM).$$
Since $F$ annihilates $\gM/\varphi^*\gM$, it also
annihilates $\pi_*(\gM/\varphi_Z^*\gM)$,
and thus we see that  $\pi_*\gM$ is naturally endowed with
the structure of a $\varphi$-module of height~$F$ with $A$-coefficients.

\subsubsection{The pushforward of the universal finite height $\varphi$-module}
To simplify notation, we write $R$ for $R^{\cC}$ from now on. We also fix a descending sequence of ideals $I_n$ of $R$
which is cofinal with the sequence $\mathfrak m_R^n$; equivalently,
which determines the $\mathfrak m_R$-adic topology on $R$.
(For the moment, the reader can just imagine that $I_n = \mathfrak m_R^n$;
however, below we will make a different choice of~$I_n$.)
We write $R_n := R/I_n,$ and let $M_n$ denote the universal
\'etale $\varphi$-module over $\Spec R_n$.
We write $\widehat{\gS}_R[1/u]$ for the $\m_R$-adic completion
of $\gS_R[1/u]$, and write $\widehat{M} = \varprojlim_n M_n;$ then $\widehat{M}$ is a 
finite projective $\widehat{\gS}_R[1/u]$-module, and we may
recover the various $M_n$ from $\widehat{M}$
via the natural isomorphism $R_n \otimes_R \widehat{M}
\iso M_n$.

We aim to prove the following proposition.

\begin{prop}\label{prop: pushforward is BK and has correct completion}
We may find 
a  $\varphi$-module $\gM'$ of  height ~$F$ with $R$-coefficients 
such that
  $\widehat{\gM'[1/u]}= \widehat{M}$ {\em (}where the hat on~$\widehat{\gM'[1/u]}$  denotes $\mathfrak m_R$-adic 
completion{\em )}.
\end{prop}

The construction of $\gM'$ involves the pushforward of the universal
finite height $\varphi$-module over~$R$, and we begin by setting up some
notation related to the study of this pushforward. Note that each of 
the morphisms \[\cC_{d,F}\times_{\cR}\Spec R_n\to \Spec R_n\] is in fact
a projective morphism of schemes $\pi_n: X_n\to \Spec R_n$. These
morphisms are compatible as~$n$ varies, and by Grothendieck existence 
give an
algebraization of
 \[\cC_{d,F}\times_{\widehat{\cR}_{d,x}}\Spf R\to \Spf R\] to a
 projective morphism of schemes $\pi: X\to\Spec R$. (In the special case of
 Breuil--Kisin modules of height 1, this
 is~ \cite[Prop.\ 2.1.10]{KisinModularity},
 and the general case can be proved in
exactly the same way (see e.g.\ \cite[Prop.\ 1.3]{MR2827797} for the
case of Breuil--Kisin modules of arbitrary height); the key point is that
$\cC_{d,F}\times_{\widehat{\cR}_{d,x}}\Spf R$ inherits a natural very
ample (formal) line bundle from the affine Grassmannian.)

\begin{lemma}
\label{lem:dominance over R}
The morphism $X \to \Spec R$  is scheme-theoretically dominant.
\end{lemma}
\begin{proof}
Recall that $\Spf R^{\square,\cC}$, which is a formal
closed subscheme of $\Spf R^{\square},$
is defined to be the scheme-theoretic image of the morphism
$$\cC_{\Spf R^{\square}} := \cC_{d,F} \times_{\cR_d} \Spf R^{\square} \to \Spf R^{\square}.$$
We may rephrase this as the conjunction of the following two statements:
the fibre product
$\cC_{\Spf R^{\square, \cC}} := \cC_{d,F} \times_{\cR_d} \Spf R^{\square,\cC}$
(which {\em a priori} is a closed formal subscheme of $\cC_{\Spf R^{\square}}$)
is equal to $\cC_{\Spf R^{\square}}$;
and the projection
\numequation
\label{eqn:C square projection}
\cC_{\Spf R^{\square, \cC}} \to \Spf R^{\square,\cC}
\end{equation}
is scheme-theoretically dominant,
in the sense that  the target of the projection is  equal to its scheme-theoretic image.
It is from this latter statement that we will deduce the lemma.

As observed prior to the statement of Definition~\ref{def:unframed C-defs},
there is a  formal group $G$ over $\Spec \cO/\varpi^a$
 acting freely on $\Spf R^{\square,\cC}$,
and (remembering  that  $R$ is  our abbreviated notation for $R^{\cC}$),
the closed formal subscheme $\Spf R$ of $\Spf R^{\square,\cC}$ is constructed so that
\numequation
\label{eqn:recalling product structure}
\Spf R \times_{\Spec \cO/\varpi^a}  G   \iso \Spf R^{\square,\cC}.
\end{equation}
(Compare Theorem~\ref{thm: existence of height F deformation ring} and its proof.)
We write $\cC_{\Spf R} := \cC_{d,F} \times_{\cR_d} \Spf R.$
What we will show is that
\numequation
\label{eqn:C projection}
\cC_{\Spf R} \to \Spf R
\end{equation}
is scheme-theoretically dominant;
the lemma follows directly from this.
(Indeed, by definition this implies that the scheme-theoretic image
of $X \to \Spec R$ contains $\Spec R_i$  for any  discrete quotient
$R_i$ of $R$, and thus  is equal to $\Spec R$; the reader
should compare the first paragraph of the proof of
Lemma~\ref{lem: formal functions implies formal image equals image formal}.
Indeed, as one sees from a consideration of the  remainder of that  proof,
the theorem on formal functions implies
that the scheme-theoretic dominance of the  morphism 
$X  \to \Spec R$ of {\em schemes}
is equivalent to the scheme-theoretic dominance
of the morphism~\eqref{eqn:C projection} of {\em formal
  schemes}.)

We first note that the $G$-action on $\Spf R^{\square,\cC}$ lifts
to a $G$-action
on~$\cC_{\Spf R^{\square,\cC}}$. 
Indeed,
the morphism $\Spf R^{\square,\cC} \to \cR_d$ which is used 
in the definition of $\cC_{\Spf R^{\square,\cC}}$  as a fibre product
is invariant under the  action of $G$  on its source,
and so this $G$-action indeed lifts to a $G$-action on the fibre product,
i.e.\ on~$\cC_{\Spf R^{\square,\cC}}$. 
The isomorphism \eqref{eqn:recalling product structure} then induces
a corresponding isomorphism
\numequation
\label{eqn:another product}
\cC_{\Spf R^{\cC}} \times_{\Spec \cO/\varpi^a} G \iso \cC_{\Spf R^{\square,\cC}}.
\end{equation}
As we observed in the proof of Theorem~\ref{thm: existence of height F deformation ring},
we have an isomorphism $G = \Spf \cO/\varpi^a[[ (V^*_{\delta})^*]]$
for a certain free $\cO/\varpi^a$-module~$V$.
From this, from~\eqref{eqn:another product}, and from the scheme-theoretic dominance
of~\eqref{eqn:C square projection},
we deduce that~\eqref{eqn:C projection} is scheme-theoretically dominant, as  required.
\end{proof}

We set $S_n := \mathcal O_{X_n}(X_n)$; since $X_n$ is proper
over the affine scheme~$\Spec R_n$, we see that $S_n$ is a finite
$R_n$-algebra.   
Let $R_n'$ denote the image of $R_n$ in $S_n$. By the scheme-theoretic
dominance of~\eqref{eqn:C projection},
the $R_n'$ are again a cofinal system of Artinian quotients 
of $R$ 
and so, replacing each $R_n$ by $R_n'$ (which doesn't
change $X_n$ 
or $S_n$), we may, and indeed do,
 assume that each morphism $R_n \to S_n$ is injective. 

We write  $Y_n := X_n \times_{\Spec \Z_p} \Spec \gS
= X_n\times_{\Spec R_n} \Spec \gS_{R_n}$ (the equality holding
because~$R_n$ is finite), 
and let $\widehat{Y}_n$ denote the completion of $Y_n$ along the 
closed subscheme ``$u = 0$''.  Equivalently,
$\widehat{Y}_n := Y_n \times_{\Spec \gS_{R_n}} \Spf \gS_{R_n},$
where  $\Spf \gS_{R_n}$ is taken with respect to the $u$-adic
topology on~$\gS_{R_n}$.
Each of the schemes $X_n$
is equipped with a tautological projection to $\cC_{d,F}$, 
and we may apply
the discussion of Subsection~\ref{subsubsec:scheme-valued points}
to this morphism 
so as to obtain
a locally free coherent sheaf $\widehat{\gM}_n$ on the formal scheme $\widehat{Y}_n$,
which by 
the Grothendieck existence theorem can be promoted
to a locally free coherent sheaf $\gM_n$ on $Y_n$ itself.

If $\pi_{Y_n}: Y_n \to \Spec \gS_{R_n}$ denotes the projection,
then the discussion of Subsection~\ref{subsubsec:scheme-valued points}
shows that $(\pi_{Y_n})_*\gM_n$ is a $\varphi$-module of height~$F$
with coefficients in $R_n$.   
Since, in particular, $(\pi_{Y_n})_* \gM_n$ is $u$-torsion free,
it  naturally embeds into the  \'etale $\varphi$-module
$\bigl((\pi_{Y_n})_* \gM_n\bigr)[1/u].$
Our next goal is to compute this \'etale $\varphi$-module,
and to compare it to~$M_n$;
we do this  in \eqref{eqn:inverting u again} below.

 
Write 
\begin{multline*}
U_n  := X_n\times_{\Spec \Z_p} \Spec \cO_{\cE} =
X_n\times_{\Spec R_n} \Spec \cO_{\cE, R_n}
\\
= Y_n\times_{\Spec \gS_{R_n}} \Spec \cO_{\cE, R_n} .
\end{multline*}
The projection  $\pi_n: X_n \to \Spec R_n$ 
induces a projection 
$\pi_{U,n}: U_n \to \Spec \cO_{\cE,R_n}.$
The open immersion
\numequation
\label{eqn:inverting u}
\Spec \cO_{\cE,R_n} \to \Spec \gS_{R_n}
\end{equation}
also  induces an open immersion $j_n: U_n \to Y_n.$

In the statement of the following lemma,
we regard the rank $d$  projective $\cO_{\cE,R_n}$-module  $M_n$  as a locally free
coherent sheaf on $\Spec \cO_{\cE,R_n}$.

\begin{lemma}
\label{lem:GAGA compatibility}
There are canonical isomorphisms $j_n^* \gM_n \iso \pi_{U,n}^* M_n,$
compatible with the  actions of~$\varphi$, and with varying~$n$.
\end{lemma}
\begin{proof}
Choose an $\cO_{\cE,R_n}$ basis for $M_n$,
and let $\gN$ denote the free $\gS_{R_n}$-submodule of $M_n$  generated
by this basis.  Then we have an isomorphism
$\varinjlim_i \frac{1}{u^i} \gN \iso M_n.$
Pulling back via $\pi_n$, we obtain
an isomorphism of quasi-coherent sheaves
\numequation
\label{eqn:inductive iso}
\varinjlim_i \frac{1}{u^i} \pi_n^* \gN \iso \pi_n^*M_n
\end{equation}
on~$Y_n$.

We pull the isomorphism~\eqref{eqn:inductive iso} back to $\widehat{Y}_n$.
On the one hand,
if we let $\widehat{\pi^*_n\gN}$  denote the $u$-adic completion of $\gN$,
then  the  left hand side pulls back to $\varinjlim_i \frac{1}{u^i} \widehat{\pi_n^*\gN}.$
On the other hand, the defining relationship between $M_n$ and the $u$-adic
completion  $\widehat{\gM}_n$  of $\gM_n$
(namely, that $M_n$  is defined by the morphism $R_n\to \cR_d$ that
was itself used to construct the morphism $X_n \to \cC_d$,
which in turn defines the coherent sheaf $\widehat{\gM}_n$ on $\widehat{Y}_n$)
shows  that $\pi_n^*M_n$ pulls back to the  sheaf
$\varinjlim_i \frac{1}{u^i} \widehat{\gM}_i$ on~$\widehat{Y}_n$.
Since $(\frac{1}{u^i} \widehat{\pi_n^* \gN})_i$
and $(\frac{1}{u^i} \widehat{\gM}_n)_i$ are inductive systems
of coherent sheaves on $\widehat{Y}_n$
with injective transition maps and having the same inductive limit
on $\widehat{Y}_n$, we see that they are mutually cofinal.
The equivalence of categories provided by Grothendieck existence then shows
that $(\frac{1}{u^i} {\pi_n^* \gN})_i$
and $(\frac{1}{u^i} \gM_n)_i$ are mutually cofinal inductive systems
of coherent  sheaves on $Y_n$ itself.
Returning to~\eqref{eqn:inductive iso},
we find that
$$\varinjlim_i \frac{1}{u^i} \gM_n \iso \pi_n^* M_n.$$
The left-hand side of this isomorphism is naturally identified with $(j_n)_* j_n^* 
\gM_n$.  The natural adjunction between $j_n^*$ and $(j_n)_*$
then shows that
$$j_n^* \gM_n \iso j_n^* \pi_n^* M_n = \pi_{U,n}^* M_n,$$
as  required.
\end{proof}

Since $j_n$ is the base change of the open immersion~\eqref{eqn:inverting u},
we find that
\numequation\label{eqn: pushing pulling open immersion}(\pi_{U,n})_* j_n^* \gM_n = \cO_{\cE,R_n}\otimes_{\gS_{R_n}}(\pi_{Y,n})_*\gM_n
= \bigl( (\pi_{Y,n})_*\gM_n\bigr) [1/u] .\end{equation}

Since $\cO_{\cE,R_n}$  is flat over~$R_n$,
we find that
$(\pi_{U,n})_* \cO_{U_n} =  S_n \otimes_{R_n} \cO_{\cE,R_n} = \cO_{\cE,S_n},$
and so by the projection formula,
\[(\pi_{U,n})_*  \pi_{U,n}^* M_n = \cO_{\cE,S_n} \otimes_{\cO_{\cE,R_n}} M_n 
= S_n \otimes_{R_n} M_n.\]
By Lemma~\ref{lem:GAGA compatibility}, we may rewrite this as
\numequation\label{eqn: pushpull zebra}(\pi_{U,n})_* j_n^* \gM_n = S_n\otimes_{R_n} M_n.\end{equation}
Comparing~\eqref{eqn: pushing pulling open immersion} and~\eqref{eqn:
  pushpull zebra}, we find that
\numequation
\label{eqn:inverting u again}
\bigl( (\pi_{Y,n})_*\gM_n\bigr) [1/u]  = S_n \otimes_{R_n} M_n.
\end{equation}

Since $R_n \to S_n$ is injective, and since $M_n$ is flat over  $R_n$ (being
projective over  the flat $R_n$-algebra~$\cO_{\cE,R_n}$), 
we find that the natural map  $M_n \to S_n\otimes_{R_n}M_n$ 
is injective. By the left exactness of pushforward, we see that
$(\pi_{Y,n})_*\gM_n $ is $u$-torsion free, so the natural map
$(\pi_{Y,n})_*\gM_n \to (\pi_{Y,n})_*\gM_n\bigr) [1/u]=
S_n\otimes_{R_n}M_n$ is also injective.
Consequently, we may define
$$\gM_n' := (\pi_{Y,n})_*\gM_n  \cap  M_n,$$
the intersection taking place in $S_n \otimes_{R_n} M_n.$
This  is a $\varphi$-stable $\gS_{R_n}$-submodule  of~$M_n$,
whose formation is compatible with change in~$n$.

\begin{lemma}
\label{lem:inverting  u}
The inclusion $\gM'_n \hookrightarrow M_n$ induces an isomorphism
$\gM'_n[1/u] = M_n.$
\end{lemma}
\begin{proof}
By virtue of the definition of $\gM_n'$, 
we  find that 
$$\gM_n'[1/u] = \bigl((\pi_{Y,n})_*\gM_n \bigr)[1/u] \cap M_n
= (S_n\otimes_{R_n} M_n) \cap M_n = M_n,$$
the second-to-last  equality following from~\eqref{eqn:inverting u again}.
\end{proof}

We now define 
$\gM' := \varprojlim_n \gM'_n;$
this is an $\gS_R$-module,
which we endow with its projective limit topology
(in which each  $\gM_n'$ is equipped with the discrete topology).

\begin{lemma}
\label{lem:sorting out topologies}
The $\gS_R$-module $\gM'$ is of finite type,
and the projective  limit topology on $\gM'$ 
coincides with the $\gm_R$-adic topology on~$\gM'$.
\end{lemma}

In order to prove this lemma,
we revisit some of the preceding constructions, with the various  $X_n$  replaced
by the projective $R$-scheme~$X$. 
To  be precise, we write $Y := X \times_{\Spec R} \Spec\gS_R$,
and write $\widehat{Y}$ to   denote
the $\mathfrak m_R$-adic completion of~$Y$.
(Note that this notation is {\em not} compatible with
the analogous notation $\widehat{Y}_n$, which denoted 
$u$-adic 
completion.
If we similarly let $\widehat{X}$ denote the $\mathfrak  m_R$-adic completion
of~$X$,
then by the construction of~$X$, there is an  isomorphism
$\widehat{X} \cong \cC_{d,F} \times_{\widehat{\cR}_{d,x}} \Spf R$,
and consequently,
an isomorphism
$\widehat{Y} \cong \cC_{d,F} \times_{\widehat{\cR}_{d,x}} \Spf \gS_R,$
where $\Spf \gS_R$  is  formed with respect to the  $\mathfrak m_R$-adic
topology on~$\gS_R$.)

We may regard $\widehat{Y}$ as the ind-scheme $\varinjlim_n  Y_n$,
and the compatible (in $n$)
collection of coherent sheaves $\gM_n$ on the various $Y_n$
give rise to a coherent  sheaf $\widehat{\gM}$ on  the formal scheme~$\widehat{Y}$.
By the Grothendieck existence  theorem, this coherent sheaf
arises as the completion of a coherent  sheaf $\gM$  on the projective $\gS_R$-scheme~$Y$.

If we write  $\pi_Y$  for  the   projection  $Y\to \Spec \gS_R$,
then the theorem  on formal  functions shows that 
$(\pi_Y)_* \gM$ is a finite type  $\gS_R$-module,
and that  there is an isomorphism
\numequation
\label{eqn:topological iso}
(\pi_Y)_* \gM  \iso \varprojlim_n 
(\pi_{Y,n})_* \gM_n,
\end{equation}
which furthermore identifies the $\gm_R$-adic topology on the
source with the projective limit topology on the target.

\begin{proof}[Proof of Lemma~\ref{lem:sorting out topologies}]
  Since each $\gM_n'$ is an $\gS_{R_n}$-submodule of
  $(\pi_{Y,n})_* \gM_n$, it follows from ~\eqref{eqn:topological iso}
  and the left exactness of inverse limits 
  that $\gM'$ ($ = \varprojlim_n \gM_n'$) is an $\gS_R$-submodule of
  $(\pi_Y)_*\gM$.  Since~$\gS_R$ is Noetherian, and~$\gM'$ is an
  $\gS_R$-submodule of~ the finite type $\gS_R$-module $(\pi_Y)_*\gM$,
  we see that~$\gM'$ is of finite type.  The projective limit topology
  on $\gM'$ coincides with the topology induced by the projective
  limit topology on $(\pi_Y)_*\gM$ induced by the right-hand side of
  the isomorphism~\eqref{eqn:topological iso}.  Thus it also coincides
  with the topology induced by the $\gm_R$-adic topology on
  $(\pi_Y)_*\gM$, which by the Artin--Rees lemma is equal to the
  $\gm_R$-adic topology on~$\gM'$.
\end{proof}

\begin{proof}[Proof  of  Proposition~\ref{prop: pushforward is BK and has correct completion}]
We have already seen in Lemma~\ref{lem:sorting out topologies} 
that $\gM'$  is a finite  type $\gS_R$-module,
and Lemmas~\ref{lem:inverting u} 
and~\ref{lem:sorting out topologies}
together show that
$$\widehat{\gM'[1/u]} = \varprojlim_n \, (\gM'/\gm_R^n \gM')[1/u] 
= \varprojlim_n \gM'_n[1/u] = \varprojlim_n M_n = \widehat{M}.$$
Since each $\gM_n'$  embeds into $M_n = \gM_n'[1/u]$,
we see that  each  $\gM_n'$ is $u$-torsion free; thus $\gM'$ is also  $u$-torsion free.

It remains to show that the $\varphi$-module structure on $\gM'$
makes it of height~$F$.
The endomorphism $\varphi$ of $\gS_R$ is finite flat,
and thus
$\varphi^*\gM'=\varprojlim_n
\varphi^*\gM_n'$. 
It hence suffices to show that each~$\gM_n'$ is of
height~$F$. 

We have already seen that each of the pushforwards  $(\pi_{Y_n})_*\gM_n$ is of height~$F$.
Now let $x\in \gM_n' = (\pi_{Y_n})_* \gM_n \cap M_n$ be arbitrary.
Then $Fx\in F(\pi_{Y,n})_*\gM_n$, so
since $(\pi_{Y,n})_*\gM_n$ has height ~$F$, we see that there is some
$y\in\varphi^*(\pi_{Y,n})_*\gM_n$ with
$\Phi_{(\pi_{Y,n})_*\gM_n}(y)=Fx$. Since
$\Phi_{M_n}:\varphi^*M_n\to M_n$ is a bijection, we have
$y\in \varphi^*M_n$, so
$y\in \varphi^*M_n\cap\varphi^*(\pi_{Y,n})_*\gM_n$. Since $\varphi$ is
flat, we have
$\varphi^*M_n\cap\varphi^*(\pi_{Y,n})_*\gM_n=\varphi^*\gM_n'$, so
$y\in\varphi^*\gM_n'$, and $\gM_n'$ is of height~$F$, as required.
\end{proof}

\subsubsection{Descent of projectivity} 

We now come to the main result of this subsection. 
  \begin{thm}
  \label{thm: effectivity of height F deformation ring}
  There is a projective \'etale $\varphi$-module $M$ over 
  $R^{\cC}$ whose $\mathfrak m_{R^{\cC}}$-adic completion 
  is isomorphic to the base change to $R^{\cC}$ of the
  (formal) universal \'etale $\varphi$-module over~$R^{\square}$.
\end{thm}
\begin{proof}
  We write $M:= \gM'[1/u]$, where $\gM'$ is the $\varphi$-module of
  height~$F$ over $R^{\cC}$ constructed in Proposition~\ref{prop:
    pushforward is BK and has correct completion}. By construction,
  $M$ is an \'etale $\varphi$-module with $R^{\cC}$-coefficients, and 
  the $\mathfrak m_{R^{\cC}}$-adic completion of $M$ is a finite
  projective module over $\widehat{\gS}_{R^{\cC}}[1/u]$, and
  corresponds to the morphism $\Spf R^{\cC}\to\cR_d$.
  Theorem~\ref{thm:effectivity of projectivity} below then shows that
  $M$ itself is a projective \'etale $\varphi$-module, as required.
\end{proof}

\begin{theorem}
	\label{thm:effectivity of projectivity}
	Let $R$ be 
a complete Noetherian local $\cO/\varpi^a$-algebra $R$ with maximal
ideal~$\m$, 
let $M$ be an
      \'etale $\varphi$-module with~$R$-coefficients,
and suppose that the $\mathfrak m$-adic completion $\widehat{M}$
is projective, or equivalently, free
{\em (}over $\widehat{R((u))}$, the $\mathfrak m$-adic completion
of $R((u))${\em )}.
Then $M$ itself is projective {\em (}over $R((u))${\em )}.
\end{theorem}

The remainder of this subsection is devoted to the proof of this theorem.
However,
before giving the proof in detail, we give an outline of it in the case
that $R = K[[T_1,\ldots,T_n]]$, where $K$ is a field of characteristic $p$.
In this case, the ring
$R((u))$ may be regarded as the ring of bounded holomorphic functions
on the open unit $n$-dimensional polydisk $\mathbb D^{n}$ over 
the complete discretely valued field $K((u))$,  
and its $\mathfrak m$-adic completion
$\widehat{R((u))}$ may be regarded as the formal completion of
the local ring $\mathcal O_{\mathbb D^{n},0}$, where $0$ denotes 
the origin of $\mathbb D^n$.  The extension
$\mathcal O_{\mathbb D^n,0} \to \widehat{R((u))}$ is then faithfully
flat, and since $\widehat{M} \iso \widehat{R((u))} \otimes_{R((u))} M$
is projective, the same is true
of $\mathcal O_{\mathbb D^n,0}\otimes_{R((u))} M.$

Writing the stalk as a direct limit of rings of holomorphic functions
on a nested sequence of polydisks centred at $0$ of shrinking radius, 
we find that the restriction of $M$ to one of these polydisks is 
projective.  The \'etale $\varphi$-module structure on $M$
then allows us to employ Frobenius amplification (``Dwork's trick'')
to show that $M$ itself gives rise to a projective module
over the ring $\mathcal O(\mathbb D^n)$ of holomorphic functions
on $\mathbb D^n$ --- i.e.\ we find that
$\mathcal O(\mathbb D^n)\otimes_{R((u))}
M$ is projective.  Finally, the inclusion $R((u)) \hookrightarrow
\mathcal O(\mathbb D^n)$ of bounded holomorphic functions in all (i.e.\
not necessarily bounded) holomorphic functions is faithfully flat,
and so $M$ itself is projective over $R((u))$.

The proof in the general case follows the same outline: we regard 
$R((u))$ as being the ring of bounded holomorphic functions on a closed 
analytic subvariety of an open unit polydisk and make the
same Frobenius amplification argument. We pass
from the case of $k$-algebras to $\cO/\varpi^a$-algebras via the usual
graded techniques.

We now fill in the details of the preceding sketch.

\begin{proof}[Proof of Theorem~\ref{thm:effectivity of projectivity}]
	Assume to begin with that $R$ is a complete Noetherian
	local $k$-algebra; we may
	then write $R = K[[T_1,\ldots,T_n]]/I,$ for some field extension
	$K$ of $k$, and some ideal $I$.
	For any integer $m\geq 0$, we write
	$$A_m := \reallywidehat{R[x_1,\dots,x_n,u]/(u^{p^m} x_1 - T_1,\dots,
		u^{p^m} x_n - T_n)} [1/u],$$
	and
	$$B_m := \reallywidehat{R[x_1,\dots,x_n,u]/(u x_1 - T_1^{p^m},\dots,
		u x_n - T_n^{p^m})}[1/u];$$
		in both cases, the hat indicates $u$-adic completion.
		There are natural morphisms
		of $R((u))$-algebras
		\numequation\label{eqn: Am to Am plus one}A_m \hookrightarrow A_{m+1},\end{equation}
		defined by mapping $x_i$ in the source to $u^{p^m(p-1)} x_i$ 
		in the target,
		and
		$$B_{m+1} \hookrightarrow B_m,$$
		defined by mapping
		$x_i$ in the source to $u^{p-1} x_i^p$ in the target;
		Lemmas~\ref{lem:more technical flatness}
		and~\ref{lem:variant flatness} below (together
		with Remark~\ref{rem:when a' is a power of a})
                show that all these morphisms are flat.
		There are also evident isomorphisms of $R((u))$-algebras
		\numequation
		\label{eqn:first frob shift}
		R((u))\otimes_{\varphi,R((u))} A_m \iso A_{m+1},
	\end{equation}
	defined by mapping $1\otimes x_i$ in the source to $x_i$ in
	the target,
		as well as isomorphisms
		\numequation
		\label{eqn:second frob shift}
		R((u))\otimes_{\varphi,R((u))} B_{m+1} \iso B_{m},
	\end{equation}
	defined by mapping $1\otimes x_i$ in the source
	to $x_i^p$ in the target.   (It is less
	evident that the morphisms~(\ref{eqn:second frob shift}) are isomorphisms,
	but this follows from Lemma~\ref{lem:isomorphism} below.)
		Note also that $A_0 = B_0$.

		We write $A := \varinjlim_m A_m.$
		Being the direct limit of a sequence
		of Noetherian rings with respect to flat transition
		morphisms, we see immediately that $A$ is a coherent ring.
		In fact, something stronger is true:
		$A$ is a Noetherian local $R((u))$-algebra,
		and there is an isomorphism
		$$\widehat{A} \iso \widehat{R((u))},$$
		where now the hats indicate $\mathfrak m$-adic
		completion. 
		To see this, it suffices to consider the case when 
		the ideal $I$ is zero.  (Standard
	        arguments with adic completions of finite type
         	modules over Noetherian rings, and with direct limits,
       	        show that the constructions
		of the various rings $A_m$, $A$, $\widehat{A}$,
		and $\widehat{R((u))}$ are compatible with
		the passage from $R$ to $R/J$, for any ideal $J$ of $R$.
		Now apply this observation to the ideal $I$ in the
		ring $K[[T_1,\ldots,T_n]]$.) 
		In this case,
		we see that $A$ is the ring of germs of holomorphic
		functions at the origin of the unit $n$-dimensional polydisk
		over 
		the complete discretely valued field $K((u))$,
		and this ring is well-known to be Noetherian,
		and to have the formal power series ring $K((u))[[T_1,
		\ldots, T_n]]$ as its $\mathfrak m$-adic completion
		(see~\cite[Prop.~7,\ \S7.3.2]{MR746961}). 

		The discussion of the preceding paragraph shows that
		the embedding
		$A \hookrightarrow \widehat{R((u))}$
		is faithfully flat, and since 
		$\widehat{M}$ (which, by the Artin--Rees lemma,
		may be identified with $\widehat{R((u))}\otimes_{R((u))} M$)
		is projective, we conclude that 
		$$A \otimes_{R((u))} M$$
		is also projective over $A$.
		This in turn implies that 
		$$A_m\otimes_{R((u))} M$$
		is projective over~$A_m$, for some sufficiently large
                value of $m$. (Indeed, this projectivity is witnessed
                by a split surjection from a finite free $A$-module,
                and both this surjection, and a splitting of the
                surjection, can be descended to some~$A_m$.)

		Composing the inclusion $A_{m-1} \to R((u))\otimes_{\varphi,
			R((u))} A_{m-1}$, defined via $a \mapsto 1\otimes a$,
		with the isomorphism ~(\ref{eqn:first frob shift}),
		we may regard $A_m$ as a faithfully flat
                $A_{m-1}$-algebra. (Note that this is not the same as
                the $R((u))$-linear map $A_{m-1}\to A_m$ considered
                in~(\ref{eqn: Am to Am plus one}).) We then compute that
		\begin{multline*}
		A_m\otimes_{A_{m-1}} (A_{m-1}\otimes_{R((u))} M)
		\iso
		A_m\otimes_{\varphi,R((u))} M
		\\
		\iso
		A_m\otimes_{R((u))}(R((u))\otimes_{\varphi,R((u))} M)
		\iso A_m\otimes_{R((u))} M
	\end{multline*}
		(the last isomorphism following from the \'etale
		$\varphi$-module structure on $M$).
		Since the property of being finitely generated projective
		can be detected after a faithfully flat base-change,
		we find that $A_{m-1}\otimes_{R((u))} M$
		is projective.
		Continuing via descending induction,
		we conclude that
		$A_0 \otimes_{R((u))} M$
		is a projective $A_0$-module, or equivalently,
		a projective $B_0$-module. 
	A completely analogous argument,
	taking into account~(\ref{eqn:second frob shift}),
        then shows that $B_m\otimes_{R((u))} M$
		is projective for each $m \geq 0.$
		
		Write $B := \varprojlim_m B_m$. 
		Since the transition morphisms in this projective
		limit are flat morphisms of Noetherian Banach
		algebras over $K((u))$, we find that 
		that $B$ is a Fr\'echet--Stein algebra,
		in the sense of~\cite{MR1990669}. As a consequence, any finitely presented $B$-module 
		$N$ is coadmissible over $B$ (by~\cite[Cor.\ 3.4]{MR1990669}), so that the natural
		morphism $N \to \varprojlim_m B_m\otimes_B N$
		is an isomorphism.
		In particular, if $X$ is a finitely generated
		(and hence finitely presented)
		$R((u))$-module, then $B\otimes_{R((u))} X$
		is a finitely presented $B$-module, and hence
		there is an isomorphism
		\numequation
		\label{eqn:f.g. iso}
		B\otimes_{R((u))} X \iso \varprojlim_m B_m\otimes_{R((u))} X.
\end{equation}

		If we take $X$ to be a finitely generated ideal $I$
		of $R((u))$, then~(\ref{eqn:f.g. iso}) shows that the morphism
		\numequation
		\label{eqn:flatness test}
		B\otimes_{R((u))} I \to B
	\end{equation}
		may be written as the
		projective limit of morphisms  
		$B_m\otimes_{R((u))} I \to B_m$.
		By Lemma~\ref{lem:flatness} below, each $B_m$ is flat over $R((u))$,
		 so we find that each of these latter morphisms is injective.
		Thus so is the morphism~(\ref{eqn:flatness test}).  It 
		follows that $B$ is flat over $R((u))$.

		In fact $B$ is faithfully flat over $R((u))$.
		To show this, it suffices (given the flatness
		that we have already proved, and the fact
		that flat maps satisfy {\em going down}) to show
		that each maximal ideal of $R((u))$ is obtained
		via restriction from a maximal ideal of $B$.
		Any such maximal ideal is the kernel of a surjection
		\numequation
		\label{eqn:surjection}
		R((u)) \to L,
	\end{equation}
		where $L$ is a finite extension
		of $K((u))$
		(see the proof of~\cite[Lem.~7.1.9]{MR1383213} for a proof
	        of this fact),
		and each $T_i$ maps to an element $a_i \in L$
		satisfying $|a_i| < 1$ (i.e.\ each $a_i$ lies
		in the maximal ideal of the ring of integers of $L$).
		If we choose $m$ so that $|a_i| \leq |u|^{1/p^m}$
		for each $a_i$,
		then the surjection~(\ref{eqn:surjection})
		extends to a surjection
		$B_m \to L$.  
		(This extension of the surjection $R((u)) \to L$
		to a surjection $B_m \to L$ for some sufficiently
		large value of $m$ is also explained carefully
		in the proof of~\cite[Lem.~7.1.9]{MR1383213}.)
		The kernel of the
		composite $B \to B_m \to L$ 
		is then a maximal ideal of $B$ that restricts to 
		the kernel of~(\ref{eqn:surjection}).

		Since $M$ is finitely generated over $R((u))$,
		we obtain, from~(\ref{eqn:f.g. iso}), an isomorphism
		$$B\otimes_{R((u))} M \iso \varprojlim_m B_m\otimes_{R((u))}M;$$
		since each $B_m\otimes_{R((u))}M$ is projective
		over $B_m$, Lemma~\ref{lem:projective limit of projectives} below
		shows that $B\otimes_{R((u))} M$ 
		is projective over $B$.
	Since $R((u)) \to B$ is faithfully flat,
		we find that $M$ is projective, as required.

		Finally, consider the general case of the theorem,
		in which $R$ is assumed to be a complete Noetherian
		local $\cO/\varpi^a$-algebra,
		without $a$ necessarily equalling~$1$. 
		Endow each of $\cO/\varpi^a$,
		$R$, $R((u))$, $\widehat{R((u))}$, $M$, and $\widehat{M}$
		with its $\varpi$-adic filtration,
		and let $\Gr^{\bullet} \cO/\varpi^a$, etc., 
		denote the corresponding associated graded object.
		We find that
	       	$\Gr^{\bullet} \cO/\varpi^a =
		k[\varepsilon]/(\varepsilon^a),$
		that $\Gr^{\bullet} R$ is a complete local Noetherian
		algebra over $k[\varepsilon]/(\varepsilon^a),$
		that $\Gr^{\bullet} R((u)) = (\Gr^{\bullet} R)((u)),$
		and that $\Gr^{\bullet} M$ is an \'etale
		$\varphi$-module over $\Gr^{\bullet} R$.
		
		Since $\mathfrak m$-adic completion of finitely
		generated modules over a 
		Noetherian local ring (such as $R((u))$) is exact, 
		we also see that the $\mathfrak m$-adic completion 
		of $\Gr^{\bullet} R((u))$ is naturally isomorphic 
		to $\Gr^{\bullet} \widehat{R((u))}$,
		and that the $\mathfrak m$-adic completion 
		of $\Gr^{\bullet} M$ is naturally isomorphic 
		to $\Gr^{\bullet} \widehat{M}$.
		The assumption that $\widehat{M}$ is finitely
		generated and projective over $\widehat{R((u))}$
		then implies that
		$\widehat{\Gr^{\bullet} M} = \Gr^{\bullet} \widehat{M}$
		is finitely generated and projective over
		$\widehat{\Gr^{\bullet} R((u))} =
		\Gr^{\bullet} \widehat{R((u))} .$
		(This is easily seen if one uses the equivalence
		between a module being finitely generated and projective
		and being a direct summand of a finite rank free module.)
		The case of the theorem already proved then shows 
		that $\Gr^{\bullet} M$ is finitely generated and projective
		over $\Gr^{\bullet} R((u))$,
		or, equivalently, finitely generated and flat 
		over $\Gr^{\bullet} R((u))$.
		Lemma~\ref{lem:graded flat}  (applied with $R = \cO/\varpi^a$
		and $A = R((u))$)
		shows that $M$ is finitely generated and flat  --- and thus 
		projective --- over $R((u)).$
\end{proof}

The following lemmas, which were used in the proof
of the preceding theorem, are presumably well-known to experts, 
but we include proofs, for lack of a reference.

\begin{lemma}
	\label{lem:tensor with coadmissible}
	If $B$ is a commutative
	Fr\'echet--Stein algebra {\em (}over a complete
	discretely valued field $K${\em )}, if $M$ is a finitely presented
	$B$-module, and if $N$ is a coadmissible $B$-module,
	then $M\otimes_B N$ is again a coadmissible $B$-module.
\end{lemma}
\begin{proof}
	If we choose a presentation $B^r \to B^s \to M \to 0,$
	then we obtain a right exact sequence
	$$ N^r \to N^s \to M\otimes_B N \to 0.$$
	Since $N$ is coadmissible, so are each of $N^r$ and $N^s$,
	and thus so is $M\otimes_B N,$ being the cokernel of a morphism
	between coadmissible $B$-modules.
\end{proof}

\begin{lemma}
	\label{lem:projective limit of projectives}
	If $B$ is a commutative Fr\'echet--Stein algebra over a 
	complete discretely valued field $K$,
	say $B \iso \varprojlim_n B_n,$ where each $B_n$
	is a Noetherian Banach $K$-algebra, with the transition morphisms
	being flat, and if $M$ is a finitely presented $B$-module
	with the property that each tensor product $B_n\otimes_B M$
	is a projective $B_n$-module, then $M$ is a projective
	$B$-module.
\end{lemma}
\begin{proof}
	Since $M$ is a finitely presented $B$-module,
	it is furthermore projective if and only it is flat.
	To show that $M$ is flat, it suffices to show
	that for each finitely generated ideal $I\subseteq B$,
	the induced morphism
	$M\otimes_B I \to M$ is injective.

	Since $I$ is finitely generated, it is the image of a morphism
	$B^r \to B$, for some $r \geq 0,$ and thus is coadmissible.
	Lemma~\ref{lem:tensor with coadmissible} then shows that $M\otimes_B I$
	is coadmissible.  The $B$-module $M$ itself is also coadmissible
	(being finitely presented).
	Thus the morphism
	\numequation
	\label{eqn:flat test}
	M\otimes_B I \to M
\end{equation}
       	may be obtained as the 
	projective limit of the morphisms
	$$B_n\otimes_B \otimes M \otimes_B I \to B_n\otimes_B M.$$
	We may rewrite each of these morphism as
	\numequation
	\label{eqn:injective map}
	(B_n\otimes_B \otimes M) \otimes_{B_n} (B_n\otimes_B I)
	\to B_n\otimes_B M.
\end{equation}
	Since $B_n$ is flat over $B$, we see that
	the inclusion of $I$ in $B$ induces a sequence of injections
	$B_n\otimes_B I \hookrightarrow B_n.$
	Since $B_n\otimes_B M$ is projective, and thus flat,
	over $B_n$, the morphisms~(\ref{eqn:injective map})
	are then also injective. 
	Hence so is their projective limit~(\ref{eqn:flat test}).
\end{proof}

\begin{lemma}
	\label{lem:completing}
	Let $A$ be a ring, and $a$ an element of $A$.
	If $M \to N$ is a morphism of $A$-modules whose kernel
	and cokernel are each
	annihilated by some power of $a$,
	then the induced morphism 
	$\widehat{M}[1/a] \to \widehat{N}[1/a]$
	{\em (}where $\widehat{\phantom{A}}$ denotes $a$-adic completion{\em )}
	is an isomorphism.
\end{lemma}
\begin{proof}
	This is standard, and follows easily from the definitions.\end{proof}

\begin{lemma}
	\label{lem:flatness}
	Let $A$ be a Noetherian ring, and let $a, b_1,
	\ldots, b_m$ be elements of~$A$.
	If we write $B :=  A[x_1,\dots,x_m]/(ax_1 - b_1,\dots,ax_m-b_m),$
	then the natural map $\widehat{A}[1/a] \to \widehat{B}[1/a]$
	{\em (}where $\widehat{\phantom{A}}$ denotes $a$-adic completion{\em )}
	is flat. \end{lemma}
\begin{proof}  We begin with the 
  case $m=1$, where we write $x,y$ for $x_1,y_1$, and $b$ for $b_1$.
	Note that each morphism in the sequence of natural morphisms
	$$\widehat{A} \to \widehat{A}[x]/(ax - b) \to \widehat{A}[1/a]$$
	becomes an isomorphism after inverting~$a$. 
	Indeed, this is evidently the case for their composite,
	and it is also evident that the first morphism becomes surjective
	after inverting~$a$.
	Similarly, each morphism in the sequence of natural morphisms
	$$B := A[x]/(ax-b) \to \widehat{A}\otimes_A B = \widehat{A}[x]/(ax-b)
	\to \widehat{B}$$
	becomes an isomorphism after passing to $a$-adic completions.
	Thus, in order to prove the lemma in the case $m=1$, it suffices to note that
	the natural morphism from
	$\widehat{A}[x]/(ax-b)$ to its $a$-adic completion is flat,
	as follows from the Artin--Rees lemma (and the fact
	that $\widehat{A}$ is Noetherian, as $A$ is). (Note that by
        the discussion about, the natural map
        $\widehat{A}[1/a]\to\widehat{B}[1/a]$ is obtained from this
        map by inverting~$a$.)

The general case follows by induction on~$m$. Indeed, writing
$$C=A[x_1,\dots,x_{m-1}]/(ax_1 - b_1,\dots,ax_{m-1}-b_{m-1}),$$
we can factor
$\widehat{A}[1/a] \to \widehat{B}[1/a]$ as $\widehat{A}[1/a] \to
  \widehat{C}[1/a] \to \widehat{B}[1/a],$ with the first map being flat by the
inductive hypothesis, and the second being flat by the case~$m=1$, as $B=C[x_m]/(ax_m-b_m)$.
\end{proof}

The following result is a somewhat technical
modification  of the preceding lemma.

\begin{lemma}
	\label{lem:more technical flatness}
	Let $A$ be a Noetherian ring, and let $a$, $a'$, $b_1,\ldots,b_m$,
	$b'_1,\ldots,b_m'$ be elements of~$A$.
	If we write $B := A[x_1,\dots,x_m]/(ax_1 - b_1b_1',\dots,ax_m-b_mb_m')$
	and $C := A[y_1,\dots,y_m]/(aa'y_1-b_1,\dots,aa'y_m-b_m),$
	then there is a morphism of $A$-algebras $B \to C$ 
	defined by mapping each $x_i$ to $a'b'_i y_i$, 
	and the induced morphism $\widehat{B}[1/aa'] \to \widehat{C}[1/aa']$
	{\em (}where $\widehat{\phantom{A}}$ denotes
	$aa'$-adic completion{\em )}
	is flat.
\end{lemma}
\begin{proof}
	We factor the morphism $B \to C$ as
	\begin{multline*}
	B \to B[y_1,\dots,y_m]/(aa'y_1-b_1,\dots,aa'y_m-b_m) \\
	\to B[y_1,\dots,y_m]/(aa'y_1-b_1,x_1-a'b_1'y_1,\dots,aa'y_m-b_m,x_m-a'b_m'y_m) = C.
\end{multline*}
	The second morphism is surjective, and its kernel, which
	is the ideal $(x_1-a'b_1'y_1,\dots,x_n-a'b_m'y_m),$
	is annihilated by~$a$.
	Lemma~\ref{lem:completing} shows that this second morphism
	becomes an isomorphism after $aa'$-adically completing and
	then inverting $aa'$.  Thus it suffices to show that
	the first morphism becomes flat after $aa'$-adically completing
	and inverting $aa'$; this follows from
	Lemma~\ref{lem:flatness}.
\end{proof}

\begin{remark}
	\label{rem:when a' is a power of a}
	We note,
	in the context of the preceding lemma,
	that if $(a')^r = a^s$ for some $r,s \geq 1$, then $aa'$-adically
	completing is the same as $a$-adically completing,
	and inverting $aa'$ is the same as inverting $a$.
\end{remark}

We also have the following variations on the preceding results.

\begin{lemma}
	\label{lem:isomorphism}
	Let $A$ be a Noetherian ring, and let $a, b_1,\ldots,b_m$
	be elements of~$A$,
	and let $n$ be a positive integer.
	If we write $B :=  A[x_1,\dots,x_n]/(a^nx_1 - b_1^n,\dots,a^nx_n-b_m^n)$
	and $C = A[y_1,\dots,y_n]/(ay_1-b_1,\dots ay_n-b_m),$
	then the morphism of $A$-algebras $B \to C$ defined by $x_i \mapsto y_i^n$
	induces an isomorphism
	$\widehat{B}[1/a] \iso \widehat{C}[1/a]$
	{\em (}where $\widehat{\phantom{A}}$ denotes $a$-adic completion{\em )}.
\end{lemma}
\begin{proof}
	The morphism $B\to C$ is finite: $C$ is generated as a $B$-module
	by the various monomials~$y_1^{e_1}\cdots y_m^{e_m}$, for $1\le e_i\le n-1$. The image of $y_i^j$ in the cokernel 
	of this morphism is annihilated by~$a^j$, so the entire cokernel
	is annihilated by $a^{m(n-1)}$.
	
	Note that the morphisms $A[1/a]\to B[1/a] \to C[1/a]$ are all
	isomorphisms, so that the kernel of the morphism $B\to C$
	is contained in the kernel of the morphism $B \to B[1/a]$.
	Each element of this kernel is annihilated by some power of $a$.
	Since this kernel is finitely generated (as $B$ is Noetherian),
	we see that this entire kernel is annihilated by some power of $a$.
	Combining this with the conclusion of the preceding paragraph,
	and with Lemma~\ref{lem:completing}, establishes
	the lemma.
\end{proof}

\begin{lemma}
	\label{lem:variant flatness}
	Let $A$ be a Noetherian ring,
	let $a, b_1,\ldots,b_m$ be elements of~$A$,
	and let $n$ be a positive integer.
	If we write $B :=  A[x_1,\dots,x_m]/\bigl(ax_1 - b_1^n,\dots,ax_m-b_m^n\bigr)$
	and $C = A[y_1,\dots,y_m]/(ay_1-b_1,\dots,ay_m-b_m),$
	then the morphism of $A$-algebras $B \to C$ defined by
	$x_i \mapsto a^{n-1}y_i^n$
	induces a flat morphism
	$\widehat{B}[1/a] \iso \widehat{C}[1/a]$
	{\em (}where $\widehat{\phantom{A}}$ denotes $a$-adic completion{\em )}.
\end{lemma}
\begin{proof}
	We factor the morphism $B\to C$ as
	\begin{multline*}
	B=A[x_1,\dots,x_m]/\bigl(ax_1 - b_1^n,\dots,ax_m-b_m^n\bigr)
	\\
	 \to A[t_1,\dots,t_m]/(a^nt_1-b_1^n,\dots,a^nt_m-b_m^n)
	 \\ \to A[y_1,\dots,y_m]/(ay_1 - b_1,\dots,ay_m-b_m)=C,
\end{multline*}
	where the first morphism is defined by $x_i \mapsto a^{n-1}t_i,$
	and the second morphism is defined by $t_i \mapsto y_i^n$.
	The present lemma then follows from
	Lemmas~\ref{lem:more technical flatness} and~\ref{lem:isomorphism}
	(taking into account Remark~\ref{rem:when a' is a power of a}).
\end{proof}

Let $R$ be an Artinian local ring, with maximal ideal $I$ and residue field~$k$. If $M$ is
an $R$-module, then we let $\Gr^{\bullet} M$ denote the graded
$k$-vector space associated to the $I$-adic filtration on $M$ (so
$\Gr^i M := I^iM/I^{i+1} M$).

\begin{lemma}
	\label{lem:graded flat}
If $A$ is an $R$-algebra, then an $A$-module $M$
is {\em (}faithfully{\em )} flat over $A$
if and only if $\Gr^{\bullet} M$ 
is {\em (}faithfully{\em )} flat over $\Gr^{\bullet} A.$
Furthermore, if any of these conditions holds,
then the natural morphism
$\Gr^{\bullet} A \otimes_{\Gr^0 A} \Gr^0 M \to
\Gr^{\bullet} M$
is an isomorphism.
\end{lemma}
\begin{proof}
	If $M$ is flat over $A$, then considering the
        result of tensoring $M$
	by the various short exact sequences
	$$0 \to I^nA \to I^mA \to I^mA/I^nA \to 0$$
	we find that the natural
	morphism $\Gr^i A \otimes_{\Gr^0 A} \Gr^0 M \to \Gr^i M$
	is an isomorphism, for each $i$.
	This proves the final assertion of the lemma.
	Furthermore, since (faithful) flatness is preserved
	under base-change, we see first that $\Gr^0 M$ is flat
	over $\Gr^0 A$ (and faithfully flat if $M$ is faithfully
	flat over $A$), and then (using the result already proved)
	that $\Gr^{\bullet} M$ is flat over $\Gr^{\bullet} A$
	(and faithfully flat if $M$ is).
	
	It is not quite as obvious that flatness of $\Gr^{\bullet} M$
	over $\Gr^{\bullet} A$ implies the flatness of $M$ over $A$,
	but this is a standard fact in commutative algebra;
	e.g.~it follows
	from~\cite[\href{http://stacks.math.columbia.edu/tag/0AS8}{Tag 0AS8}]{stacks-project}. 
	(If $i \geq 0,$ then base-changing via the map
	$\Gr^{\bullet} A \to
	\Gr^{\leq i} A := A/IA \oplus I/I^2 \oplus \cdots \oplus I^i/I^{i+1},$
	we find that $\Gr^{\leq i} M := M/IM \oplus \cdots \oplus I^iM/I^{i+1}M$
	is flat over $\Gr^{\leq i} A$.  In particular, the
	embedding $I^i/I^{i+1} A =: \Gr^i A \hookrightarrow \Gr^{\leq i} A$
	induces an embedding
	$$I^i/I^{i+1} \otimes_A M
	\iso \Gr^i A \otimes_{\Gr^{\leq i} A}
	\Gr^{\leq i} M;$$ 
	concretely, this means that the morphism 
	$$I^i/I^{i+1} \otimes_A M \to I^iM/I^{i+1}M$$
	is an embedding.  Letting $i$ vary, and recalling that $I$
	is nilpotent, we deduce 
	from~\cite[\href{http://stacks.math.columbia.edu/tag/0AS8}{Tag
          0AS8}]{stacks-project} that $M$ is flat over $A$.)

	If $\Gr^{\bullet} M$ is furthermore faithfully flat over
	$\Gr^{\bullet} A$, then (since faithful flatness is preserved
	under base-change) we see that $\Gr^0 M := M/IM$ is faithfully flat
	over $\Gr^0 A := A/IA$.  Because~$I$ is nilpotent,  an $A$-module vanishes
	if and only if its reduction mod $I$ does, and we conclude that $M$ is faithfully flat over $A$.
\end{proof}

\bibliographystyle{amsalpha}
\bibliography{universalBM}

\end{document}